\documentclass[a4,10pt]{article}
\usepackage{amsmath,amscd,amssymb}
\usepackage[margin=.8in]{geometry}

\newcommand{\nbiga}{\mathcal{A}}
\newcommand{\nbigb}{\mathcal{B}}
\newcommand{\nbigc}{\mathcal{C}}
\newcommand{\nbigd}{\mathcal{D}}
\newcommand{\nbige}{\mathcal{E}}
\newcommand{\nbigf}{\mathcal{F}}
\newcommand{\nbigg}{\mathcal{G}}
\newcommand{\nbigh}{\mathcal{H}}
\newcommand{\nbigi}{\mathcal{I}}
\newcommand{\nbigj}{\mathcal{J}}

\newcommand{\nbigl}{\mathcal{L}}
\newcommand{\nbigm}{\mathcal{M}}
\newcommand{\nbign}{\mathcal{N}}
\newcommand{\nbigo}{\mathcal{O}}
\newcommand{\nbigp}{\mathcal{P}}

\newcommand{\nbigr}{\mathcal{R}}
\newcommand{\nbigs}{\mathcal{S}}

\newcommand{\nbigu}{\mathcal{U}}
\newcommand{\nbigv}{\mathcal{V}}

\newcommand{\nbigx}{\mathcal{X}}

\newcommand{\nbigz}{\mathcal{Z}}

\newcommand{\proj}{\mathbb{P}}
\newcommand{\seisuu}{{\mathbb Z}}
\newcommand{\rnum}{{\mathbb Q}}

\newcommand{\cnum}{{\mathbb C}}
\newcommand{\real}{{\mathbb R}}

\newcommand{\hyperk}{\mathbb{K}}

\newcommand{\DD}{\mathbb{D}}
\newcommand{\EE}{\mathbb{E}}

\newcommand{\gbige}{\mathfrak E}

\newcommand{\gbigi}{\mathfrak I}

\newcommand{\gbigr}{\mathfrak R}
\newcommand{\gbigs}{\mathfrak S}

\newcommand{\gminia}{\mathfrak a}
\newcommand{\gminib}{\mathfrak b}
\newcommand{\gminic}{\mathfrak c}

\newcommand{\gminig}{\mathfrak g}

\newcommand{\gminim}{\mathfrak m}

\newcommand{\gminio}{\mathfrak o}

\newcommand{\vecalpha}{{\boldsymbol \alpha}}
\newcommand{\veca}{{\boldsymbol a}}
\newcommand{\vecb}{{\boldsymbol b}}

\newcommand{\vecdelta}{{\boldsymbol \delta}}

\newcommand{\vecm}{{\boldsymbol m}}

\newcommand{\vecI}{{\boldsymbol I}}

\newcommand{\vecW}{{\boldsymbol W}}

\newcommand{\larr}{\leftarrow}

\newcommand{\lrarr}{\longrightarrow}

\newcommand{\pf}{{\bf Proof}\hspace{.1in}}
\newcommand{\qed}{\mbox{\rule{1.2mm}{3mm}}}

\def\End{\mathop{\rm End}\nolimits}
\def\Ext{\mathop{\rm Ext}\nolimits}

\def\Cok{\mathop{\rm Cok}\nolimits}

\def\Image{\mathop{\rm Im}\nolimits}

\def\Re{\mathop{\rm Re}\nolimits}

\def\Gr{\mathop{\rm Gr}\nolimits}

\def\Cone{\mathop{\rm Cone}\nolimits}
\def\rank{\mathop{\rm rank}\nolimits}
\def\Spec{\mathop{\rm Spec}\nolimits}

\def\Ker{\mathop{\rm Ker}\nolimits}

\def\Gr{\mathop{\rm Gr}\nolimits}
\def\Sym{\mathop{\rm Sym}\nolimits}

\def\Res{\mathop{\rm Res}\nolimits}

\def\ord{\mathop{\rm ord}\nolimits}

\def\ch{\mathop{ch}\nolimits}

\def\Tr{\mathop{\rm Tr}\nolimits}

\def\id{\mathop{\rm id}\nolimits}

\def\codim{\mathop{\rm codim}\nolimits}

\def\ch{\mathop{\rm ch}\nolimits}

\def\Irr{\mathop{\rm Irr}\nolimits}

\newcommand{\del}{\partial}
\newcommand{\delbar}{\overline{\del}}

\newcommand{\nhom}{{\mathcal Hom}}

\newcommand{\sankaku}{\triangle}

\newcommand{\etabar}{\overline{\eta}}

\newcommand{\Abar}{\overline{A}}

\newcommand{\Sp}{{\mathcal Sp}}

\newcommand{\lefttop}[1]{{}^{#1}\!}

\def\reg{\mathop{\rm reg}\nolimits}

\newcommand{\Ehat}{\widehat{E}}

\newcommand{\Etilde}{\widetilde{E}}

\newcommand{\Vhat}{\widehat{V}}
\newcommand{\nablahat}{\widehat{\nabla}}

\newcommand{\nablatilde}{\widetilde{\nabla}}

\newcommand{\Chat}{\widehat{C}}

\newcommand{\ptilde}{\widetilde{p}}

\newcommand{\ftilde}{\widetilde{f}}

\newcommand{\stilde}{\widetilde{s}}

\newcommand{\Psitilde}{\widetilde{\Psi}}
\newcommand{\nbiglhat}{\widehat{\nbigl}}

\newcommand{\Rtilde}{\widetilde{R}}

\newcommand{\Fhat}{\widehat{F}}

\newcommand{\Rhat}{\widehat{R}}

\newcommand{\nbigehat}{\widehat{\nbige}}

\newcommand{\nbigvtilde}{\widetilde{\nbigv}}

\newcommand{\Utilde}{\widetilde{U}}

\newcommand{\Xtilde}{\widetilde{X}}

\newcommand{\nbigrtilde}{\widetilde{\nbigr}}

\newcommand{\Ltilde}{\widetilde{L}}

\def\ord{\mathop{\rm ord}\nolimits}

\def\Gal{\mathop{\rm Gal}\nolimits}

\def\Mod{\mathop{\rm Mod}\nolimits}

\def\Ch{\mathop{\rm Ch}\nolimits}
\def\good{\mathop{\rm good}\nolimits}

\def\Mero{\mathop{\rm Mero}\nolimits}

\def\sm{\mathop{\rm sm}\nolimits}

\def\CC{\mathop{\rm CC}\nolimits}
\def\DM{\mathop{\rm DM}\nolimits}
\def\red{\mathop{\rm red}\nolimits}
\def\Lag{\mathop{\rm Lag}\nolimits}
\def\irreg{\mathop{\rm irreg}\nolimits}

\def\RHB{\mathop{\rm RHB}\nolimits}
\def\td{\mathop{\rm td}\nolimits}

\def\app{\mathop{\rm app}\nolimits}

\newcommand{\gtilde}{\widetilde{g}}
\newcommand{\Ztilde}{\widetilde{Z}}

\newcommand{\varphitilde}{\widetilde{\varphi}}

\newcommand{\gbigetilde}{\widetilde{\gbige}}

\newcommand{\Stilde}{\widetilde{S}}
\newcommand{\tautilde}{\widetilde{\tau}}

\newcommand{\Atilde}{\widetilde{A}}

\newcommand{\gbigrtilde}{\widetilde{\gbigr}}

\newcommand{\Mtilde}{\widetilde{M}}

\newcommand{\Ybar}{\overline{Y}}

\newcommand{\Ttilde}{\widetilde{T}}

\newcommand{\nbigitilde}{\widetilde{\nbigi}}

\newcommand{\vecnbigi}{{\boldsymbol \nbigi}}

\newcommand{\nrhom}{R{\mathcal Hom}}

\newcommand{\vecH}{{\boldsymbol H}}

\newcommand{\Hhat}{\widehat{H}}
\newcommand{\Htilde}{\widetilde{H}}

\newcommand{\gminiatilde}{\widetilde{\gminia}}
\newcommand{\gminibtilde}{\widetilde{\gminib}}

\newcommand{\phitilde}{\widetilde{\phi}}

\newcommand{\gminiahat}{\widehat{\gminia}}

\newcommand{\nbigctilde}{\widetilde{\nbigc}}

\newcommand{\Sigmatilde}{\widetilde{\Sigma}}

\newcommand{\vecnbigf}{\boldsymbol \nbigf}

\newcommand{\nbigohat}{\widehat{\nbigo}}

\newcommand{\vecZ}{\boldsymbol Z}

\newcommand{\ttP}{{\tt P}}
\newcommand{\ttQ}{{\tt Q}}

\newcommand{\vecLambda}{\boldsymbol{\Lambda}}
\newcommand{\Sigmabar}{\overline{\Sigma}}

\newcommand{\Xhat}{\widehat{X}}

\newcommand{\nbigchat}{\widehat{\nbigc}}

\newcommand{\yhat}{\widehat{y}}
\newcommand{\nbigotilde}{\widetilde{\nbigo}}
\newcommand{\Locst}{\mathsf{Loc^{St}}}
\newcommand{\xhat}{\widehat{x}}

\newcommand{\nbigihat}{\widehat{\nbigi}}

\newcommand{\alphatilde}{\widetilde{\alpha}}
\newcommand{\btilde}{\widetilde{b}}
\newcommand{\vecSigma}{{\boldsymbol \Sigma}}
\newcommand{\Sigmahat}{\widehat{\Sigma}}
\newcommand{\Ihat}{\widehat{I}}
\newcommand{\zhat}{\widehat{z}}
\newcommand{\nbigjhat}{\widehat{\nbigj}}

\newtheorem{thm}{Theorem}[section]
\newtheorem{cor}[thm]{Corollary}

\newtheorem{rem}[thm]{Remark}
\newtheorem{lem}[thm]{Lemma}
\newtheorem{prop}[thm]{Proposition}
\newtheorem{df}[thm]{Definition}

\newtheorem{condition}[thm]{Condition}

\newtheorem{notation}[thm]{Notation}

\begin{document}

\title{Algebraic integrable connections with bounded irregularity}

\author{Takuro Mochizuki}
\date{}
\maketitle

\begin{abstract}
We study the boundedness of
families of algebraic flat connections
with bounded irregularity.
As an application.
we study the boundedness
of families of holonomic $\nbigd$-modules
with dominated characteristic cycles.
 
\vspace{.1in}
\end{abstract}

\section{Introduction}

\subsection{Boundedness}

\subsubsection{Boundedness for coherent sheaves}

We recall the notion of boundedness
for a class of coherent sheaves.
Let $X$ be a complex smooth and connected projective variety.
A class $\nbigc$ of coherent $\nbigo_X$-modules
is called bounded
if there exists a complex variety $S$ and
a coherent $\nbigo_{S\times X}$-module $\nbige$ flat over $S$
with the following property.
\begin{itemize}
 \item For any $E\in\nbigc$,
       there exists a closed point $s\in S$
       such that
       the restriction of $\nbige$
       to $s\times X$ is isomorphic to $E$.
\end{itemize}
For example, we may consider the class $\nbigc_1$ of
locally free $\nbigo_X$-modules
whose rank and Chern classes are fixed.
We may also consider the subclass $\nbigc_2\subset\nbigc_1$ of
semistable sheaves.
Then, it is easy to see that $\nbigc_1$ is unbounded unless $\dim X=0$.
It is well known that $\nbigc_2$ is bounded.

This kind of boundedness implies various finiteness properties.
For example,
if a class $\nbigc$ is bounded,
there exists $C>0$ such that
$\dim H^j(X,E)<C$
for any $E\in\nbigc$ and any $j\in\seisuu_{\geq 0}$.

\subsubsection{Boundedness for meromorphic flat bundles}

In this paper,
we are interested in
the boundedness of some families of meromorphic flat bundles.
Though we shall study meromorphic flat bundles
on projective varieties over any algebraically closed field
of characteristic $0$,
we restrict ourselves to the complex case
in this introduction,
and we use the complex analytic language.

Let $H$ be a hypersurface of $X$.
Let $\nbige$ be a coherent $\nbigo_X(\ast H)$-module
with an integrable connection $\nabla$.
Let $\Mero(X,H)$ denote the category of
such coherent $\nbigo_X(\ast H)$-modules with
an integrable connections.
For positive integers $r$ and $N$,
let $\Mero(X,H,r,N)$ denote the full subcategory
of $\Mero(X,H)$
consisting of 
$(V,\nabla)\in \Mero(X,H)$
satisfying $\rank(V)\leq r$
and the following condition.
\begin{itemize}
 \item There exists a lattice $\nbigl\subset V$
such that
$\nabla(\nbigl)\subset
       \nbigl\otimes\Omega^1(\log H)\otimes\nbigo_X(NH)$.
       Here, a lattice of $V$ means a coherent $\nbigo_X$-submodule
       which generates $V$ over $\nbigo_X(\ast H)$.
\end{itemize}

There are two types of boundedness results
for $\Mero(X,H,r,N)$.
First, we have the following boundedness
as meromorphic integrable connections
obtained as a special case of Corollary \ref{cor;26.2.23.2}.
\begin{thm}
\label{thm;26.2.23.21}
There exists a smooth complex variety $\nbigs_1$
and a coherent torsion-free $\nbigo_{\nbigs_1\times X}$-module
$E_{\nbigs_1}$
with a meromorphic integrable connection 
\[
 \nabla_{\nbigs_1}:
 E_{\nbigs_1}\to
 E_{\nbigs_1}\otimes\Omega^1_{\nbigs_1\times X/\nbigs_1}
 (\log (\nbigs_1\times H))
 \bigl(N(\nbigs_1\times H)\bigr)
\]
relative to $\nbigs_1$ such that the following holds.
\begin{itemize}
 \item $E_{\nbigs_1}$ is a reflexive $\nbigo_{\nbigs_1\times X}$-module
       flat over $\nbigs_1$.
       For any $s\in \nbigs_1$,
       $E_{\nbigs_1|a\times X}$ is reflexive.
 \item For any $(V,\nabla)\in \nbigc(X,H,r,\vecm)$,
       there exist $s\in\nbigs_1$
       such that
       $(E_{\nbigs_1},\nabla_{\nbigs_1})_{|s\times X}
       \otimes\nbigo_X(\ast H)
       \simeq (V,\nabla)$.
\end{itemize}
\end{thm}

We also have the boundedness as $\nbigd$-modules
obtained as a special case of Theorem \ref{thm;26.2.23.1}.
In general,
for any smooth morphism $X_T\to T$ of smooth complex algebraic varieties,
let $\nbigd_{X_T/T}$ denote the sheaf of
linear differential operators on $X_T$ relative to $T$.
For any $\nbigd_{X_T/T}$-module $\nbign$ which is flat over $T$,
and for any closed point $u\in T$,
let $\nbign_u$ denote the restriction of $\nbign$
to $X_T\times_Tu$.

\begin{thm}
\label{thm;25.11.1.2}
There exist a complex variety $\nbigs$
and a coherent $\nbigd_{\nbigs\times X/\nbigs}$-module $\nbigm$
flat over $\nbigs$
such that the following holds.
 \begin{itemize}
  \item For any $(V,\nabla)\in\Mero(X,H,r,N)$,
	there exists $s\in \nbigs(\cnum)$
	such that $\nbigm_s\simeq(V,\nabla)$.
 \end{itemize}
 \end{thm}
\begin{rem}
 $E_{\nbigs_1}(\ast(\nbigs_1\times H))$
in Theorem {\rm\ref{thm;26.2.23.21}}
is not necessarily $\nbigd_{\nbigs_1\times X/\nbigs_1}$-coherent,
in general.
\hfill\qed
\end{rem}
 
As a consequence of Theorem \ref{thm;25.11.1.2},
we obtain the following cohomological boundedness.
\begin{cor}
There exists a constant $C(X,H,r,N)>0$ such that
$\dim H^j(X,V\otimes\Omega^{\bullet})<C(X,H,r,N)$
for any $(V,\nabla)\in\Mero(X,H,r,N)$. 
\end{cor}

According to \cite{Hu-Teyssier,kedlaya3},
the boundedness of $\Mero(X,H,r,N)$ was asked by Esnault and Langer.
Hu and Teyssier \cite{Hu-Teyssier}
proved the cohomological boundedness
in the case $\dim X\leq 2$.
They obtained much more detailed upper bound for the constant $C(X,H,r,N)$.
In the case $X=\proj^1$,
Sabbah \cite{Sabbah-DM}
studied the finiteness of rigid meromorphic flat bundles
by using the Arinkin-Katz program.

\subsubsection{Boundedness for holonomic $\nbigd$-modules
with dominated characteristic cycles}

Related with the boundedness as in
Theorem \ref{thm;26.2.23.21} and Theorem \ref{thm;25.11.1.2},
we study the boundedness of
the family of holonomic $\nbigd$-modules
whose characteristic cycles are dominated.

For any holonomic $\nbigd_X$-module $M$,
there exists a good filtration $F$.
We obtain the coherent $\nbigo_{T^{\ast}X}$-module $\Gr^F(M)^{\sim}$.
The support $\Ch(M)$ is a complex Lagrangian cone in $T^{\ast}X$,
and called the characteristic variety.
Let $\Ch(M)=\bigcup_i \Ch(M)_i$
denote the irreducible decomposition.
The characteristic cycle
$\CC(M)=\sum_i \CC(M)_i\Ch(M)_i$
is a refined invariant,
where $\CC(M)_i$ denotes the multiplicity
of $\Ch(M)_i$ in $\Gr^F(M)^{\sim}$.

Let $\Lambda_i$ $(i=1,\ldots,\ell)$ be complex Lagrangian cones
in $T^{\ast}X$.
We obtain a cycle $\Lambda=\sum_{i=1}^{\ell} m_i\Lambda_i$.
We say $\CC(M)\leq \Lambda$
if $\CC(M)=\sum_{i=1}^{\ell}a_i\Lambda_i$
with $a_i\leq m_i$.
Let $\nbigc(X,\Lambda)$
denote the family of holonomic $\nbigd_X$-modules $M$
such that $\CC(M)\leq\Lambda$.

\begin{thm}[Theorem
 \ref{thm;25.12.14.20},
 Lemma \ref{lem;26.2.24.1}]
\label{thm;25.12.14.150}
The family $\nbigc(X,\Lambda)$ is bounded
in the sense that
there exists a complex algebraic variety $S$
and a coherent $\nbigd_{S\times X/S}$-module
 $\nbigm$ such that the following holds.
 \begin{itemize}
  \item $\nbigm$ is flat over $S$,
	and 
       $\nbigm_s\in\nbigc(X,\Lambda)$
	for any $s\in S$.	
  \item For any $M\in\nbigc(X,\Lambda)$,
	there exists $s\in S$
	such that
	$M\simeq \nbigm_s$.
 \end{itemize}
 \end{thm}

\subsection{Good meromorphic flat bundles}
\label{subsection;25.12.14.100}

Let $X$ be a complex projective manifold
with a normal crossing hypersurface $H$.
We recall the notion of good meromorphic flat bundles.

\subsubsection{Good set of ramified irregular values}
\label{subsection;26.2.27.1}

Let $\nbigo_X$ denote the sheaf of holomorphic functions.
Let $\nbigo_X(\ast H)$ denote the sheaf of
meromorphic functions allowing poles along $H$.
The stalks at $P\in H$ are denoted by
$\nbigo_{X,P}$ and $\nbigo_X(\ast H)_P$.
The completion of $\nbigo_{X,P}$ is denoted by $\nbigohat_{X,P}$.
We set
$\nbigohat_X(\ast H)_P=\nbigo_X(\ast H)_P\otimes\nbigohat_{X,P}$.
Let $(z_1,\ldots,z_n)$ be a local coordinate around $P$
such that $H=\bigcup_{j=1}^{\ell}\{z_j=0\}$ around $P$.
For any positive integer $k$,
we set
$\nbigo_{X,P}^{(H,k)}
=\nbigo_{X,P}[z_1^{1/k},\ldots,z_{\ell}^{1/k}]$.
It is independent of the choice of the coordinate system.
We also set
$\nbigohat^{(H,k)}_{X,P}
=\nbigo^{(H,k)}_{X,P}\otimes\nbigohat_{X,P}$.
We consider
$\nbigo^{(H,k)}_X(\ast H)_P$
and
$\nbigohat^{(H,k)}_X(\ast H)_P$
similarly.

We recall a notion of order
for an element
$\gminia$ of $\nbigohat^{(H,k)}_{X}(\ast H)_P$.
\begin{itemize}
 \item 
If $\gminia$ is contained in $\nbigohat_{X,P}^{(H,k)}$,
we set
$\ord(\gminia)=(0,\ldots,0)$
in $\seisuu^{\ell}$.
\item
     For any $\gminia\in
     \nbigohat^{(H.k)}_X(\ast H)_P\setminus\nbigohat_{X,P}^{(H,k)}$,
if there exists $\vecm\in (\frac{1}{k}\seisuu)_{\leq 0}^{\ell}$
such that
     $\prod_{i=1}^{\ell}z_i^{-m_i}\gminia$
     is an invertible element of
     $\nbigohat^{(H,k)}_{X,P}$,
     we set
     $\ord(\gminia)=\vecm$.
 \item
      Otherwise,
      $\ord(\gminia)$ is not defined.
\end{itemize}
For example, $\ord(z_1^{-1}z_2^{-1})=(-1,-1)$.
But,
$\ord(z_1^{-1}z_2)$
and $\ord(z_1^{-1}+z_2^{-1})$ are not defined.
We can also consider $\ord(\gminia)$ 
for any $\gminia\in
\nbigohat_{X}^{(H,k)}(\ast H)_P/
\nbigohat^{(H,k)}_{X,P}$.

We define the partial order on $(\frac{1}{k}\seisuu)^{\ell}$
by  $\vecm\leq\vecm'\Longleftrightarrow \forall i,m_i\leq m_i'$.
\begin{df}
Let 
$\nbigi\subset
\nbigohat^{(H,k)}_{X}(\ast H)_P/\nbigohat^{(H,k)}_{X,P}$
be a finite subset which is invariant under the Galois action.
It is called good,
if the following conditions are satisfied.
\begin{itemize}
 \item 
$\ord(\gminia)$
are defined for any $\gminia$ in $\nbigi$,
and $\ord(\gminia-\gminib)$ are defined
       for any $\gminia,\gminib$ in $\nbigi$.
\item  The sets
$\{\ord(\gminia)\,|\,\gminia\in\nbigi\}$
and
       $\{\ord(\gminia-\gminib)\,|\,
       \gminia,\gminib\in\nbigi\}$
are totally ordered
with respect to the partial order $\leq$ on
$(\frac{1}{k}\seisuu)^{\ell}$.
\hfill\qed
\end{itemize}
\end{df}

For example,
the set $\{z_1^{-1}z_2^{-1},z_1^{-1},0\}$ is good,
but
$\{z_1^{-1}z_2^{-1},z_1^{-1},z_2^{-1},0\}$ is not good.

\begin{rem}
\label{rem;26.2.27.2}
Though the concepts are explained 
for finite subsets of 
 $\cnum[\![z_1^{1/k},\ldots,z_{\ell}^{1/k},z_{\ell+1},\ldots,z_n]\!]
 [z_1^{-1},\ldots,z_{\ell}^{-1}]$,
they are obviously generalized to
the context of 
finite subsets of 
$R[\![z_1^{1/k},\ldots,z_{\ell}^{1/k},z_{\ell+1},\ldots,z_n]\!]
[z_1^{-1},\ldots,z_{\ell}^{-1}]$,
for any ring $R$.
\hfill\qed
\end{rem}

\subsubsection{Good meromorphic flat bundles}

Let $\nbige$ be a locally free $\nbigo_X(\ast H)$-module
with an integrable connection.
For $P\in H$,
we set
$\nbigehat_P^{(H,k)}=\nbige_{P}\otimes\nbigohat^{(H,k)}_{X,P}$.
We say $(\nbige,\nabla)$ is good at $P$
if there exist a good set
$\nbigi_P(\nbige,\nabla)
\subset\nbigohat^{(H,k)}_{X}(\ast P)/\nbigohat^{(H,k)}_{X,P}$
and a decomposition
\[
 (\nbigehat^{(H,k)}_{P},\nabla)
 =\bigoplus_{\gminia\in\nbigi_P(\nbige,\nabla)}
 \bigl(
 (\nbigehat^{(H,k)}_{P})_{\gminia},\nabla_{\gminia}
 \bigr)
\]
such that
each $(\nbigehat^{(H,k)}_P)_{\gminia}$ has a lattice
$\nbiglhat_{P,\gminia}$ such that
$(\nabla_{\gminia}-d\gminiatilde\id)\nbiglhat_{P,\gminia}
\subset
\nbiglhat_{P,\gminia}\otimes\Omega_X^1(\log H)$.
Here, $\gminiatilde\in\nbigohat^{(H,k)}_X(\ast H)_P$ is a lift of $\gminia$.
We say $(\nbige,\nabla)$ is good on $(X,H)$
if it is good at any $P\in H$.

\begin{rem}
 It is known that
 $\nbigi_P(\nbige,\nabla)
\subset\nbigo^{(H,k)}_{X}(\ast P)/\nbigo^{(H,k)}_{X,P}$.
\hfill\qed
\end{rem}

\subsubsection{Boundedness of
good meromorphic flat bundles as meromorphic objects}

Let $\nbigx_S\to S$ be a smooth projective morphism of
smooth complex algebraic varieties
such that each fiber is connected.
Let $\nbigh_S\subset\nbigx_S$ be
a relatively normal crossing hypersurface.

For any $s\in S$,
we set $(\nbigx_s,\nbigh_s)=(\nbigx_S,\nbigh_S)\times_Ss$.
Let $\Mero^{\good}(\nbigx_s,\nbigh_s,r,N)$
denote the full subcategory
of $\Mero(\nbigx_s,\nbigh_s,r,N)$,
consisting of good meromorphic flat bundles
contained in $\Mero(\nbigx_s,\nbigh_s,r,N)$.
Let $S(\cnum)$ denote the set of the closed points of $S$.

\begin{thm}[Theorem
\ref{thm;25.10.15.40}]
\label{thm;25.10.31.1}
The family
$\Mero^{\good}(\nbigx_s,\nbigh_s,r,N)$ $(s\in S(\cnum))$
is bounded in the following sense.
\begin{itemize}
 \item There exists a smooth complex variety
       $\nbigs$ over $S$.
       We set $(\nbigx_{\nbigs},\nbigh_{\nbigs})
       :=(\nbigx_S,\nbigh_S)\times_S\nbigs$.
 \item There exists a locally free
       $\nbigo_{\nbigx_{\nbigs}}$-module
       $\nbigl_{\nbigs}$
       equipped with
       a meromorphic integrable relative connection
       $\nabla_{\nbigs}:\nbigl_{\nbigs}
       \to\nbigl_{\nbigs}
       \otimes\Omega^1_{\nbigx_{\nbigs}/\nbigs}(\log \nbigh_{\nbigs})
       (N\nbigh_{\nbigs})$.
 \item For any $s\in S(\cnum)$
       and $(V,\nabla)\in\Mero^{\good}(\nbigx_s,\nbigh_s,r,N)$,
       there exists
       $\stilde\in\nbigs(\cnum)$ over $s$
       such that
       $(V,\nabla)\simeq
       \bigl(
       \nbigl(\ast \nbigh_{\nbigs}),
       \nabla
       \bigr)_{|\stilde\times_{\nbigs}\nbigx_{\nbigs}}$.
       Moreover,
       $\nbigl_{|\stilde\times_S\nbigx_S}$
       equals
       the good Deligne-Malgrange lattice of
       $(V,\nabla)$.
\end{itemize}
 \end{thm}

Let us mention that 
the boundedness in Theorem \ref{thm;25.10.31.1}
is reduced to a more classical boundedness result
of torsion-free coherent sheaves due to Maruyama
\cite{Maruyama} (see Proposition \ref{prop;25.12.14.110})
once we obtain the estimate for
the first and second Chern classes
of a good Deligne-Malgrange lattice
in \S\ref{section;25.12.14.120}.

\subsubsection{Resolution of turning points}

If $\dim X=1$,
any meromorphic flat bundle is good
by the classical Hukuhara-Levelt-Turrittin theorem.
However, if $\dim X\geq 2$,
there exist meromorphic flat bundles which are not good.
There is the following theorem
which was conjectured by Sabbah
and proved by Kedlaya in \cite{kedlaya,kedlaya2}
and the author in \cite{Mochizuki-surface, Mochizuki-wild}.
Note that the result of Kedlaya is more general
and available for the complex analytic case
not only for the algebraic case.

\begin{thm}
There exists a projective morphism
 $\varphi:X'\to X$
such that  
(i) $H'=\varphi^{-1}(H)$ is normal crossing,
(ii) $X'\setminus H'\simeq X\setminus H$,
(iii) $\varphi^{\ast}(\nbige,\nabla)$ is good on $(X',H')$.
\hfill\qed
\end{thm}

For the boundedness of the family of meromorphic flat bundles
which are not necessarily good,
we have to control the resolutions.

\subsubsection{Resolutions and boundedness}

Let $X$ be a smooth projective irreducible complex variety
with a simple normal crossing hypersurface $H$.
Let $r$ and $N$ be positive integers.
We have the following more useful boundedness theorem
for $\nbigc(X,H,r,N)$
as a special case
of Theorem \ref{thm;25.12.14.30} and Theorem \ref{thm;26.2.7.31}.

\begin{thm}
\label{thm;26.2.23.20}
There exist a smooth complex variety $\nbigs$,
a complex variety $\Xhat_{\nbigs}$ which is smooth projective over $\nbigs$
with a morphism $\rho_{\nbigs}:
\Xhat_{\nbigs}\to \nbigs\times X$  over $\nbigs$,
a hypersurface $\Hhat_{\nbigs}\subset \Xhat_{\nbigs}$
which is normal crossing relative to $\nbigs$,
a good meromorphic flat bundle $(\nbigv,\nabla)$ on
$(\Xhat_{\nbigs},\Hhat_{\nbigs})$ relative to $\nbigs$,
and a lattice $\Ehat_{\nbigs}\subset \nbigv$ 
such that the following holds.
\begin{itemize}
 \item $\rho^{-1}_{\nbigs}(\nbigs\times H)=\Hhat_{\nbigs}$,
       and $\Xhat_{\nbigs}\setminus\Hhat_{\nbigs}\simeq
       \nbigs\times(X\setminus H)$.
 \item For any $(V,\nabla)\in \nbigc(X,H,r,N)$,
       there exists $s\in\nbigs$
       such that
       $\rho_s^{\ast}(V,\nabla)
       \simeq
       (\nbigv,\nabla)_s$
       and that
       $\Ehat_s$
       is the Deligne-Malgrange lattice of
       $\rho_s^{\ast}(V,\nabla)$.
\end{itemize}
There also exists a coherent $\nbigd_{\Xhat_{\nbigs}/\nbigs}$-module
$\nbigm$ flat over $\nbigs$ such that the following holds.
\begin{itemize}
 \item $\nbigm(\ast\Hhat_{\nbigs})=\nbigv$.
 \item For any $(V,\nabla)\in \nbigc(X,H,r,N)$,
       there exists $s\in\nbigs$
       such that
       $\rho_s^{\ast}(V,\nabla)
       \simeq
       \nbigm_s$.
\end{itemize}
\end{thm}

We obtain Theorem \ref{thm;26.2.23.21} and Theorem \ref{thm;25.11.1.2}
from Theorem \ref{thm;26.2.23.20}
by taking the push-forward.
We also apply Theorem \ref{thm;26.2.23.20}
to study Theorem \ref{thm;25.12.14.150}.

\subsection{Outline for Theorem \ref{thm;26.2.23.20}}

\subsubsection{Meromorphic Lagrangian covers}

Let $\Sigma\subset T^{\ast}(X\setminus H)$
be a closed algebraic subset such that
each irreducible component is dominant and proper over $X\setminus H$.
We also assume that the smooth part of $\Sigma$ is Lagrangian.
Such $\Sigma$ is called a meromorphic Lagrangian cover.
We introduce
a condition for good behaviour of meromorphic Lagrangian covers
around the boundary.
\begin{df}
Let $x\in H$.
We say that $\Sigma$ is good at $x$
if there exist a good set of ramified irregular values
$\nbigi(\Sigma,x)\subset
\nbigo^{(H,k)}_{X}(\ast H)_x/\nbigo^{(H,k)}_{X,x}$,
and a covering $\rho_x: X_x'\to X_x$ branched along $H\cap X_x$
such that
\[
\rho_x^{\ast}\Sigma
=\bigsqcup_{\gminia\in\nbigi(\Sigma,x)}
\bigl(
d\gminia+
\Sigma_{\gminia}
\bigr),
\]
where $\Sigma_{\gminia}$ are logarithmic,
i.e.,
the closure in $T^{\ast}X_x'(\log H_x')$ is proper over $X_x'$.
Here, $H_x'=\rho_x^{-1}(H_x)$.

We say $\Sigma$ is good on $(X,H)$
if $\Sigma$ is good at any $x\in H$.
\hfill\qed
\end{df}
 
For positive integers $r$ and $N$,
let $\Lag(X,H,r,N)$ denote
the set of meromorphic Lagrangian covers
$\Sigma$
such that
the closure $\overline{\Sigma}$
in $T^{\ast}X(\log H)(NH)$ is proper over $X$,
and the degree of the projection $\Sigma\to X\setminus H$
is less than $r$.
This family is bounded,
and we can take a resolution in a uniform way.

\begin{prop}[Corollary
\ref{cor;26.2.23.10}]
\label{prop;25.11.1.1}
There exist a complex smooth variety $S$,
a closed subset
$\nbigz\subset S\times T^{\ast}(X\setminus H)$
flat over $S$,
and a projective morphism 
$\varphi:\nbigx_S\to S\times X$
such that the following holds.
\begin{itemize}
 \item $\nbigx_S\to S$ is smooth and projective,
       and $\nbigh_S=\varphi^{-1}(S\times H)$
       is normal crossing relative to $S$.
\item $\nbigx_S\setminus \nbigh_S\simeq S\times (X\setminus H)$,
\item  For any $\Sigma\in \Lag(X,H,r,N)$,
       there exists $s\in S(\cnum)$ such that $\nbigz\times_Ss=\Sigma$.
       Moreover,
       $\varphi_s^{\ast}(\Sigma)$ is good on
       $(\nbigx_s,\nbigh_s)$,
       where $\varphi_s:\nbigx_S\times_Ss\to X$
       is induced by $\varphi$.
       \hfill\qed
\end{itemize}
 \end{prop}

\subsubsection{Wild harmonic bundles}

Let $(E,\theta)$ be a Higgs bundle on a complex manifold $Y$.
Let $h$ be a Hermitian metric of $E$.
We obtain the Chern connection $\nabla_h$,
and the adjoint $\theta^{\dagger}_h$ of $\theta$.
We say $(E,\theta,h)$ is a harmonic bundle
if the connection $\DD_h^1=\nabla_h+\theta+\theta^{\dagger}_h$ is flat.
For the Higgs bundle $(E,\theta)$,
we have the corresponding $\nbigo_{T^{\ast}Y}$-module $\Etilde$,
whose supports are finite over $Y$.
The support $\Sigma_{\theta}$ is called the spectral cover.
By using Gabber's theorem \cite{Gabber},
if $(E,\theta)$ underlies a harmonic bundle $(E,\theta,h)$, 
$\Sigma_{\theta}$ is Lagrangian.

Let $X$ be a smooth projective complex variety
with a normal crossing hypersurface $H$.
Let $(E,\theta,h)$ be a harmonic bundle on $X\setminus H$.
We say $(E,\theta,h)$ is wild on $(X,H)$
if $\Sigma_{\theta}$ is a meromorphic Lagrangian cover,
i.e., $\Sigma_{\theta}$ is an algebraic subset of
$T^{\ast}(X\setminus H)$.
The harmonic bundle is called good on $(X,H)$
if $\Sigma_{\theta}$ is good on $(X,H)$.

\subsubsection{The associated meromorphic objects}

From a harmonic bundle on $(E,\delbar_E,\theta,h)$ on $X\setminus H$,
we have the flat bundle $(E,\DD^1_h)$ on $X\setminus H$.
We have the underlying holomorphic vector bundle
$\nbige^1=(E,\delbar_E+\theta^{\dagger}_h)$.
We regard it as an $\nbigo_{X\setminus H}$-module.
Let $\iota:X\setminus H\to X$ denote the inclusion.
We obtain the sheaf
$\iota_{\ast}(\nbige^1)$ on $X$.

For any $x\in H$,
let $(X_x,z_1,\ldots,z_n)$ be a coordinate neighbourhood
such that $X_x\cap H=\bigcup_{j=1}^{\ell}\{z_j=0\}$.
We obtain an $\nbigo_X(\ast H)_x$-module
\[
 (\nbigp^h\nbige^1)_x:=
 \Bigl\{
 s\in (\iota_{\ast}\nbige^1)_x\,\Big|\,
 |s|_h=O\Bigl(\prod_{i=1}^{\ell}|z_i|^{-C}\Bigr)
 \,\,\,
 \exists C>0
 \Bigr\}.
\]
We obtain an $\nbigo_X(\ast H)$-module
$\nbigp^h\nbige^1\subset\iota_{\ast}\nbige^1$
whose stalk at $x$ is $(\nbigp^h\nbige^1)_x$.
We recall the following proposition from \cite{Mochizuki-wild}.
\begin{prop}
\mbox{{}}
\begin{itemize}
 \item 
$\nbigp^h\nbige^1$ is a coherent $\nbigo_X(\ast H)$-module,
       and $\DD^1_h$ is an integrable connection
       of $\nbigp^h\nbige^1$.
Moreover,
$(\nbigp^h\nbige^1,\DD^1_h)\in\Mero(X,H,r,N)$
if and only if
$\Sigma_{\theta}\in\Lag(X,H,r,N)$.
\item
$(E,\delbar_E,\theta,h)$ is good on $(X,H)$
if and only if $(\nbigp^h\nbige^1,\DD^1_h)$ is good on $(X,H)$.
We have
$\nbigi(\nbigp^h\nbige^1,\DD^1_h,x)
 =\bigl\{
 2\gminia\,\big|\,\gminia\in\nbigi(\Sigma_{\theta},x)
 \bigr\}$ for any $x\in H$.
\hfill\qed
\end{itemize}
\end{prop}

We have the following theorem \cite{Mochizuki-wild}.
\begin{thm}
For any semisimple meromorphic flat bundle $(V,\nabla)$ on $(X,H)$,
there exists a wild harmonic bundle $(E,\theta,h)$ such that 
$(\nbigp^h\nbige^1,\DD^1_h)\simeq (V,\nabla)$. 
\hfill\qed
\end{thm}

Let $(V,\nabla)$ be a meromorphic flat bundle on $(X,H)$.
There exists a Jordan-H\"older filtration,
namely,
each graded piece is simple.
By the above theorem,
we have the harmonic bundle $(E,\theta)$
corresponding to
the associated graded meromorphic flat bundle.
We set $\Sigma(V,\nabla):=\Sigma_{\theta}$.
\begin{cor}
Let $\varphi:X'\to X$ be a projective morphism
such that
(i) $H'=\varphi^{-1}(H)$ is normal crossing,
(ii) $\varphi^{\ast}\Sigma(V,\nabla)$ is good on $(X',H')$.
Then, $\varphi^{\ast}(V,\nabla)$  is good on $(X',H')$.
\hfill\qed
\end{cor}

\begin{rem}
We shall introduce the notion of 
meromorphic Lagrangian irregularities,
which is equivalences 
of meromorphic Lagrangian covers.
We shall associate a meromorphic Lagrangian irregularity
with a meromorphic flat bundle.
\hfill\qed
\end{rem}

\subsubsection{Outline of a proof of Theorem \ref{thm;26.2.23.20}}

\label{subsection;25.12.14.101}

For $(X,H,r,N)$, we have
a smooth complex variety $S$
and a projective morphism
$\varphi:\nbigx_S\to S\times X$ over $S$
with a normal crossing hypersurface
$\nbigh_S$
as in Proposition \ref{prop;25.11.1.1}.
There exists $N'>0$ such that
the closure of $\nbigz$
in $T^{\ast}(\nbigx_{S}/S)(\log \nbigh_{S})(N'\nbigh_S)$
is proper over $\nbigx_S$.
There exists
$\nbigs\to S$
and
$(\nbigl_{\nbigs},\nabla)$
as in Theorem \ref{thm;25.10.31.1} with $N'$.

We have the induced morphism
$\varphi_{\nbigs}:
 \nbigx_{\nbigs}\to \nbigs\times X$,
induced by $\nbigx_S\to S\times X$.
We set
\[
(\nbigv,\nabla)=
\bigl(
 (\varphi_{\nbigs})_{\ast}
 \bigl(
 \nbigl(\ast \nbigh_{\nbigs})
 \bigr),
 \nabla
 \bigr).
\]
By taking a stratification,
we may assume that
it is flat over $\nbigs$.
We may also assume that $(\nbigv,\nabla)$ has a good lattice
$\nbigl_1$ relative to $\nbigs$.

Let $(V,\nabla)\in\Mero(X,H,r,N)$.
We have $\Sigma(V,\nabla)\in \Lag(X,H,r,N)$.
There exists $s\in S$ such that
$\varphi_s^{\ast}(\Sigma(V,\nabla))$ is good on
$(\nbigx_s,\nbigh_s)=(\nbigx_S,\nbigh_S)_s$.
It implies
$\varphi_s^{\ast}(V,\nabla)$ is good on
$(\nbigx_s,\nbigh_s)$.
There exists $\stilde\in \nbigs$ over $s$
such that
$(\nbigl\otimes
\nbigo_{\nbigx_{\nbigs}}(\ast \nbigh_{\nbigs}),\nabla)_{\stilde}
\simeq \varphi_s^{\ast}(V,\nabla)$.
We can construct a
coherent $\nbigd_{\nbigx_{\nbigs}/\nbigs}$-submodule
$\nbigm$ with the desired property
by considering the $\nbigd_{\nbigx_{\nbigs}/\nbigs}$-module
generated by $\nbigl_1(N_1\nbigh_{\nbigs})$ for a sufficiently
large $N_1$.

\paragraph{Acknowledgement}

This study was partially inspired by
the idea of arithmetic support due to Maxim Kontsevich.
The boundedness of meromorphic flat bundles 
is one of three topics in  the author's talk in the conference
``Mathematics on the Crossroad of Centuries''
held at IHES in 2024 September.
The author also gave talks about this topic
in Algebra-Geometry-Number theory held
at Institute of Mathematics, Vietnam Academy of Science and Technology
in 2025 November,
and 
in Workshop on Integral Stokes Structures and Applications
held at Mannheim University in 2026 January.
I thank the organizers for the opportunities.

I thank Claude Sabbah for discussions in various occasions
and for his kindness.

I am partially supported by
the Grant-in-Aid for Scientific Research (A) (No. 21H04429),
the Grant-in-Aid for Scientific Research (A) (No. 22H00094),
the Grant-in-Aid for Scientific Research (A) (No. 23H00083),
the Grant-in-Aid for Scientific Research (C) (No. 20K03609),
the Grant-in-Aid for Scientific Research (C) (No. 25K06973)
Japan Society for the Promotion of Science.
I am also partially supported by the Research Institute for Mathematical
Sciences, an International Joint Usage/Research Center located in Kyoto
University.

\section{Lagrangian meromorphic covers}

\subsection{Symmetric products of vector spaces}

\subsubsection{The Chow embedding}
\label{subsection;25.10.1.1}

Let $\hyperk$ be an algebraically closed field of characteristic $0$.
Let $V$ be a finite dimensional $\hyperk$-vector space.
We also regard $V$ as a $\hyperk$-variety.
Let $r$ be a positive integer.
Let $\gbigs_r$ denote the $r$-th symmetric group
acting on
$V^r=\overbrace{V\times \cdots\times V}^r$
by permutation of the components.
Let $S^r(V)$ denote the quotient space,
i.e., $S^r(V)=V^r/\gbigs_r$.
The image of $(v_1,\ldots,v_r)\in V^r$ in $S^r(V)$
is denoted by $\prod_{i=1}^r v_i$.
Let $\bigl[\prod_{i=1}^{r}v_i\bigr]\subset V$
denote the subset induced by
the tuple $v_1,\ldots,v_r$.

Let $\Sym^r(\hyperk\oplus V)$
denote the vector space obtained as
the $r$-th symmetric tensor product
of $\hyperk\oplus V$.
We also regard it as a $\hyperk$-variety.
Let $\iota:V\to \hyperk\oplus V$ be given by
$\iota(v)=1\oplus v$.
We obtain the morphism of $\hyperk$-varieties
$S^r(V) \to \Sym^r(\hyperk\oplus V)$
defined by
\[
\prod_{i=1}^rv_i
\longmapsto
\prod_{i=1}^r \iota(v_i). 
\]
According to the theorem of Chow and van der Wearden
(see \cite[\S4.1, Theorem 1.1]{GKZ}),
it is a closed immersion of algebraic $\hyperk$-varieties.

\subsubsection{Partitions}
\label{subsection;25.10.1.3}

For any closed point $\alpha\in S^r(V)$,
there exist 
$\alpha_i\in V$ and $m_i\in\seisuu_{>0}$
$(i=1,\ldots,\ell(\alpha))$
such that 
$\alpha_i\neq \alpha_j$ $(i\neq j)$,
$m_i\geq m_j$ $(i\leq j)$
and 
$\alpha=\prod_{i=1}^{\ell(\alpha)} \alpha_i^{m_i}$.
The tuple $(m_i\,|\,i=1,\ldots,\ell(\alpha))$
is a partition of $r$,
denoted by $\ttP(\alpha)$.
The number $\ell(\alpha)$
and the partition $\ttP(\alpha)$
are well defined for $\alpha\in S^r(V)$.

\subsubsection{Stratification by partitions}
\label{subsection;24.8.27.1}

Let 
$\ttP=(m_1,\ldots,m_{\ell})$
be a partition of $r$,
i.e., $\ttP$ is a decreasing sequence of
positive integers
such that $\sum m_i=r$.
We set $I_m(\ttP)=\{1\leq i\leq \ell\,|\,m_i=m\}$.
We obtain the decomposition
$\{1,\ldots,\ell\}=\bigsqcup I_m(\ttP)$
by which $V^{\ell}=\prod V^{I_m(\ttP)}$.
Let $\gbigs_{I_m(\ttP)}$ denote
the group of automorphisms of the finite subset $I_m(\ttP)$,
which is isomorphic to $\gbigs_{|I_m(\ttP)|}$.
The natural $\gbigs_{I_m(\ttP)}$
on $V^{I_m(\ttP)}$
induces
the action of
$\prod\gbigs_{I_m(\ttP)}$ on $V^{\ell}$.
The quotient is naturally isomorphic to
$\prod S^{|I_m(\ttP)|}(V)$.

We obtain the finite morphism
$\Psi_{\ttP}:V^{\ell}
\to S^r(V)$
defined by
$\Psi_{\ttP}(v_1,\ldots,v_{\ell})
=\prod_{i=1}^{\ell}v_i^{m_i}$.
Let $S^r_{\ttP}(V)$ denote the image of 
$\Psi_{\ttP}$ with the reduced structure.
It is also the image of the induced morphism
$\prod S^{|I_m(\ttP)|}(V)\to S^r(V)$.
Let 
$(V^{\ell})^{\circ}
\subset V^{\ell}$
denote the Zariski open subset of
$(v_1,\ldots,v_{\ell})$ such that
$v_i\neq v_j$ $(i\neq j)$.
Let $S^r_{\ttP}(V)^{\circ}$ denote the image of
$(V^{\ell})^{\circ}$ by $\Psi_{\ttP}$.
It is a Zariski open subset in
$S^r_{\ttP}(V)$.
We have
$(V^{\ell})^{\circ}/\prod \gbigs_{I_m(\ttP)}
\simeq
S^r_{\ttP}(V)^{\circ}$.

Let $\ttP_1=(m_i\,|\,i=1,\ldots,\ell)$ and
$\ttP_2=(n_j\,|\,j=1,\ldots,k)$ be partitions of $r$.
We say  $\ttP_1\prec\ttP_2$
if the following holds.
\begin{itemize}
 \item There exists a surjection
       $\varphi:\{1,\ldots,k\}\to\{1,\ldots,\ell\}$
       such that
       $m_i=\sum_{\varphi(j)=i}n_j$.
\end{itemize}
This defines a partial order $\prec$ on the set of partitions of $r$.
We have
$S^r_{\ttP}(V)
=\bigsqcup_{\ttP\prec\ttP'}
 S^r_{\ttP'}(V)^{\circ}$.
In particular,
we obtain the stratification of $S^r(V)$
by locally Zariski closed subsets
$S^r(V)
 =\bigsqcup_{\ttP}S^r_{\ttP}(V)^{\circ}$.

\subsubsection{Universal subsets}
\label{subsection;25.10.1.10}

For $j=1,\ldots,r$,
we set
$\Delta_{0,j}(V)
=\bigl\{
 (v_0,v_1,\ldots,v_{r})\,\big|\,
 v_0=v_j
 \bigr\}
\subset V\times V^{r}$.
We set
$Z(V^{r})=\bigcup_{j=1}^{r} \Delta_{0,j}(V)$.
The $\gbigs_{r}$-action on $V^{r}$
induces an $\gbigs_{r}$-action
on $Z(V^r)$.
Let $Z(S^r(V))$ denote the quotient space.
It is naturally a closed subvariety of
$V\times S^r(V)$.
We have the naturally defined finite morphism
$Z(S^r(V))\to S^{r}(V)$.

Let $\ttP=(m_1,\ldots,m_{\ell})$ be a partition of $r$.
We obtain
\[
 Z(S^r_{\ttP}(V)^{\circ})
 :=Z(S^r(V))\times_{S^r(V)}
 S_{\ttP}^r(V)^{\circ},
 \quad\quad
 Z(S^r_{\ttP}(V)):=
 Z(S^r(V))\times_{S^r(V)}
 S_{\ttP}^r(V).
\]
Let $Z(S^r_{\ttP}(V)^{\circ})_{\red}$
and $Z(S^r_{\ttP}(V))_{\red}$
denote the subvarieties
of $Z(S^r(V))$ with the reduced structure.
The projection
$Z(S^r_{\ttP}(V)^{\circ})_{\red}\to S^r_{\ttP}(V)^{\circ}$
is a proper etale cover of degree $\ell$.
The projection
$Z(S^r_{\ttP}(V))_{\red}\to S^r_{\ttP}(V)$
is a finite morphism.

\subsubsection{Induced finite morphisms}

Let $Y$ be an irreducible $\hyperk$-variety
with a morphism $\phi:Y\to S^r(V)$ of algebraic $\hyperk$-varieties.
Let $\Sigma(\phi)$ denote the fiber product of $\phi$
and $Z(S^r(V))\to S^r(V)$.
The natural morphism
$\Sigma(\phi)\to Y$ is finite.
For any closed point $y\in Y$,
the set of the closed points of
$\Sigma(\phi)\times_Yy$
is the tuple given by $\phi(y)\in S^r(V)$.

There exists a unique partition $\ttP(\phi)$
such that
$\phi(Y)\subset S^r_{\ttP(\phi)}(V)$
and that
$U(\phi):=\phi^{-1}(S^r_{\ttP(\phi)}(V)^{\circ})\neq\emptyset$.
We obtain
$\Sigma(\phi)_{\red}$
as the fiber product of
$Z(S^r_{\ttP}(V))_{\red}$ and $Y$.
Let $\Sigma(\phi)^{\circ}_{\red}$
denote the fiber product of
$\Sigma(\phi)_{\red}$ and $U(\phi)$.
The projection $\Sigma(\phi)^{\circ}_{\red}\to U(\phi)$
is a proper etale cover of degree $\ell$.

\subsection{Symmetric products of vector bundles}

\subsubsection{Symmetric products and locally closed subsets}

Let $X$ be an irreducible normal algebraic variety over $\hyperk$.
The product bundle $X\times\Spec\hyperk[t]$ is
also denoted by $\nbigo_X$.
Let $\nbige$ be an algebraic vector bundle on $X$.
Let $S^r(\nbige)$
be the quotient of
$\nbige^r=\overbrace{\nbige\times_X\cdots\times_X\nbige}^r$
by the natural action of $\gbigs_r$.
As in \S\ref{subsection;25.10.1.1},
there exists the Chow embedding
\[
 S^r\bigl(\nbige\bigr)
 \lrarr
 \Sym^r\bigl(
 \nbigo_X\oplus
 \nbige
 \bigr).
\]
We regard $S^r(\nbige)$
as a Zariski closed subset of
$\Sym^r\bigl(
\nbigo_X\oplus
\nbige
\bigr)$.

\begin{lem}
\label{lem;25.10.16.1}
Let $Z\subset X$ be a closed subset with $\dim Z\leq \dim X-2$.
Any section $\phi$ of $S^r(\nbige)$ on $X\setminus Z$
uniquely extends to a section of $S^r(\nbige)$ on $X$.
\end{lem}
\pf
Let $\phitilde$ be the induced section of 
$\Sym^r(\nbigo_X\oplus\nbige)$ on $X\setminus Z$.
It uniquely extends to a section of
$\Sym^r(\nbigo_X\oplus\nbige)$ on $X$.
Because $S^r(\nbige)$ is closed in
$\Sym^r(\nbigo_X\oplus\nbige)$,
it factors through $S^r(\nbige)$.
\hfill\qed

\vspace{.1in}
For any partition $\ttP$ of $r$,
we obtain the closed subvariety
$S^r_{\ttP}(\nbige)\subset S^r(\nbige)$
as in \S\ref{subsection;24.8.27.1}
such that the fiber of $S^r_{\ttP}(\nbige)$ over $x\in X$
is $S^r_{\ttP}(\nbige_{|x})$.
Similarly,
we obtain the Zariski open subset
$S^r_{\ttP}(\nbige)^{\circ}\subset S^r_{\ttP}(\nbige)$
such that
the fiber of $S^r_{\ttP}(\nbige)^{\circ}$ over $x\in X$
is $S^r_{\ttP}(\nbige_{|x})^{\circ}$.

By the construction in \S\ref{subsection;25.10.1.10},
we obtain the Zariski closed subset
$Z(S^r_{\ttP}(\nbige))_{\red}
\subset \nbige\times_XS^r_{\ttP}(\nbige)$
such that the fiber over $x\in X$
is $Z(S^r_{\ttP}(\nbige_{|x}))_{\red}$.
The induced morphism 
$Z(S^r_{\ttP}(\nbige))\to S^r_{\ttP}(\nbige)_{\red}$
is finite.
We have the closed subset
$Z(S^r_{\ttP}(\nbige)^{\circ})_{\red}
\subset \nbige\times_XS^r_{\ttP}(\nbige)^{\circ}$
such that the fiber over $x\in X$
is $Z(S^r_{\ttP}(\nbige_{|x})^{\circ})_{\red}$.
The induced morphism 
$Z(S^r_{\ttP}(\nbige)^{\circ})_{\red}
\to S^r_{\ttP}(\nbige)^{\circ}$
is proper and etale of degree $\ell$,
where $\ttP=(m_1,\ldots,m_{\ell})$.

\subsubsection{Finite covers associated with sections}

Let $\phi:X\to S^r(\nbige)$ be any section.
For each closed point $x\in X$,
we obtain the partition $\ttP(\phi(x))$ of $r$
determined by $\phi(x)\in S^r(\nbige_{|x})$.
(See \S\ref{subsection;25.10.1.3}.)
The following lemma is obvious
by the stratification in \S\ref{subsection;24.8.27.1}.

\begin{lem}
There uniquely exist a partition $\ttP(\phi)$ of $r$
and a non-empty Zariski open subset $U(\phi)\subset X$
such that
(i) $\ttP(\phi(x))=\ttP(\phi)$ $(x\in U(\phi))$,
 (ii) $\ttP(\phi)\prec\ttP(\phi(x))$ and
 $\ttP(\phi)\neq \ttP(\phi(x))$ for any $x\in X\setminus U(\phi)$.
\hfill\qed
\end{lem}
We set $\ell(\phi)=\ell(\phi(x))$ $(x\in U(\phi))$.
We obtain
$\phi(U(\phi))\subset S^r_{\ttP(\phi)}(\nbige)^{\circ}$.
Let 
$\Sigma(\phi)_{\red,U(\phi)}\subset \nbige_{|U(\phi)}$
denote the fiber product of
$\phi:U(\phi)\to S^r_{\ttP(\phi)}(\nbige)^{\circ}$
and $Z(S^r_{\ttP}(\nbige)^{\circ})_{\red}
\to S^r_{\ttP(\phi)}(\nbige)^{\circ}$.
The projection
$\Sigma(\phi)_{\red,U(\phi)}\to U(\phi)$
is a proper etale cover of degree $\ell(\phi)$.

Let $\Sigma(\phi)_{\red}$ denote the closure of
$\Sigma(\phi)_{\red,U(\phi)}$ in $\nbige$
with the reduced scheme structure.
Note that 
$\Sigma(\phi)_{\red}$ is finite over $X$
because it is contained in
the fiber product of
$\phi:X\to S^r(\nbige)$
and  $Z(S^r(\nbige))\to S^r(\nbige)$.
The natural isomorphism
$\Sigma(\phi)_{\red,U(\phi)}\to
 U(\phi)
 \times_{S^r_{\ttP(\phi)}(\nbige)^{\circ}}
Z(S_{\ttP(\phi)}(\nbige)^{\circ})_{\red}$
extends
to a morphism of schemes
$\Sigma(\phi)_{\red}
\to
 X\times_{S^r_{\ttP(\phi)}(\nbige)}Z(S_{\ttP(\phi)}^r(\nbige))_{\red}$.
It induces a bijection of the closed points.
In particular, the following holds.
\begin{lem}
For any closed point $x\in X$,
the set of the closed points of
$\Sigma(\phi)_{\red}\times_Xx$
equals
$[\phi(x)]\subset \nbige_{|x}$.
\hfill\qed
\end{lem}

\begin{df}
A finite cover $Z\to X$ is defined to be
a morphism of $\hyperk$-varieties
such that each irreducible component of $Z$
is finite and dominant over $X$.
We also say that $Z$ is a finite cover over $X$.
\hfill\qed
\end{df}

\begin{lem}
$\Sigma(\phi)_{\red}\to X$ is a finite cover.
\end{lem}
\pf
Let
$\Sigma(\phi)_{\red,U(\phi)}=
\bigsqcup_{j=1}^m \Sigma(\phi)_{\red,U(\phi),j}$
denote the decomposition into connected components.
Let $\Sigma(\phi)_{\red,j}$ denote the closure of
$\Sigma(\phi)_{\red,U(\phi),j}$,
which is irreducible,
and finite and dominant over $X$.
Because $\Sigma(\phi)_{\red}=\bigcup_{j=1}^m \Sigma(\phi)_{\red,j}$,
we obtain the claim of the lemma.
\hfill\qed

\subsubsection{Sections associated with finite covers}

Let $\Sigma\subset \nbige$ be a closed subset
such that $\Sigma\to X$ is a finite cover.
We consider the reduced structure of $\Sigma$.
Note that the properness implies the finiteness.
There exists a Zariski open subset $X_0\subset X$
such $\Sigma_0:=\Sigma_{|X_0}\to X_0$ is finite and etale.
Let $d$ denote the degree.

\begin{lem}
There exists a unique section $\phi_{\Sigma}$ of
$S^d(\nbige)$
such that
the set of the closed points of $\Sigma_{|x}$
equals $[\phi_{\Sigma}(x)]$
for any closed point $x\in X$.
\end{lem}
\pf
We set $\Sigma^d:=
\overbrace{\Sigma\times_X\cdots\times_X\Sigma}^d$
and $\Sigma_0^d:=
\overbrace{\Sigma_0\times_X\cdots\times_X\Sigma_0}^d$.
We consider
\[
A(\Sigma_0)
=\bigl\{
 (x_1,\ldots,x_d)\in
 \Sigma_0^d
 \,\big|\,
 x_i\neq x_j
 \bigr\}.
\]
Let $\gbigs_d$ denote the $d$-th symmetric group
acting on $\Sigma^d$ in a natural way.
It induces the $\gbigs_d$-action on $A(\Sigma_0)$,
and we have $A(\Sigma_0)/\gbigs_d\simeq X_0$.
Let $\Abar(\Sigma)$
denote the closure of $A(\Sigma_0)$
in $\Sigma^d$.
It is proper over $X$.
Let $Y(\Sigma)$ denote the normalization of $\Abar(\Sigma)$.
We have the natural $\gbigs_d$-action on $Y(\Sigma)$.
Note that $Y(\Sigma)/\gbigs_d$ is normal
and that $X_0=A(\Sigma_0)/\gbigs_d$ is naturally
an open subset of $Y(\Sigma)$.
There exists the natural finite morphism
$Y(\Sigma)/\gbigs_d\to X$
such that
(i) it induces the identity morphism on $X_0$,
(ii) it induces the bijection of the closed points.
Hence, $Y(\Sigma)/\gbigs_d\simeq X$.

Let $\pi:Y(\Sigma)\to X$ denote the projection.
Let $\varphi_i:Y(\Sigma)\to \Sigma$ be the morphism
induced by the projection
$\Sigma^d\to\Sigma$ onto the $i$-th component.
It induces a section $s_i$ of $\pi^{\ast}\nbige$.
We obtain the section
$\prod_{i=1}^d(1\oplus s_i)$
of $\pi^{\ast}\Sym^d\bigl(\nbigo_X\oplus\nbige\bigr)$.
The section is $\gbigs_d$-invariant.
Hence, we obtain the section
$\phitilde_{\Sigma}$ of $\Sym^d\bigl(\nbigo_X\oplus\nbige\bigr)$
such that
$\pi^{\ast}\phitilde_{\Sigma}=\prod_{i=1}^d(1\oplus s_i)$.
By the construction,
there is a section $\phi_{\Sigma_0}:X_0\to S^d(\nbige)_{|X_0}$
such that $\phitilde_{\Sigma|X_0}$ is induced by $\phi_{\Sigma_0}$.
Because $S^d(\nbige)$ is closed in
$\Sym^d(\nbigo_X\oplus\nbige)$,
there exists the section $\phi_{\Sigma}$
of $S^d(\nbige)$ which induces $\phitilde_{\Sigma}$.
By the construction,
the set of closed points of $\Sigma_x$
is given by $\phi_{\Sigma}(x)$.
The uniqueness is clear.
\hfill\qed

\vspace{.1in}
Let $Z\subset X$ be a closed subset such that
$\dim Z\leq \dim X-2$.
Let $\Sigma_1\subset \nbige_{|X\setminus Z}$
be a closed subset
such that
$\Sigma_1\to X\setminus Z$ is a finite cover.

\begin{lem}
\label{lem;25.10.15.30}
The closure $\Sigmabar_1$ of $\Sigma_1$ in $\nbige$
is a finite cover over $X$.
\end{lem}
\pf
Let $\phi_{\Sigma_1}$ be a section of $S^d(\nbige)_{|X\setminus Z}$
induced by $\Sigma_1$.
By Lemma \ref{lem;25.10.16.1},
it uniquely extends to
a section $\phi'_{\Sigma_1}$ of $S^d(\nbige)$ on $X$.
We obtain the closed subset
$\Sigma(\phi'_{\Sigma_1})_{\red}\subset S^d(\nbige)$.
Because
$\Sigma(\phi')_{\Sigma_1,\red}$
equals $\Sigma_1$ over $X\setminus Z$,
it is the closure of $\Sigma_1$,
and a finite cover over $X$.
\hfill\qed

\subsection{Stratification on the set of sections of symmetric products}

\subsubsection{The algebraic set of
sections to a closed subset in a vector bundle}
\label{subsection;25.10.17.1}

Let $Y$ be a projective normal irreducible $\hyperk$-variety.
Let $\nbige_Y$ be a coherent reflexive $\nbigo_Y$-module.
We regard the finite dimensional $\hyperk$-vector space
$H^0(Y,\nbige_Y)$ as a $\hyperk$-variety.
Let $Z$ be a closed subset of $Y$
with $\codim Z\geq 2$
such that
(i) $Y\setminus Z$ is regular,
(ii) $\nbige_{Y|Y\setminus Z}$ is a locally free sheaf.
We set $Y^{\circ}=Y\setminus Z$
and $\nbige^{\circ}_Y=\nbige_{Y|Y\setminus Z}$.
We have
$H^0(Y,\nbige_Y)=H^0(Y^{\circ},\nbige_Y^{\circ})$.
We regard $\nbige^{\circ}_Y$ as a vector bundle on $Y^{\circ}$.
In particular, it is regarded as a $\hyperk$-variety.
Let $Y^{\circ}(\hyperk)$ denote the set of the $\hyperk$-valued points
of $Y^{\circ}$.

\begin{lem}
\label{lem;25.10.1.2}
For any Zariski closed subset $C\subset \nbige^{\circ}_Y$,
there exists a Zariski closed subset
$H^0(Y^{\circ},C)\subset H^0(Y,\nbige_Y)$
determined by the following condition.
\begin{itemize}
 \item Let $\phi\in H^0(Y,\nbige_Y)$.
       Then, $\phi(Y^{\circ})\subset C$ holds
       if and only if
       $\phi\in H^0(Y^{\circ},C)$.
\end{itemize} 
\end{lem}
\pf
For any $y\in Y^{\circ}(\hyperk)$,
we consider the morphism
$\Phi_y:H^0(Y,\nbige_Y)\to \nbige^{\circ}_{Y|y}$
determined by $\Phi_y(\phi)=\phi(y)$.
We obtain the Zariski closed subset
$\Phi_y^{-1}(C)\subset H^0(Y,\nbige_Y)$.
Then,
$H^0(Y^{\circ},C)=\bigcap_{y\in Y^{\circ}(\hyperk)}\Phi_y^{-1}(C)$
has the desired property.
\hfill\qed

\begin{cor}
Sections $Y^{\circ}\to C$
bijectively correspond to
$\hyperk$-valued points of
the $\hyperk$-variety
$H^0(Y^{\circ},C)$.
\hfill\qed
\end{cor}

Let $Z_1$ be another closed subset of $Y$
with $\codim Z_1\geq 2$
such that
(i) $Y\setminus Z_1$ is regular,
(ii) $\nbige_{Y|Y\setminus Z_1}$ is a locally free sheaf.
We set $Y^{\circ}_1=Y\setminus Z_1$
and $\nbige^{\circ}_1=\nbige_{Y|Y_1^{\circ}}$.
The following lemma is clear.
\begin{lem}
Let $C\subset\nbige^{\circ}$
and $C_1\subset \nbige_1^{\circ}$
be closed subsets
such that $C=C_1$ in $\nbige^{\circ}\cap\nbige_1^{\circ}$.
Then, we have
$H^0(Y^{\circ},C)=H^0(Y_1^{\circ},C_1)$
in $H^0(Y,\nbige)$.
\hfill\qed
\end{lem}

There exists
a locally free $\nbigo_Y$-module
$\nbigg$ with an epimorphism
$\nbigg^{\lor}\to\nbige_Y^{\lor}$,
where $\nbigf^{\lor}=\nhom_{\nbigo_Y}(\nbigf,\nbigo_Y)$.
We obtain a morphism of $\nbigo_Y$-modules
$\nbige_{Y}\to \nbigg$.
Let $\overline{C}$ denote the closure of $C$
in $\nbigg$.
\begin{lem}
We have $H^0(Y^{\circ},C)=H^0(Y,\overline{C})$.
\hfill\qed
\end{lem}

\subsubsection{Set of sections of symmetric products}

Let $X$ be a projective normal $\hyperk$-variety.
Let $\nbige$ be a coherent reflexive $\nbigo_X$-module.
Let $Z$ be a closed subset $X$ with $\codim Z\geq 2$
such that (i) $X^{\circ}=X\setminus Z$ is regular,
(ii) $\nbige^{\circ}=\nbige_{|X^{\circ}}$ is locally free.
From the coherent $\nbigo_X$-module
$\Sym^r(\nbigo_X\oplus \nbige)$,
we obtain the reflexive sheaf
$(\Sym^r(\nbigo_X\oplus \nbige)^{\lor})^{\lor}$
on $X$,
where $\nbigf^{\lor}=\nhom_{\nbigo_X}(\nbigf,\nbigo_X)$.
We obtain the closed subset
$S^r(\nbige^{\circ})\subset
 \Sym^r\bigl(\nbigo_{X^{\circ}}\oplus\nbige^{\circ}
 \bigr)$.
 We obtain the $\hyperk$-variety
\[
H^0(X,S^r(\nbige)):=
H^0(X^{\circ},S^r(\nbige^{\circ}))
\]
as in \S\ref{subsection;25.10.17.1}.
It is independent of the choice of $Z$.

\begin{cor}
\mbox{{}}
\begin{itemize}
 \item Any $\phi\in H^0(X,S^r(\nbige))$
       corresponds to a section
       $X^{\circ}\to S^r(\nbige^{\circ})$,
       and hence
       we have the associated
       closed subset $\Sigma(\phi)_{\red}\subset \nbige^{\circ}$
       such that
       (i) $\Sigma(\phi)_{\red}$
       is finite over $X^{\circ}$,
       (ii) any irreducible component of $\Sigma(\phi)_{\red}$
       is dominant over $X^{\circ}$.
\item  Let $\Sigma\subset\nbige^{\circ}$ be
a closed subset such that
(i) $\Sigma\to X$ is finite,
(ii) any irreducible component of $\Sigma$ is dominant over $X$.
 Let $d\in\seisuu_{>0}$ be the degree of $\Sigma$ over $X$.
Then, there exists a closed point
$\phi\in H^0(X,S^d(\nbige))$
such that $\Sigma=\Sigma(\phi)_{\red}$.
\hfill\qed      
\end{itemize}
\end{cor}

\subsubsection{Decomposition by partitions}
\label{subsection;25.10.2.1}

\begin{lem}
Let $T\subset H^0(X,S^r(\nbige))$
be any non-empty irreducible Zariski closed subset.
Then, there exists a partition $\ttP_{T}$
and a non-empty Zariski open subset $U\subset T$
such that 
$\ttP(\phi)=\ttP_{T}$ for any closed points $\phi\in U$.
\end{lem}
\pf
There exists the universal section
$\phi_{T}:T\times X^{\circ}\to
(T\times X^{\circ})
\times_{X^{\circ}}S^r(\nbige^{\circ})$.
There exist a non-empty Zariski open subset
$U'\subset T\times X^{\circ}$
and a partition $\ttP_{T}$
such that
$\ttP(\phi)=\ttP_{T}$ for any $\phi\in U'$.
There exists a Zariski open subset $U$ of $T$
which is contained in the image of $U'\to T$.
The tuple $(\ttP_{T},U)$ satisfies the desired condition.
\hfill\qed

\vspace{.1in}
By using a Noetherian induction,
we obtain the following corollary.

\begin{cor}
\label{cor;25.10.1.11}
There exists a decomposition
$H^0(X,S^r(\nbige))=\bigsqcup_{i\in\Lambda} T_i$
by locally Zariski closed irreducible smooth affine subsets such that
$\ttP(\phi)$ $(\phi \in T_i)$ are constant.
\hfill\qed
\end{cor}

\begin{rem}
In Corollary {\rm\ref{cor;25.10.1.11}},
we do not impose $\ttP(\phi_i)\neq\ttP(\phi_j)$ 
for $\phi_i\in T_i$ and $\phi_j\in T_j$ even if $T_i\neq T_j$.
\hfill\qed
\end{rem}

\subsubsection{Refined decomposition}

There exists a locally free sheaf $\nbigg$ on $X$
with an epimorphism $\nbigg^{\lor}\to\nbige^{\lor}$.
We obtain the monomorphism $\nbige\to \nbigg$.
Each $\phi\in H^0(X,S^r(\nbige))$
induces
$\phi_{\nbigg}\in H^0(X,S^r(\nbigg))$,
and 
we obtain the closed subset
$\Sigma(\phi_{\nbigg})_{\red}\subset \nbigg$
which is finite over $X$.
It is the closure of
$\Sigma(\phi)_{\red}$,
and we have
$\Sigma(\phi)_{\red}=\Sigma(\phi_{\nbigg})_{\red}\cap
\nbigg_{|X^{\circ}}$.

\begin{lem}
\label{lem;25.10.2.2}
There exist decompositions
$T_i=\bigsqcup_{j\in\Gamma(i)}T_{i,j}$
by locally Zariski closed irreducible smooth affine subsets
such that the following holds for any $i,j$.
\begin{itemize}
 \item There exists a closed subscheme
$\Sigma(T_{i,j})_{\red}\subset T_{i,j}\times\nbigg$
       such that
       $\Sigma(T_{i,j})_{\red}\to T_{i,j}$
       is flat,
       and that
       the fiber over $\phi\in T_{i,j}$
       equals $\Sigma(\phi_{\nbigg})_{\red}\subset \nbigg$. 
 \end{itemize}
\end{lem}
\pf
Let $T_i'\subset T_i$ be any non-empty irreducible Zariski closed subset.
Let $R_i'$ denote the ring of algebraic functions on $T_i'$,
which is finitely generated over $\hyperk$.
Let $\hyperk'$ be an algebraic closure
of the fractional field of $R_i'$.
We set $X_{\hyperk'}=X\times_{\hyperk}\hyperk'$
on which we obtain the induced vector bundle
$\nbigg_{\hyperk'}$.
We have the open subset
$X^{\circ}_{\hyperk'}=X^{\circ}\times_{\hyperk}\hyperk'$
on which we obtain the induced vector bundle
$\nbige_{\hyperk'}^{\circ}$.
The universal section
$T_i'\times X^{\circ}\to T_i'\times S^r(\nbige^{\circ})$
induces a section $\phi_{\hyperk'}$ of
$S^r(\nbige_{\hyperk'}^{\circ})$ on $X^{\circ}_{\hyperk'}$.
We obtain
$\Sigma(\phi_{\hyperk',\nbigg})_{\red}\subset \nbigg_{\hyperk'}$.

Let $R_i'\subset R\subset\hyperk'$ be an extension of rings
such that $R$ is regular and finitely generated over $\hyperk$.
Let $\phi_R:X^{\circ}\times\Spec(R)\to
S^r(\nbige^{\circ})\times\Spec(R)$
denote the induced section.
Then,
$\Sigma(\phi_{\hyperk',\nbigg})_{\red}$
is induced by $\Sigma(\phi_{R,\nbigg})_{\red}$
and the extension $R\to\hyperk'$.
By the generic freeness theorem
\cite[Theorem 6.9.1]{EGA-IV-2},
we may assume that $\Sigma(\phi_{R,\nbigg})_{\red}$
is flat over $\Spec(R)$.
It is finite over $X\times \Spec(R)$.
By \cite[Theorem 12.2.4]{EGA-IV-3},
we may assume that
for any closed point $y\in \Spec(R)$
the fiber $\Sigma(\phi_{R,\nbigg})_{\red}\times_{\Spec(R)} y$ is reduced.

For any closed point $y\in\Spec(R)$,
let $\phi_y:X^{\circ}\to S^r(\nbige^{\circ})$
denote the section induced by $\phi_R$.
For any closed point $x\in X^{\circ}$,
the set of the closed points of
$\Sigma(\phi_{R,\nbigg})_{\red}\times_{X^{\circ}\times\Spec(R)}(x,y)$
equals
$[\phi_{y,\nbigg}(x)]$.
Hence, the closed points of
$\Sigma(\phi_{R,\nbigg})_{\red}\times_{\Spec(R)} y$
equals
the closed points of
$\Sigma(\phi_{y,\nbigg})_{\red}$.
Because both 
$\Sigma(\phi_{R,\nbigg})\times_{\Spec(R)} y$
and
$\Sigma(\phi_{y,\nbigg})_{\red}$
are equipped with the reduced structure,
we obtain
$\Sigma(\phi_{R,\nbigg})_{\red}\times_{\Spec(R)} y
=\Sigma(\phi_{y,\nbigg})_{\red}$.

There exists a non-empty affine open subset $U'$ of $T'_i$
contained in the image of $\Spec(R)\to T'_i$.
By shrinking $\Spec(R)$,
we may assume that $U'$ equals the image.
We have the universal section
$\phi_{U'}:U'\times X\to U'\times S^r(\nbige^{\circ})$.
We obtain
$\Sigma(\phi_{U',\nbigg})_{\red}\subset U'\times\nbigg$.
We may assume that it is flat over $U'$.
Because $\Spec(R)\to U'$ is smooth,
$\Sigma(\phi_{U',\nbigg})_{\red}\times_{U'}\Spec(R)$
is reduced,
and hence
we obtain
$\Sigma(\phi_{U',\nbigg})_{\red}\times_{U'}\Spec(R)
=\Sigma(\phi_{R,\nbigg})_{\red}$.
For any $\phi\in U'$,
there exists $y\in \Spec(R)$ which is mapped to $\phi$.
We have $\phi_y=\phi$.
In $\nbigg$,
we have
$\Sigma(\phi_{U',\nbigg})_{\red}\times_{U'}z
=\Sigma(\phi_{R,\nbigg})_{\red}\times_{\Spec(R)}y
=\Sigma(\phi_{\nbigg})_{\red}$.

By using this consideration with the Noetherian induction,
we can construct desired decompositions.
\hfill\qed

\vspace{.1in}
There exist a Zariski dense open subset
$\Sigma(T_{i,j})^{\sm}_{\red}\subset
\Sigma(T_{i,j})_{\red}\cap(\nbige^{\circ}\times T_{i,j})$
such that
(i) $\Sigma(T_{i,j})_{\red}^{\sm}\to T_{i,j}$
is smooth,
(ii)  for each $y\in T_{i,j}$,
$\Sigma(T_{i,j})^{\sm}_{\red}\times_{T_{i,j}}y$
intersects with each irreducible component of
$\Sigma(\phi_{y,\nbigg})_{\red}$,
where $\phi_y$ denotes the section corresponding to $y$.
	Let $\iota_{0}:\nbige^{\circ}\times T_{i,j}
	\to\nbigg\times T_{i,j}$,
	$\iota_1:\Sigma(T_{i,j})_{\red}\to \nbigg\times T_{i,j}$
	and
	$\iota_2:\Sigma(T_{i,j})_{\red}^{\sm}
	\to \nbigg\times T_{i,j}$
	denote the inclusions.

\begin{lem}
\label{lem;25.10.17.10}
 There exist torsion-free sheaves
       $\nbigf^k$ $(k\in\seisuu_{\geq 0})$ on
       $\Sigma(T_{i,j})_{\red}$
       such that the following holds.
 \begin{itemize}
  \item
       The restriction of 
       $\nbigf^k$
       to $\Sigma(T_{i,j})_{\red}^{\sm}$
       equal the sheaves of algebraic relative $k$-forms
       $\Omega^k_{\Sigma(T_{i,j})^{\sm}_{\red}/T_{i,j}}$.
       In particular,
       there exist the natural morphisms
       $\iota_{1\ast}(\nbigf^k)
       \to
       \iota_{2\ast}\Omega^k_{\Sigma(T_{i,j})_{\red}^{\sm}/T_{i,j}}$.
  \item 
	There exists
	a morphism
	$\iota_{0\ast}\bigl(
	\Omega^k_{(\nbige^{\circ}\times T_{i,j})/T_{i,j}}
	\bigr)
	\to
	\iota_{1\ast}(\nbigf^k)$
	such that the induced morphism
	$\iota_{0\ast}\bigl(
	\Omega^k_{(\nbige^{\circ}\times T_{i,j})/T_{i,j}}
	\bigr)
	\to
	\iota_{2\ast}\Omega^k_{\Sigma(T_{i,j})_{\red}^{\sm}/T_{i,j}}$
	equals
	the morphism induced by the pull back of relative $k$-forms.	
 \end{itemize}
By refining the decomposition,
we may also assume that
$\nbigf^k$ are flat over $T_{i,j}$.
\end{lem}
\pf
We set
$\Sigma(T_{i,j})_{\red}^{\circ}
=
\Sigma(T_{i,j})_{\red}
\cap(\nbige^{\circ}\times T_{i,j})$.
Let
$j_1:\Sigma(T_{i,j})_{\red}^{\sm}
\to
\Sigma(T_{i,j})_{\red}^{\circ}$
and
$j_2:\Sigma(T_{i,j})_{\red}^{\circ}\to
\nbige^{\circ}\times T_{i,j}$
denote the inclusions.
We obtain the morphism 
of quasi-coherent sheaves
$j_2^{\ast}
 \Omega^k_{(\nbige^{\circ}\times T_{i,j})/T_{i,j}}
 \to
 j_{1\ast}
 \Omega^k_{\Sigma(T_{i,j})_{\red}^{\sm}}$
on
$\Sigma(T_{i,j})_{\red}^{\circ}$.
There exists a coherent sheaf $\nbigf^k_1$
on $\Sigma(T_{i,j})_{\red}^{\circ}$
whose restriction to 
$\Sigma(T_{i,j})_{\red}^{\sm}$
equals
$\Omega^k_{\Sigma(T_{i,j})_{\red}^{\sm}/T_{i,j}}$
(see \cite[Problem II.5.15]{Hartshorne}, for example).
By taking the sum of
$\nbigf^k_1$ and
the images of 
$j_2^{\ast}
\Omega^k_{(\nbige^{\circ}\times T_{i,j})/T_{i,j}}$
in
$j_{1\ast}
 \Omega^k_{\Sigma(T_{i,j})_{\red}^{\sm}/T_{i,j}}$,
we obtain a coherent sheaf $\nbigf^k_2$ on
$\Sigma(T_{i,j})_{\red}^{\circ}$
with a morphism
$j_2^{\ast}
 \Omega^k_{\nbige^{\circ}\times T_{i,j}/T_{i,j}}
 \to \nbigf_2^k$.
We extend $\nbigf_2^k$ a coherent sheaf $\nbigf^k$
on $\Sigma(T_{i,j})_{\red}$.
\hfill\qed

\subsubsection{Meromorphic $k$-forms on $\nbige$}

Let $H$ be a hypersurface of $X$.
Let $\tau$ be a meromorphic $k$-form on $\nbige^{\circ}$
whose restriction to
$\nbige^{\circ}_{|X^{\circ}\setminus H}$ is a regular $k$-form.
Let $T_{i,j}(\tau)$ denote the set of
the closed points $\phi$ of $T_{i,j}$
such that
the restriction of $\tau$
to $\Sigma(\phi)_{U(\phi),\red}$ is $0$.

\begin{lem}
\label{lem;25.10.2.3}
$T_{i,j}(\tau)$ are Zariski closed subsets of
$T_{i,j}$.
\end{lem}
\pf
We may assume that $H$ is ample
and that $X\setminus X^{\circ}\subset H$.
Let $\pi:\Sigma(T_{i,j})\to X^{\circ}$
and $p:\Sigma(T_{i,j})\to T_{i,j}$
denote the projections.
Let $\nbigf^k$ be a torsion-free sheaf
on $\Sigma(T_{i,j})$
as in Lemma \ref{lem;25.10.17.10}.
Then,
for a sufficiently large integer $N$,
$\tau$ induces a section
$\Phi_{\tau}$
of the locally free sheaf
$p_{\ast}
\Bigl(
\nbigf^k
\otimes
\pi^{\ast}\nbigo_X(NH)
\Bigr)$.
Because $T_{i,j}(\tau)$
equals the $0$-set of $\Phi_{\tau}$,
we obtain the claim of the lemma.
\hfill\qed

\subsection{Conditions for meromorphic finite covers}
\label{subsection;25.12.10.2}

\subsubsection{General case}

Let $X$ be a normal algebraic variety over $\hyperk$.
Let $H$ be a hypersurface of $X$ such that
$X\setminus H$ is smooth.
\begin{df}
A meromorphic finite cover over $(X,H)$
is a closed subset $\Sigma\subset T^{\ast}(X\setminus H)$
such that
$\Sigma\to X\setminus H$ is a finite cover.
\hfill\qed
\end{df}

For later use,
we also introduce the following condition.
\begin{df}
\label{df;25.12.7.15}
A meromorphic finite cover $\Sigma$ over $(X,H)$
is called Lagrangian 
if the smooth part of $\Sigma$ is Lagrangian
with respect to the natural symplectic structure of
$T^{\ast}(X\setminus H)$.
\end{df}

\subsubsection{Formal conditions}
\label{subsection;25.12.7.30}

We assume that $X$ is smooth, and that $H$ is normal crossing.
Let $x\in H$.
Let $\gminim_x$ denote the maximal ideal of
the local ring $\nbigo_{X,x}$.
Let $\nbigohat_{X,x}$ denote the completion of $\nbigo_{X,x}$
with respect to $\gminim_x$.
For any $\nbigohat_{X,x}$-module $\nbigm$,
we set
$\nbigm(\ast H):=
\nbigm\otimes_{\nbigo_{X,x}}\nbigo_{X}(\ast H)_x$.
For any affine morphism $Z\to X$,
let $\nbigo_{X,x}(Z)$ denote the ring of
algebraic functions on $Z\times_{X}\Spec(\nbigo_{X,x})$.
We set
$\nbigohat_{X,x}(Z):=\nbigo_{X,x}(Z)\otimes_{\nbigo_{X,x}}\nbigohat_{X,x}$.
In particular,
we obtain the $\nbigohat_{X,x}(\ast H)$-algebra
$\nbigohat_{X,x}(T^{\ast}(X\setminus H))$.
We have
$\nbigohat_{X,x}(T^{\ast}(X\setminus H))
=\nbigohat_{X,x}\bigl(T^{\ast}X(\log H)(NH)\bigr)(\ast H)$
for any $N\in\seisuu$.

\begin{df}
A surjective morphism of 
$\nbigohat_{X,x}(\ast H)$-algebras
$\nbigohat_{X,x}(T^{\ast}(X\setminus H))\to \nbigb$
is called a formal meromorphic cover at $x$
if there exist
$N\in\seisuu_{>0}$
and 
a surjection of
$\nbigohat_{X,x}$-algebras
$\nbigohat_{X,x}(T^{\ast}X(\log H)(NH))\to
 \nbigb_{N,H}$
such that 
(i) $\nbigb_{N,H}$ is finite over $\nbigohat_{X,x}$,
(ii) $\nbigb=\nbigb_{N,H}(\ast H)$,
 (iii) $\nbigb_{N,H}$ is torsion-free as an $\nbigohat_{X,x}$-module.
\hfill\qed
\end{df}

\begin{df}
A formal meromorphic cover
$\nbigohat_{X,x}(T^{\ast}(X\setminus H))\to \nbigb$
is called logarithmic
if there exists a surjection
of $\nbigohat_{X,x}$-algebras
$\nbigohat_{X,x}(T^{\ast}X(\log H))\to \nbigb_0$
such that (i) $\nbigb_0$ is finite over $\nbigohat_{X,x}$,
(ii) $\nbigb_0(\ast H)=\nbigb$
as $\nbigohat_{X,x}(T^{\ast}X(\log H))$-algebras.
\hfill\qed
\end{df}

Any $\gminia\in\nbigohat_{X,x}(\ast H)$
induces
$d\gminia\in
(T^{\ast}X(\log H)(NH))_x\otimes_{\nbigo_{X,x}}\nbigohat_{X,x}$
for some $N\in\seisuu_{>0}$.
The addition of $d\gminia$
induces an automorphism of
$T^{\ast}X(\log H)(NH)
\times_{X}\Spec(\nbigohat_{X,x})$.
It induces an automorphism of the $\nbigohat_{X,x}$-algebra
$\nbigohat_{X,x}(T^{\ast}X(\log H)(NH))$,
and an automorphism 
of the $\nbigohat_{X,x}(\ast H)$-algebra
$\nbigohat_{X,x}(T^{\ast}(X\setminus H))$.
The automorphisms of the algebras are denoted by
$\kappa^{\ast}_{d\gminia}$.

Let
$\nbigohat_{X,x}(T^{\ast}(X\setminus H))\to \nbigb$
be a formal meromorphic cover at $x$.
We set
$\kappa_{d\gminia}^{\ast}\nbigb:=\nbigb$
as an $\nbigohat_{X,x}(\ast H)$-algebra,
and we consider the morphism
$\nbigohat_{X,x}(T^{\ast}(X\setminus H))
\to \kappa_{d\gminia}^{\ast}\nbigb$
obtained as the composition of
$\nbigohat_{X,x}(T^{\ast}(X\setminus H))
 \stackrel{\kappa^{\ast}_{d\gminia}}{\lrarr}
  \nbigohat_{X,x}(T^{\ast}(X\setminus H))
  \lrarr
  \nbigb$.
  
\begin{df}
A formal meromorphic cover
$\nbigohat_{X,x}(T^{\ast}(X\setminus H))\to \nbigb$
is unramifiedly good
if there exist
a good set of irregular values 
$\nbigi(\nbigb)\subset \nbigohat_{X,x}(\ast H)/\nbigohat_{X,x}$
and
the decomposition of
$\nbigohat_{X,x}(\ast H)$-algebras
$\nbigb
=\prod_{\gminia\in\nbigi(\nbigb)}
 \nbigb_{\gminia}$
such that the following holds for any $\gminia\in\nbigi(\nbigb)$.
\begin{itemize}
 \item Let $\gminiatilde\in\nbigohat_{X,x}(\ast H)$
       be a lift of $\gminia$.
       Then, $\kappa_{d\gminiatilde}^{\ast}\nbigb_{\gminia}$
       is logarithmic. 
\end{itemize}
See {\rm\S\ref{subsection;26.2.27.1}},
in particular Remark {\rm\ref{rem;26.2.27.2}} 
for the notion of good set of irregular values.
\hfill\qed
\end{df}

Let $n=\dim X$.
Let $(z_1,\ldots,z_n)$ be a generator of
the maximal ideal $\gminim_x$ of $\nbigo_{X,x}$
such that the ideal sheaf of $H$ at $x$
is generated by $\prod_{i=1}^{\ell}z_i$
for some $1\leq \ell\leq n$.
For any $e\in\seisuu_{>0}$,
we obtain the flat extension of rings
\[
 \nbigohat^{(e,H)}_{X,x}
 =\nbigohat_{X,x}[z_j^{1/e}\,|\,j=1,\ldots,\ell].
\]
Let $\nbigohat_{X,x}(T^{\ast}(X\setminus H))
\to \nbigb$
be a formal meromorphic cover.
We obtain
the $\nbigohat^{(e,H)}_{X,x}$-algebras
\[
\nbigohat_{X,x}^{(e,H)}(T^{\ast}(X\setminus H))=
\nbigohat_{X,x}(T^{\ast}(X\setminus H))\otimes_{\nbigohat_{X,x}}
 \nbigohat^{(e,H)}_{X,x},
\quad\quad
\nbigb^{(e,H)}=
\nbigb\otimes_{\nbigohat_{X,x}}
 \nbigohat^{(e,H)}_{X,x},
\]
and the natural surjective morphism
$\nbigohat_{X,x}^{(e,H)}(T^{\ast}(X\setminus H))
\to
\nbigb^{(e,H)}$.
Note that
there exists a morphism of $\hyperk$-varieties
$f:X'\to X$ with a point $x'\in X'$ satisfying $f(x')=x$
such that the following holds.
\begin{itemize}
 \item  $f$ induces
	an isomorphism
	$\nbigohat^{(e,H)}_{X,x}\simeq \nbigohat_{X',x'}$.
 \item 	$T^{\ast}X'(\log H')\simeq f^{\ast}T^{\ast}X(\log H)$
	and
	$T^{\ast}X'(\log H')(eNH')\simeq f^{\ast}T^{\ast}X(\log H)(NH)$
	around $x'$.
\end{itemize}
Hence,
$\nbigohat_{X,x}^{(e,H)}(T^{\ast}(X\setminus H))\to
 \nbigb^{(e,H)}$
is regarded as a formal meromorphic cover.

\begin{df}
\label{df;25.12.6.2}
The formal meromorphic cover
$\nbigohat_{X,x}(T^{\ast}(X\setminus H))\to\nbigb$ is called good
if there exists $e\in\seisuu_{>0}$
such that
 $\nbigohat_{X,x}^{(e,H)}(T^{\ast}(X\setminus H))
 \to\nbigb^{(e,H)}$
is unramifiedly good.
 \hfill\qed 
\end{df}
In Definition \ref{df;25.12.6.2},
the index set 
$\nbigi(\nbigb^{(e,H)})\subset
 \nbigohat^{(e,H)}_{X}(\ast H)_x\big/
 \nbigohat^{(e,H)}_{X,x}$
is also denoted as $\nbigi(\nbigb)$.
It is invariant under the Galois action.
We set
\[
 \nbigohat^{(H)}_{X,x}=
 \varinjlim
 \nbigohat^{(e,H)}_{X,x}.
\]
If the formal meromorphic cover
$\nbigohat_{X,x}(T^{\ast}(X\setminus H))\to\nbigb$ is good,
$\nbigi(\nbigb)$ is naturally regarded as a subset of
$\nbigohat^{(H)}_{X,x}(\ast H)/\nbigohat^{(H)}_{X,x}$.

\subsubsection{Algebraic good set of irregular values}

Let $U=\Spec R$ be an $n$-dimensional
smooth connected affine $\hyperk$-variety
with a normal crossing hypersurface $H_U$.
The set of closed points are denoted by $U(\hyperk)$
and $H_U(\hyperk)$, respectively.
Let $x\in H_U(\hyperk)$.
Assume that there exists
an etale morphism
$\psi:U\to \hyperk^n$
such that the following holds.
\begin{itemize}
 \item $H_U=\psi^{-1}\bigl(
 \bigcup_{i=1}^{\ell}\{z_i=0\}
\bigr)$ 
       for some $1\leq \ell\leq n$.
 \item $\psi(x)=(0,\ldots,0)$.
 \item $\psi^{-1}\{z_i=0\}$ are connected.
\end{itemize}
The functions $\psi^{\ast}(z_j)$ are also denoted by $z_j$.
We set $R(\ast H_U):=H^0(U,\nbigo_U(\ast H_U))$.
\begin{df}
Such a tuple of functions
$(z_1,\ldots,z_n)$ is called an etale coordinate system 
on $U$.
\hfill\qed
\end{df}

For any closed point $y\in H_U$,
let $\Rhat_y$ denote the completion of $R$ at $y$.
Let $\Rhat_y(\ast H_U)$ denote the formal completion of
$R(\ast H_U)$ at $y$.
We have the naturally induced map
$R(\ast H_U)/R\to \Rhat_y/\Rhat_y(\ast H_U)$.

\begin{df}
A finite subset $\nbigi$ of $R(\ast H_U)/R$
is called good if
the induced tuple 
$\nbigi_{|\yhat}\subset\Rhat_y(\ast H_U)/\Rhat_y$
is good for any $y\in H_U$.
\hfill\qed
 \end{df}

\begin{lem}
\label{lem;25.10.6.20}
Let $\nbigi\subset R(\ast H_U)/R$
be a finite subset.
Suppose that
$\nbigi_{|\xhat}\subset \Rhat_x(\ast H_U)/\Rhat_x$ is good.
Let $Z$  be the set of the closed points
$y\in H_U$
such that $\nbigi$ is not good in 
$\Rhat_y(\ast H_U)/\Rhat_y$.
Then,
$Z$ is a Zariski closed subset
of $H_U$,
and  $H_U\setminus Z$ contains $x$.
In particular, $H_U\setminus Z$ is a Zariski dense
open subset of $H_U$. 
\end{lem}
\pf
Let $\nbigitilde\subset R(\ast H_U)$
be a lift of $\nbigi$.
Because $\nbigi$ induces
a good set of irregular values at $x$,
we have
$\ord(\gminia)\in\seisuu_{\leq 0}^{\ell}$
and $\ord(\gminia-\gminib)\in\seisuu_{\leq 0}^{\ell}$.
We have 
$z^{-\ord(\gminia)}\gminia\in \Rhat_x$
and
$z^{-\ord(\gminia-\gminib)}(\gminia-\gminib)\in \Rhat_x$.
Let $R_x$ denote the local ring at $x$.
Because $\Rhat_x$ is faithfully flat over $R_x$,
$R_x(\ast H_U)/R_x\to \Rhat_x(\ast H_U)/\Rhat_x$
is injective.
We obtain that
$z^{-\ord(\gminia)}\gminia\in R_x$
and 
$z^{-\ord(\gminia-\gminib)}(\gminia-\gminib)\in R_x$.
We obtain that
$z^{-\ord(\gminia)}\gminia$
and $z^{-\ord(\gminia-\gminib)}(\gminia-\gminib)$
are contained in the local rings
of the generic points of $H$.
Hence, we obtain that
$z^{-\ord(\gminia)}\gminia$
and $z^{-\ord(\gminia-\gminib)}(\gminia-\gminib)$
are contained in $R$.

Let $\ord(\gminia)_i$ $(1\leq i\leq \ell)$ denote the $i$-th components of
$\ord(\gminia)\in\seisuu_{\leq 0}^{\ell}$.
We obtain the following closed subset:
\[
 Z(z^{-\ord(\gminia)}\gminia)
 :=
 \{z^{-\ord(\gminia)}\gminia=0\}
 \cap 
 \bigcup_{\ord(\gminia)_i<0}
 \{z_i=0\}.
\]
Similarly,
let $\ord(\gminia-\gminib)_i$ denote the $i$-th component of
$\ord(\gminia-\gminib)\in\seisuu_{\leq 0}^{\ell}$.
We set
\[
 Z(z^{-\ord(\gminia-\gminib)}(\gminia-\gminib))
  :=\{z^{-\ord(\gminia-\gminib)}(\gminia-\gminib)=0\}
  \cap
  \bigcup_{\ord(\gminia-\gminib)_i<0}
  \{z_i=0\}.
\]
We obtain
\[
 Z(\nbigi)
 :=\bigcup_{\gminia\in\nbigi\setminus\{0\}}
  Z(z^{-\ord(\gminia)}\gminia)
\cup
 \bigcup_{\gminia\neq \gminib}
 Z(z^{-\ord(\gminia-\gminib)}(\gminia-\gminib)).
\]
Then,
$\nbigi$ induces a good set of irregular values at $y$
if and only if 
$y\in H_U\setminus Z(\nbigi)\neq\emptyset$.
\hfill\qed

\subsubsection{Algebraic local conditions}
\label{subsection;25.12.7.31}

Let 
$\Sigma_U\subset T^{\ast}(U\setminus H_U)$
be a meromorphic finite cover over $(U,H_U)$.

\begin{df}
We say that $\Sigma_U$ is logarithmic on $(U,H_U)$
if the closure of $\Sigma_U$
in $T^{\ast}U(\log H_U)$  is proper over $U$.

\hfill\qed
\end{df}

For any $\gminia\in\nbigi$,
let $\gminiatilde\in R(\ast H_U)$ be a lift.
We obtain the section
$d\gminiatilde$
of $\Omega^1_{U}(\log H_U)(\ast H_U)$.
We obtain the automorphism
$\kappa_{d\gminiatilde}$ of
$T^{\ast}(U\setminus H_U)$
by the addition.
For any closed subset $Z\subset T^{\ast}(U\setminus H_U)$,
we obtain
$\kappa_{d\gminiatilde}^{\ast}Z
\subset T^{\ast}(U\setminus H_U)$.

\begin{df}
$\Sigma_U$ is unramifiedly good on $(U,H_U)$ 
if there exits
a good finite set $\nbigi(\Sigma_U)\subset R(\ast H_U)/R$
and a decomposition of $\hyperk$-varieties
\[
 \Sigma_U
 =\bigsqcup_{\gminia\in\nbigi(\Sigma_U)}
 \Sigma_{U,\gminia}
\] 
such that
$\kappa_{d\gminiatilde}^{\ast}
\Sigma_{U,\gminia}$
are logarithmic.
\hfill\qed
\end{df}

For a positive integer $e$,
we set $R^{(e)}=R[z_i^{1/e}\,|\,i=1,\ldots,\ell]$.
We obtain
$U^{(e)}=\Spec R^{(e)}$
with the finite morphism
$\varphi_{e,H_U}:U^{(e)}\to U$.

\begin{df}
We say that 
$\Sigma_U$ is good on $(U,H_U)$
if $\varphi_{e,H_U}^{\ast}\Sigma_U$
is unramifiedly good for some positive integer $e$.
In case, 
we obtain
$\nbigi(\Sigma_U):=
\nbigi(\varphi_{e,H_U}^{\ast}(\Sigma_U))$.
\hfill\qed
\end{df}

\subsubsection{Comparison of the algebraic local condition
and the formal condition}
\label{subsection;25.12.7.1}

Let $X$ be an $n$-dimensional smooth algebraic variety over $\hyperk$.
Let $H$ be a normal crossing hypersurface of $X$.

An etale coordinate neighbourhood
of $x\in X$
is a tuple
$(U,x',\varphi,\psi)$
of a smooth affine $\hyperk$-variety $U$,
a closed point $x'\in U$,
an etale morphism $\varphi:(U,x')\to (X,x)$,
and an etale morphism
$\psi:U\to \hyperk^n$
satisfying the following conditions.
\begin{itemize}
 \item $H_U:=\varphi^{-1}(H)$ equals
       $\psi^{-1}\bigl(
       \bigcup_{i=1}^{\ell}
       \{z_i=0\}
       \bigr)$
       for some $1\leq \ell\leq n$.
 \item $\psi(x')=(0,\ldots,0)$.
 \item $\psi^{-1}(\{z_i=0\})$
       are connected.
\end{itemize}

\begin{prop}
\label{prop;25.10.6.1}
A meromorphic finite cover
$\Sigma\subset T^{\ast}(X\setminus H)$
is formally good at $x\in X$
if and only if
there exists an etale neighborhood
$(U,x',\varphi,\psi)$ of $x$
such that 
 $\varphi^{\ast}\Sigma\subset T^{\ast}(U\setminus H_U)$
 is good on $(U,H_U)$.
\end{prop}
\pf
The ``only if'' part is clear.
Suppose that $\Sigma$ is formally good at $x$.
We may assume that $X$ is affine
and that there exists an etale morphism
$\psi:X\to \hyperk^n$
such that
(i) $H=\psi^{-1}\bigl(
 \bigcup_{j=1}^{\ell}\{z_j=0\}
\bigr)$,
(ii) $\psi(x)=(0,\ldots,0)$,
(iii) $\psi^{-1}(\{z_j=0\})$ are connected.
We set 
$R^{(e)}=R[z_i^{1/e}\,|\,i=1,\ldots,\ell]$
and $X^{(e)}=\Spec R^{(e)}$.
We have the ramified cover
$\rho_e:X^{(e)}\to X$.
We set
$H^{(e)}=\rho_e^{-1}(H)$
and $\Sigma^{(e)}=\rho_e^{\ast}\Sigma
\subset T^{\ast}(X^{(e)}\setminus H^{(e)})$.
We have $x^{(e)}\in X^{(e)}$
such that $\rho_e(x^{(e)})=x$.
Then,
$\Sigma^{(e)}$ is formally unramifiedly good at $x^{(e)}$.

For any $N\in\seisuu_{>0}$,
we set
$\nbige^{(e)}_{N}:=T^{\ast}X^{(e)}(\log H^{(e)})(NH^{(e)})$.
There exists $N>0$ such that
the closure of
$\Sigmabar^{(e)}_{N}\subset\nbige^{(e)}_{N}$ of $\Sigma^{(e)}$
is proper over $X^{(e)}$.
Let $d$ be the number of the general fiber of
$\Sigma^{(e)}\to X^{(e)}\setminus H^{(e)}$.
We set
\[
 A(\Sigma^{(e)}):=
 \bigl\{
 (x_1,\ldots,x_d)\in
 \overbrace{\Sigma^{(e)}
 \times_{X^{(e)}\setminus H^{(e)}}
 \cdots\times_{X^{(e)}\setminus H^{(e)}}
 \Sigma^{(e)}}^d
 \,\big|\,
 x_i\neq x_j\,\,(i\neq j)
 \bigr\}.
\]
Let $\Abar(\Sigmabar^{(e)}_{N})$ denote 
the closure of $A(\Sigma^{(e)})$
in
$\overbrace{\nbige^{(e)}_{N}\times_{X^{(e)}}
\cdots\times_{X^{(e)}}\nbige^{(e)}_{N}}^d$.
Let $Y(\Sigmabar^{(e)}_{N})$
denote the normalization of $\Abar(\Sigmabar^{(e)}_{N})$.
The projection
$Y(\Sigmabar^{(e)}_{N})\to X^{(e)}$ is an affine morphism.
Let $\gbigs_d$ denote the $d$-th symmetric group.
There exists the natural $\gbigs_d$-action on
$Y(\Sigmabar^{(e)}_{N})$.
We have
$Y(\Sigmabar^{(e)}_{N})/\gbigs_d\simeq X^{(e)}$.
Let $\pi:Y(\Sigmabar^{(e)}_{N})\to X^{(e)}$
denote the projection,
which is a finite morphism.

Let $\nbigo_{X,x}^h$ be a Henselization of $\nbigo_{X,x}$.
By \cite[Corollary 4.3]{Milne},
$\nbigo^h_{X^{(e)},x^{(e)}}:=
\nbigo^h_{X,x}[z_i^{1/e}\,|\,i=1,\ldots,\ell]$
is a Henselian ring,
which is a Henselization of $\nbigo_{X^{(e)},x^{(e)}}$.
By \cite[Theorem I.4.2]{Milne},
we have the decomposition into local rings:
\[
 \nbigo_{X^{(e)},x^{(e)}}(Y(\Sigmabar^{(e)}_{N}))
 \otimes
 \nbigo^h_{X^{(e)},x^{(e)}}
 =\bigoplus_{y\in\pi^{-1}(x^{(e)})}
 \Bigl(
  \nbigo_{X^{(e)},x^{(e)}}(Y(\Sigmabar^{(e)}_{N}))
 \otimes
 \nbigo^h_{X^{(e)},x^{(e)}}
 \Bigr)_y.
\]
By replacing $X$ with an etale coordinate neighbourhood of $x$,
we may assume that
there exists the decomposition of $\hyperk$-varieties
\[
\Ybar(\Sigmabar^{(e)}_{N})
=\bigsqcup_{y\in \pi^{-1}(x^{(e)})}
\Ybar(\Sigmabar^{(e)}_{N})_y
\]
where
$\Ybar(\Sigmabar_{N})_y\times_{X^{(e)}} x^{(e)}=\{y\}$.
Let $G_y\subset \gbigs_d$ denote the stabilizer.
We have the natural action $G_y$ on $Y(\Sigmabar^{(e)}_{N})_y$,
and $Y(\Sigmabar^{(e)}_{N})_y/G_y\simeq X^{(e)}$.

Let $p_j:Y(\Sigmabar^{(e)}_{N})\to
\Sigmabar^{(e)}_{N}\subset \nbige^{(e)}_{N}$
denote the morphism induced by the projection onto
the $j$-th component.
It induces a section
$s_j$ of $\pi^{\ast}\nbige^{(e)}_{N}$
on $Y(\Sigmabar^{(e)}_{N})$.

Let $y\in \pi^{-1}(x^{(e)})$.
Let $\nbigz(s_i-s_j)_y$ $(i\neq j)$
denote the $0$-set of $s_i-s_j$
in $Y(\Sigmabar^{(e)}_{N})_{y|X^{(e)}\setminus H^{(e)}}$.
Let $\overline{\nbigz(s_i-s_j)_y}$
denote the closure of $\nbigz(s_i-s_j)_y$
in $Y(\Sigmabar^{(e)}_{N})_{y}$.
We consider the equivalence relation $\sim_y$ on
$\{1,\ldots,d\}$
defined by the following condition.
\begin{itemize}
 \item We have $i\sim_y j$
       if there exists a sequence
       $i=p(0),p(1),\ldots,p(m)=j$
       such that
       $y\in \overline{\nbigz(s_{p(k)}-s_{p(k-1)})}$
       for $k=1,\ldots,m$.
\end{itemize}
Let $S_y$ denote the quotient set.
We obtain the decomposition
$\{1,\ldots,d\}=\bigsqcup_{I\in S_y} I$.

Let $\Psi_{y,i}:Y(\Sigmabar^{(e)}_{N})_y\to \Sigmabar^{(e)}_{N}$
be the finite morphism induced by $s_i$.
Because $\dim\Image\Psi_{y,i}=\dim X$,
\[
 (\Sigmabar^{(e)}_{N})_{y,I}:=
\bigcup_{i\in I}\Image\Psi_{y,i}
\]
is a union of some irreducible components of
$\Sigmabar^{(e)}_{N}$.
We obtain the decomposition
$\Sigmabar^{(e)}_{N}
=\bigcup_{I\in S_y}
 (\Sigmabar^{(e)}_{N})_{y,I}$.
We set
$\Sigma^{(e)}_{y,I}
:=(\Sigmabar^{(e)}_{N})_{y,I|X^{(e)}\setminus H^{(e)}}$.
By shrinking $X$,
we may assume
$\Sigma^{(e)}_{y,I_1}
 \cap
 \Sigma^{(e)}_{y,I_2}
=\emptyset$.
We obtain the decomposition of $\hyperk$-varieties:
\[
 \Sigma^{(e)}=\bigsqcup_{I\in S_y} \Sigma^{(e)}_{y,I}.
\]
Let $A(\Sigma^{(e)})_y\subset A(\Sigma^{(e)})$
denote the closed subset of
$(x_1,\ldots,x_d)\in A(\Sigma^{(e)})$ satisfying the following condition.
\begin{itemize}
 \item For any $I\in S_y$ and $i\in I$,
       we have $x_i\in \Sigma^{(e)}_{y,I}$.
\end{itemize}
We obtain the decomposition
of $\hyperk$-varieties:
\[
 A(\Sigma^{(e)})
 =\bigsqcup_{y\in \pi^{-1}(x^{(e)})}
 A(\Sigma^{(e)})_y.
\]
Let $\Ybar(\Sigmabar^{(e)}_N)'_y$
denote the normalization of
the closure of 
$A(\Sigma^{(e)})_y$
in $\Abar(\Sigmabar^{(e)}_N)$.
We obtain the decomposition
\[
\Ybar(\Sigmabar^{(e)}_N)
=\bigsqcup_{y\in \pi^{-1}(x^{(e)})}
\Ybar(\Sigmabar^{(e)}_N)'_y.
\]
Because $y\in \Ybar(\Sigmabar^{(e)}_N)$,
we obtain
$\Ybar(\Sigmabar^{(e)}_N)'_y
=
\Ybar(\Sigmabar^{(e)}_N)_y$.
It implies the following lemma.
\begin{lem}
For $g\in \gbigs_d$,
we have
$g\Ybar(\Sigmabar^{(e)}_N)_y
=\Ybar(\Sigmabar^{(e)}_N)_{g(y)}$.
As a result,
we obtain 
$g(I)=I$ for any $g\in G_y$
and $I\in S_y$.
\hfill\qed
\end{lem}

Because
the formal cover
$\nbigohat_{X^{(e)},x^{(e)}}(T^{\ast}(X^{(e)}\setminus H^{(e)}))
\to
\nbigohat_{X^{(e)},x^{(e)}}(\Sigma^{(e)})$
of $\nbigohat_{X^{(e)},x^{(e)}}$
is unramifiedly good,
there exists
$\gminia_j\in\nbigohat_{X^{(e)},x^{(e)}}(\ast H^{(e)})$
such that
$s_j-\pi^{\ast}d\gminia_j
\in
\pi^{\ast}(\nbige_{H^{(e)},0})\otimes\nbigohat_{X^{(e)},x^{(e)}}$.
Because $\{\gminia_j\}$ is good,
if $i,j\in I$ for some $I$,
then we may assume $\gminia_i=\gminia_j$.
Hence,
by setting
\[
 \ttQ_1(s_i)
=\frac{1}{|G_y|}
 \sum_{g\in G_y} g^{\ast}(s_i),
\]
we obtain that
$\ttQ_1(s_i)-\pi^{\ast}d\gminia_i
 \in
 \pi^{\ast}(\nbige_{H^{(e)},0})\otimes\nbigohat_{X^{(e)},x^{(e)}}$.
Because $\ttQ_1(s_i)$ is $G_y$-invariant,
there exists
\[
 \ttQ(s_i)\in H^0(X^{(e)}\setminus H^{(e)},
 \Omega^1_{X^{(e)}\setminus H^{(e)}})
\]
such that
$\pi^{\ast}(\ttQ(s_i))=\ttQ_1(s_i)$.
We obtain the induced section
$[\ttQ(s_i)]$
of
$\Omega_{X^{(e)}}^1(\ast H^{(e)})/\Omega_{X^{(e)}}^1(\log H^{(e)})$.
Note that
$\ttQ(s_i)-d\gminia_i\in
\Omega_{X^{(e)}}^1(\log H^{(e)})_{x^{(e)}}
\otimes\nbigohat_{X^{(e)},x^{(e)}}$.
Hence, we obtain that
$d[\ttQ(s_i)]=0$
as a section of
\[
\Omega^2_{X^{(e)}}(\ast H^{(e)})/\Omega^2_{X^{(e)}}(\log H^{(e)}).
\]
Note that
$\Omega^{\bullet}_{X^{(e)}}(\log H^{(e)})_{x^{(e)}}
\to
\Omega^{\bullet}_{X^{(e)}}(\ast H^{(e)})_{x^{(e)}}$
is a quasi-isomorphism.
Hence,
there exists
\[
 \gminib_i\in
\nbigo_{X^{(e)},x^{(e)}}(\ast H^{(e)})/\nbigo_{X^{(e)},x^{(e)}}
\]
such that
$d\gminib_i=[\ttQ(s_i)]$.
For each $i,j\in I$,
we have $\gminib_i=\gminib_j$.
We set $\gminib_I:=\gminib_i$ $(i\in I)$.
Let $\gminibtilde_I\in \nbigo_{X,x}$ be a lift of $\gminib_I$.
By the previous consideration,
we obtain that
$\kappa_{d\gminibtilde_I}^{\ast}\Sigma^{(e)}_I$ is logarithmic.
Thus, we obtain Proposition \ref{prop;25.10.6.1}.
\hfill\qed

\subsubsection{Etale locally good meromorphic finite covers}

Let $\Sigma\subset T^{\ast}(X\setminus H)$
be a meromorphic finite cover over $(X,H)$.
For any $x\in X$,
we obtain a local meromorphic finite cover
$\nbigohat_{X,x}(T^{\ast}(X\setminus H))
\to
\nbigohat_{X,x}(\Sigma)$
over $\nbigohat_{X,x}$.

\begin{df}
\label{df;25.10.4.1}\mbox{{}}
\begin{itemize}
 \item 
$\Sigma$ is called formally good at $x$
 if
the induced formal cover
 $\nbigohat_{X,x}(T^{\ast}(X\setminus H))
\to
\nbigohat_{X,x}(\Sigma)$ is good.
In case, let 
$\nbigi(\Sigma,x)\subset
\nbigohat^{(H)}_{X,x}(\ast H)/\nbigohat^{(H)}_{X,x}$
denote the set of ramified irregular values
appearing as the index set of the decomposition.
(See {\rm\S\ref{subsection;25.12.7.30}}.)
\item
$\Sigma$ is called formally good on $(X,H)$
if it is formally good at any $x\in H$.
In case,
let 
$\vecnbigi(\Sigma)=(\nbigi(\Sigma,x)\,|\,x\in H)$.     
denote the tuple of the set of ramified irregular values.
\hfill\qed
\end{itemize}
\end{df}

\begin{df}
$\Sigma$ is called etale locally good at $x$
if
there exists an etale coordinate neighbourhood
$(U,x',\varphi,\psi)$ of $x\in X$
such that 
$\varphi^{\ast}\Sigma$ is good on $(U,H_U)$.
We say that $\Sigma$ is etale locally good on $(X,H)$
if it is etale locally good at any $x\in H$.
\hfill\qed
\end{df}

By Proposition \ref{prop;25.10.6.1},
we obtain the following proposition.
\begin{prop}
\label{prop;25.10.6.2}
$\Sigma$ is formally good on $(X,H)$
if and only if 
 it is etale locally good on $(X,H)$.
In the case,
there exist a finite etale coordinate cover
$\{(U_i,\varphi_i,\psi_i)\,|\,i=1,\ldots,k\}$
such that
$\varphi_i^{\ast}\Sigma$ are
good on $(U_i,H_{U_i})$. 
\hfill\qed
\end{prop}

\subsubsection{Pull back and push-forward via birational projective morphisms}

Let $\varphi:X'\to X$ be a birational projective morphism
such that $H'=\varphi^{-1}(H)$ is normal crossing
and that $X'\setminus H'\simeq X\setminus H$.
We have the natural isomorphism
$T^{\ast}(X'\setminus H')\simeq T^{\ast}(X\setminus H)$.
For any meromorphic finite cover $\Sigma\subset T^{\ast}(X\setminus H)$
over $(X,H)$,
we obtain the induced meromorphic finite cover
$\varphi^{\ast}(\Sigma)\subset T^{\ast}(X'\setminus H')$
over $(X',H')$.
Conversely,
for any meromorphic finite cover
$\Sigma'\subset T^{\ast}(X'\setminus H')$
over $(X',H')$,
we obtain the induced meromorphic finite cover
$\varphi_{\ast}(\Sigma')\subset T^{\ast}(X\setminus H)$
over $(X,H)$.

\begin{prop}
\mbox{{}}\label{prop;25.10.14.20}
\begin{itemize}
 \item Let 
       $\Sigma\subset T^{\ast}(X\setminus H)$
       be a formally good meromorphic finite cover over $(X,H)$.
       Then, $\varphi^{\ast}(\Sigma)$ is formally good on $(X',H')$,
       and
       $\vecnbigi(\varphi^{\ast}\Sigma)
	=\varphi^{\ast}(\vecnbigi(\Sigma))$ holds.
 \item Let $\Sigma'\subset T^{\ast}(X'\setminus H')$ be
       a formally good meromorphic finite cover over $(X',H')$.
       Suppose that there exists
       a tuple of the set of ramified irregular values
       $\vecnbigi=(\nbigi_x\,|\,x\in H)$
       such that
       $\vecnbigi(\Sigma')\subset\varphi^{\ast}\vecnbigi$.
       Then, $\varphi_{\ast}(\Sigma')\subset T^{\ast}(X\setminus H)$
       is formally good on $(X,H)$,
       and we have
       $\vecnbigi(\varphi_{\ast}\Sigma)\subset\vecnbigi$.
\end{itemize}
\end{prop}
\pf
Let $U=\Spec(R)$ be an affine $\hyperk$-variety
with a closed point $x$.
Let $\psi:U\to \Spec\hyperk[z_1,\ldots,z_n]$ be an etale morphism
such that $\psi(x)=(0,\ldots,0)$,
and that $\{\psi^{\ast}(z_i)=0\}$ are connected.
We set $H=\psi^{-1}(\bigcup_{i=1}^{\ell}\{z_i=0\})$ for some
$1\leq \ell\leq n$.
Let $U'=\Spec(R')$ be an affine $\hyperk$-variety
with a closed point $x'$.
Let $\psi':U'\to\Spec \hyperk[z_1,\ldots,z_n]$
be an etale morphism
such that
$\psi'(x')=(0,\ldots,0)$
and that $\{(\psi')^{\ast}z_i=0\}$ are connected.
We set $H'=(\psi')^{-1}\bigl(
\bigcup_{i=1}^{m}\{z_i=0\}\bigr)$.
We set $\psi_i=\psi^{\ast}(z_i)$
and $\psi'_i=(\psi')^{\ast}(z_i)$.
Let $\varphi:U'\to U$ be a birational morphism
such that
(i) $\varphi(x')=x$,
(ii) $\varphi^{\ast}(\psi_i)
=\prod_{j=1}^m(\psi'_j)^{k_{i,j}}\gtilde_i$
such that $\gtilde_i(x')\neq 0$.

Let $R^{(e)}=R[\psi_1^{1/e},\ldots,\psi_{\ell}^{1/e}]$
and $(R')^{(e)}=R[(\psi'_1)^{1/e},\ldots,(\psi'_{m})^{1/e}]$
for any $e\in\seisuu_{>0}$.
Let $U^{(e)}=\Spec R^{(e)}$
and $(U')^{(e)}=\Spec (R')^{(e)}$.
For any $e\in\seisuu_{>0}$,
there exists $e'\in\seisuu_{>0}$
with the following commutative diagram:
\[
\begin{CD}
 (U')^{(e')} @>>> U^{(e)} \\
 @VVV @VVV \\
 U' @>>> U.
\end{CD}
\]
A good finite subset of
$R^{(e)}(\ast H^{(e)})/R^{(e)}$
induces
a good finite subset of
$(R')^{(e')}(\ast (H')^{(e')})/(R')^{(e')}$.

Let $\gminia\in R^{(e)}(\ast H^{(e)})$.
If $\Sigma=\kappa_{d\gminia}^{\ast}\Sigma_1$ holds
for logarithmic $\Sigma_1$,
then
$\varphi^{\ast}\Sigma
=\kappa_{d\varphi^{\ast}\gminiatilde}^{\ast}\bigl(
\varphi^{\ast}\Sigma_1\bigr)$
holds for logarithmic $\varphi^{\ast}\Sigma_1$.
If $\Sigma'=\kappa_{d\gminiatilde}^{\ast}\Sigma_1'$ holds
for logarithmic $\Sigma_1'$,
then
$\varphi_{\ast}(\Sigma')
=\kappa_{d\gminiatilde}^{\ast}(\varphi_{\ast}\Sigma_1')$ holds
for logarithmic $\varphi_{\ast}\Sigma'_1$.
Then, we obtain the claim of
Proposition \ref{prop;25.10.14.20}
from Proposition \ref{prop;25.10.6.2}.
\hfill\qed

\subsection{Family case}

\subsubsection{Local algebraic conditions}

Let $S$ be a regular integral domain over $\hyperk$.
Let $R$ be an $S[x_1,\ldots,x_n]$-algebra
which is finite and etale over $S[x_1,\ldots,x_n]$.
Let $U=\Spec(R)$.
We obtain the normal crossing hypersurface 
$H=\bigcup_{i=1}^{\ell}\{x_i=0\}$ of $U$.
We set $U^{\circ}=U\setminus H$.
Let $T^{\ast}(U^{\circ}/S)$ denote the relative cotangent bundle
of $U^{\circ}$ over $\Spec(S)$.
Let $T^{\ast}(U/S)(\log H)$
denote the relative logarithmic cotangent bundle
of $(U,H)$ over $\Spec(S)$.

Let $\Sigma\subset T^{\ast}(U^{\circ}/S)$
be a finite cover of $U^{\circ}$,
which is flat over $S$.

\begin{df}
We say that $\Sigma$ is logarithmic relative to $S$
if the closure of $\Sigma$
in $T^{\ast}(U/S)(\log H)$ is proper over $U$.
\hfill\qed
\end{df}

For any $f\in R(\ast H)$,
we obtain the section $df$ of $T^{\ast}(U^{\circ}/S)$
by taking the exterior derivative relative to $\Spec(S)$.
We obtain the automorphism
$\kappa_{df}$ of $T^{\ast}(U^{\circ}/S)$
by adding $df$.

\begin{df}
\label{df;26.1.21.1}
We say that $\Sigma$ is unramifiedly good on $(U,H)$ relative to $S$
if there exist a good set of irregular values
$\nbigi\subset R(\ast H)/R$
and a decomposition
\[
 \Sigma=\bigsqcup_{\gminia\in\nbigi}
 \Sigma_{\gminia}
\]
such that 
each $\kappa_{d\gminiatilde}^{\ast}\Sigma_{\gminia}$ 
is logarithmic.
Here, $\gminiatilde\in R(\ast H)$ denotes a lift of $\gminia$.
\hfill\qed
\end{df}

Let $\Spec(\Stilde)\to \Spec(S)$ be an etale covering
of affine schemes.
We set $\Rtilde=\Stilde\otimes_SR$
and $\Utilde=\Rtilde$.
For $e\in\seisuu_{>0}$,
we set $\Rtilde^{(e)}=\Rtilde[x_1^{1/e},\ldots,x_{n}^{1/e}]$
and $\Utilde^{(e)}=\Spec(\Rtilde^{(e)})$.
We obtain the morphism
$\rho_{e}:\Utilde^{(e)}\to U$.
We set $\Htilde^{(e)}=\rho_{e}^{-1}(H)$
and $\Sigmatilde^{(e)}=\rho_{e}^{\ast}(\Sigma)$.
\begin{df}
We say that $\Sigma$ is good on $(U,H)$ relative to $S$
if there exist
an etale covering of affine schemes $\Spec(\Stilde)\to\Spec(S)$
and $e\in\seisuu_{>0}$ such that 
$\rho_{e}^{\ast}(\Sigma)$ is
unramifiedly good on $(\Utilde^{(e)},\Htilde^{(e)})$ relative to
$\Stilde$.
\hfill\qed
\end{df}

Let $\hyperk'$ be any algebraically closed field of characteristic $0$.
Let $S\to \hyperk'$ be a morphism of $\hyperk$-algebras.
We obtain
$(U_{\hyperk'},H_{\hyperk'},\Sigma_{\hyperk'})$
from $(U,H,\Sigma)$
and $S\to\hyperk'$.
The following lemma is clear.
\begin{lem}
\label{lem;25.12.7.11}
If $\Sigma$ is good on $(U,H)$ relative to $S$,
then $\Sigma_{\hyperk'}$ is good on $(U_{\hyperk'},H_{\hyperk'})$.
\hfill\qed
\end{lem}

\begin{lem}
\label{lem;25.12.7.10}
Let $\hyperk_1$ be an algebraically closed field
which contains $S$.
Suppose that $\Sigma_{\hyperk_1}$ is
good (resp. logarithmic, unramifiedly good)
on $(U_{\hyperk_1},H_{\hyperk_1})$.
Then, there exists an extension 
of regular $\hyperk$-algebras $S\subset S_1$ in $\hyperk_1$
such that 
$\Sigma_{S_1}$ is good
(resp. logarithmic, unramifiedly good)
on $(U_{S_1},H_{S_1})$ relative to $S_1$,
where $(U_{S_1},H_{S_1},\Sigma_{S_1})$
is induced by $(U,H,\Sigma)$ and the extension $S\subset S_1$.
\end{lem}
\pf
It is enough to consider the case where
$\Sigma_{\hyperk_1}$ is unramifiedly good.
We set
$R_{\hyperk_1}=R\otimes_S\hyperk_1$.
There exist a good set of irregular values
$\nbigi\subset
R_{\hyperk_1}(\ast H_{\hyperk_1})
\Big/
R_{\hyperk_1}$
and a decomposition
$\Sigma_{\hyperk_1}
=\bigsqcup_{\gminia\in\nbigi}
\Sigma_{\hyperk_1,\gminia}$
such that
$\kappa_{d\gminiatilde}^{\ast}
\Sigma_{\hyperk_1,\gminia}$
are logarithmic.
By enlarging $S$ in $\hyperk_1$,
we may assume that
there exists a finite subset
$\nbigi\subset
R(\ast H)
\Big/R$
which induces $\nbigi$.
By Lemma \ref{lem;25.10.6.20},
we may assume that $\nbigi_S$
is a good set of irregular values
because $\nbigi$ is a good set of irregular values.

There exists a decomposition
$\Sigma
 =\bigcup_{\gminia\in\nbigi}
 \Sigma_{\gminia}$
such that
$\Sigma_{\hyperk_1,\gminia}=\Sigma_{\gminia}\times_S\hyperk_1$.
For $\gminia\neq\gminib$,
because 
$\Sigma_{\hyperk_1,\gminia}\cap\Sigma_{\hyperk_1,\gminib}
=\emptyset$,
$\Sigma_{\gminia}\cap\Sigma_{\gminib}$
is not dominant over $\Spec(S)$.
By enlarging $S$,
we may assume that
$\Sigma_{\gminia}\cap\Sigma_{\gminib}=\emptyset$.

For $\gminia\in\nbigi$,
the corresponding element of $\nbigi_S$
is also denoted by $\gminia$.
Let $\gminiatilde\in R(\ast H)$ be a lift of $\gminia$.
Let us consider
$\kappa_{d\gminiatilde}^{\ast}\Sigma_\gminia$.
Let $\proj$ be a projective completion of
$T^{\ast}(U/S)(\log H)$.
Let $\proj_{\infty}$ denote
the complement of $T^{\ast}(U/S)(\log H)$ in $\proj$.
Let $\overline{\kappa_{d\gminiatilde}^{\ast}\Sigma_{\gminia}}$
denote the closure of
$\kappa_{d\gminiatilde}^{\ast}\Sigma_{\gminia}$ in $\proj$.
Because 
$\kappa_{d\gminiatilde}^{\ast}\Sigma_{\hyperk_1,\gminia}$
is logarithmic,
$\overline{\kappa_{d\gminiatilde}^{\ast}\Sigma_{\gminia}}
\cap \proj_{\infty}$
is not dominant over $\Spec(S)$.
By enlarging $S$,
we may assume that
$\overline{\kappa_{d\gminiatilde}^{\ast}\Sigma_{\gminia}}
\cap \proj_{\infty}
=\emptyset$.
It implies that
$\kappa_{d\gminiatilde}^{\ast}\Sigma_{\gminia}$
is logarithmic.
\hfill\qed

\vspace{.1in}
Let $\Sigma^{(i)}\subset T^{\ast}(U^{\circ}/S)$ $(i=1,2)$
be good covers of $U^{\circ}$.
For any closed point $z$ of $\Spec(S)$,
we set
$(U_z,H_z,\Sigma^{(i)}_z)=(U,H,\Sigma^{(i)})\times_{\Spec(S)}z$.
We obtain 
good meromorphic covers
$\Sigma^{(i)}_z$ on $(U_z,H_z)$,
and the associated tuples
$\vecnbigi(\Sigma^{(i)}_z)
=(\nbigi(\Sigma^{(i)}_z,x)\,|\,x\in H_z)$.

\begin{lem} 
\label{lem;26.1.21.2}
There exists a closed subset $Z\subset \Spec(S)$
such that
$\vecnbigi(\Sigma^{(1)}_z)
=\vecnbigi(\Sigma^{(2)}_z)$ holds
if and only if  
$z\in Z$.
\end{lem}
\pf
It is enough to consider the case where 
both $\Sigma^{(i)}$ are unramifiedly good.
Let $\nbigi^{(i)}\subset R(\ast H)/R$
be the index set in Definition \ref{df;26.1.21.1}.
For each $\gminia\in\nbigi^{(i)}$
and $z\in\Spec(S)$,
let $\gminia_z\in\nbigi^{(i)}_z\subset R_z(\ast H_z)/R_z$
denote the induced element.
It is enough to consider the case
where the hypersurfaces $\{z_i=0\}\subset U$
are connected,
and there exists a section
$\alpha:\Spec(S)\to \bigcap\{z_i=0\}\subset U$.
Let $\Rhat_{\alpha}$ denote the completion of $R$
along the image of $\alpha$.
We obtain the induced subsets
$\nbigihat^{(i)}_{\alpha}\subset
\Rhat_{\alpha}(\ast H)/\Rhat_{\alpha}$.
For $\vecm\in\seisuu^{\ell}$,
we set $x^{\vecm}=\prod_{i=1}^{\ell}x_i^{m_i}$.

For $\gminia\in\nbigi^{(i)}$,
we have the expansion
$\gminia=\sum \gminia_{\vecm}x^{\vecm}$,
where $\gminia_{\vecm}=0$ for $\vecm\in\seisuu_{\geq 0}^{\ell}$.
We may regard $\gminia_{\vecm}$
as functions on $\Spec(S)$.

For $\gminia\in\nbigihat^{(1)}_{\alpha}$
and $\gminib\in\nbigihat^{(2)}_{\alpha}$,
we obtain the following closed subsets
$Z(\gminia,\gminib)
 =\bigcap_{\vecm}
 \{\gminia_{\vecm}=\gminib_{\vecm}\}
 \subset
 \Spec(S)$.
Then,
$\gminia_z=\gminib_z$ holds
if and only if $z\in Z(\gminia,\gminib)$.
We have $\gminia_z\in\nbigi^{(2)}_z$
if and only if
$z\in \bigcup_{\gminib\in\nbigi^{(2)}}Z(\gminia,\gminib)$.
Hence, 
$\nbigi^{(1)}_z\subset\nbigi^{(2)}_z$
if and only if
$z\in\bigcap_{\gminia\in\nbigi^{(1)}}
 \bigcup_{\gminib\in\nbigi^{(2)}}Z(\gminia,\gminib)$.
Similarly,
$\nbigi^{(2)}_z\subset\nbigi^{(1)}_z$
if and only if
$z\in\bigcap_{\gminia\in\nbigi^{(2)}}
 \bigcup_{\gminib\in\nbigi^{(1)}}Z(\gminia,\gminib)$.
Then, we obtain the claim of the lemma.
\hfill\qed

\subsubsection{Families of etale locally good meromorphic finite covers}
\label{subsection;25.12.7.40}

Let $X$ be a $\hyperk$-variety smooth over $\Spec(S)$.
Let $H\subset X$ be a hypersurface of $X$
which is normal crossing relative to $S$
in the following sense.
\begin{itemize}
 \item $H$ is normal crossing.
 \item Let $H=\bigcup_{j\in\Lambda} H_{j}$
       denote the irreducible decomposition.
       For any finite subset $I\subset \Lambda$,
       $H_{I}:=\bigcap_{j\in I}H_{j}$
       is smooth over $\Spec S$.
\end{itemize}
We also assume that the monodromy of $H$ is trivial
in the following sense.
\begin{itemize}
 \item Let $H_{I}=\bigcup_{j\in\Lambda(I)} H_{I,j}$
       denote the irreducible decomposition.
       Then, any geometric fiber of
       $H_{I,j}\to \Spec(S)$
       is irreducible.
\end{itemize}
Let $T^{\ast}(X/S)$ denote the relative cotangent bundle
of $X$ over $\Spec(S)$.
Let $T^{\ast}(X/S)(\log H)$
denote the relative cotangent bundle of $(X,H)$ over $S$.

Let $\Sigma\subset T^{\ast}((X\setminus H)/S)$
be a finite cover of $X\setminus H$,
which is flat over $S$.
Let $x$ be any closed point of $H$.
For any etale coordinate neighbourhood
$(U,x',\varphi,\psi)$
of $x$ in $X$
in the sense of \S\ref{subsection;25.12.7.1},
we set
$H_{U}=\varphi^{-1}(H)$
and 
$\Sigma_U=\varphi^{\ast}\Sigma\subset T^{\ast}((U\setminus H_U)/S)$.

\begin{df}
We say that $\Sigma$ is etale locally good
(resp. logarithmic, unramifiedly good)
at $x$ relative to $S$
if there exists an etale coordinate neighbourhood 
$(U,x',\varphi,\psi)$ of $x$
such that $\varphi^{\ast}\Sigma_U$
is good
(logarithmic, unramifiedly good)
on $(U,H_U)$.
\hfill\qed
\end{df}

\begin{df}
\label{df;25.12.7.2}
We say that $\Sigma$ is etale locally good
(resp. logarithmic, unramifiedly good) on $(X,H)$ relative to $S$
if it is etale locally good
(logarithmic, unramifiedly good)
at any $x\in X$ relative to $S$.
\hfill\qed
\end{df}

Let $\hyperk'$ be any algebraically closed field of characteristic $0$.
Let $S\to \hyperk'$ be a morphism of $\hyperk$-algebras.
We obtain $(X_{\hyperk'},H_{\hyperk'},\Sigma_{\hyperk'})$
from $(X,H,\Sigma)$ by $S\to\hyperk'$.
We obtain the following lemma from Lemma \ref{lem;25.12.7.11}.
\begin{lem}
\label{lem;25.12.7.31}
If $\Sigma$ is good
(resp. logarithmic, unramifiedly good)
on $(X,H)$ relative to $S$,
 then  $\Sigma_{\hyperk'}$ is good
(resp. logarithmic, unramifiedly good)
 on $(X_{\hyperk'},H_{\hyperk'})$.
\hfill\qed
\end{lem}

\begin{lem}
\label{lem;25.12.7.12}
Let $\hyperk_1$ be an algebraically closed field
which contains $S$.
Suppose that $\Sigma_{\hyperk_1}$ is
good (resp. logarithmic, unramifiedly good)
on $(X_{\hyperk_1},H_{\hyperk_1})$.
Then, there exists an extension 
of regular $\hyperk$-algebras $S\subset S_1$ in $\hyperk_1$
such that 
$\Sigma_{S_1}$ is good
(resp. logarithmic, unramifiedly good)
on $(X_{S_1},H_{S_1})$ relative to $S_1$,
where $(X_{S_1},H_{S_1},\Sigma_{S_1})$
is induced by $(X,H,\Sigma)$ and the extension $S\subset S_1$.
\end{lem}
\pf
There exists
a finite tuple
of etale morphisms
$\varphi_j:U_{\hyperk_1,j}\to X_{\hyperk_1}$
with an etale morphism
$\psi_j:U_{\hyperk_1,j}\to \Spec \hyperk_1[x_1,\ldots,x_n]$
such that the following holds.
\begin{itemize}
 \item $X_{\hyperk_1}=\bigcup_{j=1}^m \Image \varphi_j$.
 \item $H_{U_{\hyperk_1,j}}:=\varphi_j^{-1}(H_{\hyperk_1})
       =\bigcup_{i=1}^{\ell(j)}\{\psi_j^{\ast}(x_i)=0\}$.
       Moreover,
       $\{\psi_j^{\ast}(x_i)=0\}$ is connected.
 \item $\Sigma_{U_{\hyperk_1,j}}=\varphi_j^{\ast}\Sigma_{\hyperk_1}$
       is good on
       $(U_{\hyperk_1,j},H_{U_{\hyperk_1,j}})$.
\end{itemize}
By enlarging $S$,
we may assume that
there exist
a tuple of etale morphisms
$\varphi_{S,j}:U_{S,j}\to X_{S}$
with an etale morphism
$\psi_{S,j}:U_{S,j}\to \Spec S[x_1,\ldots,x_n]$
which induce $(\varphi_{j},\psi_j)$ $(j=1,\ldots,m)$.

We set $H_{U_{S,j}}=\varphi_{S,j}^{-1}(H_S)$
and
$U_{S,j}^{\circ}=U_{S,j}\setminus H_{U_{S,j}}$.
We obtain the closed subsets 
$\Sigma_{U_{S,j}}=
\varphi_{S,j}^{\ast}(\Sigma_{\hyperk_1})$
of $T^{\ast}(U_{S,j}^{\circ}/S)$.
According to \cite[Theorem 6.9.1]{EGA-IV-2},
by enlarging $S$,
we may assume that
$\Sigma_{U_{S,j}}$ are flat over $S$.
According to Lemma \ref{lem;25.12.7.10},
by enlarging $S$,
we may assume that
$\Sigma_{U_{S,j}}$
are good (resp. logarithmic, unramifiedly good)
on $(U_{S,j},H_{U_{S,j}})$ relative to $S$.
Hence, we obtain the claim of Lemma \ref{lem;25.12.7.12}.
\hfill\qed

\subsubsection{Families of Lagrangian covers}

Let $\Sigma\subset T^{\ast}((X\setminus H)/S)$
be a finite cover of $X\setminus H$,
which is flat over $\Spec(S)$.
Let $\Sigma^{\sm}$ denote the smooth part of $\Sigma$.
We also assume the following.
\begin{itemize}
 \item Each geometric fiber of $\Sigma\to\Spec(S)$
is reduced.
 \item For any closed point $s$ of $\Spec(S)$,
       any irreducible component of
       $\Sigma_{S,s}=\Sigma_S\times_{\Spec(S)}s$
       has hon-empty intersection
       with $\Sigma^{\sm}$.
\end{itemize}

\begin{df}
$\Sigma$ is called Lagrangian if
$\Sigma_{S,s}$ is Lagrangian 
in the sense of Definition {\rm\ref{df;25.12.7.15}}
 for any closed point $s$ of $\Spec(S)$.
 \hfill\qed
\end{df}

Let $\omega$ denote the relative symplectic structure of
$T^{\ast}(X/S)$.
Let $T^{\ast}(\Sigma^{\sm}/S)$ denote
the relative cotangent bundle of $\Sigma^{\sm}$ over $S$.
We obtain the induced section
$\Psi_{\omega}$
of $\bigwedge^2T^{\ast}(\Sigma^{\sm}/S)$
induced by $\omega$.

\begin{lem}
$\Sigma$ is Lagrangian if and only if 
$\Psi_{\omega}=0$.
\hfill\qed
\end{lem}

The following lemmas are clear.
\begin{lem}
\label{lem;25.12.7.30}
If $\Sigma$ is Lagrangian,
then any geometric fiber of $\Sigma_S\to\Spec(S)$
is Lagrangian.
\hfill\qed
\end{lem}

\begin{lem}
Suppose that $\Sigma_{\hyperk_1}$ is Lagrangian
for an injection of $S$ to an algebraically closed field $\hyperk_1$.
Then, $\Sigma$ is Lagrangian.
\hfill\qed
\end{lem}

\subsubsection{Comparison of good meromorphic covers}

Let $\Sigma^{(i)}\subset T^{\ast}((X\setminus H)/S)$ $(i=1,2)$ be 
finite covers of $X\setminus H$ which are etale locally good on $(X,H)$.
For each closed point $z$ of $\Spec(S)$,
we obtain good meromorphic covers $\Sigma^{(i)}_z$
on $(X_z,H_z)$.
We obtain the tuples
$\vecnbigi(\Sigma^{(i)}_z)$.
We obtain the following proposition from Lemma \ref{lem;26.1.21.2}.

\begin{prop}
\label{prop;26.1.21.3}
There exists a closed subset $Z\subset \Spec(S)$
such that
$\vecnbigi(\Sigma^{(1)}_z)
=\vecnbigi(\Sigma^{(2)}_z)$
if and only if
$z\in Z$. 
\hfill\qed
\end{prop}

\subsection{Models of meromorphic finite covers}

\subsubsection{Preliminary}

Let $\hyperk_1$ be an algebraically closed field
which contains $\hyperk$.
Let $X_{\hyperk_1}$ be
an irreducible projective normal $\hyperk_1$-variety.
Let
$\vecZ_{\hyperk_1}=(Z_{\hyperk_1,j}\,|\,j=1,\ldots,m_1)$
be a tuple of closed subsets of $X_{\hyperk_1}$
such that
$X^{\circ}_{\hyperk_1}
:=X_{\hyperk_1}\setminus \bigcup_{j=1}^{m_1}Z_{\hyperk_1,j}$
is smooth.
Let $\vecW_{\hyperk_1}=(W_{\hyperk_1,k}\,|\,k=1,\ldots,m_2)$
be a tuple of closed subsets of
$T^{\ast}X_{\hyperk_1}^{\circ}$
such that each $W_{\hyperk_1,k}$ is a finite and dominant over
$X^{\circ}_{\hyperk_1}$.
For each $J\subset \{1,\ldots,m_1\}$,
we set
$Z_{\hyperk_1,J}=\bigcap_{j\in J}Z_{\hyperk_1,j}$.
We have the irreducible decompositions
$Z_{\hyperk_1,J}=\bigcup_{p\in\Lambda(J,1)} Z_{\hyperk_1,J,p}$.

Let $S\subset\hyperk_1$
be a regular subring finitely generated over $\hyperk$.

\begin{df}
A model of $(X_{\hyperk_1},\vecZ_{\hyperk_1},\vecW_{\hyperk_1})$
over $S$
consists of the data
$(X_S,\vecZ_{S},\vecW_S)$ as follows.
\begin{itemize}
 \item $X_S$ denotes an irreducible $\hyperk$-variety
       flat over $\Spec(S)$
       such that
       (i) $X_S$ is projective over $\Spec(S)$,
       (ii) each geometric fiber of $X_S\to\Spec(S)$
       is normal and irreducible,
       (iii) $X_S\times\hyperk_1=X_{\hyperk_1}$.
 \item $\vecZ_{S}=(Z_{S,j}\,|\,j=1,\ldots,m_1)$ denotes
      a tuple of irreducible closed subset
	of $X_S$
       such that
       (i) $Z_{S,j}$ are flat over $\Spec S$,
       (ii) each geometric fiber of
       $Z_{S,j}\to\Spec(S)$
       is irreducible and reduced,
       (iii) $Z_{S,j}\times\hyperk_1=Z_{\hyperk_1,j}$,
       (iv)
       $X^{\circ}_S=X_S\setminus\bigcup_{j=1}^{m_1} Z_{S,j}$
       is smooth over $S$.
 \item For any non-empty subset $J\subset \{1,\ldots,m_1\}$,
       we have the decomposition
       $\bigcap_{j\in J}Z_{S,j}
       =\bigcup_{p\in\Lambda(J,1)}Z_{S,J,p}$ such that
       (i) $Z_{S,J,p}$ are flat over $\Spec(S)$,
       (ii) each geometric fiber of
       $Z_{S,J,p}\to\Spec(S)$
       is irreducible and reduced,
       (iii) $Z_{S,J,p}\times\hyperk_1=Z_{\hyperk_1,J,p}$.
 \item $\vecW_S=(W_{S,k}\,|\,k=1,\ldots,m_2)$ denotes
       a tuple of irreducible closed subsets of
       $T^{\ast}(X_S^{\circ}/S)$
       such that
       (i) $W_{S,k}$ are flat over $\Spec S$,
       (ii) each geometric fiber of
       $W_{S,k}\to\Spec(S)$
       is irreducible and reduced,
       (iii)
       $W_{S,k}\times \hyperk_1=W_{\hyperk_1,k}$,
       (iv) $W_{S,k}$ are finite and dominant over
       $X^{\circ}_{S}$.
       Here,
       $T^{\ast}(X_S^{\circ}/S)$
       denotes the relative cotangent bundle of
       $X^{\circ}_S$ over $\Spec S$.
\end{itemize}
We shall impose the following additional conditions.
\begin{itemize}
 \item
      Let $I\subset \{1,\ldots,m_1\}$ be any subset.
      If
      $\bigcap_{j\in I} Z_{\hyperk_1,j}$
      is smooth over $\hyperk_1$,
      then $\bigcap_{j\in I}Z_{S,j}$ is smooth over $\Spec S$.
      (We formally set $Z_{\emptyset}=X$.)
       \hfill\qed
\end{itemize}
\end{df}

The following lemma is standard.
\begin{lem}
There exist a regular subring $S\subset \hyperk_1$
finitely generated over $\hyperk$
and a model $(X_S,\vecZ_S,\vecW_S)$
of $(X_{\hyperk_1},\vecZ_{\hyperk_1},\vecW_{\hyperk_1})$
over $S$.
\end{lem}
\pf
By taking an embedding
$X_{\hyperk_1}\subset \proj^N_{\hyperk_1}$ for a large $N$,
we regard
$X_{\hyperk_1}$ and $Z_{\hyperk_1,j}$
are subvarieties of $\proj^N_{\hyperk_1}$.
Let $\gbigi$, $\gbigi_{j}$, $\gbigi_{J,p}$
denote the homogeneous ideals
of $\hyperk_1[y_0,\ldots,y_N]$
corresponding to $X_{\hyperk_1}$, $Z_{\hyperk_1,j}$
and $Z_{\hyperk_1,I,p}$,
respectively.
There exists a regular subring $S_1\subset \hyperk_1$
finitely generated over $\hyperk$
such that
$\gbigi$, $\gbigi_j$ $\gbigi_{J,p}$ are induced by
homogeneous ideals
$\gbigi_{S_1}$, $\gbigi_{S_1,j}$ and $\gbigi_{S_1,J.p}$
of $S_1[y_0,\ldots,y_N]$,
respectively.
Let $X_{S_1}$,  $Z_{S_1,j}$, $Z_{S_1,J,p}$
denote the corresponding closed subvarieties of $\proj^N_{S_1}$.
We may assume that they are dominant over $\Spec S_1$.
We may assume that $Z_{S_1,j}\subset X_{S_1}$
and $Z_{S_1,J,p}\subset \bigcap_{j\in J}Z_{S_1,j}$.
By the generic freeness theorem
\cite[Theorem 6.9.1]{EGA-IV-2},
we may assume that
$X_{S_1}$, $Z_{S_1,j}$
and $Z_{S_1,J,p}$ are flat over $\Spec(S_1)$.
By \cite[Theorem 12.2.4]{EGA-IV-3},
we may assume the following.
\begin{itemize}
 \item Each geometric fiber
of $X_{S_1}\to \Spec(S_1)$ is normal and irreducible.
We may also assume that
$X^{\circ}_{S_1}$ is smooth over $\Spec(S_1)$.
 \item Each geometric fiber of
       $Z_{S_1,j}\to\Spec(S_1)$ is irreducible and reduced.
 \item Each geometric fiber of
       $Z_{S_1,J,p}\to\Spec(S_1)$ is irreducible and reduced.
\end{itemize}
In particular,
$X_{S_1}$,
$Z_{S_1,j}$
and $Z_{S_1,J,p}$ are irreducible.
We have
$\bigcup_{p\in\Lambda(J,1)}Z_{S_1,J.p}
\subset
\bigcap_{j\in J}Z_{S_1,j}$.
Let $Z'_{S_1}$ be an irreducible component of
$\bigcap_{j\in J}Z_{S_1,j}$
which is not contained in
$\bigcup_{p\in\Lambda(J,1)}Z_{S_1,J.p}$
Because 
$\bigcup_{p\in\Lambda(J,1)}Z_{\hyperk_1,J,p}
=\bigcap_{j\in J}Z_{\hyperk_1,j}$,
$Z'_{S_1}$ is not dominant over $\Spec(S_1)$.
Hence, we may assume that
$\bigcup_{p\in\Lambda(J,1)}Z_{S_1,J.p}
=
\bigcap_{j\in J}Z_{S_1,j}$.

Let $U_{S_1}=\Spec R_{S_1}$ be an affine open subset
of $X_{S_1}^{\circ}$.
Let $T^{\ast}(U_{S_1}/S_1)$
denote the relative cotangent bundle of $U_{S_1}$ over $\Spec S_1$.
Let $\Rtilde_{S_1}$ denote the ring of algebraic functions
of $T^{\ast}(U_{S_1}/S_1)$.
We have
$T^{\ast}(U_{S_1}/S_1)=\Spec \Rtilde_{S_1}$.
Let $R_{\hyperk_1}=R_{S_1}\otimes\hyperk_1$
and $\Rtilde_{\hyperk_1}=\Rtilde_{S_1}\otimes\hyperk_1$.
We obtain the induced affine open subset
$U_{\hyperk_1}=\Spec R_{\hyperk_1}$
of $X_{\hyperk_1}^{\circ}$.
We have
$T^{\ast}U_{\hyperk_1}=\Spec \Rtilde_{\hyperk_1}$.
We have the ideal $\gbigi'_{k}$
of $\Rtilde_{\hyperk_1}$
corresponding to
$T^{\ast}U_{\hyperk_1}\cap W_{\hyperk_1,k}$.
We may assume that
there exists an ideal $\gbigi'_{S_1,k}$ of $\Rtilde_{S_1}$
which induces $\gbigi'_{k}$.
Let $W_{U_{S_1},k}\subset T^{\ast}(U_{S_1}/S_1)$
denote the corresponding closed subset.
By \cite[Theorem 6.9.1]{EGA-IV-2} and
\cite[Theorem 12.2.4]{EGA-IV-3},
we may assume the following.
\begin{description}
 \item[(a1)] $W_{U_{S_1},k}$ is flat over $\Spec S_1$.
 \item[(a2)] Each geometric fiber of $W_{U_{S_1},k}\to \Spec (S_1)$
       is irreducible and reduced.
 \item[(a3)] $W_{U_{S_1},k}$ induces
	    $W_{\hyperk_1,k}\cap T^{\ast}(U_{\hyperk_1})$.
\end{description}
Because $W_{U_{S_1},k}$ is flat over $\Spec(S_1)$,
the following morphism is injective.
\[
 \Rtilde_{S_1}/\gbigi'_{S_1,k}
 \to
 \Bigl(
 \Rtilde_{S_1}/\gbigi'_{S_1,k}
 \Bigr)\otimes_{S_1}\hyperk_1
=
 \Rtilde_{\hyperk_1}/\gbigi'_k.
\]
Hence, the above conditions
(a1,a2,a3) characterize $W_{U_{S_1},k}$.

Let $X_{S_1}=\bigcup_{q=1}^{M} U_{S_1,q}$
be a covering by affine open subsets.
We obtain the closed subsets
$W_{U_{S_1,q},k}\subset T^{\ast}(U_{S_1,q}/S_1)$
$(k=1,\ldots,m_2,\,\,q=1,\ldots,M)$.
Because of the characterization by
the conditions (a1,a2,a3),
we obtain
\[
W_{U_{S_1,q(1)},k}
\cap
T^{\ast}(U_{S_1,q(1)}/S_1)\cap T^{\ast}(U_{S_1,q(2)}/S_1)
=
W_{U_{S_1,q(2)},k}
\cap
T^{\ast}(U_{S_1,q(1)}/S_1)\cap T^{\ast}(U_{S_1,q(2)}/S_1)
\]
in $T^{\ast}(U_{S_1,q(1)}/S_1)\cap T^{\ast}(U_{S_1,q(2)}/S_1)$.
Hence, there exist closed subsets
$W_{S_1,k}$ such that
\[
W_{S_1,k}\cap T^{\ast}(U_{S_1,q}/S_1)
=W_{U_{S_1,q},k}.
\]
By the construction,
the following holds.
\begin{itemize}
 \item $W_{S_1,k}$ are flat over $\Spec(S_1)$.
 \item Each geometric fiber of
       $W_{S_1,k}\to \Spec(S_1)$
       is irreducible and reduced.
 \item $W_{S_1,k}$ induces $W_{\hyperk_1,k}$.
\end{itemize}

Clearly, $W_{S_1,k}\to X^{\circ}_{S_1}$ is dominant.
We may assume that
$W_{S_1,k}\to X^{\circ}_{S_1,k}$ is proper.
Indeed, let $T^{\ast}(X_{S_1}^{\circ}/S_1)\subset \proj$
denote the projective completion.
Let $\proj_{\infty}$ denote the complement of
$T^{\ast}(X_{S_1}^{\circ}/S_1)$.
Let $\overline{W}_{S_1,k}$ denote the closure of
$W_{S_1,k}$ in $\proj$.
Let $T\subset X_{S_1}^{\circ}$
denote the image of
$\overline{W}_{S_1,k}\cap \proj_{\infty}$
by the projection $\proj\to X_{S_1}^{\circ}$.
It is a closed subset.
Because $W_{\hyperk_1,k}\to X^{\circ}_{\hyperk_1}$
is proper,
we obtain that $T\to \Spec(S_1)$ is not dominant.
Hence, we may assume that $T=\emptyset$,
i.e.,
$W_{S_1,k}\to X^{\circ}_{S_1}$ is proper.
\hfill\qed

\vspace{.1in}

Let $(X_S,\vecZ_S,\vecW_S)$ be a model over $S$.
For any closed point $s$ of $S$,
by taking the fiber over $s$,
we obtain
the $\hyperk$-variety $X_{S,s}$
with the tuple of closed subsets
$\vecZ_{S,s}=(Z_{S,s,j}\,|\,j=1,\ldots,m_1)$ of $X_{S,s}$
and the tuple of closed subsets
$\vecW_{S,s}=(W_{S,s,k}\,|\,k=1,\ldots,m_2)$ of
$T^{\ast}X^{\circ}_{S,s}$,
where $X^{\circ}_{S,s}=X_{S,s}\setminus \bigcup_{j=1}^{m_1} Z_{S,s,j}$.

\begin{lem}
$W_{S,s,k}\to X^{\circ}_{S,s}$
is finite and dominant over $X^{\circ}_{S,s}$.
\end{lem}
\pf
Because $W_{S,k}\to X^{\circ}_{S}$ is proper,
we obtain that
$W_{S,s,k}\to X^{\circ}_{S,s}$ is proper,
which implies the finiteness.
By using the valuative criterion,
we obtain that the morphism
$W_{S,k}\to X^{\circ}_{S}$ 
induces 
a surjection of the closed points.
In particular, $W_{S,s,k}\to X^{\circ}_{S,s}$ is dominant.
\hfill\qed

\subsubsection{Models for meromorphic covers}
\label{subsection;25.10.13.20}

Let $X_{\hyperk_1}$ be a normal irreducible projective $\hyperk_1$-variety
with a hypersurface
$H_{\hyperk_1}$.
Let
$H_{\hyperk_1}=\bigcup_{j=1}^{m_1} H_{\hyperk_1,j}$
denote the irreducible decomposition.
Let $\Sigma_{\hyperk_1}\subset
T^{\ast}(X_{\hyperk_1}\setminus H_{\hyperk_1})$
be a meromorphic finite cover over $(X_{\hyperk_1},H_{\hyperk_1})$.
Let $\Sigma_{\hyperk_1}=\bigcup_{k=1}^{m_2} \Sigma_{\hyperk_1,k}$
denote the irreducible decomposition.
We obtain the tuples
$\vecH_{\hyperk_1}=(H_{\hyperk_1,j})$
and
$\vecSigma_{\hyperk_1}=(\Sigma_{\hyperk_1,k})$.

Let $S\subset \hyperk_1$ be
a regular ring finitely generated over $\hyperk$.
Let $(X_S,\vecH_{S},\vecSigma_S)$
be a model of
$(X_{\hyperk_1},\vecH_{\hyperk_1},\vecSigma_{\hyperk_1})$.
We obtain the hypersurface
$H_{S}=\bigcup_{j=1}^{m_1}H_{S,j}$ of $X_{S}$,
and the closed subset
$\Sigma_{S}=\bigcup_{k=1}^{m_2} \Sigma_{S,k}$
of $T^{\ast}(X^{\circ}_{S}/S)$,
where
$X^{\circ}_{S}=X_{S}\setminus H_S$.

\begin{df}
Such $(X_S,H_S,\Sigma_S)$
is called a model of $(X_{\hyperk_1},H_{\hyperk_1},\Sigma)$
over $S$.
\hfill\qed
\end{df}

We obtain the following proposition
by Lemma \ref{lem;25.12.7.12}.
\begin{prop}
\label{prop;25.12.7.20}
Suppose that $X_{\hyperk_1}$ is smooth,
that $H_{\hyperk_1}$ is simply normal crossing,
and that $\Sigma$ is good over $(X_{\hyperk_1},H_{\hyperk_1})$.
Then, 
there exists a regular subring $S\subset \hyperk_1$
finitely generated over $\hyperk_1$
and a model $(X_S,H_S,\Sigma_S)$ of
$(X_{\hyperk_1},H_{\hyperk_1},\Sigma_{\hyperk_1})$
over $S$
such that the following holds. 
 \begin{itemize}
  \item $X_S$ is smooth over $S$,
	and $H_S$ is normal crossing relative to $S$
	with trivial monodromy.
	(See {\rm\S\ref{subsection;25.12.7.40}}.)
\item
$\Sigma_S$ is good on $(X_S,H_S)$ relative to $S$
in the sense of Definition {\rm\ref{df;25.12.7.2}}.
 \item If $\Sigma$ is Lagrangian,
       we can impose the additional condition
       that $\Sigma_S$ is also Lagrangian.
\hfill\qed
 \end{itemize}
\end{prop}

\subsection{Change of base fields}

\subsubsection{Statements}

Let $X_{\hyperk}$ be a smooth $\hyperk$-variety
with a normal crossing hypersurface $H_{\hyperk}$.
Let $H_{\hyperk}=
\bigcup_{i\in \Lambda}H_{\hyperk,i}$ denote the irreducible decomposition.
For any $I\subset\Lambda$,
we set
$H_{\hyperk,I}=\bigcap_{i\in I}H_{\hyperk,i}$
and
$H_{\hyperk,I}^{\circ}=
H_{\hyperk,I}\setminus
\bigcup_{i\not\in I}H_{\hyperk,i}$.

Let $\hyperk\subset\hyperk_1$ be an extension of
algebraically closed fields of characteristic $0$.
Let $(X_{\hyperk_1},H_{\hyperk_1})$
denote the pair of varieties induced by
$(X_{\hyperk},H_{\hyperk})$
and the extension $\hyperk\to\hyperk_1$.
The following proposition is obvious.
\begin{prop}
\label{prop;25.12.8.40}
A meromorphic cover
$\Sigma_{\hyperk}\subset T^{\ast}(X_{\hyperk}\setminus H_{\hyperk})$
over $(X_{\hyperk},H_{\hyperk})$
is etale locally good if and only if
the induced meromorphic cover
$\Sigma_{\hyperk}\times\hyperk_1$  over $(X_{\hyperk_1},H_{\hyperk_1})$
is etale locally good.
\end{prop}
\pf
The ``only if'' part is clear.
Suppose that $\Sigma_{\hyperk}\times\hyperk_1$ is etale locally good.
Let $S\subset\hyperk_1$ be a regular $\hyperk$-algebra
which is finitely generated over $\hyperk$.
Let $(X_S,H_S,\Sigma_S)$
be induced by
$(X_{\hyperk},H_{\hyperk},\Sigma_{\hyperk})$
by $\hyperk\to S$.
According to Lemma \ref{lem;25.12.7.12},
if $S$ is sufficiently large,
we may assume that
$\Sigma_S$ is etale locally good.
By taking the fiber over any closed point of $\Spec(S)$,
we obtain that
$\Sigma_{\hyperk}$ is etale locally good.
\hfill\qed

\vspace{.1in}

Let
$\Sigma_{\hyperk_1}\subset
T^{\ast}(X_{\hyperk_1}\setminus H_{\hyperk_1})$
be a formally good meromorphic cover over $(X_{\hyperk_1},H_{\hyperk_1})$.
Any closed point $x$ of $X_{\hyperk}$
induces a closed point $x_1$ of $X_{\hyperk_1}$.
We have the natural inclusion
\[
 \nbigohat^{(H_{\hyperk})}_{X_{\hyperk},x}(\ast H_{\hyperk})\big/
       \nbigohat^{(H_{\hyperk})}_{X_{\hyperk},x}
       \subset
\nbigohat^{(H_{\hyperk_1})}_{X_{\hyperk_1},x_1}
(\ast H_{\hyperk_1})/\nbigohat^{(H_{\hyperk_1})}_{X_{\hyperk_1},x_1}.
\]

We shall prove the following proposition in \S\ref{subsection;25.12.8.32}.

\begin{prop}
\label{prop;25.10.14.10}
We assume the following condition.
\begin{itemize}
 \item If $H_{\hyperk,I}\neq\emptyset$,
       any irreducible component of
       $H_{\hyperk,I}^{\circ}$
       contains a closed point $x$
       such that
\[
       \nbigi(\Sigma_{\hyperk_1},x_1)\subset
       \nbigohat^{(H_{\hyperk})}_{X_{\hyperk},x}(\ast H_{\hyperk})\big/
       \nbigohat^{(H_{\hyperk})}_{X_{\hyperk},x}.
\]
\end{itemize}
Then, there exists 
$\Sigma'_{\hyperk}\subset T^{\ast}(X_{\hyperk}\setminus H_{\hyperk})$
such that the following holds.
 \begin{itemize}
 \item Each irreducible component of $\Sigma'_{\hyperk}$ is
	proper and dominant over $X_{\hyperk}\setminus H_{\hyperk}$.
 \item $\Sigma'_{\hyperk}$ is good on $(X_{\hyperk},H_{\hyperk})$.
       In particular,
       $\Sigma'_{\hyperk}\times \hyperk_1$ is good on
       $(X_{\hyperk_1},H_{\hyperk_1})$.
 \item $\vecnbigi(\Sigma'_{\hyperk}\times \hyperk_1)
	     =\vecnbigi(\Sigma_{\hyperk_1})$.
 \end{itemize} 
If $\Sigma_{\hyperk_1}$ is Lagrangian,
we can impose the additional condition that
$\Sigma'_{\hyperk}$ is Lagrangian.
\end{prop}

\begin{rem}
$\Sigma_{\hyperk}$ does not necessarily induce $\Sigma_{\hyperk_1}$.
\hfill\qed
\end{rem}

\subsubsection{Injectivity}

Let $U_{\hyperk}=\Spec R_{\hyperk}$ be
an $n$-dimensional smooth connected affine $\hyperk$-variety
with a normal crossing hypersurface $H_{\hyperk}$ of $U_{\hyperk}$.
Let $x\in U_{\hyperk}$ be any closed point.
Assume that there exists an etale coordinate
$\psi_{\hyperk}:U_{\hyperk}\to \Spec\hyperk[z_1,\ldots,z_n]$
such that
(i) $H_{\hyperk}=\psi_{\hyperk}^{-1}\bigl(
\bigcup_{j=1}^{\ell}\{z_j=0\}
\bigr)$,
(ii)
$\psi_{\hyperk}^{-1}(\{z_j=0\})$
are connected,
(iii) $\psi_{\hyperk}(x)=(0,\ldots,0)$.
Let $\Rhat_{\hyperk,x}$ denote the completion of $R_{\hyperk}$ at $x$.
 
\begin{lem}
$R_{\hyperk}(\ast H_{\hyperk})/R_{\hyperk}
\to \Rhat_{\hyperk,x}(\ast H_{\hyperk})/\Rhat_{\hyperk,x}$
is injective.
\end{lem}
\pf
Let $a\in R_{\hyperk,x}(\ast H_{\hyperk})/R_{\hyperk,x}$
be any element.
It generates the submodule
$R_{\hyperk,x}a\subset R_{\hyperk}(\ast H_{\hyperk})/R_{\hyperk,x}$.
Because $\Rhat_{\hyperk,x}$
is faithfully flat over $R_{\hyperk,x}$,
we obtain the injection
$0\neq (R_{\hyperk,x}a)\otimes
\Rhat_{\hyperk,x}
\to \Rhat_{\hyperk,x}(\ast H_{\hyperk})/\Rhat_{\hyperk,x}$.
Hence,
$R_{\hyperk,x}(\ast H_{\hyperk})/R_{\hyperk,c}
\to \Rhat_{\hyperk,x}(\ast H_{\hyperk})/\Rhat_{\hyperk,x}$
is injective.
Let $a\in R_{\hyperk,x}\cap R_{\hyperk}(\ast H_{\hyperk})$.
There exists $f\in R_{\hyperk}$ such that
$a\in R_{\hyperk}(\ast f)\cap R_{\hyperk}(\ast H_{\hyperk})$
and that $f\not\in \gminim_x$.
Because $R_{\hyperk}$ is regular,
we have $R_{\hyperk}(\ast f)\cap R_{\hyperk}(\ast H)=R_{\hyperk}$.
\hfill\qed

\begin{cor}
Let $\nbigi_1,\nbigi_2\subset R_{\hyperk}(\ast H_{\hyperk})/R_{\hyperk}$
be good sets of irregular values. 
Let $\nbigi_{1|\xhat}$ and $\nbigi_{2|\xhat}$
denote the induced subsets of
$\Rhat_{x}(\ast H)/\Rhat_x$.
If $\nbigi_{1|\xhat}=\nbigi_{2|\xhat}$,
then $\nbigi_1=\nbigi_2$.
\hfill\qed
\end{cor}

\subsubsection{Good sets of irregular values}

Let $\hyperk_1$ be an algebraically closed field which contains $\hyperk$.
Let $U_{\hyperk_1}=\Spec R_{\hyperk_1}$ and
$H_{\hyperk_1}$ denote the $\hyperk_1$-varieties
induced by $U_{\hyperk}$ and $H_{\hyperk}$,
respectively.
Let 
$\nbigi_{\hyperk_1}\subset
R_{\hyperk_1}(\ast H_{\hyperk_1})/R_{\hyperk_1}$
be a good finite subset.
Let $x_{1}$ be the closed point of $U_{\hyperk_1}$
induced by $x$.

\begin{lem}
\label{lem;25.12.6.1}
Suppose that
$\nbigi_{\hyperk_1|\xhat_{1}}\subset
\Rhat_{\hyperk_1,x_{1}}(\ast H_{\hyperk_1})
/\Rhat_{\hyperk_1,x_{1}}$
is contained in
$\Rhat_{\hyperk,x}(\ast H_{\hyperk})/\Rhat_{\hyperk,x}$.
Then, 
$\nbigi_{\hyperk_1}$ is contained in
$R_{\hyperk}(\ast H)/R_{\hyperk}$. 
\end{lem}
\pf
Let $S\subset\hyperk_1$
be a regular subring finitely generated over $\hyperk$.
We set $U_S=\Spec(R_S):=U_{\hyperk}\times\Spec S$
and $H_S:=H_{\hyperk}\times\Spec S$.
We set $R_S=R\otimes S$.
If $S$ is sufficiently large,
there exists 
$\nbigi_{S}\subset R_{S}(\ast H_S)/R_S$
which induces 
$\nbigi_{\hyperk_1}$
by the extension $S\to \hyperk_1$.
For each closed point $s\in \Spec S$,
we obtain $\nbigi_{S,s}\subset
R_{\hyperk}(\ast H_{\hyperk})/R_{\hyperk}$
as the pull back of $\nbigi_S$
by $U_{\hyperk}\simeq U_S\times s\to U_S$.

Let $x_S:\Spec S\to \Spec R_S$ denote the section
induced by $x:\Spec\hyperk\to \Spec R_{\hyperk}$.
Let $I\subset R_S$ denote the ideal corresponding to
the morphism $R_S\to S$.
Let $\Rhat_{S,I}$ denote the completion of $R_S$
with respect to $I$.
Let $\nbigi_{S|\xhat_S}\subset
\Rhat_{S,I}(\ast H_S)/\Rhat_{S}$
denote the tuple induced by
$\nbigi_S$.
Because
\[
 \Rhat_{\hyperk,x}
 =\hyperk[\![z_1,\ldots,z_{\ell}]\!],\quad
 \Rhat_{\hyperk_1,x_1}
 =\hyperk_1[\![z_1,\ldots,z_{\ell}]\!],
 \quad\quad
 \Rhat_{S,x_S}
 =S[\![z_1,\ldots,z_{\ell}]\!],
\]
$\Rhat_{\hyperk,x}(\ast H_{\hyperk})/\Rhat_{\hyperk,x}
\to
\Rhat_{S,I}(\ast H_S)/\Rhat_{S,I}
\to
\Rhat_{\hyperk_1,x_1}(\ast H_{\hyperk_1})/\Rhat_{\hyperk_1,x_1}$
are injective.
Because 
$\nbigi_{\hyperk_1|\xhat_1}\subset
\Rhat_{\hyperk,x}(\ast H)/\Rhat_{\hyperk,x}$,
we obtain that
$\nbigi_{S|\xhat_S}
\subset
\Rhat_{\hyperk,x}(\ast H)/\Rhat_{\hyperk,x}$.

Let $s$ be any closed point of $S$.
By the previous consideration,
we obtain
$\nbigi_{S,s|\xhat_S}=\nbigi_{S|\xhat_S}=\nbigi_{S|\xhat}$.
We obtain
$\nbigi_{S,s|\xhat}=\nbigi_{\hyperk_1|\xhat_1}$,
and hence
$\nbigi_{S,s}=\nbigi_{\hyperk_1}$.
\hfill\qed

\subsubsection{Unramifiedly good meromorphic finite cover}

Let
$\Sigma_{\hyperk_1}\subset T^{\ast}(U_{\hyperk_1}\setminus H_{\hyperk_1})$
be a meromorphic finite cover over $(U_{\hyperk_1},H_{\hyperk_1})$.
Assume that $\Sigma_{\hyperk_1}$
is unramifiedly good on $(U_{\hyperk_1},H_{\hyperk_1})$.
We set $\nbigi_{\hyperk_1}:=\nbigi(\Sigma_{\hyperk_1})$.
(See \S\ref{subsection;25.12.7.31}.)
Suppose that 
$\nbigi_{\hyperk_1|\xhat_1}
\subset \Rhat_{\hyperk,x}(\ast H)/\Rhat_{\hyperk,x}$.
By Lemma \ref{lem;25.12.6.1},
$\nbigi_{\hyperk_1}$ is induced by
a good set
$\nbigi_{\hyperk}\subset R_{\hyperk}(\ast H)/R_{\hyperk}$.

\begin{lem}
\label{lem;25.10.14.1}
There exists 
an unramifiedly good meromorphic finite cover
$\Sigma_{\hyperk}$
over $(U_{\hyperk},H_{\hyperk})$
such that 
$\nbigi(\Sigma_{\hyperk})=\nbigi_{\hyperk}$. 
It implies 
$\vecnbigi(\Sigma_{\hyperk}\times\hyperk_1)
=\vecnbigi(\Sigma_{\hyperk_1})$.
\end{lem}
\pf
Let $S\subset\hyperk_1$ be a regular subring finitely generated over
$\hyperk$
such that there exists an unramifiedly good model
$\Sigma_{S}\subset T^{\ast}((U_S\setminus H_S)/S)$
of $\Sigma_{\hyperk_1}$ as in Proposition \ref{prop;25.12.7.20}.
Let $s\in \Spec(S)$ be any closed point.
We obtain
$\Sigma_{S,s}\subset T^{\ast}(U_{\hyperk}\setminus H_{\hyperk})$
which is unramifiedly good
such that
$\nbigi(\Sigma_{S,s}\times\hyperk_1)
=\vecnbigi(\Sigma_{\hyperk_1})$.
\hfill\qed

\subsubsection{Proof of Proposition \ref{prop;25.10.14.10}}
\label{subsection;25.12.8.32}

Let $(S,X_S,H_S,\Sigma_S)$
be a model as in Proposition \ref{prop;25.12.7.20}.
For any closed point $s$ of $\Spec(S)$,
we obtain a good meromorphic finite cover
$\Sigma_{S,s}$ over $(X_{\hyperk},H_{\hyperk})$.
By Lemma \ref{lem;25.10.14.1},
we obtain that
$\vecnbigi(\Sigma_{S,s}\times \hyperk_1)=\vecnbigi(\Sigma_{\hyperk_1})$.
\hfill\qed

\subsection{Resolutions of meromorphic Lagrangian covers}

\subsubsection{The absolute case}

Let $X$ be a projective normal $\hyperk$-variety
with a hypersurface $H$
such that $X\setminus H$ is smooth.

\begin{prop}
\label{prop;25.10.6.10}
For any meromorphic Lagrangian cover
$\Sigma\subset T^{\ast}(X\setminus H)$ over $(X,H)$,
there exists a birational projective morphism
$\rho:X'\to X$ 
such that the following holds.
\begin{itemize}
 \item $X'$ is smooth,
       and
       $H'=\rho^{-1}(H)$ is normal crossing.
 \item $\rho$ induces an isomorphism
       $X'\setminus H'\simeq X\setminus H$.
 \item $\rho^{\ast}\Sigma$ is etale locally good on $(X',H')$.
\end{itemize}
\end{prop}
\pf
It is enough to consider the case where
$X$ is smooth and $H$ is normal crossing.
It is enough to consider the case where
$\hyperk$ is an algebraic closure of
a field which is finitely generated over $\rnum$.
Therefore, it is enough to consider the case
where $\hyperk$ is a subfield of
the complex number field $\cnum$.
We set $X_{\cnum}=X\times_{\hyperk}\cnum$,
$H_{\cnum}=H\times_{\hyperk}\cnum$
and
$\Sigma_{\cnum}=\Sigma\times_{\hyperk}\cnum$.
By \cite[Theorem 15.2.7]{Mochizuki-MTM},
there exists
a projective birational morphism of smooth $\cnum$-varieties
$\rho_{\cnum}:X'_{\cnum}\to X_{\cnum}$
such that
(i) $H'_{\cnum}=\rho^{-1}(H_{\cnum})$ is normal crossing,
(ii) $\rho_{\cnum}$ induces an isomorphism
$X'_{\cnum}\setminus H'_{\cnum}
\simeq X_{\cnum}\setminus H_{\cnum}$,
(iii) $\Sigma_{\cnum}'=\rho_{\cnum}^{\ast}\Sigma_{\cnum}$
is formally good on $(X'_{\cnum},H'_{\cnum})$.
By Proposition \ref{prop;25.10.6.2},
it is etale locally good on $(X'_{\cnum},H'_{\cnum})$.

Let $(S,X'_S,H'_S,\Sigma'_S)$ be a model as in
Proposition \ref{prop;25.12.7.20}.
We may assume that
there exists a morphism
$\rho_S:X'_S\to X\times S$ over $\Spec S$
which induces $\rho_{\cnum}$.
Let $s$ be any closed point of $\Spec(S)$.
We obtain a birational projective morphism
$\rho_{S,s}:X'_{S,s}\to X$
such that
$H'_{S,s}=\rho_{S,s}^{-1}(H)$
and
$\Sigma'_{S,s}=\rho_{S,s}^{\ast}\Sigma$.
By Lemma \ref{lem;25.12.7.31} and Lemma \ref{lem;25.12.7.30},
we obtain that
$\Sigma'_{S,s}$ is good on $(X'_{S,s},H'_{S,s})$.
\hfill\qed

\subsubsection{The family case}

Let $X$ be a projective normal $\hyperk$-variety.
There exists a closed subset $Z\subset X$
with $\codim Z\geq 2$
such that $X\setminus Z$ is regular.
We set $X^{\circ}=X\setminus Z$.
There exists the reflexive $\nbigo_X$-coherent sheaf
$T^{\ast}X$ whose restriction to $X\setminus Z$
equals the cotangent bundle $T^{\ast}(X^{\circ})$.
Let $H$ be a hypersurface of $X$ such that
$Z\subset H$.
We set
$\nbige_{H,N}=T^{\ast}X(NH)$
and
$\nbige_{H,N}^{\circ}=T^{\ast}X^{\circ}(NH)$.
As in \S\ref{subsection;25.10.2.1},
we obtain the algebraic set 
$H^0(X^{\circ},S^r(\nbige^{\circ}_{H,N}))$
of sections $X^{\circ}\to S^r(\nbige_{H,N})$.
Let $\nbigg_{H,N}$ be a locally free sheaf on $X$
with an epimorphism $\nbigg_{H,N}^{\lor}\to \nbige^{\lor}_{H,N}$.
By Corollary \ref{cor;25.10.1.11}
and Lemma \ref{lem;25.10.2.2},
there exists a stratification
$H^0(X^{\circ},S^r(\nbige^{\circ}_{N,H}))
=\bigsqcup_{i\in\Lambda} T_i$
by locally Zariski closed subsets
such that the following holds.
\begin{itemize}
 \item $\ttP(\phi)$ $(\phi\in T_i)$ are constant.
 \item There exist closed subschemes
       $\Sigma(T_i)_{\red}\subset T_i\times\nbigg_{H,N}$
       such that
       $\Sigma(T_i)_{\red}\to T_i$ are flat,
       and that the fibers over
       $\phi\in T_i$ are $\Sigma(\phi_{\nbigg_{H,N}})_{\red}$.
 \item There exist torsion-free sheaves $\nbigf^k$
       on
       $\Sigma(T_i)_{\red}$
       as in 
       Lemma \ref{lem;25.10.17.10},
       which are flat over $T_i$.
\end{itemize}

Let $T_i^{\Lag}$ denote the set of
the closed points $\phi\in T_i$
such that
\[
\Sigma(\phi)_{U(\phi),\red}\times_{X^{\circ}}(U(\phi)\setminus H)
\subset T^{\ast}(U(\phi)\setminus H)
\]
are Lagrangian.
By Lemma \ref{lem;25.10.2.3},
we obtain the following lemma.

\begin{lem}
\label{lem;25.10.18.11}
$T_i^{\Lag}$ are Zariski closed subsets of $T_i$.
\hfill\qed
\end{lem}

Let $\Sigma(T_i^{\Lag})_{\red}$
denote the fiber product of $\Sigma(T_i)_{\red}$
and $T_i^{\Lag}$ over $T_i$.

\begin{prop}
\label{prop;25.10.18.10}
There exists a morphism
$\hyperk$-varieties
$\psi:\Ttilde_i^{\Lag}\to T_i^{\Lag}$
and 
       a projective birational morphism
\[
       \rho_{\Ttilde_i^{\Lag}}:
       \widehat{\Ttilde_i^{\Lag}\times X}\to \Ttilde_i^{\Lag}\times X
\]
such that the following holds.
\begin{itemize}
 \item[(i)] $\Ttilde_i^{\Lag}$ is smooth and affine over $\hyperk$,
	    and $\psi$ induces a surjection of
	    closed points.
 \item[(ii)]
	   $\rho_{\Ttilde_i^{\Lag}}$ induces
	   an isomorphism
      $\rho_{\Ttilde_i^{\Lag}}^{-1}(\Ttilde_i^{\Lag}\times(X\setminus H))
       \simeq
	   \Ttilde_i^{\Lag}\times(X\setminus H)$.
 \item[(iii)]
	    The induced morphism
	    $\widehat{\Ttilde_i^{\Lag}\times X}\to \Ttilde_i^{\Lag}$
	    is smooth,
	    and
	    $\widehat{\Ttilde_i^{\Lag}\times H}=
	    \rho_{\Ttilde_i^{\Lag}}^{-1}(\Ttilde_i^{\Lag}\times H)$
	    is a normal crossing hypersurface of
	     $\widehat{\Ttilde_i^{\Lag}\times X}\to \Ttilde_i^{\Lag}$
	     relative to $\Ttilde^{\Lag}$
	     without monodromy.
	     (See {\rm\S\ref{subsection;25.12.7.40}}.)
 \item[(iv)] Let $\Sigmahat(T^{\Lag}_i)_{\red}$
	     denote the fiber product of
	     $\Sigma(T_i^{\Lag})_{\red}$
	     and
	     $\widehat{\Ttilde_i^{\Lag}\times X}$
	     over
	     $\Ttilde_i^{\Lag}\times X$.
	     Then,
	     $\Sigmahat(T^{\Lag}_i)_{\red}$
	     is etale locally good
	     over
	     $\bigl(
	     \widehat{\Ttilde_i^{\Lag}\times X},
	     \widehat{\Ttilde_i^{\Lag}\times H}
	     \bigr)$.
\end{itemize} 
\end{prop}
\pf
Let $T_i'\subset T_i$ be any Zariski closed irreducible subset.
Let $R_i'$ be the ring of algebraic functions on $T_i'$.
Let $\hyperk'$ be an algebraic closure of
the fractional field of $R_i'$.
The universal section
$T_i'\times X^{\circ}\to T_i'\times S^r(\nbige^{\circ}_{H,N})$
induces the section
$\phi_{\hyperk'}:X^{\circ}_{\hyperk'}\to
S^r(\nbige^{\circ}_{H,N})_{\hyperk'}$.

We obtain
the meromorphic Lagrangian cover
$\Sigma(\phi_{\hyperk'})\subset
T^{\ast}(X_{\hyperk'}\setminus H_{\hyperk'})$
over $X_{\hyperk'}\setminus H_{\hyperk'}$.
By Proposition \ref{prop;25.10.6.10},
there exists a projective birational morphism
$\rho_{\hyperk'}:
\Xhat_{\hyperk'}\to X_{\hyperk'}$
such that the following holds.
\begin{itemize}
 \item $\Xhat_{\hyperk'}$ is smooth over $\hyperk'$,
       and
       $\Hhat_{\hyperk'}=\rho_{\hyperk'}^{-1}(H_{\hyperk'})$
       is a normal crossing hypersurface of
       $\Xhat_{\hyperk'}$.
 \item $\rho_{\hyperk'}$  induces an isomorphism
       $\Xhat_{\hyperk'}\setminus \Hhat_{\hyperk'}
       \simeq
       X_{\hyperk'}\setminus H_{\hyperk'}$.
 \item $\Sigmahat(\phi_{\hyperk'})=
       \rho_{\hyperk'}^{\ast}\Sigma(\phi_{\hyperk'})$
       is etale locally good on
       $(\Xhat_{\hyperk'},\Hhat_{\hyperk'})$.
\end{itemize}

There exists a model
Let $(S,\Xhat_S,\Hhat_S,\Sigmahat_S)$
of $(\Xhat_{\hyperk'},\Hhat_{\hyperk'},\Sigmahat(\phi_{\hyperk'}))$
as in Proposition \ref{prop;25.12.7.20}.
We may assume that $R_i'\subset S\subset \hyperk'$
and that
there exists a morphism
$\rho_{S}:\Xhat_S\to X\times\Spec(S)$
which induces $\rho_{\hyperk'}$.
The image of the morphism $\Spec(S)\to T_i'$ contains
a Zariski open subset $U_i'$.
We may assume that
the image is $U_i'$,
and that $\Spec(S)\to U_i'$ is smooth.

By this procedure and a Noetherian induction,
we can construct $\psi$ and $\rho_{\Ttilde_i^{\Lag}}$
as desired.
\hfill\qed

\begin{cor}
\label{cor;26.2.23.10}
There exist a smooth $\hyperk$-variety $S$,
a meromorphic Lagrangian cover
$\Sigma_{S}\subset \nbige^{\circ}_{H,N}\times S$,
and a smooth projective variety $\Xhat_{\nbigs}$ over $\nbigs$
with a morphism $\rho:\Xhat_{\nbigs}\to \nbigs\times X$
such that the following holds.
\begin{itemize}
 \item $\Hhat_{S}=\rho^{-1}(H\times S)$ is normal crossing relative to $S$.
 \item $\rho^{\ast}\Sigma_S$ is good on
       $(\Xhat_{S},\Hhat_S)$.
 \item For any meromorphic cover $\Sigma\subset \nbige_{H,N}$,
       there exists $a\in S$
       such that
       $\Sigma_{S}\times S_a=\Sigma$.
       Moreover,
       $\rho_s^{\ast}(\Sigma)$ is good on
       $(\Xhat_s,\Hhat_s)$.
       \hfill\qed
\end{itemize} 
\end{cor}

\subsection{Meromorphic Lagrangian irregularity}

\subsubsection{Preliminary}

Let $\Sigma_1,\Sigma_2$ be meromorphic Lagrangian covers
on $(X,H)$.
Suppose that
there exists a projective morphism
of smooth $\hyperk$-varieties
$\rho:X'\to X$
such that
(i) $H'=\rho^{-1}(H)$ is normal crossing,
(ii) $X'\setminus H'\simeq X\setminus H$,
(iii) $\rho^{\ast}\Sigma_i$ are
etale locally good on $(X',H')$,
(iv)
$\vecnbigi(\rho^{\ast}\Sigma_1)
=\vecnbigi(\rho^{\ast}\Sigma_2)$.

\begin{lem}
\label{lem;25.12.8.20}
Let $\rho_1:X_1\to X$
be a projective morphism of smooth $\hyperk$-varieties
such that
(i) $H_1=\rho_1^{-1}(H)$ is normal crossing,
(ii) $X_1\setminus H_1\simeq X\setminus H$,
(iii) $\rho_1^{\ast}\Sigma_1$ are
etale locally good on $(X_1,H_1)$.
Then, $\rho_1^{\ast}\Sigma_2$ is good on $(X_1,H_1)$,
and we have
$\vecnbigi(\rho_1^{\ast}\Sigma_1)
=\vecnbigi(\rho_1^{\ast}\Sigma_2)$.
\end{lem}
\pf
It follows from Proposition \ref{prop;25.10.14.20}.
\hfill\qed

\subsubsection{Meromorphic Lagrangian irregularity}

We introduce an equivalence relation
for meromorphic Lagrangian covers.

\begin{df}
For meromorphic Lagrangian covers $\Sigma_1,\Sigma_2$ on $(X,H)$,
we say $\Sigma_1\sim\Sigma_2$ 
if the following holds.
\begin{itemize}
 \item There exists a projective morphism
       of smooth $\hyperk$-varieties
       $\rho:X'\to X$
       such that
       (i) $H'=\rho^{-1}(H)$ is normal crossing,
       (ii) $X'\setminus H'\simeq X\setminus H$,
       (iii) $\rho^{\ast}\Sigma_i$ are
       etale locally good on $(X',H')$,
       (iv)
       $\vecnbigi(\rho^{\ast}\Sigma_1)
       =\vecnbigi(\rho^{\ast}\Sigma_2)$.
\end{itemize} 
An equivalence class is called
a meromorphic Lagrangian irregularity. 
An equivalence class of a meromorphic Lagrangian cover $\Sigma$
is denoted by $\vecI(\Sigma)$
and called the irregularity of $\Sigma$.
\hfill\qed
\end{df}

\begin{rem}
The degrees of $\Sigma_i\to X\setminus H$
are not necessarily equal
because we do not consider the multiplicity
in this paper. 
\hfill\qed
\end{rem}

We obtain the following proposition from 
Lemma \ref{lem;25.12.8.20}.
\begin{prop}
Let $\Sigma$ be a meromorphic Lagrangian cover over $(X,H)$.
Let $\vecI(\Sigma)$ denote the irregularity of $\Sigma$.
Let $\rho:X'\to X$ be any morphism of smooth $\hyperk$-varieties
such that
(i) $H'=\rho^{-1}(H)$ is normal crossing,
(ii) $X'\setminus H'\simeq X\setminus H$,
(iii) $\rho^{\ast}\Sigma$ is etale locally good on $(X',H')$.
Then, for any $\Sigma_1\in \vecI(\Sigma)$,
$\rho^{\ast}(\Sigma_1)$ is etale locally good on $(X',H')$,
and $\vecnbigi(\rho^{\ast}\Sigma_1)=\vecnbigi(\rho^{\ast}\Sigma)$ holds.
\hfill\qed 
\end{prop}

\begin{df}
A meromorphic Lagrangian regularity $\vecI$
is called good (resp. logarithmic, unramifiedly good)
on $(X,H)$
if any $\Sigma\in\vecI$ is 
good (resp. logarithmic, unramifiedly good)
on $(X,H)$.
It is equivalent to the condition 
that there exists $\Sigma\in\vecI$
which is
good (resp. logarithmic, unramifiedly good)
on $(X,H)$.
\hfill\qed
\end{df}

\subsubsection{Pull back and push-forward via birational projective morphisms}

Let $\varphi:X'\to X$ be a birational projective morphism 
such that $H'=\varphi^{-1}(H)$ is normal crossing
and that $X'\setminus H'\simeq X\setminus H$.
We have the natural isomorphism
$T^{\ast}(X'\setminus H')\simeq T^{\ast}(X\setminus H)$.

\begin{lem}\mbox{{}}
\begin{itemize}
 \item Let $\vecI$ be a meromorphic Lagrangian irregularity over $(X,H)$.
For any $\Sigma_1,\Sigma_2\in\vecI$,
we have $\vecI(\varphi^{\ast}\Sigma_1)=\vecI(\varphi^{\ast}\Sigma_2)$.
 \item Let $\vecI'$ be a meromorphic Lagrangian irregularity over $(X',H')$.
       For any $\Sigma'_1,\Sigma'_2\in\vecI'$,
       we have $\vecI(\varphi_{\ast}\Sigma_1')
       =\vecI(\varphi_{\ast}\Sigma_2')$.
\end{itemize}
\end{lem}
\pf
The first claim is clear.
The second claim follows from Proposition \ref{prop;25.10.14.20}.
\hfill\qed

\vspace{.1in}

Hence, we define
$\varphi^{\ast}(\vecI):=\vecI(\varphi^{\ast}\Sigma)$
for any $\Sigma\in\vecI$.
We also define
$\varphi_{\ast}(\vecI'):=\vecI(\varphi_{\ast}\Sigma')$
for any $\Sigma'\in\vecI'$.

\subsubsection{Comparison}

Let $S$, $X$, $H$ be as in \S\ref{subsection;25.12.7.40}.
Let $\Sigma^{(i)}\subset T^{\ast}((X\setminus H)/S)$ $(i=1,2)$ be 
Lagrangian covers of $X\setminus H$ which are etale locally good on $(X,H)$.
For each closed point $z$ of $\Spec(S)$,
we obtain good meromorphic Lagrangian covers $\Sigma^{(i)}_z$
on $(X_z,H_z)$.
We obtain the following proposition from 
Proposition \ref{prop;26.1.21.3}.
\begin{prop}
\label{prop;26.1.21.4}
There exists a closed subset $Z\subset \Spec(S)$
such that
$\vecI(\Sigma^{(1)}_z)
=\vecI(\Sigma^{(2)}_z)$
if and only if
$z\in Z$. 
\hfill\qed
\end{prop}

\section{Complex analytic meromorphic flat bundles}

\subsection{Good meromorphic flat bundles}

\subsubsection{Formal good meromorphic flat bundles}

Let $X$ be any $n$-dimensional complex manifold.
Let $H$ be a simple normal crossing hypersurface of $X$.
Let $\nbigohat_{X,x}$ denote the completion
of the local ring $\nbigo_{X,x}$
with respect to the maximal ideal $\gminim_x$.
There exists a parameter system $(x_1,\ldots,x_n)$ of $\gminim_x$
such that the ideal of $H$ at $x$ is generated by
$\prod_{i=1}^{\ell}x_i$.
For any positive integer $e$,
we set
$\nbigohat^{(e)}_{X,x}
=\nbigohat_{X,x}[x_i^{1/e}\,|\,i=1,\ldots,\ell]$.
For any $\nbigohat_{X,x}$-module $\nbigf$,
we set
$\nbigf(\ast H)=\nbigf\otimes_{\nbigo_{X,x}}\nbigo_{X}(\ast H)_x$
and
$\varphi_e^{\ast}(\nbigf)
=\nbigf\otimes\nbigohat^{(e)}_{X,x}$.
For any free $\nbigohat_{X,x}(\ast H)$-module $\nbigf$,
a lattice of $\nbigf$ means
a free $\nbigo_{X,x}$-submodule $\nbigf_0\subset\nbigf$
such that $\nbigf_0(\ast H)=\nbigf$.

Let $\nbigv$ be a free $\nbigohat_{X,x}(\ast H)$-module of finite rank
with an integrable connection $\nabla:\nbigv\to\nbigv\otimes\Omega^1_X$.
Any $\varphi_e^{\ast}(\nbigv)$ is equipped with
the induced integrable connection $\nabla$.
We recall some definitions.
\begin{df}
\mbox{{}}\label{df;26.1.21.10}
\begin{itemize}
 \item 
$(\nbigv,\nabla)$ is called regular
if there exists a free $\nbigohat_{X,x}$-module $\nbigv_0\subset\nbigv$
such that 
(i)  $\nbigv_0(\ast H)=\nbigv$,
(ii) $\nabla(\nbigv_0)\subset\nbigv_0\otimes\Omega^1(\log H)$.
\item $(\nbigv,\nabla)$ is called unramifiedly good
      if there exist a good set of irregular values
      $\nbigi\subset\nbigohat_{X,x}(\ast H)/\nbigo_{X,x}$
      and a decomposition
\begin{equation}
\label{eq;25.10.10.2}
      (\nbigv,\nabla)
      =\bigoplus_{\gminia\in\nbigi}
      (\nbigv_{\gminia},\nabla_{\gminia})
\end{equation}
      such that
      $\nabla_{\gminia}-d\gminiatilde\id_{\nbigv_{\gminia}}$
      are regular.
      Here, $\gminiatilde\in\nbigohat_{X,x}$
      denote lifts of $\gminia$.
      We allow $\nbigv_{\gminia}=0$
      in {\rm(\ref{eq;25.10.10.2})}.
      The set
      $\nbigi(\nbigv)=
      \bigl\{\gminia\in\nbigi\,|\,\nbigv_{\gminia}\neq 0\bigr\}$
      is well defined.
 \item $(\nbigv,\nabla)$ is called good
       if there exists $e\in\seisuu_{>0}$
       such that
       $\varphi_e^{\ast}(\nbigv,\nabla)$
       is unramifiedly good.
       We obtain the set
       $\nbigi(\nbigv)\subset
       \nbigohat^{(e)}_{X,x}(\ast h)/\nbigohat^{(e)}_{X,x}$.
      \hfill\qed
\end{itemize}
\end{df}

\subsubsection{The associated Deligne-Malgrange good filtered bundles
in the formal case}

Let $H_{x,i}=\{x_i=0\}$.
For any $\nbigohat_{X,x}$-module $\nbigf$,
we set $\nbigf_{|H_{x,i}}=\nbigf/x_i\nbigf$.

Suppose that $(\nbigv,\nabla)$ is regular.
For each lattice $\nbigv_0\subset\nbigv$ such that
$\nabla(\nbigv_0)\subset\nbigv_0\otimes\Omega^1_X(\log H)$,
we obtain the residue endomorphism
$\Res_{i}(\nabla)$ on $\nbigv_{|H_{x,i}}$
by taking the residue with respect to $dx_i/x_i$.

\begin{lem}
For any $\veca=(a_i)\in\real^{\ell}$,
there exists a unique lattice 
$\nbigp^{\DM}_{\veca}\nbigv$
such that the following holds.
\begin{itemize}
 \item Any eigenvalue $\alpha$
       of $\Res_{i}(\nabla)$
       satisfies
       $a_i-1<-\Re(\alpha)\leq a_i$.
\end{itemize}
We obtain the filtered bundle
 $\nbigp_{\ast}^{\DM}(\nbigv)
 =(\nbigp_{\veca}^{\DM}(\nbigv)\,|\,\veca\in\real^{\ell})$
over $\nbigv$.
\hfill\qed
\end{lem}

Suppose that $(\nbigv,\nabla)$ is unramifiedly good
with the decomposition (\ref{eq;25.10.10.2}).
We obtain
the filtered bundles
$\nbigp^{\DM}_{\ast}(\nbigv_{\gminia})$
over $\nbigv_{\gminia}$
by applying the previous construction
to $(\nbigv_{\gminia},\nabla_{\gminia}-d\gminia\id_{\nbigv_{\gminia}})$
which are regular.
We obtain
\[
\nbigp^{\DM}_{\ast}(\nbigv)
=\bigoplus\nbigp^{\DM}_{\ast}(\nbigv_{\gminia}).
\]

Let us consider the case
$(\nbigv,\nabla)$ is good but not necessarily unramifiedly good.
By applying the previous construction
to the unramifiedly good 
$\varphi_e^{\ast}(\nbigv,\nabla)$,
we obtain the filtered bundle
$\nbigp_{\ast}^{\DM}(\varphi_e^{\ast}\nbigv)$
over $\varphi_e^{\ast}(\nbigv)$.
Let $\Gal_e$ denote the Galois group
of the extension
$\nbigohat^{(e)}_{X,x}/\nbigohat_{X,x}$.
We have the natural action of $\Gal_e$
on $\varphi_e^{\ast}(\nbigv)$.
It induces the action on
the filtered bundles
$\nbigp_{\ast}(\varphi^{\ast}(\nbigv))$.
As the descent,
we obtain the filtered bundle
$\nbigp^{\DM}_{\ast}(\nbigv)$ over $\nbigv$.

In particular, $\nbigp^{\DM}_0(\nbige)$
is called the good Deligne-Malgrange lattice.

\subsubsection{Good meromorphic flat bundles in the complex analytic case}

Let $X$ be any complex manifold.
Let $H$ be a simple normal crossing hypersurface of $X$.
For each $x\in H$,
there exists a holomorphic coordinate system
$(x_1,\ldots,x_n)$ such that
$H=\bigcup_{i=1}^{\ell}\{x_i=0\}$.
We set
$\nbigo_{X,x}^{(e)}
=\nbigo_{X,x}[x_i^{1/e}\,|\,i=1,\ldots,\ell]$.

Let $(\nbigv,\nabla)$ be a meromorphic flat bundle on $(X,H)$.
For any $x\in X$,
we obtain
$\nbigv_{|\xhat}=\nbigv_x\otimes_{\nbigo_{X,x}}\nbigohat_{X,x}$
with the induced integrable connection.
We recall some definitions.
\begin{df}
\label{df;25.10.10.3}
$(\nbigv,\nabla)$ is called regular
(resp. unramifiedly good, good)
if $(\nbigv_{|\xhat},\nabla)$
are regular (resp. unramifiedly good, good)
for any $x\in X$.
\hfill\qed
\end{df}

\begin{prop}[\mbox{\cite{Mochizuki-Stokes-good}}]
If $(\nbigv,\nabla)$ is good,
there exists
a filtered bundle $\nbigp^{\DM}_{\ast}(\nbigv_0)$
over $\nbigv_0$
such that
$\nbigp^{\DM}_{\ast}(\nbigv_0)_{|\xhat}
=\nbigp^{\DM}_{\ast}(\nbigv_{0|\xhat})$.
In particular,
there exists a locally free $\nbigo_X$-submodule $\nbigv_0$
such that $\nbigv_{0|\xhat}$
are good Deligne-Malgrange lattice of
$\nbigv_{|\xhat}$.
\hfill\qed
\end{prop}

If $(\nbigv,\nabla)$ is good,
we obtain the set
$\nbigi(\nbigv,x)
\subset
\nbigohat^{(e)}_{X}(\ast H)_x/\nbigohat^{(e)}_{X,x}$
for any $x\in H$
as the index set of the decomposition of
$\varphi_e^{\ast}(\nbigv,\nabla)_{|\xhat}$.
(See Definition \ref{df;26.1.21.10}.)

\begin{prop}
\mbox{{}}
\begin{itemize}
 \item
We have
$\nbigi(\nbigv,x)
\subset
\nbigo^{(e)}_{X}(\ast H)_x/\nbigo^{(e)}_{X,x}$
for any $x\in H$.
 \item
      $\vecnbigi(\nbigv)=(\nbigi(\nbigv,x)\,|\,x\in H)$
      is a good system of ramified irregular values.
\end{itemize}
\end{prop}
\pf
In the unramified case,
it is proved in \cite[Proposition 2.19]{Mochizuki-Stokes-good}.
The ramified case is reduced to the unramified case.
\hfill\qed

\subsection{Meromorphic Lagrangian covers in the complex analytic case}

\subsubsection{Meromorphic finite covers}

Let $X$ be a complex manifold.
Let $H$ be a complex analytic closed subset of $X$.
We set $X^{\circ}=X\setminus H$.
Let $\Sigma\subset T^{\ast}X^{\circ}$ be a closed subset.
Let $\proj_{T^{\ast}X}$ denote
the projective completion of $T^{\ast}X$.

\begin{df}
$\Sigma$ is called a meromorphic cover over $(X,H)$
if the closure of $\Sigma$ in $\proj_{T^{\ast}X}$
is a complex analytic closed subset. 
If the smooth part of $\Sigma$ is Lagrangian with respect to 
the natural symplectic structure of $T^{\ast}(X\setminus H)$,
it is called a meromorphic Lagrangian cover.
\hfill\qed
\end{df}

\subsubsection{Some conditions in the case where $H$ is normal crossing}

Let us consider the case where $H$ is a normal crossing hypersurface.
\begin{df}
$\Sigma$ is called logarithmic
if the closure of $\Sigma$ in $T^{\ast}X(\log H)$
is proper over $X$.
\hfill\qed
\end{df}

For any section $\tau$ of $T^{\ast}(X\setminus H)$,
let $\kappa_{\tau}$ denote the automorphism of
$T^{\ast}(X\setminus H)$
obtained as the addition of $\tau$.

\begin{df}
\mbox{{}}
\begin{itemize}
 \item 
Let $x\in H$.
We say that $\Sigma$ is unramifiedly good at $x$
if there exists
a good set of irregular values
$\nbigi(\Sigma,x)\subset \nbigo_{X}(\ast H)_x/\nbigo_{X,x}$,
a neighbourhood $X_x$,
and a decomposition
\[
 \Sigma_{|X_x\setminus H}
 =\bigsqcup_{\gminia\in\nbigi(\Sigma,x)}
 \Sigma_{\gminia}
\] 
such that
$\kappa_{d\gminiatilde}^{\ast}\Sigma_{\gminia}$
is logarithmic on $(X_x,H\cap X_x)$.
 \item We say that $\Sigma$ is unramifiedly good on $(X,H)$
       if it is unramifiedly good at any $x\in H$.
       \hfill\qed
\end{itemize}
\end{df}

Let $x\in H$.
There exists a holomorphic coordinate neighbourhood
$(X_x,z_1,\ldots,z_n)$ around $x$
such that
$H\cap X_x=\bigcup_{j=1}^{\ell}\{z_j=0\}$.
By the coordinate, we may regard $X_x\subset\cnum^n$.
Let $\varphi_{e}:\cnum^n\to\cnum^n$
be defined by
$\varphi_{e}(\zeta_1,\ldots,\zeta_n)
=(\zeta_1^e,\ldots,\zeta_{\ell}^e,\zeta_{\ell+1},\ldots,\zeta_n)$.
We set
$X_x^{(e)}=\varphi_{x,e}^{-1}(X_x)$
and $H_x^{(e)}=\varphi_{x,e}^{-1}(H\cap X_x)$.
The induced map
$X_x^{(e)}\to X_x$ is also denoted by $\varphi_{x,e}$.

\begin{df}\mbox{{}}
\begin{itemize}
 \item $\Sigma$ is called good at $x$
       if there exist $e\in\seisuu_{>0}$
       and
       a neighbourhood
       $X_x$ of $x$
       such that $\varphi_{x,e}^{\ast}(\Sigma)$
       is unramifiedly good
       on $(X_x^{(e)},H_x^{(e)})$.
       In this case,
       we obtain
       $\nbigi(\Sigma,x)\subset
       \nbigo^{(e)}_X(\ast H)_x/\nbigo^{(e)}_{X,x}$.
 \item $\Sigma$ is called good on $(X,H)$
       if it is good at any $x\in H$.
       In this case,
       we obtain the tuple
       $\vecnbigi(\Sigma)
       =\bigl(
       \nbigi(\Sigma,x)\,|\,x\in H
       \bigr)$.
       \hfill\qed
\end{itemize}
\end{df}

\subsubsection{Resolution}

Let us consider the case where $H$ is not necessarily
a normal crossing hypersurface.
The following proposition is already mentioned in 
Proposition \ref{prop;25.10.6.10}.

\begin{prop}[\mbox{\cite[Theorem 15.2.7]{Mochizuki-MTM}}]
\label{prop;25.12.8.10}
There exists a projective morphism of complex manifolds $\rho:X'\to X$
such that
(i) $H'=\rho^{-1}(H)$ is a normal crossing hypersurface,
(ii) $\rho$ induces an isomorphism $X'\setminus H'\simeq X\setminus H$,
(iii) $\rho^{\ast}(\Sigma)$ is good on $(X',H')$. 
\end{prop}
\pf
It is enough to consider the case where $H$ is a normal crossing hypersurface,
and it is studied in \cite[Theorem 15.2.7]{Mochizuki-MTM}.
\hfill\qed

\subsubsection{Residues}

Suppose that $\Sigma$ is a good Lagrangian cover.
Let $x\in H$.
Let $\varphi_{x,e}:X^{(e)}_x\to X_x$ be as before.
We have the decomposition
\[
 \varphi_{x,e}^{\ast}(\Sigma)
 =\bigsqcup_{\gminia\in\nbigi(\Sigma,x)}
 \Sigma_{\gminia}
\]
such that
$\kappa_{d\gminia}^{\ast}\Sigma_{\gminia}$
is logarithmic.
Moreover,
if $X_x$ is sufficiently small,
there exist
$\EE(\Sigma,x,\gminia)\subset\cnum^{\ell}$
and a decomposition
\[
  \Sigma_{\gminia}
 =\bigsqcup_{\vecalpha\in\EE(\Sigma,x,\gminia)}
 \Sigma_{\gminia,\vecalpha}
\]
such that the following holds.
\begin{itemize}
 \item For $\vecalpha\in\EE(\Sigma,x,\gminia)$,
       we set
       $\tau(\vecalpha)
       =\sum_{i=1}^{\ell}\alpha_i\varphi_{x,e}^{\ast}(dz_i/z_i)$.
       Let
       $\overline{
       \kappa_{d\gminia+\tau(\vecalpha)}^{\ast}\Sigma_{\gminia,\vecalpha}}$
       denote the closure of
       $\kappa_{d\gminia+\tau(\vecalpha)}^{\ast}\Sigma_{\gminia,\vecalpha}$
       in
       $T^{\ast}X_x^{(e)}(\log H_x^{(e)})$.
       Let $\pi_{x,e}:T^{\ast}X^{(e)}_x\to X^{(e)}_x$
       denote the projection.
       Then,
       $\overline{
       \kappa_{d\gminia+\tau(\vecalpha)}^{\ast}\Sigma_{\gminia,\vecalpha}}
       \cap
       \pi_{x,e}^{-1}(H^{(e)}_{x})$
       is contained in the $0$-section.     
\end{itemize}

\subsection{Wild harmonic bundles}

\subsubsection{Harmonic bundles}

Let $(E,\theta)$ be a Higgs bundle on a complex manifold $Y$.
Let $h$ be a Hermitian metric of $E$.
We obtain the Chern connection $\nabla_h$,
and the adjoint $\theta^{\dagger}_h$ of $\theta$.
If the connection $\DD_h^1=\nabla_h+\theta+\theta^{\dagger}_h$ is flat,
$(E,\theta,h)$ is called a harmonic bundle.
For the Higgs bundle $(E,\theta)$,
we have the corresponding $\nbigo_{T^{\ast}Y}$-module $\Etilde$,
whose supports are finite over $Y$.
The support $\Sigma_{\theta}$ is called the spectral cover.
We recall the following lemma,
which follows from Gabber's theorem \cite{Gabber}.
\begin{lem}
If $(E,\theta)$ underlies a harmonic bundle $(E,\theta,h)$, 
$\Sigma_{\theta}$ is Lagrangian.
\hfill\qed
\end{lem}

\subsubsection{Wild harmonic bundles and good wild harmonic bundles}

Let $X$ be a complex manifold
with a closed complex analytic subspace $H$.
Let $(E,\theta,h)$ be a harmonic bundle on $X\setminus H$.
\begin{df}
$(E,\theta,h)$ is called wild on $(X,H)$
if $\Sigma_{\theta}$ is a meromorphic cover over $(X,H)$.
\hfill\qed
\end{df}

\begin{df}
Suppose that $H$ is normal crossing hypersurface.
\begin{itemize}
 \item The harmonic bundle $(E,\theta,h)$ is called tame
       if $\Sigma_{\theta}$ is logarithmic.
 \item  The harmonic bundle $(E,\theta,h)$ is called good wild
	(resp. unramifiedly good wild) on $(X,H)$
	if $\Sigma_{\theta}$ is good
	(resp. unramifiedly good) on $(X,H)$.
\hfill\qed	
\end{itemize}
\end{df}

\begin{df}
The harmonic bundle $(E,\theta,h)$ is called
$\sqrt{-1}\real$-wild harmonic bundle
if the following holds. 
\begin{itemize}
 \item Let $\rho:X'\to X$ be a resolution for
       $\Sigma_{\theta}$
       as in Proposition {\rm\ref{prop;25.12.8.10}}.
       Then, for any $x\in H'$ and
       $\gminia\in \nbigi(\rho^{\ast}\Sigma)$,
       $\EE(\rho^{\ast}\Sigma,x',\gminia)$
       consists of
       tuples of purely imaginary numbers.
       \hfill\qed
\end{itemize}
\end{df}

\subsubsection{The associated meromorphic integrable connections}

Let $X$ be a complex manifold
with a hypersurface $H$.
Let $\iota:X\setminus H\to X$ denote the inclusion.

Let $(E,\theta,h)$ be a wild harmonic bundle on $(X,H)$.
From a harmonic bundle on $(E,\delbar_E,\theta,h)$ on $X\setminus H$,
we have the flat bundle $(E,\DD^1_h)$ on $X\setminus H$.
We have the underlying holomorphic vector bundle
$\nbige^1=(E,\delbar_E+\theta^{\dagger}_h)$.
We regard it as an $\nbigo_{X\setminus H}$-module.
We obtain the sheaf
$\iota_{\ast}(\nbige^1)$ on $X$.

For any $x\in H$,
let $X_x$ be a neighbourhood of $x$
with a holomorphic function $f\in\nbigo(X_x)$
such that $H=f^{-1}(0)$.
We obtain an $\nbigo_X(\ast H)_x$-module
\[
 (\nbigp^h\nbige^1)_x:=
 \Bigl\{
 s\in (\iota_{\ast}\nbige^1)_x\,\Big|\,
 |s|_h=O\Bigl(|f|^{-C}\Bigr)
 \,\,\,
 \exists C>0
 \Bigr\}.
\]
We obtain the $\nbigo_X(\ast H)$-module
$\nbigp^h\nbige^1\subset\iota_{\ast}\nbige^1$
whose stalk at $x$ is $(\nbigp^h\nbige^1)_x$.

\begin{prop}
\mbox{{}}
\begin{itemize}
 \item 
$\nbigp^h\nbige^1$ is a coherent $\nbigo_X(\ast H)$-module,
       and $\DD^1_h$ is an integrable connection
       of $\nbigp^h\nbige^1$.
 \item
      Suppose that $H$ is normal crossing.
      Then, $(E,\theta,h)$ is good on $(X,H)$
if and only if
$(\nbigp^h\nbige^1,\DD^1_h)$ is a good meromorphic flat bundle on $(X,H)$.
In the case, we have
$\nbigi(\nbigp^h\nbige^1,\DD^1_h,x)
 =\bigl\{
 2\gminia\,\big|\,\gminia\in\nbigi(\Sigma_{\theta},x)
 \bigr\}$.
\hfill\qed
\end{itemize}
\end{prop}

\subsubsection{The correspondence in the projective case}

Suppose that $X$ is projective.

\begin{thm}
\label{thm;25.12.9.1}
For any semisimple meromorphic flat bundle $(V,\nabla)$ on $(X,H)$,
there exists a $\sqrt{-1}\real$-wild harmonic bundle $(E,\theta,h)$ such that 
$(\nbigp^h\nbige^1,\DD^1_h)\simeq (V,\nabla)$,
which is unique up to isomorphism.
\hfill\qed
\end{thm}

\subsection{The associated meromorphic Lagrangian irregularity
in the algebraic case}
\label{subsection;25.12.8.30}

Let $X$ be a normal complex projective variety.
Let $H$ be any hypersurface of $X$ such that $X\setminus H$ is smooth.
Let $(V,\nabla)$ be any meromorphic flat connection on $(X,H)$.
There exists a Jordan-H\"older filtration,
namely,
each graded piece is simple.
By Theorem \ref{thm;25.12.9.1},
we have the $\sqrt{-1}\real$-wild harmonic bundle $(E,\theta)$
corresponding to
the associated graded meromorphic flat bundle.
Let $\Psi_2:T^{\ast}(X\setminus H)\to T^{\ast}(X\setminus H)$
be the automorphism obtained as $v\mapsto 2v$
in the fiber direction.
We set
$\Sigma(V,\nabla):=\Psi_2(\Sigma_{\theta})$.

\begin{cor}
\label{cor;25.12.8.31}
Let $\varphi:X'\to X$ be a projective morphism
such that $H'=\varphi^{-1}(H)$ is normal crossing.
Then, 
$\varphi^{\ast}(V,\nabla)$ is good on $(X',H')$
if and only if
$\varphi^{\ast}\Sigma(V,\nabla)$ is good on $(X',H')$.
In the case,
we have
$\vecnbigi(\varphi^{\ast}(V,\nabla))
=\vecnbigi(\varphi^{\ast}\Sigma(V,\nabla))$. 
 \hfill\qed
\end{cor}

\subsection{Local systems with Stokes structure}

\subsubsection{Partial orders on the index sets}

Let $\varpi:\Xtilde(H)\to X$ denote the oriented real blow up
along $H$.
Let $\nbigo_{\Xtilde(H)}$ denote the sheaf of
holomorphic functions on $\Xtilde(H)$.
For any $\nbigotilde_X$-module $\nbigf$,
we set
$\nbigf(\ast H)=
\nbigf\otimes_{\varpi^{-1}(\nbigo_X)}\varpi^{-1}(\nbigo_X(\ast H))$.

Let $x\in H$.
We have a holomorphic coordinate neighbourhood
$(X_x;x_1,\ldots,x_n)$ around $x$
such that $H\cap X_x=\bigcup_{i=1}^{\ell}\{x_i=0\}$.
We define
$\nbigotilde^{(e)}_{\Xtilde_x(H_x)}
=\nbigotilde_{\Xtilde_x(H_x)}[x_i^{1/e}\,|\,i=1,\ldots,\ell]$.

Let $\nbigi_x\subset \nbigo^{(e)}_{X,x}(\ast H)/\nbigo^{(e)}_{X,x}$
be a good set of ramified irregular values.
By choosing branches $x_i^{1/e}$ $(i=1,\ldots,\ell)$,
we obtain
$\nbigi\subset
\nbigo^{(e)}_{\Xtilde_x(H)}(\ast H)_x
\big/\nbigo^{(e)}_{\Xtilde_x,x}$.
The image is independent of the choices
of $x^{1/e}_i$.

For each
$\gminia\in\nbigi$,
we choose a lift
$\gminiatilde\in \nbigo_{\Xtilde_x(H)}(\ast H)_y$.
We define $\gminia\leq_y\gminib$
if there exists a neighbourhood $\nbigu_y$ of $y$
in $\Xtilde(H)$ such that
$-\Re(\gminia)\leq -\Re(\gminib)$
on $\nbigu_y\setminus\varpi^{-1}(H)$.

Let $\vecnbigi=(\nbigi_x\,|\,x\in H)$ be
a good system of ramified irregular values.
If $x'$ is sufficiently close to $x$,
we have the natural surjection
$\nbigi_{x}\to\nbigi_{x'}$.
If $y'\in\varpi^{-1}(x')$ is sufficiently close to
$y\in \varpi^{-1}(x)$,
the map
$\varphi_{x',x}:\nbigi_x\to\nbigi_{x'}$
preserves the partially orders
$\leq_x$ and $\leq_{x'}$,
i.e.,
if $\gminia\leq_y\gminib$,
then $\varphi_{x',x}(\gminia)\leq_{y'}\varphi_{x',x}(\gminib)$.

\subsubsection{Stokes structures}

Let $L$ be a local system on $X\setminus H$.
Let $\Ltilde$ denote the local system on $\Xtilde(H)$
induced by $L$.
A Stokes structure of $\Ltilde$
is a tuple of filtrations
$\vecnbigf=\bigl(
\nbigf^y(\Ltilde_y)\,\big|\,
y\in \varpi^{-1}(H)
\bigr)$
indexed by
the partially ordered set $(\nbigi_{\varpi(y)},\leq_y)$
satisfying the following conditions.
\begin{itemize}
 \item There exists a decomposition
       $\Ltilde_y=
       \bigoplus_{\gminia\in \nbigi_{\varpi(y)}}
       G_{\gminia}$
       such that
       $\nbigf^y_{\gminia}(\Ltilde_y)=
       \bigoplus_{\gminib\leq_y\gminia}
       G_{\gminib}$.
 \item If $y'$ is sufficiently close to $y$,
       we have
\[
       \nbigf^{y'}_{\gminib}(\Ltilde_{y'})
       =\sum_{\varphi_{x',x}(\gminia)\leq \gminib}
       \nbigf^{y}_{\gminia}(\Ltilde_{y})
\]
       under the isomorphism
       $\Ltilde_{y}\simeq\Ltilde_{y'}$.
       Here, $x'=\varpi(y')$ and $x=\varpi(y)$.
\end{itemize}

\subsubsection{Deformation (1)}
\label{subsection;25.10.10.40}

Let $p:X\to B$ be a smooth fibration of
complex manifolds.
Let $H$ be a normal crossing hypersurface of $X$.
Let $b\in B$.
The fibers over $b$ are denoted by $X_b$ and $H_b$.
We assume the following conditions.
\begin{itemize}
 \item $B$ is simply connected.
 \item There exists $b_0\in B$
       and a diffeomorphism
       $(X,H)\simeq B\times (X_{b_0},H_{b_0})$
       such that
       the restriction to $(X_{b_0},H_{b_0})$
       is the identity.
\end{itemize}
Let $\vecnbigi$ be a good system of ramified irregular values
on $(X,H)$.
We obtain the induced good system of ramified irregular values
$\vecnbigi_{b}$ on $(X_b,H_b)$ by the restriction.
We have the natural functor
\begin{equation}
\label{eq;25.10.10.10}
 \Mero(X,H,\vecnbigi)
 \to
 \Mero(X_b,H_b,\vecnbigi_b).
\end{equation}

\begin{prop}
\label{prop;25.10.10.30}
The functor {\rm(\ref{eq;25.10.10.10})}
is an equivalence.
\end{prop}
\pf
If $\vecnbigi$ is unramified,
it is explained in
\cite[Theorem 3.9]{Mochizuki-Stokes-good}.
Let us study the case where
there exists a ramified Galois covering
$\varphi:(X',H')\to (X,H)$ over $B$
such that
$\varphi^{\ast}\vecnbigi$ is unramified.
We have the equivalence
$\Locst(X',H',\varphi^{\ast}\vecnbigi)
\simeq
\Locst(X'_b,H'_b,\varphi^{\ast}\vecnbigi_b)$.
It is easy to obtain the equivalence for
$\Gal$-equivariant objects.
So, we obtain
$\Locst(X,H,\vecnbigi)
\simeq
\Locst(X_b,H_b,\vecnbigi_b)$.

Let us study the general case.
Let $x\in H_b$.
There exists a diffeomorphism
$F:B\times (X_b,H_b)\simeq(X,H)$.
We have a $C^{\infty}$-submanifold
$F(B\times x)\subset (X,H)$.
There exists an open tubular neighbourhood $N$
of $F(B\times x)$
such that the following holds.
\begin{itemize}
 \item $H\cap N=\bigcup_{j=1}^{\ell} H_j$,
       where $H_j$ are smooth hypersurfaces of $N$.
       Let $H_{j,b}$ denote the fiber over $b$.
 \item There exists a diffeomorphism
       $(N,H_{1},\ldots,H_{\ell})
       \simeq B\times
       (N_b,H_{1,b},\ldots,H_{\ell,b})$.      
\end{itemize}
There exists a complex manifold
$N'$
with a ramified Galois covering
$\varphi:(N',H')\to (N,H\cap N)$
ramified along each $H_j'$
such that
$\varphi^{\ast}\vecnbigi_{|N}$
is unramified.
By the previous consideration,
$\Locst(N,H\cap N,\vecnbigi_{|N})
\simeq
\Locst(N_b,(H\cap N)_b,\vecnbigi_{|N_b})$.

Let $(\Ltilde_i,\vecnbigf)\in\Locst(X,H,\vecnbigi)$ $(i=1,2)$.
Let $f:(\Ltilde_1,\vecnbigf)_{|X_b}\to (\Ltilde_2,\vecnbigf)_{|X_b}$
be a morphism.
Let $X=\bigcup N^{(i)}$ be an open covering
such that
$(N^{(i)},H\cap N^{(i)})\simeq B\times(N^{(i)}_b,(H\cap N^{(i)})_b)$
where $N^{(i)}$ are as above.
There exist unique morphisms
$g^{(i)}:
(\Ltilde_1,\vecnbigf)_{|N^{(i)}}\to (\Ltilde_2,\vecnbigf)_{|N^{(i)}}$
such that
$g^{(i)}|_{N^{(i)}_b}=f_{|N^{(i)}_b}$.
By using the uniqueness,
we obtain $g^{(i)}_{|N^{(i)}\cap N^{(j)}}=g^{(j)}_{|N^{(i)}\cap N^{(j)}}$.
Hence, we obtain
$g:(\Ltilde_1,\vecnbigf)\to (\Ltilde_2,\vecnbigf)$
such that $g_{|X_b}=f$.

Let $(\Ltilde_b,\vecnbigf)\in \Locst(X_b,H_b,\vecnbigf)$.
There exist
$(\Ltilde^{(i)},\vecnbigf)
\in\Locst(N^{(i)},(H\cap N^{(i)}),\vecnbigf_{N^{(i)}})$
such that
$(\Ltilde^{(i)},\vecnbigf)_{|N^{(i)}_b}
\simeq
 (\Ltilde_b,\vecnbigf)_{|N^{(i)}_b}$.
By using the fully faithfulness,
we can patch $(\Ltilde^{(i)},\vecnbigf)$
to construct $(\Ltilde,\vecnbigf)$
such that $(\Ltilde,\vecnbigf)_{|N^{(i)}}\simeq (\Ltilde^{(i)},\vecnbigf)$.
We obtain $(\Ltilde,\vecnbigf)_{|X_b}\simeq
(\Ltilde_b,\vecnbigf)$.
\hfill\qed

\subsubsection{Deformation (2)}
\label{subsection;25.10.10.41}

Let $p:X\to B$ be a smooth fibration of
complex manifolds.
Let $H$ be a normal crossing hypersurface of $X$.
Let $B_0\subset B$ be a complex submanifold.
\begin{itemize}
 \item $B$ and $B_0$ are connected,
       and $\pi_1(B_0,b_0)\to \pi_1(B,b)$
       is surjective for any $b_0\in B_0$.
 \item For any $b\in B$,
       there exists a neighbourhood $B_b\subset B$
       of $b$ such that
       there exists a diffeomorphism
       $(X,H)_{|B_b}\simeq B_b\times (X_{b},H_{b})$.
\end{itemize}
Let $\vecnbigi$ be a good system of ramified irregular values
on $(X,H)$.
We obtain the induced good system of ramified irregular values
$\vecnbigi_{b}$ on $(X_b,H_b)$ by the restriction.

\begin{prop}
\label{prop;25.10.10.50}
Let $b_0\in B$ and
Let $(\Ltilde_{b_0},\vecnbigf)\in\Locst(X_{b_0},H_{b_0},\vecnbigi_{b_0})$.
Suppose that there exists a local system
$\Ltilde$ on $\Xtilde(\Htilde)$ such that
$\Ltilde_{|\Xtilde_{b_0}(H_{b_0})}\simeq \Ltilde_{b_0}$.
Then, there exists a unique 
Stokes structure $\vecnbigf$ of $\Ltilde$ indexed by $\vecnbigi$
such that
 $(\Ltilde,\vecnbigf)_{|\Xtilde_{b_0}(H_{b_0})}
 \simeq (\Ltilde_{b_0},\vecnbigf)$.
\end{prop}
\pf
Let us consider the case $B_0=\{b_0\}$.
Let $\ptilde:\Xtilde(\Htilde)\to B$ denote the projection.
For any $b\in B$,
let $\gamma:[0,1]\to B$ be a path
such that $\gamma(0)=b_0$ and $\gamma(1)=b$
and that $\gamma$ is injective.
There exist
a simply connected neighbourhood $W$ of $\gamma([0,1])$
and a diffeomorphism
$(X,H)_{|W}\simeq W\times (X_{b_0},H_{b_0})$.
By Proposition \ref{prop;25.10.10.30},
there exists the Stokes structure
$\vecnbigf$ of $\Ltilde_{|\pi^{-1}(W)}$
whose restriction to $\Xtilde_{b_0}(H_{b_0})$
equals $(\Ltilde_{b_0},\vecnbigf)$.
We obtain
$(\Ltilde_{|\Xtilde_{b}(H_{b})},\vecnbigf)$.
Because $B$ is simply connected,
by using Proposition \ref{prop;25.10.10.30},
we can check that $(\Ltilde_{|\Xtilde_{b}(H_{b})},\vecnbigf)$
is independent of the choice of $\gamma$.
For any $b\in B$,
there exist a simply connected neighbourhood $B_{b}$
and a diffeomorphism
$(X,H)_{|B_b}\simeq B_b\times (X_b,H_b)$.
By Proposition \ref{prop;25.10.10.30},
there exists the Stokes structure
$\vecnbigf$ of $\Ltilde_{|\pi^{-1}(B_b)}$
whose restriction to
$\Xtilde_b(H_b)$ equals
$(\Ltilde_{|\Xtilde_b(H_b)},\vecnbigf)$.
For $b,b'\in B$,
we can check that
the induced Stokes structures of
$\Ltilde_{|\pi^{-1}(B_b\cap B_{b'})}$
are the same.
Hence, we obtain the Stokes structure
$\vecnbigf$ of $\Ltilde$
such that
$(\Ltilde,\vecnbigf)_{|\Xtilde_{b_0}(H_{b_0})}
\simeq
 (\Ltilde_{b_0},\vecnbigf)$.

Let us study the general case.
Let $\Psi:B'\to B$ be a universal covering.
We obtain
$X'=B'\times_BX$ and
$H'=B'\times_BH$.
We obtain the induced good system of ramified irregular values
$\vecnbigi'$ over $(X',H')$.
We have the map
$\Psitilde:\Xtilde'(H')\to \Xtilde(H)$.
We set $\Gamma=\pi_1(B,b_0)$.
We have the natural Galois action of $\Gamma$ on $B'$.
It lifts to the action on $(X',H')$
and $\Xtilde'(H')$.

Let $b'_0\in \Psi^{-1}(b_0)$.
By the previous consideration,
there exists
a unique Stokes structure $\vecnbigf'$ of
$\Psi^{\ast}(\Ltilde)$
such that
\[
\Psi^{\ast}(\Ltilde,\vecnbigf')_{|\Xtilde'_{b'_0}(H_{b'_0})}
=(\Ltilde_{\Xtilde_{b_0}(H_{b_0})},\vecnbigf_{|\Xtilde_{b_0}(H_{b_0})}).
\]

Let $B_0'=B'\times_BB_0$.
Because $\pi_1(B_0,b_0)\to \pi_1(B,b_0)$ is surjective,
$B'_0$ is connected.
We set
$(X_{B_0},H_{B_0}):=(X,H)_{|B_0}$
and
$(X'_{B_0'},H'_{B_0'})=(X',H')_{|B_0'}$.
We have the induced map
$\Psi_1:(X'_{B_0'},H'_{B_0'})\to (X_{B_0},H_{B_0})$
and
$\Psitilde_1:
\Xtilde'_{B'_0}(H'_{B_0})\to \Xtilde_{B_0}(H_{B_0})$.
We obtain the Stokes structure
$\Psitilde_1^{\ast}(\vecnbigf_{B_0})$
of
$\Psi^{\ast}(\Ltilde)_{|\Xtilde'_{B_0'}(H'_{B_0'})}$
which is $\Gamma$-equivariant.
For any $\gamma\in\Gamma$,
the isomorphism
$\gamma^{\ast}:
\Psitilde^{\ast}(\Ltilde)_{\gamma(b'_0)}
\simeq
\Psitilde^{\ast}(\Ltilde)_{b'_0}$
is enhanced to the isomorphism of
the local systems with Stokes structure
$\gamma^{\ast}:
 (\Psitilde^{\ast}(\Ltilde),\vecnbigf')_{\gamma(b'_0)}
 \simeq
 (\Psitilde^{\ast}(\Ltilde),\vecnbigf')_{b'_0}$.
Hence,
$\vecnbigf'$ is $\Gamma$-equivariant,
and there exists 
a Stokes structure $\vecnbigf$ of
$\Ltilde$ indexed by $\vecnbigi$
such that
$\Psitilde^{\ast}(\Ltilde,\vecnbigf)
=(\Psitilde^{\ast}(\Ltilde),\vecnbigf')$.
Then, $(\Ltilde,\vecnbigf)$
satisfies the desired condition.
\hfill\qed

\subsection{Riemann-Hilbert-Birkhoff correspondence}

Let $\vecnbigi$ be a good system of ramified irregular values.
Let $\Mero(X,H,\vecnbigi)$ denote the category of 
good meromorphic flat bundles $(V,\nabla)$ on $(X,H)$
such that $\nbigi(\nbigv,x)\subset\nbigi_x$ for any $x\in H$.
Let $\Locst(X,H,\vecnbigi)$ denote the category of
local systems on $\Xtilde(H)$
with Stokes structure indexed by $\vecnbigi$.

Let $(V,\nabla)\in \Mero(X,H,\vecnbigi)$.
Let $L_V$ be the local system on $X\setminus H$
obtained as the sheaf of flat sections of $(V,\nabla)$.
It induces the local system
$\Ltilde_V$ on $\Xtilde(H)$.
Let $x\in X$.
There exists a neighbourhood $X_x\subset X$ of $x$
and a ramified Galois covering
$\varphi_x:(X_x',H_x')\to (X_x,H_x)$
such that
$\varphi_x^{\ast}(V,\nabla)$ is unramifiedly good
on $(X_x',H_x')$,
and that $\varphi_x^{-1}(x)$ consists one point $x'$.
Let $\varpi':\Xtilde'_x(H'_x)\to \Xtilde_x$
denote the oriented real blow up.
It induces a morphism
$\varphitilde_x:\Xtilde_x'(H_x')\to \Xtilde(H)$.
As explained in \cite[\S4]{Mochizuki-Stokes-good},
we obtain the Stokes structure
$\vecnbigf_{x'}=\bigl(
\vecnbigf^{y'}(\varphi_x^{\ast}(\Ltilde_V)_{y'})\,|\,
y'\in(\varpi')^{-1}(H_x')
\bigr)$
of $\varphitilde_x^{\ast}(\Ltilde_V)$
associated with $\varphi_y^{\ast}(V,\nabla)$.
Because $\varphi_y^{\ast}(V,\nabla)$
is equivariant with respect to the action
of the Galois group $\Gal$ of the covering $\varphi_x$,
$(\varphitilde^{\ast}(\Ltilde_V),\vecnbigf_{x'})$
is also equivariant with respect to $\Gal$-action.
We obtain the tuple of the filtrations
$\vecnbigf_x=(\nbigf^y(\Ltilde_{V,y})\,|\,y\in\varpi^{-1}(H_x))$.
It is the Stokes structure of $(V,\nabla)_{|X_x}$.
By varying $x$,
we obtain
the Stokes structure $\vecnbigf$ of $\Ltilde_V$.
In this way,
we obtain
\[
 \RHB(V,\nabla)=(\Ltilde_V,\vecnbigf)
 \in\Locst(X,H,\vecnbigi).
\]
It is a functor of abelian categories.

\begin{prop}
The functor $\RHB$ is an equivalence.
\end{prop}
\pf
If $\vecnbigi$ is unramified,
it is already explained in \cite[\S4.2]{Mochizuki-Stokes-good}.
Let us study the case where there exists
a ramified Galois covering $\varphi:(X',H')\to (X,H)$
such that $\varphi^{\ast}\vecnbigi$ is unramified.
Let $\Gal$ denote the Galois group.
For any $(\Ltilde,\vecnbigf)\in\Locst(X,H,\vecnbigi)$,
we obtain
$\Gal$-equivariant object
$\varphi^{\ast}(\Ltilde,\vecnbigf)\in\Locst(X',H',\varphi^{\ast}\vecnbigi)$.
There exits
$(V',\nabla')\in \Mero(X',H',\varphi^{\ast}\vecnbigi)$
such that
$\RHB(V',\nabla')\simeq \varphi^{\ast}(\Ltilde,\vecnbigf)$.
Because
$\RHB:\Mero(X',H',\varphi^{\ast}\vecnbigi)
\simeq \Locst(X',H',\varphi^{\ast}\vecnbigi)$,
$(V',\nabla')$ is $\Gal$-equivariant,
and there exists $(V,\nabla)\in\Mero(X,H,\vecnbigi)$
such that
$\varphi^{\ast}(V,,\nabla)\simeq (V',\nabla')$.
The $G$-equivariant isomorphism
$\varphi^{\ast}\RHB(V,\nabla)\simeq
\varphi^{\ast}(\Ltilde,\vecnbigf)$
induces an isomorphism
$\RHB(V,\nabla)\simeq (\Ltilde,\vecnbigf)$.
It implies that $\RHB$ is essentially surjective.
Clearly, $\RHB$ is faithful.
Let $(V_i,\nabla)\in\Mero(X,H,\vecnbigi)$
with a morphism $f:\RHB(V_1,\nabla)\to\RHB(V_2,\nabla)$.
Because 
$\RHB:\Mero(X',H',\varphi^{\ast}\vecnbigi)
\simeq \Locst(X',H',\varphi^{\ast}\vecnbigi)$,
we obtain
$g':\varphi^{\ast}(V_1,\nabla)\to \varphi^{\ast}(V_2,\nabla)$
such that 
$\RHB(g')=\varphi^{\ast}(f)$.
We obtain that $g'$ is $\Gal$-equivariant,
and hence there exists
$g:(V_1,\nabla)\to (V_2,\nabla)$
such that $\varphi^{\ast}(g)=g'$.
We can check  that $\RHB(g)=f$,
i.e.,
$\RHB$ is full in this case.

Let us study the general case.
By using the equivalence in the local case above,
we can check that
$\RHB$ is fully faithful.
Let $(\Ltilde,\vecnbigf)\in\Locst(X,D,\vecnbigi)$.
By the equivalence in the local case above,
there exist a open covering $X=\bigcup\nbigu_i$
and $(V_i,\nabla)\in \Mero(\nbigu_i,H\cap\nbigu_i,\vecnbigi_{|\nbigu_i})$
such that
$\RHB(V_i)\simeq (\Ltilde,\vecnbigf)_{|\varpi^{-1}(\nbigu_i)}$.
By the local equivalence above,
we have the isomorphisms
$\psi_{j,i}:
\RHB(V_i)_{|\nbigu_i\cap\nbigu_j}
\simeq
\RHB(V_j)_{|\nbigu_i\cap\nbigu_j}$
which satisfies the cocycle condition.
By patching them,
there exists $(V,\nabla)$
such that
$(V,\nabla)_{|\nbigu_i}=(V_i,\nabla)$.
Hence, $\RHB$ is essentially surjective.
\hfill\qed

\begin{cor}
In the setting of {\rm\S\ref{subsection;25.10.10.40}},
we obtain the equivalence
$\Mero(X,H,\vecnbigi)\simeq
\Mero(X_b,H_b,\vecnbigi_b)$.
\hfill\qed
\end{cor}

\begin{cor}
Let $X,B,H,B_0,\vecnbigi$
be as in {\rm\S\ref{subsection;25.10.10.41}}.
Let $(V_0,\nabla)\in\Mero(X_{B_0,}H_{B_0},\vecnbigi_{B_0})$.
Suppose that the associated local system $L_{V_0}$
extends to a local system on $X\setminus H$.
Then, there exists $(V,\nabla)\in \Mero(X,H,\vecnbigi)$
such that $(V,\nabla)_{|X_{B_0}}\simeq (V_0,\nabla)$.
Such $(V,\nabla)$ is unique up to canonical isomorphism.
\hfill\qed
\end{cor}

\subsubsection{Restriction to ample hypersurfaces}

Let $X$ be a smooth complex projective variety with $\dim X\geq 3$.
Let $H$ be a simple normal crossing hypersurface
with the irreducible decomposition $H=\bigcup_{i\in \Lambda}H_i$.
Let $\vecnbigi$ be a good system of ramified  irregular values on $(X,H)$.
Let $Z$ be an ample smooth hypersurface of $X$ such that
(i) $Z\subsetneq H$,
(ii) $Z\cup H$ is also simple normal crossing hypersurface.
Let $\vecnbigi_Z$ be the good system of ramified irregular values
on $(Z,Z\cap H)$ induced by $\vecnbigi$.

By the restriction,
we obtain the functors
\begin{equation}
\label{eq;25.10.10.1} 
\Mero(X,H,\vecnbigi)\to \Mero(Z,Z\cap H,\vecnbigi_Z),
\end{equation}
\begin{equation}
\label{eq;25.10.10.50}
 \Locst(X,H,\vecnbigf)  \to
 \Locst(Z,Z\cap H,\vecnbigi_Z).
\end{equation}

\begin{thm}
\label{thm;25.10.13.11}
The restriction functors
{\rm(\ref{eq;25.10.10.1})}
and {\rm(\ref{eq;25.10.10.50})}
are equivalences.
\end{thm}
\pf
For $I\subset \Lambda$,
we set $H_I=\bigcap_{i\in I}H_i$
and $H_I^{\circ}=H_I\setminus \bigcup_{j\not\in I}(H_I\cap H_j)$.
We also set
$H_{[j]}=\bigcup_{|I|=j}H_I$.

We have $H^{\circ}_{i,j}\cap Z\neq\emptyset$
for any $i,j\in\Lambda$ with $i\neq j$.
By \cite{Hamm-Le},
the following holds.
\begin{itemize}
 \item $\pi_1(Z\setminus H)\simeq \pi_1(X\setminus H)$.
 \item $\pi_1(H_i^{\circ}\cap Z)\to \pi_1(H_i^{\circ})$
       is onto for any $i\in I$.
\end{itemize}

Let $(V_i,\nabla)\in\Mero(X,H,\vecnbigi)$.
Suppose there exists
$f:(V_1,\nabla)_{|Z}\to (V_2,\nabla)_{|Z}$.
We obtain that
$(\Ltilde_1,\vecnbigf)_{|\Ztilde(H\cap Z)}
\to
(\Ltilde_2,\vecnbigf)_{|\Ztilde(H\cap Z)}$.
It extends to
$L_1\to L_2$
and
$\Ltilde_1\to\Ltilde_2$.
Moreover,
it is compatible with the Stoke structures
on $\varpi^{-1}(H_i^{\circ})$ for any $i$.
Hence, there exists
$f':(V_1,\nabla)_{|X\setminus H_{[2]}}
\to (V_2,\nabla)_{|X\setminus H_{[2]}}$.
Because $f'$ preserves the good Deligne-Malgrange lattices
on $X\setminus H_{[2]}$,
by the Hartogs property,
it extends to a morphism
$\ftilde:(V_1,\nabla)\to (V_2,\nabla)$
such that $\ftilde_{|Z}=f$.
Therefore, the functors are fully faithful.

Let $(\Ltilde_Z,\vecnbigf)\in \Locst(Z,H\cap Z,\vecnbigi_Z)$.
Because $\pi_1(Z\setminus H)\simeq \pi_1(X\setminus H)$,
there exists a local system $L$ on $X\setminus H$
such that $L_{|Z\setminus H}=L_Z$.
By Proposition \ref{prop;25.10.10.50},
there exists Stokes structure $\vecnbigf$
of $\Ltilde_{|\varpi^{-1}(X\setminus H_{[2]})}$
over $\vecnbigi_{|X\setminus H_{[2]}}$.
Because $H^{\circ}_{i,j}\cap Z\neq\emptyset$,
it extends to a Stokes structure
of $\Ltilde_{|\varpi^{-1}(X\setminus H_{[3]})}$
over $\vecnbigi_{|X\setminus H_{[3]}}$.
Hence, there exists
$(V',\nabla')\in \Mero(X\setminus H_{[3]},H\setminus H_{[3]},\vecnbigi)$
such that
$\RHB(V',\nabla')=(\Ltilde_{|\varpi^{-1}(X\setminus H_{[3]})},\vecnbigf)$.
By \cite{Malgrange-Lattice, Siu-extension},
$(V',\nabla)$ extends to a meromorphic connection
$(V,\nabla)$ on $(X,H)$.
By using \cite[Lemma 2.37, Lemma 2.38]{Mochizuki-Stokes-good},
we obtain that $(V,\nabla)\in\Mero(X,D,\vecnbigi)$.
\hfill\qed

\section{Estimate of characteristic numbers in the complex surface case}
\label{section;25.12.14.120}

\subsection{Statement}

Let $X$ be a smooth connected complex projective surface.
Let $H$ be a simple normal crossing hypersurface of $X$.
Let $H=\bigcup_{i\in\Lambda}H_i$ be the irreducible decomposition.
Let $[H_i]\in H^2(X,\real)$
denote the cohomology classes induced by $H_i$.

Let $(V,\nabla)$ be a good meromorphic flat bundle on $(X,H)$.
We obtain the good Deligne-Malgrange filtered bundle
$\nbigp^{\DM}_{\ast}V$.
We obtain the locally free $\nbigo_X$-module
$V^{\DM}=\nbigp_0^{\DM}V$.
For $-1<a\leq 0$,
we obtain the locally free $\nbigo_{H_i}$-module
$\lefttop{i}\Gr^{\nbigp^{\DM}}_{a}(V^{\DM})$.
We shall prove the following proposition
in \S\ref{subsection;25.12.9.50}
after preliminaries
in \S\ref{subsection;25.12.9.100} and 
\S\ref{subsection;25.12.9.101}.

\begin{prop}
\label{prop;25.10.9.40}
The class $c_1(V^{\DM})$
is contained in
\[
 \Bigl\{
 \sum_{i\in\Lambda} a_i[H_i]\in H^2(X,\real)\,\Big|\,
 0\leq a_i<\rank (V)
 \Bigr\}. 
\]
There exists $C_1>0$
depending only on $([H_i]^2,\rank(V))$
such that the following holds for any $-1<a\leq 0$:
\[
 \left|
 \deg_{H_i}\bigl(
 \lefttop{i}\Gr^{\nbigp^{\DM}}_a(V^{\DM})
 \bigr)
 \right|
 <C_1. 
\]
There exists $C_2>0$
depending only on
$\rank(V)$
and the intersection numbers
$[H_i]\cdot[H_j]$ $(i,j\in\Lambda)$
such that
\[
 |c_2(V^{\DM})|<C_2.
\]
\end{prop}

\subsection{Ramified coverings and good system ramified irregular values}
\label{subsection;25.12.9.100}

\subsubsection{Oriented real blow up}
\label{subsection;25.10.7.40}

Let $X$ be a complex surface.
Let $C\subset X$ be
a smooth $1$-dimensional compact connected complex submanifold.
We assume the following:
\begin{itemize}
 \item The inclusion map $C\to X$ is homotopy equivalent.       
\end{itemize}
Let $H$ be a hypersurface of $X$
such that
(i) $C\not\subset H$,
(ii) $C\cup H$ is a simple normal crossing hypersurface,
(iii) each connected component of $H$ is isomorphic to a disc.
We set $H_C:=C\cap H\subset C$.
Let $\varpi:\Xtilde(H\cup C)\to X$ denote the oriented real blow up
along $C\cup H$.
Let $\nbigo_{\Xtilde(H\cup C)}$ denote the sheaf of
holomorphic functions on $\Xtilde(H\cup C)$.
We set
$\nbigo_{\Xtilde(H\cup C)}(\ast (H\cup C)):=
\nbigo_{\Xtilde(H\cup C)}\otimes_{\varpi^{-1}\nbigo_X}
\varpi^{-1}(\nbigo_X(\ast(H\cup C)))$.

\subsubsection{Some sheaves on the oriented real blow up}
\label{subsection;25.12.9.30}

Let $\nbigo^{(e)}_{\Xtilde(H\cup C)}$
denote the sheaf on $\Xtilde(H\cup C)$
determined as follows.
\begin{itemize}
 \item $\nbigo^{(e)}_{\Xtilde(H\cup C)|X\setminus H}=\nbigo_{X\setminus H}$.
 \item For each $x\in C\setminus H$,
       there exists a holomorphic coordinate neighbourhood
       $(X_x,z_1,z_2)$
       such that $C\cap X_x=\{z_1=0\}$.
       Around $y\in\varpi^{-1}(x)$,
       by taking a branch of $z_1^{1/e}$,
       we set
       $\nbigo^{(e)}_{\Xtilde(H\cup C)}
       =\nbigo^{(e)}_{\Xtilde(H\cup C)}[z_1^{1/e}]$,
       which is independent of the choices of
       a coordinate system and a branch.
 \item For each $x\in H\setminus C$,
       there exists a holomorphic coordinate system
       $(X_x,z_1,z_2)$
       such that $H\cap X_x=\{z_2=0\}$.
       Around $y\in\varpi^{-1}(x)$,
       by taking a branch of $z_2^{1/e}$,
       we set
       $\nbigo^{(e)}_{\Xtilde(H\cup C)}
       =\nbigo^{(e)}_{\Xtilde(H\cup C)}[z_2^{1/e}]$,
       which is independent of the choices of
       a coordinate system and a branch.
 \item For each $x\in H\cap C$,
       there exists a holomorphic coordinate system
       $(X_x,z_1,z_2)$
       such that $H\cap X_x=\{z_2=0\}$
       and $C\cap X_x=\{z_1=0\}$.
       Around $y\in \varpi^{-1}(x)$,
       by taking branches of $z_i^{1/e}$,
       we set
       $\nbigo^{(e)}_{\Xtilde(H\cup C)}
       =\nbigo^{(e)}_{\Xtilde(H\cup C)}[z_1^{1/e},z_2^{1/e}]$,
       which is independent of
       the choices of a coordinate system
       and branches.
\end{itemize}
We set
$\nbigotilde_{\Xtilde(H\cup C)}
=\varinjlim \nbigo^{(e)}_{\Xtilde(H\cup C)}$.
We also set
\[
 \nbigotilde_{\Xtilde(H\cup C)}(\ast (H\cup C))
 = \nbigotilde_{\Xtilde(H\cup C)}
 \otimes_{\varpi^{-1}\nbigo_X}
 \varpi^{-1}\nbigo_X(\ast (H\cup C)).
\]
We obtain the sheaf of sets
\[
 \nbigotilde_{\Xtilde(H\cup C)}(\ast (H\cup C))
 \big/
 \nbigotilde_{\Xtilde(H\cup C)}
\]
on $\varpi^{-1}(H\cup C)$.

\subsubsection{Local systems induced by
good system of ramified irregular values on $(X,H)$}

Let $\vecnbigi=(\nbigi_x\,|\,x\in C\cup H)$
be a good system of ramified irregular values on $(X,C\cup H)$
in the sense of \cite[\S15.1]{Mochizuki-MTM}.
For $x\in H$ and $y\in \varpi^{-1}(x)$,
$\nbigi_x$ determines the subset
\[
 L(\vecnbigi)_y
 \subset
 \bigl(
 \nbigotilde_{\Xtilde(H\cup C)}(\ast (H\cup C))
 \big/
 \nbigotilde_{\Xtilde(H\cup C)}
 \bigr)_y.
\]
This procedure induces a subsheaf of sets
\[
 L(\vecnbigi)
 \subset
\nbigotilde_{\Xtilde(H\cup C)}(\ast (H\cup C))
\big/
\nbigotilde_{\Xtilde}(H\cup C) 
\]
on $\varpi^{-1}(C\cup H)$.
By restricting to
$\varpi^{-1}(C\setminus H)$,
we obtain a local system of sets
\[
 L(\vecnbigi)_{C\setminus H}:=
 L(\vecnbigi)_{|\varpi^{-1}(C\setminus H)}
\]
By restricting to $\varpi^{-1}(C\cap H)$,
we obtain a local system of sets
\[
 L(\vecnbigi)_{C\cap H}:=
 L(\vecnbigi)_{|\varpi^{-1}(C\cap H)}.
\]

\subsubsection{Statement}

We shall prove the following proposition
in \S\ref{subsection;25.12.9.10} after preliminaries.

\begin{prop}
\label{prop;25.10.7.30}
There exist complex manifolds
$X^{(j)}$ $(j=1,2,3,4)$
and holomorphic maps
\[
 X^{(4)}
 \stackrel{\varphi_4}{\lrarr}
 X^{(3)}
 \stackrel{\varphi_3}{\lrarr}
 X^{(2)}
 \stackrel{\varphi_2}{\lrarr}
 X^{(1)}\stackrel{\varphi_1}{\lrarr} X
\]
such that the following holds.
\begin{itemize}
 \item $X'=\varphi_1(X^{(1)})$ is
       a neighbourhood of $C$ in $X$.
       The induced map
       $\varphi_1:X^{(1)}\to X'$
       is a proper and finite map.
 \item There exists a hypersurface $H'$ of $X'$
       such that
       (i) any connected component of $H'$ is isomorphic to a disc,
       (ii) $H\cap H'=\emptyset$,
       (iii) $C\cup H'\cup H$ is normal crossing,
       (iv) $\varphi_1:X^{(1)}\to X'$ is a finite proper map
       ramified along $\varphi_1^{-1}(H\cup H')$.
       We set
       $C^{(1)}=\varphi_1^{-1}(C)$ and
       $H^{(1)}=\varphi_1^{-1}(H)$.
 \item $\varphi_2$ is a finite Galois covering map.
       We set
       $C^{(2)}=\varphi_2^{-1}(C^{(1)})$
       and
       $H^{(2)}=\varphi_2^{-1}(H^{(1)})$.              
 \item
       $\varphi_3$ is a finite proper map
       ramified along $C^{(3)}=\varphi_3^{-1}(C^{(2)})$,
       and
       the induced map
       $X^{(3)}\setminus C^{(3)}\to X^{(2)}\setminus C^{(2)}$
       is a cyclic covering map.
       We set
       $H^{(3)}=\varphi_3^{-1}(H^{(2)})$.
 \item $\varphi_4$ is a finite Galois covering map.
       We set
       $C^{(4)}=\varphi_4^{-1}(C^{(3)})$
       and
       $H^{(4)}=\varphi_4^{-1}(H^{(3)})$.              
 \item Let $\vecnbigi^{(4)}$ denote the induced
       good system of ramified irregular values
       on $(X^{(4)},C^{(4)}\cup H^{(4)})$.
       Then,
       the local systems
       $L(\vecnbigi_4)_{C^{(4)}\setminus H^{(4)}}$
       and
       $L(\vecnbigi_4)_{C^{(4)}\cap H^{(4)}}$
       are globally constant sheaves
       on each connected component.
\end{itemize}
\end{prop}

We have the following complement.
Let $G_3$ denote the Galois group of $\varphi_3$,
and let $G_4$ denote the Galois group of $\varphi_4$.
We have the natural actions of $G_j$ on $X^{(j)}$.
We shall prove the following lemma in
\S\ref{subsection;25.12.9.11}.

\begin{lem}
\label{lem;25.10.9.10}
There exists a natural action of 
$G_3\times G_4$ on $X^{(4)}$
which induces the $G_3$-action on $X^{(3)}$.
\end{lem}

\subsubsection{Preliminary}

Let $Y$ be a complex manifold.
Let $L$ be a holomorphic line bundle on $Y$.
\begin{lem}
Assume that there exists $\alpha\in H^2(Y,\seisuu)$
such that $c_1(L)=m\cdot \alpha$
for some integer $m$.
Then, there exists a holomorphic line bundle $L^{1/m}$
with an isomorphism
$(L^{1/m})^{\otimes m}\simeq L$. 
\end{lem}
\pf
From the exact sequence
$0\lrarr \seisuu\lrarr\nbigo
 \lrarr \nbigo^{\ast}\lrarr 0$,
we obtain the exact sequence
\[
 H^1(Y,\nbigo_Y)
 \lrarr
 H^1(Y,\nbigo_Y^{\ast})
 \stackrel{\psi_1}{\lrarr}
 H^2(Y,\seisuu)
 \stackrel{\psi_2}{\lrarr}
 H^2(Y,\nbigo_Y).
\]
We note that
$H^j(Y,\nbigo_Y)$ are $\cnum$-vector spaces.
Note that $\psi_2(c_1(L))=0$
in $H^2(Y,\nbigo_Y)$.
Hence, $m\cdot \psi_2(\alpha)=0$,
which implies that
$\psi_2(\alpha)=0$.
Hence, there exists a holomorphic line bundle $L_0$
such that
$c_1(L_0)=\alpha$.
Because
$c_1(L^{-1}\otimes L_0^{m})=0$,
there exists $\beta\in H^1(Y,\nbigo_Y)$
such that
$\psi_1(\beta)=L^{-1}\otimes L_0^{m}$.
We set
$L^{1/m}:=L_0\otimes \psi_1(\beta/m)$.
Then, it satisfies
$(L^{1/m})^{\otimes m}=L$.
\hfill\qed

\subsubsection{Covering branched along $H\cup H'$}
\label{subsection;25.12.9.20}

Let $H'$ be a hypersurface of $X$
such that
$H\cap H'=\emptyset$ and that
$H\cup H'\cup C$ is normal crossing.

\begin{lem}
Suppose that 
each connected component of $H\cup H'$ is homeomorphic to a disc.
Suppose that
$|H\cap C|=|H'\cap C|$.
Then,
for any $m\in\seisuu_{>0}$,
there exists
a complex manifold $X'$
with a holomorphic map 
$\varphi_m:X'\to X$ 
such that the following holds. 
\begin{itemize}
 \item $X'\setminus\varphi_m^{-1}(H\cup H')\to X\setminus(H\cup H')$
       is a cyclic covering of order $m$.
 \item $\varphi_m^{-1}(H\cup H')\to H\cup H'$
       is an isomorphism.
\end{itemize}
\end{lem}
\pf
We consider the line bundle
$L=\nbigo_X(H-H')$.
We have $c_1(L)=0$ in $H^2(X,\seisuu)$.
There exists a holomorphic line bundle
$L^{1/m}$
with an isomorphism
$(L^{1/m})^{\otimes m}\simeq L$.
There exists the morphism
$L^{1/m}\to L$
given by $v\longmapsto v^{m}$.
It extends to a morphism
$\Phi:
\proj(\nbigo_X\oplus L^{1/m})
\to
\proj(\nbigo_X\oplus L)$.
We set
$X_0=\proj(\nbigo_X\oplus 0)$
and
$X_{\infty}=\proj(0\oplus L)$.
We also set
$X_0^{(m)}=\proj(\nbigo_X\oplus 0)$
and
$X_{\infty}^{(m)}=\proj(0\oplus L^{1/m})$.
Then,
$\proj(\nbigo_X\oplus L^{1/m})
\setminus
(X_0^{(m)}\cup X_{\infty}^{(m)})
\to 
\proj(\nbigo_X\oplus L)
\setminus
(X_0\cup X_{\infty})$
is a cyclic covering of order $m$,
and $\Phi$ induces an isomorphism
$X_0^{(m)}\simeq X_0$
and
$X_{\infty}^{(m)}\simeq X_{\infty}$.

The morphism
$\nbigo_X\to \nbigo_X(H)$
induces
$X\to \proj^1(\nbigo_X\oplus L)$.
Let $X'$ denote the fiber product of
$X$ and $\proj(\nbigo_X\oplus L^{1/m})$
over $\proj(\nbigo_X\oplus L)$.
Then, the induced morphism $X'\to X$
satisfies the desired conditions.
\hfill\qed

\subsubsection{Covering branched along $C$}
\label{subsection;25.12.9.21}

Let $\chi$ denote the degree of
the normal bundle of $C$ in $X$.

\begin{lem}
Let $m$ be a positive integer
such that $\chi/m$ is an integer.
Then,
there exists a complex manifold $X'$
with a morphism
$\varphi_m:X'\to X$
such that
the following holds.
\begin{itemize}
 \item $X'\setminus \varphi_m^{-1}(C)\to X\setminus C$
       is a cyclic covering of order $m$.
 \item $\varphi_m^{-1}(C)\simeq C$.
\end{itemize} 
The Euler class of the normal bundle of $\varphi_m^{-1}(C)$
in $X'$ is $\chi/m$.
\end{lem}
\pf
We consider the holomorphic line bundle
$L=\nbigo_X(C)$.
Because $H^2(X,\seisuu)\simeq H^2(C,\seisuu)$,
there exists a $\alpha\in H^2(X,\seisuu)$
such that $c_1(L)=m\alpha$.
There exists a holomorphic line bundle $L^{1/m}$
with an isomorphism
$(L^{1/m})^{\otimes m}\simeq L$.
There exists the natural section $X\to L$.
Then, we can construct $X'\to X$
as in the previous lemma.
\hfill\qed

\subsubsection{Proof of Proposition \ref{prop;25.10.7.30}}
\label{subsection;25.12.9.10}

By shrinking $X$,
we may assume that there exists a hypersurface $H'$ of $X$
such that
(i) each connected component of $H'$ is isomorphic to disc,
(ii) $H\cap H'=\emptyset$,
(iii) $X\cup H\cup H'$ is normal crossing,
(iv) $|H\cap C|=|H'\cap C|$.
We also assume that $C\to X$ is homotopy equivalent.
There exists $m_1\in\seisuu_{>0}$
such that the following holds.
\begin{itemize}
 \item Let $x\in H\setminus C$.
       Then,
       $\nbigi_x\subset
       \nbigo^{(m_1)}_{X}(\ast H)_x/\nbigo^{(m_1)}_{X,x}$.
\end{itemize}
Let $\varphi^{(1)}_{m_1}:X^{(1)}\to X$
be the cyclic covering of order $m_1$
ramified along $H\cup H'$,
as in \S\ref{subsection;25.12.9.20}.
We set $C^{(1)}=(\varphi_{m_1}^{(1)})^{-1}(C)$,
$H^{(1)}=\varphi_{m_1}^{-1}(H)$
and
$H^{\prime(1)}=\varphi_{m_1}^{-1}(H')$.
We obtain
a good system of ramified irregular values
$\vecnbigi^{(1)}$ on $(X^{(1)},C^{(1)}\cup H^{(1)})$.

For each $x\in C^{(1)}\setminus H^{(1)}$,
let $[\nbigi^{(1)}_x]$ denote the quotient set of
$\nbigi^{(1)}_x$
by the natural action of the Galois group of
the extension
$\nbigo_{X^{(1)},x}\to \nbigo^{(e)}_{X^{(1)},x}$.
The tuple
$[\vecnbigi^{(1)}]=([\nbigi^{(1)}_x]\,|\,x\in C^{(1)}\setminus H^{(1)})$
gives a local system of finite sets on $C^{(1)}\setminus H^{(1)}$.
Note that the monodromy around any point of $C^{(1)}\cap H^{(1)}$
is trivial by the construction of $\varphi_{m_1}^{(1)}$.
Hence, $[\vecnbigi^{(1)}]$ extends to a local system
on $C^{(1)}$,
which is denoted by $[\vecnbigi^{(1)}]'$.

Let $x_0\in C^{(1)}$.
We obtain a homomorphism $\kappa^{(1)}$
of $\pi_1(C^{(1)},x_0)$
to the group of automorphisms of
the finite set $[\vecnbigi^{(1)}]'_{x_0}$.
Let $\varphi^{(2)}:X^{(2)}\to X^{(1)}$
denote the finite Galois covering map
corresponding to
the normal subgroup $\Ker(\kappa^{(1)})$
of $\pi_1(X^{(1)},x_0)=\pi_1(C^{(1)},x_0)$.
We set
$C^{(2)}=(\varphi^{(2)})^{-1}(C^{(1)})$
and
$H^{(2)}=(\varphi^{(2)})^{-1}(H^{(1)})$.
Note that
$(\varphi_{|C^{(2)}})^{-1}([\vecnbigi^{(1)}]')$
is constant.

Let $\chi$ be the Euler number of the normal bundle of
$C^{(2)}$ in $X^{(2)}$.
There exists $e_2\in\seisuu_{>0}$
such that the following holds.
\begin{itemize}
 \item Let $x\in C^{(2)}\setminus H^{(2)}$.
       Then,
       $\nbigi^{(2)}_x\subset
       \nbigo^{(e_2)}_{X^{(2)}}(\ast C^{(2)})_x/\nbigo^{(e_2)}_{X^{(2)},x}$.
\end{itemize}
We set $m_2=|\chi|$ if $\chi\neq 0$
or $m_2=e_2$ if $\chi=0$.
Let $\varphi^{(3)}_{m_2}:X^{(3)}\to X^{(2)}$
be the cyclic covering of order $m_2$
ramified along $C^{(2)}$
as in \S\ref{subsection;25.12.9.21}.
We set $C^{(3)}=(\varphi_{m_2}^{(3)})^{-1}(C^{(2)})$,
$H^{(3)}=(\varphi^{(2)}_{m_2})^{-1}(H^{(2)})$.
We obtain
a good system of ramified irregular values
$\vecnbigi^{(3)}$ on $(X^{(3)},C^{(3)}\cup H^{(3)})$.

Let
$\varpi^{(3)}_0:\Xtilde^{(3)}(C^{(3)})
\to X^{(3)}$
denote the oriented real blow up
of $X^{(3)}$ along $C^{(3)}$.
\begin{lem}
The local system of sets
$L(\vecnbigi^{(3)})_{C^{(3)}\setminus H^{(3)}}$
on $C^{(3)}\setminus H^{(3)}$
naturally extends to
a local system of sets
$L(\vecnbigi^{(3)})'_{C^{(3)}}$
on $(\varpi^{(3)}_0)^{-1}(C^{(3)})$.
\end{lem}
\pf
Let $x\in C^{(3)}\cap H^{(3)}$.
Let $\gamma$ be any loop around $(\varpi^{(3)}_0)^{-1}(x)$.
By the construction of $\varphi^{(1)}$,
the monodromy
of $L(\vecnbigi^{(3)})'_{C^{(3)}}$
along $\gamma$ is trivial.
Then, we obtain the claim of the lemma.
\hfill\qed

\vspace{.1in}
The orbit decompositions of $\nbigi^{(3)}_x$
$(x\in C^{(3)}\setminus H^{(3)})$
induce the decomposition
\[
 L(\vecnbigi^{(3)})'_{C^{(3)}}
 =\bigsqcup _{\gminio\in\Lambda}
  L(\vecnbigi^{(3)})'_{C^{(3)},\gminio},
\]
where 
$L(\vecnbigi^{(3)})'_{C^{(3)},\gminio}$
are principal $\seisuu/\ell(\gminio)$-bundles
on $(\varpi^{(3)}_0)^{-1}(C^{(3)})$.

\begin{lem}
\label{lem;25.12.9.20}
For any $x\in C^{(3)}\setminus H^{(3)}$,
the monodromy of $L(\vecnbigi^{(3)})'_{C^{(3)}}$
along $(\varpi_0^{(3)})^{-1}(x)$ is trivial. 
\end{lem}
\pf
If $\chi=0$, it follows from the construction of
$\varphi^{(3)}_{e_2}$.
If $\chi\neq 0$,
we note that Euler number of
the normal bundle of $C^{(3)}$ in $X^{(3)}$
is $1$ or $-1$.
Hence, the natural morphism
\[
 H_1\bigl(
 (\varpi_0^{(3)})^{-1}(C^{(3)}),\seisuu
 \bigr)
 \lrarr
 H_1(C^{(3)},\seisuu)
\]
is an isomorphism.
Hence, the natural morphism
\[
 H^1(C^{(3)},\seisuu/\ell\seisuu)
 \lrarr
 H^1\bigl(
 (\varpi_0^{(3)})^{-1}(C^{(3)}),\seisuu/\ell\seisuu
 \bigr)
\]
is an isomorphism for any $\ell\in\seisuu_{>0}$.
Then, we obtain the claim of the lemma.
\hfill\qed

\vspace{.1in}

Let
$\varpi^{(3)}_1:\Xtilde^{(3)}(C^{(3)}\cup H^{(3)})
\to X^{(3)}$
denote the oriented real blow up
of $X^{(3)}$ along $C^{(3)}\cup H^{(3)}$.
We obtain the following lemma
by the construction of $\varphi^{(1)}_{m_1}$
and Lemma \ref{lem;25.12.9.20}.

\begin{lem}
\label{lem;25.10.7.31}
The local system
$L(\vecnbigi^{(3)})_{C^{(3)}\cap H^{(3)}}$
is constant on each connected component of
$(\varpi_1^{(3)})^{-1}(C^{(3)}\cap H^{(3)})$.

\hfill\qed
\end{lem}

By Lemma \ref{lem;25.12.9.20},
$\vecnbigi^{(3)}=(\nbigi^{(3)}_x\,|\,x\in C^{(3)}\setminus H^{(3)})$
determines a local system on $C^{(3)}\setminus H^{(3)}$,
which extends to a local system of finite sets
on $C^{(3)}$.
Let $x_0\in C^{(3)}$.
We have the natural homomorphism $\kappa^{(3)}$
of $\pi_1(C^{(3)},x_0)$ to the automorphisms of $\nbigi^{(3)}_{x_0}$.
Let $X^{(4)}\to X^{(3)}$ be the covering space
corresponding to the normal subgroup
$\Ker(\kappa^{(3)})$
of $\pi_1(X^{(3)},x_0)=\pi_1(C^{(3)},x_0)$.
Then, the induced map $X^{(4)}\to X$
satisfies the desired conditions.
\hfill\qed

\subsubsection{Proof of Lemma \ref{lem;25.10.9.10}}
\label{subsection;25.12.9.11}

Let $g_{X^{(4)}}$ be a $G_4$-invariant Riemannian metric
of $X^{(4)}$
such that
$g_{X^{(4)}}(Jv,Jv)=g_{X^{(4)}}(v,v)$
for any tangent vector $v$,
where $J$ is the automorphism of $TX^{(4)}$
induced by the multiplication of $\sqrt{-1}$.
It induces a Riemannian metric $g_{X^{(3)}}$ of $X^{(3)}$
such that the covering map $\varphi_4$ is locally isometric.
For $j=3,4$,
let $N^{\real}_{C^{(j)}}X^{(j)}$ denote the normal bundles
of $C^{(j)}$ in $X^{(j)}$
obtained as the orthogonal complement of
$TC^{(j)}\subset TX^{(j)}_{|C^{(j)}}$.
They are naturally complex vector bundles.
The restriction $\varphi_{4|C^{(4)}}$ induces
an isomorphism
\begin{equation}
\label{eq;25.10.9.20}
N^{\real}_{C^{(4)}}X^{(4)}\simeq
\varphi_{4|C^{(4)}}^{\ast}\bigl(
N^{\real}_{C^{(3)}}X^{(3)}\bigr).
\end{equation}
Let $0_{C^{(j)}}$ denote the image of the zero section of
$N^{\real}_{C^{(j)}}X^{(j)}$.

We regard
$G_3\subset S^1=\bigl\{a\in\cnum\,\big|\,|a|=1\bigr\}$.
There exists the natural $G_3$-action on
$N^{\real}_{C^{(4)}}X^{(4)}$.
The isomorphism (\ref{eq;25.10.9.20})
is $G_3$-equivariant.
The $G_3$-action on $N^{\real}_{C^{(4)}}X^{(4)}$
extends to the $G_3\times G_4$-action on $N^{\real}_{C^{(4)}}X^{(4)}$.

There exist $G_3$-invariant neighbourhoods
$C^{(3)}\subset U^{(3)}\subset X^{(3)}$
and 
$0_{C^{(3)}}\subset V^{(3)}\subset N^{\real}_{C^{(3)}}X^{(3)}$
such that
the exponential map with respect to $g_{X^{(3)}}$
induces a $G_3$-equivariant diffeomorphism
$V^{(3)}\simeq U^{(3)}$.
We obtain the neighbourhood
$0_{C^{(4)}}\subset V^{(4)}\subset N^{\real}_{C^{(4)}}X^{(4)}$
as the pull back of $V^{(3)}$.
It is $G_3\times G_4$-invariant.
Let $C^{(4)}\subset U^{(4)}\subset X^{(4)}$
denote the neighbourhood obtained as the image of $V^{(4)}$
by the exponential map
with respect to $g_{X^{(4)}}$.
We have the $G_3\times G_4$-action on $U^{(4)}$,
which induces
the $G_3$-action on $U^{(3)}=U^{(4)}/G_4$.

Note that
the inclusion maps
$U^{(4)}\to X^{(4)}$ are homotopy equivalent.
For any $g\in G_3$,
the map $g:U^{(4)}\simeq U^{(4)}$
uniquely extends to $X^{(4)}\simeq X^{(4)}$
because
$g:U^{(3)}\simeq U^{(3)}$ extends to
$X^{(3)}\simeq X^{(3)}$,
and $X^{(4)}\to X^{(4)}$ is a covering map.
We obtain the $G_3$-action on $X^{(4)}$.
Let $g_3\in G_3$ and $g_4\in G_4$.
Note that $g_4^{-1}g_3g_4=g_3$ on $U^{(4)}$,
and that $g_4^{-1}g_3g_4$ and $g_3$ induce
the same map on $X^{(4)}$.
Hence, we obtain that
$g_3^{-1}g_2g_3=g_2$ on $X^{(4)}$,
i.e.,
the $G_3$-action and $G_4$-action are commuting.
We note that
the actions are independent of the choice of $g_{X^{(4)}}$.
\hfill\qed

\subsection{Degree of the graded pieces}
\label{subsection;25.12.9.101}

\subsubsection{Statement}

Let $X$, $C$ and $H$ be as in \S\ref{subsection;25.10.7.40}.
Let $(V,\nabla)$ be a good meromorphic flat bundle
on $(X,H\cup C)$
of rank $r$.
We obtain the good Deligne-Malgrange filtered bundle
$\nbigp^{\DM}_{\ast}(V)$ on $(X,C\cup H)$
indexed by $\real\times\real^H$.
For $a\in\real$,
by taking the graduation with respect to the filtration along $C$,
we obtain
the locally free $\nbigo_C(\ast H_C)$-module
$\lefttop{C}\Gr^{\nbigp^{\DM}}_a(V)$.
Moreover,
for $a\in\real$ and $\vecb\in\real^H$,
we obtain locally free sheaves
\[
 \nbigp^{\DM}_{\vecb}\bigl(
 \lefttop{C}\Gr^{\nbigp^{\DM}}_a(V)
 \bigr)
:= 
\lefttop{C}\Gr^{\nbigp^{\DM}}_a\bigl(
 \nbigp^{\DM}_{a,\vecb}(V)
 \bigr)
\]
on $C$.
Thus, we obtain a filtered bundle
$\nbigp_{\ast}^{\DM}\bigl(
 \lefttop{C}\Gr^{\nbigp^{\DM}}_a(V)
\bigr)$ on $(C,H_C)$.
Let $N_CX$ denote the normal bundle of $C$ in $X$.
We shall prove the following proposition in 
\S\ref{subsection;25.12.9.40}.

\begin{prop}
\label{prop;25.10.7.3}
We have
\begin{equation}
\label{eq;25.10.7.50}
 \deg\left(
  \nbigp^{\DM}_{\ast}
 \bigl(
 \lefttop{C}
 \Gr^{\nbigp^{\DM}}_{a}(V)
 \bigr)
 \right)
 =
a\cdot \deg(N_CX)
 \cdot
 \rank
 \lefttop{C}
 \Gr^{\nbigp^{\DM}}_{a}(V).
\end{equation}
\end{prop}

\subsubsection{Refinement to the equivariant case}
\label{subsection;25.10.9.1}

Suppose that
$X$ is equipped with an action of a finite cyclic group $G$
satisfying the following condition.
\begin{itemize}
 \item $G\cdot x=x$ for any $x\in C$.
 \item The $G$-action preserves
       any connected component of $H$.
 \item $(V,\nabla)$ is $G$-equivariant.
 \item The quotient space $X/G$ is a complex manifold
       such that
       the projection $X\to X/G$
       is a ramified covering
       of complex manifolds.
\end{itemize}
The Deligne-Malgrange filtered bundle
$\nbigp^{\DM}_{\ast}(V)$ is $G$-equivariant.
Let $\Irr(G)$ denote the set of
the irreducible representations of $G$.
We obtain the decomposition
\[
 \nbigp^{\DM}_{\ast}
 \bigl(
 \lefttop{C}\Gr^{\nbigp^{\DM}}_a(V)
 \bigr)
=\bigoplus_{\chi\in\Irr(G)}
 \nbigp^{\DM}_{\ast}
 \bigl(
 \lefttop{C}\Gr^{\nbigp^{\DM}}_a(V)
 \bigr)_{\chi},
\]
where
we have $g^{\ast}v=\chi(g)v$
for any $g\in G$ and a section $v$ of
$\nbigp^{\DM}
 \bigl(
 \lefttop{C}\Gr^{\nbigp^{\DM}}_a(V)
 \bigr)$.
The following proposition is a refinement of
Proposition \ref{prop;25.10.7.3},
which we shall prove in \S\ref{subsection;25.12.9.32}.

\begin{prop}
\label{prop;25.12.9.31}
We have
\begin{equation}
\label{eq;25.10.9.2}
 \deg\left(
  \nbigp^{\DM}_{\ast}
 \bigl(
 \lefttop{C}
 \Gr^{\nbigp^{\DM}}_{a}(V)
 \bigr)_{\chi}
 \right)
 =
a\cdot \deg(N_CX)
 \cdot
 \rank
 \lefttop{C}
 \Gr^{\nbigp^{\DM}}_{a}(V)_{\chi}.
\end{equation}
\end{prop}

\subsubsection{Global decomposition by irregular values in the special case}

Suppose that $X$ and $(V,\nabla)$ is $G$-equivariant
as in \S\ref{subsection;25.10.9.1}.
We study the special case.
 
\begin{condition}
\label{condition;25.10.8.11}
Suppose that
$\vecnbigi$ satisfies the following conditions.
 \begin{itemize}
  \item $\nbigi_x\subset \nbigo_{X}(\ast (H\cup C))_x/\nbigo_{X,x}$
	for any $x\in C\cup H$.
       We obtain the local system
       $L_{C\setminus H}(\vecnbigi)=(\nbigi_x\,|\,x\in C\setminus H)$.
 \item The local system
       $L_{C\setminus H}(\vecnbigi)$
       is constant.
       \hfill\qed
 \end{itemize}
\end{condition}

We set $\Gamma:=H^0(C\setminus H,L_{C\setminus H}(\vecnbigi))$.
We have the action of $G$ on $\Gamma$ by the pull back.
We have the decomposition
\[
 (V,\nabla)_{|\Chat}
 =\bigoplus_{\gminia\in \Gamma}
 (\Vhat_{\gminia},\nablahat_{\gminia})
\]
such that the following holds.
\begin{itemize}
 \item Let $x$ be any point of $C\setminus H$.
       Let $X_x$ be a neighbourhood of $x$,
       and we set $C_x:=X_x\cap C$.
       For $\gminia\in\Gamma$,
       let $\gminiatilde_x\in \nbigo_{X_x}(\ast C_x)$
       be a lift of $\gminia$.
       Then, 
       $\nablahat_{\gminia}-d\gminiatilde_P\id_{\Vhat_{\gminia}}$
       are regular singular along $C_P$.
 \item Let $x\in H\cap C$.
       Let $X_x$ be a neighbourhood of $x$,
       and $C_x:=X_x\cap C$ and $H_x:=X_x\cap H$.
       For $\gminia\in\Gamma$,
       let $\gminiatilde_x\in \nbigo_{X_x}(\ast (C_x\cup H_x))$
       be a lift of $\gminia$.
       Then,
       there exists a lattice
       $\Vhat^{(0)}_{\gminia}\subset \Vhat_{\gminia}$
       such that 
\[
       (\nablahat_{\gminia}
       -d\gminiatilde_x\id_{\Vhat_{\gminia}})
       \bigl(
       \Vhat^{(0)}_{\gminia}
       \bigr)
       \subset
       \Vhat^{(0)}_{\gminia}
       \otimes
       \Omega^1_{X_x}(\log C_x)(\ast H_x).
\]
 \item For any $g\in G$,
       we have
       $g^{\ast}\Vhat_{\gminia}=\Vhat_{g^{\ast}\gminia}$.
\end{itemize}
Let $\nbigp_{\ast}^{\DM}(V)$
be the unramifiedly good Deligne-Malgrange filtered bundle of $(V,\nabla)$.
We have the decomposition of filtered bundles
\[
 \nbigp_{\ast}^{\DM}(V)_{|\Chat}
 =\bigoplus_{\gminia\in\Gamma}
\nbigp_{\ast}^{\DM}\Vhat_{\gminia}.
\]

We set
$V^{\DM}=\nbigp^{\DM}_0(V)$
and $\Vhat^{\DM}_{\gminia}=\nbigp^{\DM}_0(\Vhat_{\gminia})$.
We obtain the following locally free sheaf on $C$
\[
 \lefttop{C}\Gr^{\nbigp^{\DM}}_a(V^{\DM})
 =\bigoplus_{\gminia\in\Gamma}
 \lefttop{C}\Gr^{\nbigp^{\DM}}_a(\Vhat^{\DM}_{\gminia}).
\]
It is equipped with the endomorphisms $\Res_C(\nabla)$.
The eigenvalues $\alpha$ are constant on $C$,
and satisfy $\Re(\alpha)=-a$.

For each $\gminia\in\Gamma$,
let $G(\gminia)=\{g\in G\,|\,g^{\ast}\gminia=\gminia\}$.
Let $\Irr(G(\gminia))$
denote the set of irreducible representations of $G(\gminia)$.
Each $\lefttop{C}\Gr^{\nbigp^{\DM}}_a(\Vhat^{\DM}_{\gminia})$
with $\Res_C(\nabla)$
is $G(\gminia)$-equivariant.
We obtain the decomposition
\[
 \lefttop{C}\Gr^{\nbigp^{\DM}}_a(\Vhat^{\DM}_{\gminia})
 =\bigoplus_{\chi\in\Irr(G(\gminia))}
 \lefttop{C}\Gr^{\nbigp^{\DM}}_a(\Vhat^{\DM}_{\gminia})_{\chi},
\]
where $g^{\ast}v=\chi(g)v$
for any section $v$ of $\Vhat_{\gminig,\chi}$
and $g\in G(\gminia)$.
We obtain the generalized eigen decomposition
with respect to $\Res_C(\nabla)$
\[
 \lefttop{C}\Gr^{\nbigp^{\DM}}_a(\Vhat_{\gminia}^{\DM})_{\chi}
 =\bigoplus_{\Re(\alpha)=-a}
 \lefttop{C}\EE_{\alpha}\Gr^{\nbigp^{\DM}}_a(\Vhat_{\gminia}^{\DM})_{\chi}.
\]
The nilpotent part of $\Res_C(\nabla)$
induces the weight filtration $W$.
We obtain the locally free $\nbigo_C$-modules
\[
 \lefttop{C}\Gr^W_j\Gr^{\nbigp^{\DM}}_a(\Vhat_{\gminia}^{\DM})_{\chi}
 =\bigoplus_{\Re(\alpha)=-a}
 \lefttop{C}
 \Gr^W_j\EE_{\alpha}\Gr^{\nbigp^{\DM}}_a(\Vhat_{\gminia}^{\DM})_{\chi}.
\]
We have the induced filtered bundles
$\nbigp^{\DM}_{\ast}\Bigl(
\lefttop{C}\Gr^W_j
\EE_{\alpha}
\Gr^{\nbigp^{\DM}}_a(\Vhat_{\gminia}^{\DM})_{\chi}
(\ast H_C)\Bigr)$
on $(C,H_C)$.

\subsubsection{An auxiliary function}

We may assume that 
$H^0(C,\nbigo_C(-C+H'))$ has a non-trivial section $\tau$
satisfying the following conditions.
\begin{itemize}
 \item $\tau_{|x}\neq 0$ in $\nbigo_C(-C+H')_{|x}$
       for any $x\in H'\cap C$.
 \item $\tau$ has only simple zeroes.
 \item There exists $H''$
       such that
       (i) $H''\cap C=\tau^{-1}(0)$,
       (ii)
       $H''\cap (H'\cup H)=\emptyset$,
       (iii) $C\cup H\cup H'\cup H''$ is normal crossing.
\end{itemize}
We set $H'\cap C=H'_C$
and $H''\cap C=H''_C$.

\begin{lem}
$\deg(N_CX)=|H'_C|-|H''_C|$.
\end{lem}
\pf
It follows from
$N_C(X)\simeq\nbigo_C(C)\simeq \nbigo_C(H'_C-H''_C)$.
\hfill\qed

\subsubsection{Connections on the graded pieces}

Let $x\in C$.
Let $(X_x,z_1,z_2)$ be a holomorphic coordinate neighbourhood
as in \S\ref{subsection;25.12.9.30}.
There exists a $G$-invariant section $\tautilde_x$
of $\nbigo_{X}(-C+H'-H'')$ on $X_x$
which induces $\tau_{|C_x}$.
For any $\alpha\in\cnum$,
we obtain the $G$-equivariant meromorphic flat bundle on
$(X_x,C_x\cup H_x'\cup H_x'')$:
\[
 \nbigl_{x,\tautilde_x}(\alpha)
 =\Bigl(
 \nbigo_{X_x}(\ast (C_x\cup H'_x\cup H''_x)),
 d+\alpha \tautilde_x^{-1}d\tautilde_x
 \Bigr).
\]
Let $\gminia_x$ be a $G(\gminia)$-invariant section of
$\nbigo_X(\ast (C\cup H\cup H'))$ on $X_x$
which is a lift of $\gminia$.
We obtain the $G(\gminia)$-equivariant meromorphic flat bundle
on $(X_x,C_x\cup H_x\cup H_x')$:
\[
 \nbigl_{x}(\gminia_x)
 =\Bigl(
 \nbigo_{X_x}(\ast (C_x\cup H'_x\cup H_x)),
 d+d\gminia_x 
 \Bigr).
\]
We obtain the $G(\gminia)$-equivariant meromorphic flat bundle
on $(X_x,C_x\cup H_x\cup H_x'\cup H_x'')$:
\[
 \bigl(
 V_x(-\alpha,-\gminia_x),\nabla^{(-\alpha,-\gminia_x)}
 \bigr)
 :=
 (V,\nabla)_{|X_x}
 \otimes
 \nbigl_{x,\tautilde_x}(\alpha)^{-1}
 \otimes
 \nbigl_{x}(\gminia_x)^{-1}.
\]
We obtain the $G(\gminia)$-equivariant filtered bundle
\[
 \nbigp^{\DM}_{\ast}\bigl(
  V_x(-\alpha,-\gminia_x)
  \bigr)
  =\nbigp^{\DM}_{\ast}\bigl(
 V
  \bigr)_{|X_x}
  \otimes
  \nbigp^{\DM}_{\ast}\bigl(
  \nbigl_{x,\tautilde_x}(\alpha)^{-1}
  \bigr)
\otimes
  \nbigp^{\DM}_{\ast}\bigl(
  \nbigl_x(\gminia_x)^{-1}
  \bigr).
\]
For $\gminib\in\nbigi_x$,
we obtain
\[
 \nbigp^{\DM}_{\ast}\bigl(
 \widehat{
 V_x(-\alpha,-\gminia_x)}
 _{\gminib-\gminia}
 \bigr)
=\nbigp^{\DM}_{\ast}\bigl(
 \Vhat_{\gminib}
 \bigr)_{|X_x}
 \otimes
   \nbigp^{\DM}_{\ast}\bigl(
  \nbigl_{x,\tautilde_x}(\alpha)^{-1}
  \bigr)
\otimes
  \nbigp^{\DM}_{\ast}\bigl(
  \nbigl_x(\gminia_x)^{-1}
  \bigr).
\]
We set $H_{C,x}=C_x\cap H_x$,
$H'_{C,x}=C_x\cap H'_x$,
and
$H''_{C,x}=C_x\cap H''_x$.
We set $G(\gminia,\gminib)=\{g\in G(\gminia)\,|\,g^{\ast}\gminib=\gminib\}$
for $\chi\in \Irr(G(\gminia,\gminib))$,
there exists the natural isomorphism of
\begin{multline}
 \lefttop{C}
 \Gr^W_j
 \EE_{\beta-\alpha}
 \Gr^{\nbigp^{\DM}}_{b+\Re(\alpha)}
 \Bigl(
 \widehat{
 V_x(-\alpha,-\gminia_x)}^{\DM}
 _{\gminib-\gminia}
 \Bigr)_{\chi}(\ast (H'_{C,x}\cup H_{C,x}\cup H''_{C,x}))
 \simeq
 \\
 \lefttop{C}
 \Gr^W_j
 \EE_{\beta}
 \Gr^{\nbigp^{\DM}}_{b}
 \Bigl(
 \Vhat_{\gminib}^{\DM}
 \Bigr)_{\chi}(\ast (H'_x\cup H_x\cup H''_x)).
\end{multline}
Note that $G(\gminia,\gminia)=G(\gminia)$.
For $(\gminib,\beta)=(\gminia,\alpha)$,
we have the induced connection
$\nabla(\gminia,\chi,\tau,\alpha,j)$
on
\[
 \lefttop{C}
 \Gr^W_j
 \EE_{0}
 \Gr^{\nbigp^{\DM}}_{0}
 \Bigl(
 \widehat{
 V_x(-\alpha,-\gminiahat_x)}^{\DM}
 _{0}
 \Bigr)_{\chi}(\ast (H'_x\cup H_x\cup H''_x))
 \simeq
 \lefttop{C}
 \Gr^W_j
 \EE_{\alpha}
 \Gr^{\nbigp^{\DM}}_{-\Re(\alpha)}
 \Bigl(
 \Vhat_{\gminia}^{\DM}
 \Bigr)_{\chi}(\ast (H'_x\cup H_x\cup H''_x)).
\]
The connection
$\nabla(\gminia,\chi,\tau,\alpha,j)$
is independent of the choice of
$\gminia_x$ and $\tautilde_x$,
i.e.,
it depends only on
$\gminia$ and $\tau$.

By varying $x\in C$,
we obtain
a connection 
$\nabla(\gminia,\chi,\tau,\alpha,j)$
on
the locally free $\nbigo_C(\ast (H'_C\cup H_C\cup H''_C))$-module
\[
 \lefttop{C}
 \Gr^W_j
 \EE_{\alpha}
 \Gr^{\nbigp^{\DM}}_{-\Re(\alpha)}
 \Bigl(
 \Vhat_{\gminia}^{\DM}
 \Bigr)_{\chi}(\ast (H'_C\cup H_C\cup H''_C)),
\]
depending only on
the choice of $\gminia$ and $\tau$.

\subsubsection{Degree in the special case}

\begin{lem}
We have
\begin{multline}
\label{eq;24.8.21.1}
\deg\left(
  \nbigp^{\DM}_{\ast}
 \Bigl(
 \lefttop{C}
 \Gr^W_j
 \EE_{\alpha}
 \Gr^{\nbigp^{\DM}}_{-\Re(\alpha)}
 \bigl(
 \Vhat_{\gminia}^{\DM}
 \bigr)_{\chi}(\ast (H'_C\cup H_C\cup H''_C))
 \Bigr)
 \right)
 = \\
-(|H'_C|-|H''_C|)
 \times
 \Re(\alpha)
 \times
 \rank
 \lefttop{C}
 \Gr^W_j
 \EE_{\alpha}
 \Gr^{\nbigp^{\DM}}_{-\Re(\alpha)}
 \bigl(
 \Vhat_{\gminia}^{\DM}
 \bigr)_{\chi}.
\end{multline}
\end{lem}
\pf
We obtain the filtered bundle
$\nbigp^{\prime\DM}_{\ast}
 \Bigl(
 \lefttop{C}
 \Gr^W_j
 \EE_{\alpha}
 \Gr^{\nbigp^{\DM}}_{-\Re(\alpha)}
 \bigl(
 \Vhat_{\gminia}^{\DM}
 \bigr)_{\chi}(\ast (H'_C\cup H_C\cup H''_C))
 \Bigr)$
with respect to the connection
$\nabla(\gminia,\chi,\tau,\alpha,j)$.
By the general property of
the Deligne-Malgrange filtration,
we have
\begin{equation}
\label{eq;25.10.8.10}
 \deg\left(
  \nbigp^{\prime\DM}_{\ast}
 \Bigl(
 \lefttop{C}
 \Gr^W_j
 \EE_{\alpha}
 \Gr^{\nbigp^{\DM}}_{-\Re(\alpha)}
 \bigl(
 \Vhat_{\gminia}^{\DM}
 \bigr)_{\chi}(\ast (H'_C\cup H_C\cup H''_C))
 \Bigr)
 \right)
=0.
\end{equation}
Around $H_C$,
we have 
\[
 \nbigp^{\prime\DM}_{\ast}
 \Bigl(
 \lefttop{C}
 \Gr^W_j
 \EE_{\alpha}
 \Gr^{\nbigp^{\DM}}_{-\Re(\alpha)}
 \bigl(
 \Vhat_{\gminia}^{\DM}
 \bigr)_{\chi}(\ast (H'_C\cup H_C\cup H''_C))
 \Bigr)
=\nbigp^{\DM}_{\ast}
 \Bigl(
 \lefttop{C}
 \Gr^W_j
 \EE_{\alpha}
 \Gr^{\nbigp^{\DM}}_{-\Re(\alpha)}
 \bigl(
 \Vhat_{\gminia}^{\DM}
 \bigr)_{\chi}(\ast (H'_C\cup H_C\cup H''_C))
 \Bigr).
\]
Around $P\in H'_C$,
we have
\[
 \nbigp^{\prime\DM}_{-\Re(\alpha)}
 \Bigl(
 \lefttop{C}
 \Gr^W_j
 \EE_{\alpha}
 \Gr^{\nbigp^{\DM}}_{-\Re(\alpha)}
 \bigl(
 \Vhat_{\gminia}^{\DM}
 \bigr)_{\chi}(\ast (H'_C\cup H_C\cup H''_C))
 \Bigr)
=\lefttop{C}
 \Gr^W_j
 \EE_{\alpha}
 \Gr^{\nbigp^{\DM}}_{-\Re(\alpha)}
 \bigl(
 \Vhat_{\gminia}^{\DM}
 \bigr)_{\chi}(\ast (H'_C\cup H_C\cup H''_C)),
\]
\[
 \nbigp^{\prime\DM}_{<-\Re(\alpha)}
 \Bigl(
 \lefttop{C}
 \Gr^W_j
 \EE_{\alpha}
 \Gr^{\nbigp^{\DM}}_{-\Re(\alpha)}
 \bigl(
 \Vhat_{\gminia}^{\DM}
 \bigr)_{\chi}(\ast (H'_C\cup H_C\cup H''_C))
 \Bigr)=0.
\]
Around $P\in H''_C$,
we have
\[
 \nbigp^{\prime\DM}_{\Re(\alpha)}
 \Bigl(
 \lefttop{C}
 \Gr^W_j
 \EE_{\alpha}
 \Gr^{\nbigp^{\DM}}_{-\Re(\alpha)}
 \bigl(
 \Vhat_{\gminia}^{\DM}
 \bigr)_{\chi}(\ast (H'_C\cup H_C\cup H''_C))
 \Bigr)
=\lefttop{C}
 \Gr^W_j
 \EE_{\alpha}
 \Gr^{\nbigp^{\DM}}_{-\Re(\alpha)}
 \bigl(
 \Vhat_{\gminia}^{\DM}
 \bigr)_{\chi}(\ast (H'_C\cup H_C\cup H''_C)),
\]
\[
 \nbigp^{\prime\DM}_{<\Re(\alpha)}
 \Bigl(
 \lefttop{C}
 \Gr^W_j
 \EE_{\alpha}
 \Gr^{\nbigp^{\DM}}_{-\Re(\alpha)}
 \bigl(
 \Vhat_{\gminia}^{\DM}
 \bigr)_{\chi}(\ast (H'_C\cup H_C\cup H''_C))
 \Bigr)=0.
\]
Then, we obtain (\ref{eq;24.8.21.1})
from (\ref{eq;25.10.8.10}).
\hfill\qed

\begin{cor}
\label{cor;24.8.22.1}
If Condition {\rm\ref{condition;25.10.8.11}} is satisfied,
the equality {\rm(\ref{eq;25.10.9.2})} holds.
\end{cor}
\pf
On $\lefttop{C}\Gr^{\nbigp^{\DM}}_a(V)$,
any eigenvalue $\alpha$ of $\Res_C(\nabla)$
satisfies $\Re(\alpha)=-a$.
We also have
$\deg(N_CX)=|H'_C|-|H''_C|$.
We obtain (\ref{eq;25.10.9.2})
from (\ref{eq;24.8.21.1})
by taking the sum with respect to
$(\gminia,\alpha,j)$.
\hfill\qed

\subsubsection{Proof of
Proposition \ref{prop;25.10.7.3}}
\label{subsection;25.12.9.40}

Let $\vecnbigi$ be a good system of ramified irregular values.
We use the notation in Proposition \ref{prop;25.10.7.30}
and Lemma \ref{lem;25.10.9.10}.
We obtain good meromorphic flat bundle
$(V^{(j)},\nabla)$ on $(X^{(j)},C^{(j)}\cup H^{(j)})$
$(j=1,2,3,4)$
as the pull back of $(V,\nabla)$.
Note that
$(V^{(4)},\nabla)$ and $(V^{(3)},\nabla)$
are naturally $G_3$-equivariant.
Because $(V^{(4)},\nabla)$
satisfies Condition \ref{condition;25.10.8.11},
we obtain the following equality
for any $a\in\real$ and $\chi\in\Irr(G_3)$:
\[
 \deg\Bigl(
 \nbigp^{\DM}_{\ast}\bigl(
 \lefttop{C^{(4)}}\Gr^{\nbigp^{\DM}}_a(V^{(4)})_{\chi}
 \bigr)
 \Bigr)
 =a\deg(N_{C^{(4)}}X^{(4)})
 \rank \lefttop{C^{(4)}}\Gr^{\nbigp^{\DM}}_a(V^{(4)})_{\chi}.
\]
We obtain
\begin{multline}
 \deg\Bigl(
 \nbigp^{\DM}_{\ast}\bigl(
 \lefttop{C^{(3)}}\Gr^{\nbigp^{\DM}}_a(V^{(3)})_{\chi}
 \bigr)
 \Bigr)
=
\frac{1}{|G_4|}
 \deg\Bigl(
 \nbigp^{\DM}_{\ast}\bigl(
 \lefttop{C^{(4)}}\Gr^{\nbigp^{\DM}}_a(V^{(4)})_{\chi}
 \bigr)
 \Bigr) 
\\
=a\frac{\deg(N_{C^{(4)}}X^{(4)})}{|G_4|}
 \rank \lefttop{C^{(4)}}\Gr^{\nbigp^{\DM}}_a(V^{(4)})_{\chi}
=a\deg(N_{C^{(3)}}X^{(3)})
 \rank \lefttop{C^{(3)}}\Gr^{\nbigp^{\DM}}_a(V^{(3)})_{\chi}.
\end{multline}
Let $1_{G_3}$ denote the trivial representation of $G_3$.
Let $m_3$ denote the ramification index of $\varphi_3$ along $C^{(3)}$.
We have
\[
 \nbigp^{\DM}_{\ast}\bigl(
 \lefttop{C^{(3)}}\Gr^{\nbigp^{\DM}}_{m_3a}(V^{(3)})_{1_{G_3}}
 \bigr)
\simeq
 \varphi_{3|C^{(3)}}^{\ast}\Bigl(
 \nbigp^{\DM}_{\ast}\bigl(
 \lefttop{C^{(2)}}\Gr^{\nbigp^{\DM}}_{a}(V^{(2)})
 \bigr)
 \Bigr).
\]
Note that $\varphi_{3|C^{(3)}}=\id_{C^{(3)}}$.
We obtain
\begin{multline}
 \deg\Bigl(
 \nbigp^{\DM}_{\ast}\bigl(
 \lefttop{C^{(2)}}\Gr^{\nbigp^{\DM}}_a(V^{(2)})
 \bigr)
 \Bigr)
=\deg\Bigl(
 \nbigp^{\DM}_{\ast}\bigl(
 \lefttop{C^{(3)}}\Gr^{\nbigp^{\DM}}_{m_3a}(V^{(3)})
 \bigr)_{1_{G_3}}
 \Bigr)
\\
=am_3\deg(N_{C^{(3)}}X^{(3)})
 \rank \lefttop{C^{(3)}}\Gr^{\nbigp^{\DM}}_{m_3a}(V^{(3)})
 =a\deg(N_{C^{(2)}}X^{(2)})
 \rank \lefttop{C^{(2)}}\Gr^{\nbigp^{\DM}}_{a}(V^{(2)}).
\end{multline}

Let $G_2$ denote the Galois group of $\varphi_2$.
We have
\begin{multline}
 \deg\Bigl(
 \nbigp^{\DM}_{\ast}\bigl(
 \lefttop{C^{(1)}}\Gr^{\nbigp^{\DM}}_a(V^{(1)})
 \bigr)
 \Bigr)
=\frac{1}{|G_2|}
 \deg\Bigl(
 \nbigp^{\DM}_{\ast}\bigl(
 \lefttop{C^{(2)}}\Gr^{\nbigp^{\DM}}_a(V^{(2)})
 \bigr)
 \Bigr)
\\
=a\frac{\deg(N_{C^{(2)}}X^{(2)})}{|G_2|}
 \rank
 \lefttop{C^{(2)}}
 \Gr^{\nbigp^{\DM}}_a(V^{(2)})
 =a\deg(N_{C^{(1)}}X^{(1)})
  \rank
 \lefttop{C^{(1)}}
 \Gr^{\nbigp^{\DM}}_a(V^{(1)}).
\end{multline}
We have
\[
 \nbigp^{\DM}_{\ast}\bigl(
 \lefttop{C^{(1)}}\Gr^{\nbigp^{\DM}}_{a}(V^{(1)})
 \bigr)
\simeq
 \varphi_{1|C}^{\ast}\Bigl(
 \nbigp^{\DM}_{\ast}\bigl(
 \lefttop{C}\Gr^{\nbigp^{\DM}}_{a}(V)
 \bigr)
 \Bigr).
\]
We also have
$N_{C^{(1)}}X^{(1)}
 \simeq
 \varphi_{1|C^{(1)}}^{\ast}\Bigl(
 N_{C}X
 \Bigr)$.
Let $m_1$ be the degree of $\varphi_1$.
By using Lemma \ref{lem;25.10.9.30} below,
we obtain
\begin{multline}
 \deg\Bigl(
 \nbigp^{\DM}_{\ast}\bigl(
 \lefttop{C}\Gr^{\nbigp^{\DM}}_a(V)
 \bigr)
 \Bigr)
 =\frac{1}{m_1}
 \deg\Bigl(
 \nbigp^{\DM}_{\ast}\bigl(
 \lefttop{C^{(1)}}\Gr^{\nbigp^{\DM}}_a(V^{(1)})
 \bigr)
 \Bigr)
\\
 =
a\frac{\deg(N_{C^{(1)}}X^{(1)})}{m_1}
 \rank \lefttop{C^{(1)}}\Gr^{\nbigp^{\DM}}_{a}(V^{(1)})
=a\deg(N_{C}X)
 \rank \lefttop{C}\Gr^{\nbigp^{\DM}}_{a}(V).
\end{multline}
Thus, we obtain Proposition \ref{prop;25.10.7.3}.
\hfill\qed

\begin{lem}
\label{lem;25.10.9.30}
Let $\varphi:C_1\to C_2$ be a surjective morphism of compact Riemann surfaces.
Let $D_2\subset C_2$ be a finite subset.
We set $D_1:=\varphi^{-1}(D_2)$.
Let $\nbigp_{\ast}V$ be a filtered bundle on $(C_2,D_2)$.
Let $d(\varphi)$ denote the degree of $\varphi$.
Then, we have
$\deg(\varphi^{\ast}(\nbigp_{\ast}V))
=d(\varphi)\cdot \deg(\nbigp_{\ast}V)$.
\end{lem}
\pf
Let $U_2$ be a neighbourhood of $D_2$ in $C_2$.
Let $h$ be a Hermitian metric of $V_{|C_2\setminus D_2}$
such that
(i) the Chern connection of $h$ is flat on $U_2\setminus D_2$,
(ii) $h$ is adapted to $\nbigp_{\ast}V$.
Let $R(h)$ denote the curvature of the Chern connection of $h$.
As well known, we have
\[
\deg(\nbigp_{\ast}V)
=\frac{\sqrt{-1}}{2\pi}
\int_{C_2\setminus D_2}
\Tr R(h).
\]
We have the induced metric $\varphi^{-1}(h)$
of $\varphi^{\ast}(V)_{|C_1\setminus D_1}$.
We obtain
\[
 \deg(\varphi^{\ast}\nbigp_{\ast}V)
= \frac{\sqrt{-1}}{2\pi}
\int_{C_1\setminus D_1}
 \Tr R(\varphi^{\ast}h)
 =d(\varphi)
 \frac{\sqrt{-1}}{2\pi}
 \int_{C_2\setminus D_2}
 \Tr R(\varphi^{\ast}h)
=d(\varphi)\deg(\nbigp_{\ast}V).
\]
\hfill\qed

\subsubsection{Proof of Proposition \ref{prop;25.12.9.31}}
\label{subsection;25.12.9.32}

We set $X'=X/G$ and $H'=H/G$.
Let $\pi:X\to X'$ denote the quotient map.
We set $C'=\pi(C)\subset X'$.
There exists a good meromorphic flat bundle $(V',\nabla)$
on $(X',H')$
such that $\pi^{\ast}(V',\nabla)\simeq(V,\nabla)$.
We have
$\pi^{\ast}(\nbigp_{\ast}^{\DM}V')\simeq
 \nbigp_{\ast}^{\DM}(V)$.
By Proposition \ref{prop;25.10.7.3},
we obtain
\[
  \deg\left(
  \nbigp^{\DM}_{\ast}
 \bigl(
 \lefttop{C'}
 \Gr^{\nbigp^{\DM}}_{b}(V')
 \bigr)
 \right)
 =
b\cdot \deg(N_{C'}X')
 \cdot
 \rank
 \lefttop{C'}
 \Gr^{\nbigp^{\DM}}_{b}(V').
\]
We identify $G=\seisuu/m\seisuu$ for
a positive integer $m$.
Let $\chi:G\to \cnum^{\ast}$ be the character
such that
\[
 \chi(1)=\exp(2\pi\sqrt{-1}(\ell/m))
\]
where $0\leq \ell<m$.
Let $a\in\real$.
We set $b=(a-\ell)/m$.
There exists a natural isomorphism
\[
  \nbigp^{\DM}_{\ast}
 \bigl(
 \lefttop{C'}
 \Gr^{\nbigp^{\DM}}_{b}(V')
 \bigr)
\simeq
 \nbigp^{\DM}_{\ast}
 \bigl(
 \lefttop{C}
 \Gr^{\nbigp^{\DM}}_{a}(V)
 \bigr)_{\chi}
 \otimes
 N_C(X)^{-\ell}.
\]
Note that
$\deg(N_{C'}X')=m\deg(N_CX)$.
We obtain
\begin{multline}
 \deg\left(
 \nbigp^{\DM}_{\ast}
 \bigl(
 \lefttop{C}
 \Gr^{\nbigp^{\DM}}_{a}(V)
 \bigr)
 \right)
=
 \deg\left(
 \nbigp^{\DM}_{\ast}
 \bigl(
 \lefttop{C'}
 \Gr^{\nbigp^{\DM}}_{b}(V')
 \bigr)
 \right)
+\ell\deg(N_CX)
 \\
=b\cdot \deg(N_{C'}X')
 \cdot
 \rank
 \lefttop{C'}
 \Gr^{\nbigp^{\DM}}_{b}(V')
 +\ell\deg(N_CX)
  \rank
 \lefttop{C'}
 \Gr^{\nbigp^{\DM}}_{b}(V')
\\
 =a\deg(N_CX)  \rank
 \lefttop{C}
 \Gr^{\nbigp^{\DM}}_{a}(V)_{\chi}.
\end{multline}
Thus, we obtain Proposition \ref{prop;25.12.9.31}.
\hfill\qed

\subsection{Proof of Proposition \ref{prop;25.10.9.40}}
\label{subsection;25.12.9.50}

Note there exists $0\leq \alpha_i<\rank(V)$ such that
\[
c_1(\nbigp_{\ast}^{\DM}V)
=c_1(\nbigp_0^{\DM}V)
+\sum_{i\in\Lambda}\alpha_i[H_i].
\]
Because $c_1(\nbigp_{\ast}^{\DM}V)=0$,
we obtain the first claim.
By Proposition \ref{prop;25.10.7.3},
we have
\[
 \deg_{H_i}
 \Bigl(
 \nbigp_{\ast}^{\DM}\bigl(
 \lefttop{i}\Gr^{\nbigp^{\DM}}_a(V^{\DM})\bigr)
  \Bigr)=
  a\cdot \deg(N_{H_i}X)
  \rank \lefttop{i}\Gr^{\nbigp^{\DM}}_a(V^{\DM})
= a\cdot [H_i]^2
  \rank \lefttop{i}\Gr^{\nbigp^{\DM}}_a(V^{\DM}).
\]
Hence, we obtain the second claim.
Because
$\int\ch_2(\nbigp^{\DM}_{\ast}V)=0$,
we obtain
\begin{multline}
\int_X\ch_2(\nbigp_0^{\DM}V)
-\sum_{i\in\Lambda}
 \sum_{-1<a\leq 0}
 a\deg_{H_i}\bigl(
 \lefttop{i}\Gr^{\nbigp^{\DM}}_a(\nbigp_0^{\DM}V)
 \bigr)
 =
\\
 -\frac{1}{2}
 \sum_{i\in\Lambda}\sum_{-1<a\leq 0}
 a^2\rank
 \bigl(
 \lefttop{i}\Gr^{\nbigp^{\DM}}_a(\nbigp_0^{\DM}V)
 \bigr)\int_X[H_i]^2
\\
 -\frac{1}{2}
 \sum_{\substack{i\neq j\\ P\in H_i\cap H_j}}
 \sum_{-1<a_i,a_j\leq 0}
 a_i\cdot a_j
 \cdot
 \rank
 \lefttop{P}\Gr^{\nbigp^{\DM}}_{(a_i,a_j)}
 \bigl(
 \nbigp^{\DM}_0(V)
 \bigr).
\end{multline}
We obtain the third claim.
\hfill\qed

\section{Algebraic meromorphic flat bundles}

\subsection{Good meromorphic flat bundles}

Let $\hyperk$ be an algebraically closed field of characteristic $0$.
Let $X$ be any $n$-dimensional smooth $\hyperk$-variety.
Let $H$ be a simple normal crossing hypersurface of $X$.
Let $\nbigohat_{X,x}$ denote the completion
of the local ring $\nbigo_{X,x}$
with respect to the maximal ideal $\gminim_x$.
There exists a parameter system $(x_1,\ldots,x_n)$ of $\gminim_x$
such that the ideal of $H$ at $x$ is generated by
$\prod_{i=1}^{\ell}x_i$.
For any positive integer $e$,
we set
$\nbigohat^{(e)}_{X,x}
=\nbigohat_{X,x}[x_i^{1/e}\,|\,i=1,\ldots,\ell]$.
For any $\nbigohat_{X,x}$-module $\nbigf$,
we set
$\nbigf(\ast H)=\nbigf\otimes_{\nbigo_{X,x}}\nbigo_{X}(\ast H)_x$
and
$\varphi_e^{\ast}(\nbigf)
=\nbigf\otimes\nbigohat^{(e)}_{X,x}$.
For any free $\nbigohat_{X,x}(\ast H)$-module $\nbigf$,
a lattice of $\nbigf$ means
a free $\nbigo_{X,x}$-submodule $\nbigf_0\subset\nbigf$
such that $\nbigf_0(\ast H)=\nbigf$.

Let $V$ be a free $\nbigo_{X}(\ast H)$-module of finite rank
with an integrable connection $\nabla:V\to V\otimes\Omega^1_X$.
We set $V_{|\xhat}=V_x\otimes\nbigohat_{X,x}$
which is equipped with the induced integrable connection $\nabla$.
Any $\varphi_e^{\ast}(V_{|\xhat})$ is equipped with
the induced integrable connection $\nabla$.
We recall some definitions.
\begin{df}
\mbox{{}}  
\begin{itemize}
 \item 
$(V,\nabla)$ is called regular at $x$
if there exists a free $\nbigohat_{X,x}$-module $V_0\subset V$
such that 
(i)  $V_0(\ast H)=V$,
(ii) $\nabla(V_0)\subset V_0\otimes\Omega^1(\log H)$.
\item $(V,\nabla)$ is called unramifiedly good at $x$
      if there exist a good set of irregular values
      $\nbigi\subset\nbigohat_{X,x}(\ast H)/\nbigo_{X,x}$
      and a decomposition
\begin{equation}
\label{eq;25.10.13.1}
      (V,\nabla)
      =\bigoplus_{\gminia\in\nbigi}
      (V_{\gminia},\nabla_{\gminia})
\end{equation}
      such that
      $\nabla_{\gminia}-d\gminiatilde\id_{V_{\gminia}}$
      are regular.
      Here, $\gminiatilde\in\nbigohat_{X,x}$
      denote lifts of $\gminia$.
      We allow $V_{\gminia}=0$
      in {\rm(\ref{eq;25.10.13.1})}.
      The set
      $\nbigi(V,\nabla,x)=
      \bigl\{\gminia\in\nbigi\,|\,V_{\gminia}\neq 0\bigr\}$
      is well defined.
 \item $(V,\nabla)$ is called good at $x$
       if there exists $e\in\seisuu_{>0}$
       such that
       $\varphi_e^{\ast}(V,\nabla)$
       is unramifiedly good.
       We obtain the set
       $\nbigi(V,\nabla,x)\subset
       \nbigohat^{(e)}_{X,x}(\ast H)/\nbigohat^{(e)}_{X,x}$,
       which is invariant under
       the action of the Galois group
       of the extension
       $\nbigohat_{X,x}\subset\nbigohat^{(H,e)}_{X,x}$.
\end{itemize}
We say that $(V,\nabla)$ is
good (resp. regular, unramifiedly good) on $(X,H)$
if $(V,\nabla)_{|\xhat}$ is
good (resp. regular, unramifiedly good) at any $x\in H$. 
\hfill\qed
\end{df}

\subsubsection{Extension of fields}

Let $X$ be a smooth $\hyperk$-variety with 
normal crossing hypersurface $H$.
Let $\hyperk'$ be an algebraically closed field which contains $\hyperk$.
We set $(X_{\hyperk'},H_{\hyperk'})=(X,H)\times_{\hyperk}\hyperk'$.
Let $(V,\nabla)$ be a meromorphic flat bundle on $(X,H)$.
Let $(V_{\hyperk'},\nabla)$ denote the induced meromorphic
flat bundle on $(X_{\hyperk'},H_{\hyperk'})$.

\begin{prop}
\label{prop;25.10.13.10}
$(V,\nabla)$ is good if and only if
$(V_{\hyperk'},\nabla)$ is good.
\end{prop}
\pf
Suppose that
$(V_{\hyperk'},\nabla)$ is good.
Let $x$ be any closed point of $X$.
The induced closed point of $X_{\hyperk'}$
is also denoted by $x$.
Note that
$(V_{\hyperk'|\xhat},\nabla)$
is induced by
$(V_{|\xhat},\nabla)$
and the extension
$\nbigohat_{X,x}\to \nbigohat_{X_{\hyperk'},x}$.
By \cite[Lemma 1.6.5, Definition 3.2.1, Theorem 4.4.2]{kedlaya},
$(V_{|\xhat},\nabla)$ is good
if and only if
$(V_{\hyperk'|\xhat},\nabla)$ is good.
Hence,
$(V,\nabla)$ is good.

Let us study the converse.
Let $H=\bigcup_{i\in\Lambda}H_i$ be the irreducible decomposition.
For any non-empty $I\subset \Lambda$,
we set
$H_I=\bigcap_{j\in I}H_j$
and 
$H_I^{\circ}=H_I\setminus
\bigcup_{j\in \Lambda\setminus I}H_j$.
Let $y$ be a closed point of
$(H_{\hyperk'})_I^{\circ}$.
Let $\eta$ denote the point of $H_I^{\circ}$
obtained as the image of $y$
by the natural morphism
$X_{\hyperk'}\to X$.
There exist
an affine neighbourhood $U$ of $\eta$ in $X$
with sections $x_1,\ldots,x_m\in \nbigo_X(U)$
such that
the ideal of $H\cap U$ is generated by
$\prod_{i=1}^mx_i$.
For any $e\in\seisuu_{>0}$,
let $U^{(e)}$ denote the affine $\hyperk$-variety
associated with
$\nbigo_X(U)[x_1^{1/e},\ldots,x_m^{1/e}]$
with the natural morphism
$\varphi_e:U^{(e)}\to U$.
There exists $e\in\seisuu_{>0}$
such that
$\varphi_e^{\ast}V$
is unramifiedly good at any point of
$\varphi_e^{-1}(H)$.
Let $R$ denote the ring
corresponding to $\varphi_e^{-1}(U\cap H_I)$.
We set $z_i=x_i^{1/e}$ $(i=1,\ldots,m)$.
We obtain the locally free
$R[\![z_1,\ldots,z_m]\!][z_1^{-1},\ldots,z_m^{-1}]$-module
$V_I$
with the integrable connection
from $\varphi_e^{\ast}(V,\nabla)$.
By \cite[Proposition 4.4.1]{kedlaya2}
and a remark after \cite[Definition 5.1.1]{kedlaya},
there exists a finite etale extension
$R\subset R'$
such that the following holds.
\begin{itemize}
 \item $\varphi_e^{-1}(y)$ is contained
       in the image of $\Spec(R')\to\Spec(R)$.
 \item There exist a good set of irregular values
       $\nbigi\subset
       R'[\![z_1,\ldots,z_m]\!][z_1^{-1},\ldots,z_{m}^{-1}]
       \big/
       R'[\![z_1,\ldots,z_m]\!]$
       and a good decomposition
\begin{equation}
\label{eq;26.2.19.1}
       V_I\otimes_{R[\![z_1,\ldots,z_m]\!]} R'[\![z_1,\ldots,z_m]\!]
       =\bigoplus_{\gminia\in\nbigi}
       (V_{I,\gminia},\nabla_{\gminia})       
\end{equation}
       such that $\nabla_{\gminia}-d\gminiatilde$
       are regular singular.
\end{itemize}
We obtain the ramified good decomposition
of $(V_{\hyperk'|\yhat},\nabla)$
from (\ref{eq;26.2.19.1}).
\hfill\qed

\subsubsection{Restriction}

Let $X$ be a smooth complex $\hyperk$-variety with $\dim X\geq 3$.
Let $H$ be a simple normal crossing hypersurface
with the irreducible decomposition $H=\bigcup_{i\in \Lambda}H_i$.
Let $\Sigma\subset T^{\ast}(X\setminus H)$
be good on $(X,H)$.
Let $\vecnbigi=\vecnbigi(\Sigma)$
be a good system of ramified  irregular values on $(X,H)$.

Let $Z$ be an ample smooth hypersurface of $X$ such that
(i) $Z\subsetneq H$,
(ii) $Z\cup H$ is also simple normal crossing hypersurface.
Let $\vecnbigi_Z$ be the good system of ramified irregular values
on $(Z,Z\cap H)$ induced by $\vecnbigi$.
By the restriction,
we obtain the functors
\begin{equation}
\label{eq;25.10.11.1} 
\Mero(X,H,\vecnbigi)\to \Mero(Z,Z\cap H,\vecnbigi_Z),
\end{equation}
We obtain the following proposition
from Proposition \ref{prop;25.10.13.10}
and Theorem \ref{thm;25.10.13.11}.

\begin{thm}
The restriction functor
{\rm(\ref{eq;25.10.11.1})}
is an equivalence.
\hfill\qed
\end{thm}

\subsection{Good lattices and the eigenvalues of the residues}

\subsubsection{Unramifiedly good lattices}

Let $X$ be a smooth projective $\hyperk$-variety
with a simple normal crossing hypersurface $H$
with the irreducible decomposition $H=\bigcup_{k\in\Lambda}H_k$.
Let $(V,\nabla)$ be an unramifiedly good meromorphic flat bundle on $(X,H)$.
For any $x\in H$,
there exists a decomposition
\[
 (V_{|\xhat},\nabla)
 =\bigoplus_{\gminia\in\nbigi(V,\nabla,x)}
 (V_{\xhat,\gminia},\nabla_{\gminia})
\]
such that $\nabla_{\gminia}-d\gminiatilde$
are regular.

Let $E$ be a lattice of $(V,\nabla)$,
i.e., $E$ be a locally free $\nbigo_X$-submodule of $V$
such that $E(\ast H)=V$.
We obtain
$E_{|\xhat}\subset V_{|\xhat}$.

\begin{df}
We say that $E$ is unramifiedly good at $x$
if
\begin{equation}
\label{eq;26.2.22.1}
 E_{|\xhat}=\bigoplus_{\gminia\in\nbigi(V,\nabla,x)}
(E_{|\xhat}\cap V_{\xhat,\gminia}),
\end{equation}
and $\nabla_{\gminia}-d\gminia$ are logarithmic
with respect to $E_{|\xhat}\cap V_{\xhat,\gminia}$.
We say that $E$ is unramifiedly good on $(X,H)$
if it is unramifiedly good at any $x\in H$.
\hfill\qed
\end{df}

\subsubsection{Good lattices}

Let $(V,\nabla)$ be a good meromorphic flat bundle on $(X,H)$.
Let $x\in H$.
Let $U=\Spec(R)$ be an affine open neighbourhood of $x$
with an etale coordinate system $(x_1,\ldots,x_n)$
such that $U\cap H=\bigcup_{i=1}^{\ell}\{x_i=0\}$.
For $e\in\seisuu_{>0}$,
$R^{(e)}=R[x_1^{1/e},\ldots,x_{\ell}^{1/e}]$
and $U^{(e)}=\Spec(R^{(e)})$.
Let $H^{(e)}_U$ denote the inverse image of $H$ by $U^{(e)}\to X$.
We obtain the unramifiedly good meromorphic flat bundle
$(V^{(e)}_U,\nabla)$ on $(U^{(e)},H^{(e)}_U)$
as the pull back of $(V,\nabla)$.
Let $\Gal^{(H,e)}_U$ denote the Galois group of
the ramified covering $U^{(e)}\to U$.

\begin{df}
Let $E$ be a lattice of $V$.
We say that $E$ is good at $x$
if there exist an etale coordinate neighbourhood
$(U,x_1,\ldots,x_n)$ of $x$ as above
and a $\Gal^{(H,e)}_U$-equivariant
unramifiedly good lattice $E_U^{(e)}$
of $(V^{(e)}_U,\nabla)$
such that $E_{|U}$ is the descent of $E_U^{(e)}$.
We say that $E$ is good on $(X,H)$
if it is good at any $x\in H$.
\hfill\qed
\end{df}

\subsubsection{Decompositions into the regular part and the irregular part}

Let $(V,\nabla)$ be a good meromorphic flat bundle on $(X,H)$.
Let $E$ be a good lattice.
Let $x\in H_k$.
We have the natural map
\[
 \rho_k:
 \nbigi(V,\nabla,x)
 \to
 \nbigohat^{(H,e)}_X(\ast H)_x
 \big/
 \nbigohat^{(H,e)}_X(\ast H_{\neq k})_x
\]
Let $\nbigi(V,\nabla,x)_{k,0}=\rho_k^{-1}(0)$.
We set
\[
(E_{|\xhat})^{\reg,k}
 =\bigoplus_{\gminia\in \nbigi(V,\nabla,x)_{k,0}}
  (E_{|\xhat}\cap V^{(H,e)}_{\xhat,\gminia}),
  \quad\quad
(E_{|\xhat})^{\irreg,k}
=\bigoplus_{\gminia\not\in \nbigi(V,\nabla,x)_{k,0}}
 (E_{|\xhat}\cap V^{(H,e)}_{\xhat,\gminia}).
\]
We obtain the decomposition
\begin{equation}
\label{eq;26.2.21.12}
 E_{|\xhat}=
 (E_{|\xhat})^{\reg,k}
 \oplus
 (E_{|\xhat})^{\irreg,k}.
\end{equation}
It is easy to check the following
\begin{equation}
\label{eq;26.2.21.14}
 \nabla (E_{|\xhat})^{\reg,k}
  \subset
  (E_{|\xhat})^{\reg,k}
  \otimes
  \Bigl(
  \Omega^1(\log H_k)
  +
  \Omega^1(\ast H_{\neq k})
  \Bigr),
  \quad\quad
  \nabla (E_{|\xhat})^{\irreg,k}
  \subset
  (E_{|\xhat})^{\irreg,k}
 \otimes
  \Omega^1(\ast H).
\end{equation}

Let $E_{|\Hhat_k}$ denote the completion of $E$ along $H_k$,
i.e.,
$E_{|\Hhat_k}=\varprojlim E/E(-mH_k)$.

\begin{prop}
\label{prop;26.2.21.20}
There exists the decomposition
\begin{equation}
\label{eq;26.2.21.13}
 E_{|\Hhat_k}
=E_{\Hhat_k}^{\reg}\oplus
E_{\Hhat_k}^{\irreg}
\end{equation}
which induces {\rm(\ref{eq;26.2.21.12})}
for any $x\in H_k$.
 Moreover, we have
 $\nabla(E_{\Hhat_k}^{\reg})
 \subset
 E_{\Hhat_k}^{\reg}
 \otimes
 \Bigl(
 \Omega^1_X(\log H_k)
+\Omega^1_X(\ast H_{\neq k})\Bigr)$
and 
$\nabla(E_{\Hhat_k}^{\irreg})
\subset
E_{\Hhat_k}^{\irreg}
\otimes
\Omega^1_X(\ast H)$. 
\end{prop}
The second claim follows from (\ref{eq;26.2.21.14}).
We shall study the first claim
for the unramifiedly good case
in \S\ref{subsection;26.2.21.21},
and for the general case
in \S\ref{subsection;26.2.21.22},

\subsubsection{Proof of Proposition \ref{prop;26.2.21.20}
in the unramified case}
\label{subsection;26.2.21.21}

Suppose that $(V,\nabla)$ and $E$ are unramifiedly good.
Let $x\in H_k$.
Let $U\subset X$ be a Zariski open neighbourhood of $x$
with an etale coordinate system
$(x_1,\ldots,x_n)$
such that
the ideal of $H$ at $x$ is generated by
$\prod_{i=1}^{\ell}x_i$
and that the ideal of $H_k$ at $x$ is generated by $x_1$.
We set $H_{U,k}=U\cap H_k$.
We also assume that $E_{|U}$ is $\nbigo_U$-free.

We set
$\nbiga^{(0)}=H^0(H_{U,k},\nbigo_{H_k}(\ast \del H_k))$,
where $\del H_k=H_k\cap H_{\neq k}$.
We set
\[
 \nbigr^{(0)}
 =\varprojlim \nbiga^{(0)}[x_1]
 /x_1^m\nbiga^{(0)}[x_1]
\simeq\varprojlim
\nbigo_{X}(\ast H_{\neq k})(U)/x_1^m\nbigo^h_{X}(\ast H_{\neq k})(U).
\]
We also set
$\nbigm^{(0)}=\nbigr^{(0)}\otimes_{\nbigo(U)} E(U)$.
It is equipped with the induced integrable connection
$\nabla:\nbigm^{(0)}\to\nbigm^{(0)}\otimes\Omega^1_X(\ast H)$.
There exists a positive integer $N$ such that
$x_1^{N+1}\nabla_{x_1}\nbigm^{(0)}\subset\nbigm^{(0)}$.
We obtain the $\cnum$-linear map
$G^{(0)}:\nbigm^{(0)}\to\nbigm^{(0)}$ induced by $x_1^{N+1}\nabla_{x_1}$
satisfying
$G^{(0)}(fs)=
x_1^{N+1}\del_{x_1}(f)s+fG^{(0)}(s)$
for $f\in\nbigr^{(0)}$ and $s\in\nbigm^{(0)}$.
It induces the $\nbiga^{(0)}$-endomorphism $F^{(0)}$ of
$\nbigm^{(0)}/x_1\nbigm^{(0)}$.
\begin{lem}
\label{lem;26.2.22.4}
There exist a finitely generated
Galois extension $\nbiga^{(0)}\subset\nbiga^{(1)}$,
a finite subset
$\nbigj_1\subset \nbiga^{(1)}$
and the decomposition
\begin{equation}
\label{eq;26.2.22.2}
 \nbiga^{(1)}\otimes_{\nbiga^{(0)}}
 (\nbigm^{(0)}/x_1\nbigm^{(0)})
 =\bigoplus_{\gminic\in\nbigj_1}
 \Bigl(
  \nbiga^{(1)}\otimes_{\nbiga^{(0)}}
 (\nbigm^{(0)}/x_1\nbigm^{(0)})
 \Bigr)_{\gminic}
\end{equation}
such that (i) $F^{(0)}-\gminic=0$ on 
$\Bigl(
\nbiga^{(1)}\otimes_{\nbiga^{(0)}}
(\nbigm^{(0)}/x_1\nbigm^{(0)})
\Bigr)_{\gminic}$,
(ii) for any two distinct elements $\gminic_1,\gminic_2\in\nbigj_1$,
 the difference $\gminic_1-\gminic_2$ is invertible
 in $\nbiga^{(1)}$,
(iii) there exists an open neighbourhood $x\in U'\subset H_k$
such that
 $\Spec(\nbiga^{(1)})\times_UU'
 \to
 \Spec(\nbiga^{(0)})\times_UU'$ is etale.
\end{lem}
\pf
We have the decomposition (\ref{eq;26.2.22.1})
of $E_{|\xhat}$.
We have the map
$\nu:\nbigi\to
\nbigohat_{X,x}(\ast H_{\neq k})\big/
x_1\nbigohat_{X,x}(\ast H_{\neq k})
=\nbigohat_{H_k,x}(\ast \del H_k)$
by considering 
$\gminia\mapsto x_1^{N+1}\del_{x_1}\gminia$.
Let $\nbigjhat_1$ denote the image of $\nu$.
For $\gminic\in\nbigjhat_1$, we set
$E_{\xhat,\gminic}
 =\bigoplus_{\nu(\gminia)=\gminic}
 (E_{|\xhat}\cap V_{\xhat,\gminia})$.
We have the decomposition
$E_{|\xhat}=\bigoplus _{\gminic\in\nbigjhat_1}
 E_{\xhat,\gminic}$,
and
$x_1^{N+1}\nabla_{x_1}E_{\xhat,\gminic}
\subset E_{\xhat,\gminic}\otimes \nbigo_X(\ast H_{\neq k})_x$.
For two distinct $\gminic_1,\gminic_2\in\nbigjhat_1$,
$\gminic_1-\gminic_2$ is invertible.
Let $\Fhat^{(0)}$ be the endomorphism of
$(E_{|\xhat}/x_1 E_{|\xhat})\otimes\nbigo_{H_k}(\ast \del H_k)_x$.
We have the decomposition
\begin{equation}
\label{eq;26.2.22.3}
 (E_{|\xhat}/x_1 E_{|\xhat})\otimes\nbigo_{H_k}(\ast \del H_k)_x
 =\bigoplus_{\gminic\in\nbigjhat_1}
 (E_{\xhat,\gminic}/x_1)\otimes\nbigo_{H_k}(\ast \del H_k)_x,
\end{equation}
and 
$\Fhat^{(0)}-\gminic=0$
on 
$(E_{\xhat,\gminic}/x_1)\otimes\nbigo_{H_k}(\ast \del H_k)_x$.

Let $P^{(0)}$ denote the characteristic polynomial of $F^{(0)}$.
We fix an algebraically closed field which contains $\nbiga^{(0)}$.
Let $\Sp(P^{(0)})$ denote the set of the solutions of $P^{(0)}$.
Let $\nbiga^{(0)}\subset\nbiga^{(1)}$
denote a finitely generated Galois extension such that
(i) $\Sp(P^{(0)})\subset \nbiga^{(1)}$,
(ii) for any two distinct elements
$\gminic_1,\gminic_2\in \Sp(P^{(0)})$,
$\gminic_1-\gminic_2\in\nbiga^{(1)}$.
We obtain the decomposition (\ref{eq;26.2.22.2}),
which is preserved by $F^{(0)}$,
such that
$F^{(0)}-\gminic$ are nilpotent on
$\Bigl(
\nbiga^{(1)}\otimes_{\nbiga^{(0)}}
(\nbigm^{(0)}/x_1\nbigm^{(0)})
\Bigr)_{\gminic}$.
By the comparison with (\ref{eq;26.2.22.3}),
we obtain that
the conditions (i) and (iii) are also satisfied.
\hfill\qed

\vspace{.1in}
We set
$\nbigr^{(1)}=
\varprojlim \nbiga^{(1)}[x_1]/x^m\nbiga^{(1)}[x_1]$
and
$\nbigm^{(1)}=\nbiga^{(1)}\otimes_{\nbigo(U)} E(U)$.
It is equipped with the induced integrable connection
$\nabla:\nbigm^{(1)}\to \nbigm^{(1)}\otimes\Omega^1_X(\ast H)$.
It is standard to obtain the unique decomposition
$\nbigm^{(1)}
=\bigoplus_{\gminic\in\nbigj_1}
 \nbigm^{(1)}_{\gminic}$
such that the following holds.
\begin{itemize}
 \item The decomposition is preserved by $G^{(0)}$.
 \item $\nbigm^{(1)}_{\gminic}/x_1\nbigm^{(1)}_{\gminic}
       =\Bigl(
       \nbiga^{(1)}\otimes_{\nbiga^{(0)}}(\nbigm^{(0)}/x_1\nbigm^{(0)})
       \Bigr)_{\gminic}$.
 \item $x_1^{-1}(x_1^{N+1}\nabla_{x_1}-\gminic)\nbigm^{(1)}_{\gminic}
       \subset
        \nbigm^{(1)}_{\gminic}$.
 \item $\nabla(\nbigm^{(1)}_{\gminic})
       \subset
        \nbigm^{(1)}_{\gminic}\otimes\Omega^1(\ast H)$.
\end{itemize}

\begin{lem}
Inductively,
for $j=1,\ldots,N$,
we can construct the following tuples.
\begin{itemize}
 \item  Finitely generated Galois extensions
	$\nbiga^{(j)}$ of $\nbiga^{(0)}$.
	We set
	$\nbigr^{(j)}=
	\varprojlim \nbiga^{(j)}[x_1]/x^m\nbiga^{(j)}[x_1]$
	and
	$\nbigm^{(j)}=\nbiga^{(j)}\otimes_{\nbigo(U)} E(U)$.
	It is equipped with the induced integrable connection
	$\nabla:\nbigm^{(j)}\to \nbigm^{(j)}\otimes\Omega^1_X(\ast H)$.
 \item 	Finite subsets
	$\nbigj_j\subset \nbigr^{(j)}/x_1^{j}\nbigr^{(j)}$,
	and decompositions
       $\nbigm^{(j)}=\bigoplus_{\gminic\in\nbigj_j}
	\nbigm^{(j)}_{\gminic}$.
\end{itemize}
They satisfy the following conditions.
\begin{itemize}
 \item The projection
       $\nbigr_x/x_1^{j+1}\nbigr_x\to
       \nbigr_x/x_1^{j}\nbigr_x$
       induces the surjection
       $\nbigj_j\to\nbigj_{j-1}$.
 \item 
       If two distinct elements $\gminic_1,\gminic_2\in\nbigj_j$
       are mapped to the same element in $\nbigj_{j-1}$,
       then
       $x_1^{-j}(\gminic_1-\gminic_2)$ 
       is invertible in
       $\nbiga^{(j)}$.
 \item $x_1^{-j}(x_1^{N+1}\nabla_{x_1}-\gminic)
       \nbigm^{(j+1)}_{\gminic}\subset\nbigm^{(j+1)}_{\gminic}$.
 \item $\nabla(\nbigm^{(j)}_{\gminic})
       \subset
       \nbigm^{(j)}_{\gminic}\otimes\Omega^1(\ast H)$.
 \item There exists an open neighbourhood $x\in U'\subset U$
       such that
       $\Spec\nbiga^{(j)}\times_UU'
       \to
       \Spec\nbiga^{(0)}\times_UU'$
       is etale.      
\end{itemize}
In particular,
we have
$(x_1\nabla_{x_1}-x_1^{-N}\gminic)\nbigm^{(N)}_{\gminic}
\subset
\nbigm^{(N)}_{\gminic}$
for any $\gminic\in\nbigj_{N-1}$.
\end{lem}
\pf
Suppose we have already constructed
$\nbiga^{(j)}$, $\nbigj_j$ and
$\nbigm^{(j)}=\bigoplus\nbigm^{(j)}_{\gminic}$.
We obtain the $\hyperk$-linear map
$G_{\gminic}^{(j)}:
\nbigm^{(j)}_{\gminic}
\to \nbigm^{(j)}_{\gminic}$
induced by $x_1^{-j}(x_1^{N+1}\nabla_{x_1}-\gminic)$.
It satisfies
$G_{\gminic}^{(j)}(fs)=
fG_{\gminic}^{(j)}(s)+x_1^{N+1-j}\del_{x_1}(f)s$
for $f\in \nbigr^{(j)}$ and $s\in\nbigm^{(j)}_{\gminic}$.
It induces the $\nbiga^{(j)}$-endomorphism $F^{(j)}_{\gminic}$
of $\nbigm^{(j)}_{\gminic}/x_1\nbigm^{(j)}_{\gminic}$.
As in Lemma \ref{lem;26.2.22.4},
there exist a finitely generated
Galois extension $\nbiga^{(j+1)}_{\gminic}$ of $\nbiga^{(j)}$,
a finite subset
$\nbigj_{\gminic,j+1}\subset \nbiga^{(j+1)}_{\gminic}$
and the decomposition
\begin{equation}
 \nbiga_{\gminic}^{(j+1)}\otimes_{\nbiga^{(j)}}
 (\nbigm^{(j)}_{\gminic}/x_1\nbigm^{(j)}_{\gminic})
 =\bigoplus_{\gminib\in\nbigj_{\gminic,j+1}}
 \Bigl(
  \nbiga^{(j+1)}_{\gminic}\otimes_{\nbiga^{(j)}}
 (\nbigm^{(j)}_{\gminic}/x_1\nbigm^{(j)}_{\gminic})
 \Bigr)_{\gminib}
\end{equation}
such that (i) $F_{\gminic}^{(j)}-\gminib=0$ on 
$\Bigl(
\nbiga_{\gminic}^{(j+1)}\otimes_{\nbiga^{(j)}}
(\nbigm^{(j)}_{\gminic}/x_1\nbigm^{(j)}_{\gminic})
\Bigr)_{\gminic}$,
(ii)for any two distinct elements
$\gminib_1,\gminib_2\in\nbigj_{\gminic,j+1}$,
the difference $\gminib_1-\gminib_2$ is invertible
in $\nbiga_{\gminic}^{(j+1)}$,
(iii) there exists a neighbourhood $x\in U'\subset U$
such that
$\Spec \nbiga^{(j+1)}_{\gminic}\times_UU'
\to \Spec\nbiga^{(0)}\times_UU'$
is etale.
Let $\nbiga^{(j+1)}$ be 
a finitely generated
Galois extension of $\nbiga^{(0)}$,
which contains all $\nbiga^{(j+1)}_{\gminic}$.
We set
$\nbigm^{(j+1)}_{\gminic}=\nbiga^{(j+1)}\otimes_{\nbigo(U)} E(U)$,
equipped with the induced integrable connection
$\nabla:\nbigm^{(j+1)}\to \nbigm^{(j+1)}\otimes\Omega^1(\ast H)$.
It is standard to obtain the unique decomposition
$\nbigm^{(j+1)}_{\gminic}
=\bigoplus_{\gminib\in\nbigj_{\gminic,j+1}}
 \nbigm^{(j+1)}_{\gminic+x_1^j\gminib}$
such that the following holds.
\begin{itemize}
 \item The decomposition is preserved by $G^{(j)}$.
 \item $\nbigm^{(j+1)}_{\gminic+x_1^j\gminib}
       /x_1\nbigm^{(j+1)}_{\gminic+x_1^j\gminib}
       =\Bigl(
       \nbiga^{(j+1)}\otimes_{\nbiga^{(j)}}
       (\nbigm^{(j)}_{\gminic}/x_1\nbigm^{(j)}_{\gminic})
       \Bigr)_{\gminib}$.
 \item $x_1^{-j}(x_1^{N+1}\nabla_{x_1}-(\gminic+x_1^{j+1}\gminib))
       \nbigm^{(1)}_{\gminic+x_1^j\gminib}
       \subset
        \nbigm^{(1)}_{\gminic+x_1^j\gminib}$.
 \item $\nabla(\nbigm^{(j+1)}_{\gminic+x_1^j\gminib})
       \subset
        \nbigm^{(j+1)}_{\gminic+x_1^j\gminib}\otimes\Omega^1(\ast H)$.
\end{itemize}
We set
$\nbigj_{j+1}=\bigsqcup_{\gminic\in \nbigj_j}
\{\gminic+x_1^j\gminib\,|\,\gminib\in\nbigj_{\gminic,j+1}\}$.
Thus, the inductive construction can proceed.
\hfill\qed

\vspace{.1in}

We set $\nbigm^{(N)\reg}=\nbigm^{(N)}_0$
and $\nbigm^{(N)\irreg}=
\bigoplus_{\gminic\neq 0}
\nbigm^{(N)}_{\gminic}$.
We obtain the decomposition
\begin{equation}
\label{eq;26.2.20.10}
\nbigm^{(N)}=
\nbigm^{(N)\reg}\oplus\nbigm^{(N)\irreg}.
\end{equation}
The decomposition is preserved by $\nabla$.

By shrinking $U$,
we may assume that
$\Spec(\nbiga^{(N)})\to (U\cap H_k)\setminus\del H_k$
is finite and etale.
Because the decomposition (\ref{eq;26.2.20.10})
is invariant under the Galois extension
$\nbiga^{(0)}\to\nbiga^{(N)}$,
there exists the decomposition
\begin{equation}
\label{eq;26.2.22.5}
\nbigm^{(0)}=\nbigm^{(0)\reg}\oplus\nbigm^{(0)\irreg}
\end{equation}
as the descent of (\ref{eq;26.2.20.10}).
We can check the following.
\begin{equation}
\label{eq;26.3.2.1}
 \nabla \nbigm^{(0)\reg}
  \subset
  \nbigm^{(0)\reg}
  \otimes
  \Bigl(
  \Omega^1(\log H_k)
  +
  \Omega^1(\ast H_{\neq k})
  \Bigr),
  \quad\quad
  \nabla \nbigm^{(0)\irreg}
  \subset
  \nbigm^{(0)\irreg}
 \otimes
  \Omega^1(\ast H).
\end{equation}
We can also check that
\begin{equation}
\label{eq;26.3.2.2}
\sum_{j\geq 0}
(x_1\nabla_{x_1})^j
(\nbigm^{(0)\irreg})_{\yhat}
=(\nbigm^{(0)\irreg})_{\yhat}(\ast H_k),
\quad
\sum_{j\geq 0}
(x_1\nabla_{x_1})^j
 (E_{|\xhat})^{\irreg,k}
=(E_{|\xhat})^{\irreg,k}.
\end{equation}

\begin{lem}
There exists the decomposition
$E_{|\Hhat_k}(U)
=(E_{|\Hhat_k})^{\reg}(U)
 \oplus
 (E_{|\Hhat_k})^{\irreg}(U)$,
 which induces
 {\rm(\ref{eq;26.2.21.12})} and
{\rm(\ref{eq;26.2.22.5})}.
It also induces the decompositions 
{\rm(\ref{eq;26.2.21.12})}
for any $y\in U\cap H_k$.
\end{lem}
\pf
We set
$\nbigrtilde=\nbigohat_{H_k,x}(\ast \del H_k)[\![x_1]\!]$
which contains
both 
$\nbigr^{(0)}\cap \nbigohat_{X,x}$ and
$\nbigo_{H_k}(H_{U,k})[\![x_1]\!]$.
There exists the natural isomorphism
$\nbigrtilde\otimes  E_{|\yhat}
\simeq
\nbigrtilde\otimes \nbigm^{(0)}$.
By the properties (\ref{eq;26.2.21.14}),
(\ref{eq;26.3.2.1})
and (\ref{eq;26.3.2.2}),
we obtain
\[
\nbigrtilde\otimes 
(E_{|\xhat})^{\reg,k}
=
\nbigrtilde\otimes 
\nbigm^{(0)\reg},
\quad
\nbigrtilde\otimes 
(E_{|\xhat})^{\irreg,k}
=
\nbigrtilde\otimes 
\nbigm^{(0)\irreg}.
\]
Because
$\nbigr^{(0)}\cap \nbigohat_{X,x}
=\nbigo_{H_k}(H_{U,k})[\![x_1]\!]$
in $\nbigrtilde$,
we obtain the first claim.
Applying the same argument to $y\in H_{U,k}$,
we obtain the second claim.
\hfill\qed

\vspace{.1in}
Thus, we obtain the decomposition (\ref{eq;26.2.21.13}).

\subsubsection{Proof of Proposition \ref{prop;26.2.21.20} in the general case} 
\label{subsection;26.2.21.22}

Let $H_{U,k}^{(e)}$ denote the inverse image of $H_k$
by $U^{(e)}\to X$.
We have the $\Gal^{(H,e)}_U$-equivariant decomposition
\begin{equation}
\label{eq;26.2.20.20}
 (E_U^{(e)})_{|\Hhat_{U,k}^{(e)}}
=(E_U^{(e)})^{\reg}_{\Hhat_{U,k}^{(e)}}
 \oplus
 (E_U^{(e)})^{\irreg}_{\Hhat_{U,k}^{(e)}}.
\end{equation}
Let $H_{U,k}=H_k\cap U$.
Because $E_{|U}$ is the descent of $E^{(e)}_U$,
we obtain the following decomposition
as the descent of (\ref{eq;26.2.20.20}):
\[
 (E_U)_{|\Hhat_{U,k}}
=(E_{|U})^{\reg}_{\Hhat_{U,k}}
 \oplus
 (E_{|U})^{\irreg}_{\Hhat_{U,k}}.
\]
Then, we obtain
(\ref{eq;26.2.21.12}).
\hfill\qed

\subsubsection{Residues}

If $E$ is unramifiedly good,
for $k\in\Lambda$ such that $x\in H_k$,
we obtain the endomorphism 
$\Res_{H_k}(\nabla)$ of $(E_{|H_k})_{|\xhat}$
as the residue of
$\bigoplus (\nabla_{\gminia}-d\gminia)$
along $H_k$.

\begin{prop}
\label{prop;26.2.21.11}
If $E$ is unramifiedly good,
there exists the endomorphism $\Res_{H_k}(\nabla)$ of $E_{|H_k}$
which induces
$\Res_{H_k}(\nabla)$ of $(E_{|H_k})_{\xhat}$
for any $x\in H_k$. 
\end{prop}
\pf
We continue to use the notation in the proof of
Proposition \ref{prop;26.2.21.20}.
We obtain the endomorphism $\Res_{H_k}(\nabla)$
of $\nbigm^{(N)}/x_1\nbigm^{(N)}$
as 
$\bigoplus_{\gminic}
\Res_{x_1}
\bigl(
x_1\nabla_{x_1}-x_1^{-N}\gminic\id_{\nbigm^{(N)}_{\gminic}}
\bigr)$.
Because it is invariant under the action of the Galois group,
there exists
$\Res_{H_k}(\nabla)$ on $\nbigm^{(0)}/x_1\nbigm^{(0)}$
which induces
$\Res_{H_k}(\nabla)$ on $\nbigm^{(N)}/x_1\nbigm^{(N)}$.
We obtain the claim of Proposition \ref{prop;26.2.21.11}
from the following lemma.

\begin{lem}
\label{lem;26.2.22.6}
There exists the $\nbigo_{H_k}(H_{U,k})$-endomorphism $\Res_{H_k}(\nabla)$
of $E(U)/x_1E(U)$
which induces $\Res_{H_k}(\nabla)$
on $\nbigm^{(0)}/x_1\nbigm^{(0)}$.
\end{lem}
\pf
By using the decomposition (\ref{eq;26.2.21.12}),
we obtain
the $\nbigohat_{X,x}/x_1\nbigohat_{X,x}$-endomorphisms
$\Res_{x_1}(\nabla_{\gminia}-d\gminia)$
of $(E_{|\xhat}\cap V_{\xhat,\gminia})\big/
x_1(E_{|\xhat}\cap V_{\xhat,\gminia})$.
We obtain the $\nbigohat_{X,x}/x_1\nbigohat_{X,x}$-endomorphisms
$\Res_{x_1}(\nabla)
 =\bigoplus \Res_{x_1}(\nabla_{\gminia}-d\gminia)$
of $E_{|\xhat}/x_1E_{|\xhat}$.
By comparing two residue endomorphisms,
we obtain that
there exists the endomorphism $\Res_{H_k}(\nabla)$
on $E_{|H_k}(U)$
which induces 
$\Res_{H_k}(\nabla)$
of $\nbigm^{(0)}/x_1\nbigm^{(0)}$.
Thus, we obtain Lemma \ref{lem;26.2.22.6}
and Proposition \ref{prop;26.2.21.11}.
\hfill\qed

\vspace{.1in}

Let us consider the case where $E$ is good but not necessarily
unramifiedly good.
By (\ref{eq;26.2.21.14}),
we obtain the endomorphism
$\Res_{H_k}(\nabla)$
of $(E^{\reg}_{\Hhat_k|H_k})_{|\xhat}$.

\begin{prop}
Let $E$ be good lattice.
We obtain the residue endomorphism
$\Res_{H_k}(\nabla)$
on $(E_{\Hhat_k}^{\reg})_{|H_k}$,
which induces $\Res_{H_k}(\nabla)$
on $(E^{\reg}_{\Hhat_k|H_k})_{|\xhat}$
for any $x\in H_k$.
\hfill\qed
\end{prop}

Note that the set of the eigenvalues of $\Res_{H_k}(\nabla)_x$
are independent of $x\in H_k$.

\begin{notation}
Let $\Sp(E,\nabla,H_k)^{\reg}\subset\hyperk$
denote the set of the eigenvalues of
$\Res_{H_k}(\nabla)$
on $(E_{\Hhat_k}^{\reg})_{|H_k}$. 
\end{notation}

\subsection{Good $\cnum$-Deligne-Malgrange filtered bundles}

\subsubsection{The case $\hyperk\subset\cnum$}

Let $\hyperk\subset\cnum$ be an algebraically closed field.
Let $X$ be a smooth projective $\hyperk$-variety
with a simple normal crossing hypersurface $H$.
Let $(V,\nabla)$ be a good meromorphic flat bundle on $(X,H)$.

Let $(U=\Spec R,\varphi,\psi)$
be an etale coordinate system as in \S\ref{subsection;25.12.7.1}.
For $e\in \seisuu_{>0}$,
let $R^{(e)}=R[x_1^{1/e},\ldots,x_{\ell}^{1/e}]$
and $U^{(e)}=\Spec R^{(e)}$.
Let $\rho^{(e)}:U^{(e)}\to U$ denote the induced morphism.
Let $\Gal_e$ denote the Galois group of the ramified covering.
For any $e\in (\rank V)!\seisuu$,
we set
$(V^{(e)},\nabla)=\rho_e^{\ast}(V,\nabla)$,
which is $\Gal_e$-equivariant
and unramifiedly good
on $(U^{(e)},H^{(e)})$.

By the extension $\hyperk\subset\cnum$,
we obtain the unramifiedly good meromorphic flat bundle
$(V^{(e)}_{\cnum},\nabla)$
on $(U^{(e)}_{\cnum},H^{(e)}_{\cnum})$.
Let $E$ be any unramifiedly good lattice of $(V^{(e)}_{\cnum},\nabla)$.
There exists the residue endomorphisms
$\Res_i(\nabla)$ $(i=1,\ldots,\ell)$ of $E_{|H^{(e)}_{\cnum,i}}$.
The eigenvalues of $\Res_i(\nabla)_{|y}$
are independent of $y\in H_{\cnum,i}^{(e)}$.

\begin{lem}
The eigenvalues of $\Res_i(\nabla)_{|y}$
are contained in $\hyperk$.
\end{lem}
\pf
Let $z\in
H^{(e)}_{i}\setminus \bigcup_{j\neq i}H^{(e)}_{j}$.
There exists an unramifiedly good lattice
$L\subset V_{|\zhat}$
for which we have the decomposition
\[
 (L,\nabla)_{|\zhat}
 =\bigoplus_{\gminia\in\nbigi(V,z)}
 (L_{\gminia},\nabla_{\gminia})
\]
such that $\nabla_{\gminia}-d\gminiatilde$
are logarithmic.
We obtain the residues
$\Res_i(\nabla)=
\bigoplus_{\gminia\in\nbigi}
\Res_i(\nabla_{\gminia}-d\gminiatilde)$
of $L_{|z}$.

For any eigenvalue $\alpha$ of $\Res_{i}(\nabla)$ on $E_{|z}$,
there exists an eigenvalue $\beta$ of $\Res_{i}(\nabla)$ on $L_{|z}$
such that $\alpha-\beta\in\seisuu$.
Conversely,
for any eigenvalue $\beta$ of $\Res_{i}(\nabla)$ on $L_{|z}$,
there exists an eigenvalue $\alpha$ of $\Res_{i}(\nabla)$ on $E_{|z}$
such that $\alpha-\beta\in\seisuu$.
\hfill\qed

\vspace{.1in}

For any $\veca\in\real^{\ell}$,
there exists the unramifiedly good lattice
$\nbigp_{e\veca}^{\DM}(V^{(e)}_{\cnum})$
determined by the following condition.
\begin{itemize}
 \item Any eigenvalue $\alpha$
       of $\Res_i(\nabla)$
       satisfies
       $ea_i-1<\Re(\alpha)\leq ea_i$.
\end{itemize}
Because it is $\Gal_e$-equivariant,
we obtain the good lattice
$\nbigp_{\veca}^{\DM}(V_{\cnum})$ of $V_{\cnum|U_{\cnum}}$
as the descent of
$\nbigp_{e\veca}^{\DM}(V^{(e)}_{\cnum})$.
We also have the invariance of
$\nbigp_{e\veca}^{\DM}(V^{(e)}_{\cnum})$
under the action of the Galois group of $\cnum/\hyperk$,
and hence
there exists a lattice
$\nbigp_{e\veca}^{\DM}(V^{(e)})$ of $V^{(e)}$
which induces
$\nbigp_{e\veca}^{\DM}(V^{(e)}_{\cnum})$.
We also have the lattice
$\nbigp^{\DM}_{\veca}(V)_U$ of $V_{|U}$.

Let $H=\bigcup_{i\in\Lambda}H_i$ be the irreducible decomposition.
By patching $\nbigp^{\DM}_{\veca}(V)_U$ for any $\veca\in\real^{\Lambda}$,
we obtain the lattices $\nbigp^{\DM}_{\veca}(V)$.

\subsubsection{General case}
\label{subsection;26.2.4.1}

Let $\hyperk$ be an algebraically closed field of characteristic $0$,
which is not necessarily contained in $\cnum$.
Let $X$ be a smooth projective $\hyperk$-variety
with a simple normal crossing hypersurface $H$.
Let $(V,\nabla)$ be a good meromorphic flat bundle on $(X,H)$.

There exists a subfield  $\hyperk_0\subset\hyperk$
such that
(i) $\hyperk_0$ is the algebraic closure of a field 
finitely generated over $\rnum$,
(ii) $(X,H)$ and $(V,\nabla)$ are defined over $\hyperk_0$,
i.e.,
there exist $(X_{\hyperk_0},H_{\hyperk_0})$ over $\hyperk_0$
and $(V_{\hyperk_0},\nabla)$ on $(X_{\hyperk_0},H_{\hyperk_0})$
such that they induce $(X,H)$ and $(V,\nabla)$.

There exists an embedding
$\hyperk_0\subset \cnum$.
Depending on the embedding,
there exists a filtered bundle
$\nbigp^{\DM}_{\ast}(V_{\hyperk_0})$
indexed by $\real^{\Lambda}$.
It induces a filtered bundle
$\nbigp^{\DM}_{\ast}(V)$
indexed by $\real^{\Lambda}$.

\subsubsection{Pull back}

Let $Y$ be a smooth projective variety over $\hyperk$
with a simple normal crossing hypersurface $H_Y$.
Let $f:Y\to X$ be a morphism such that $f^{-1}(H)=H_Y$.
Recall the following.
\begin{lem}
$(V_Y,\nabla)=f^{\ast}(V,\nabla)$
is a good meromorphic flat bundle on $(Y,H_Y)$. 
\end{lem}
\pf
Let $y\in Y$ be a closed point.
We set $x=f(y)$.
For some $e\in\seisuu_{>0}$,
there exist
$\nbigi(V,x)\subset
\nbigohat^{(H,e)}_{X}(\ast H)_x/
\nbigohat^{(H,e)}_{X,x}$
and a decomposition
\[
 (V,\nabla)\otimes\nbigohat^{(H,e)}_x
=\bigoplus_{\gminia\in\nbigi(V,x)}
 (\Vhat_{\gminia},\nabla_{\gminia})
\]
such that
each $(\Vhat_{\gminia},\nabla_{\gminia}-d\gminiatilde)$
has a logarithmic lattice.
There exists $e'$
such that the morphism
$\nbigo_{X,x}\to \nbigo_{Y,y}$
induces
$\nbigo^{(H,e)}_{X,x}\to \nbigo^{(H_Y,e')}_{Y,y}$.
We obtain
$\nbigi(V_Y,y)\subset
\nbigo^{(H_Y,e')}_{Y}(\ast H_Y)_y/
\nbigo^{(H_Y,e')}_{Y,y}$
and the decomposition
\[
(V_Y,\nabla)\otimes\nbigohat^{(H_Y,e')}_{Y,y}
=\bigoplus_{\gminia\in\nbigi(V,x)}
 (\Vhat_{Y,\gminia},\nabla_{\gminia})
\]
such that
each 
$(\Vhat_{Y,\gminia},\nabla_{\gminia}-d\gminiatilde)$
has a logarithmic lattice.
\hfill\qed

\vspace{.1in}

Let $H_Y=\bigcup_{i\in\Lambda_Y}H_{Y,i}$
denote the irreducible decomposition.
For any $i\in \Lambda$,
let $\vecb^{(i)}\in\seisuu_{\geq 0}^{\Lambda_Y}$
be determined by the pull back
$f^{\ast}(H_i)=\sum_{j\in\Lambda_Y} b^{(i)}_jH_{Y,j}$
as the divisor.
The matrix $(b^{(i)}_j)$ 
determines the map
$f^{\ast}:\real^{\Lambda}\to \real^{\Lambda_Y}$.
Let $\vecdelta_{\Lambda}=(1,\ldots,1)\in\real^{\Lambda}$.

\begin{lem}
For any $\veca\in\real^{\Lambda}$,
we have
$\nbigp^{\DM}_{f^{\ast}(\veca-\vecdelta_{\Lambda})}(f^{\ast}V)
 \subset
 f^{\ast}
 \nbigp^{\DM}_{\veca}(V)
\subset
\nbigp^{\DM}_{f^{\ast}(\veca)}(f^{\ast}V)$.
\end{lem}
\pf
It is enough to consider the unramified case
and $\hyperk=\cnum$.
Let $y\in H_Y$ and $x=f(y)$.
Let $(y_1,\ldots,y_m)$ be an etale coordinate around $y$
such that
$H_Y=\bigcup_{i=1}^k\{y_i=0\}$.
Let $(x_1,\ldots,x_n)$ be an etale coordinate neighbourhood around $x$
such that
$H=\bigcup_{i=1}^{\ell}\{x_i=0\}$.
Let $\alpha_{i,j}$ $(1\leq j\leq n,\,\,1\leq i\leq m)$
be the non-negative integers determined by
\[
 f^{\ast}(dx_j/x_j)
=\sum \alpha_{i,j}dy_i/y_i. 
\]
Let $\beta_i\in\hyperk$
such that $a_j-1<-\Re(\beta_j)\leq a_i$.
We obtain
\[
 \sum \alpha_{i,j}(a_j-1)
 <-\sum \alpha_{i,j}\Re(\beta_j)
 \leq
 \sum \alpha_{i,j}a_j.
\]
If $\veca\in\real_{<0}^{\Lambda}$,
we have
$-\sum \alpha_{i,j}\Re(\beta_j)<0$.
If $\veca\in\real_{\geq 1}^{\Lambda}$,
we have 
$0<-\sum \alpha_{i,j}\Re(\beta_j)$.
\hfill\qed

\begin{cor}
For $\veca\in\real_{<0}^{\Lambda}$,
we have
$f^{\ast}\nbigp^{\DM}_{\veca}V
\subset \nbigp^{\DM}_{<0}(f^{\ast}V)$.
For $\veca\in\real_{\geq 1}^{\Lambda}$,
we have
$\nbigp^{\DM}_{0}(f^{\ast}V)\subset f^{\ast}\nbigp^{\DM}_{\veca}V$.
\end{cor}

\subsection{The associated meromorphic Lagrangian irregularity}
\label{subsection;25.10.15.20}

Let $X_{\hyperk}$ be a normal projective $\hyperk$-variety
with a hypersurface $H_{\hyperk}$
such that $X_{\hyperk}\setminus H_{\hyperk}$
is smooth.
Let $(V,\nabla)$ be a meromorphic integrable connection on
$(X_{\hyperk},H_{\hyperk})$,
i.e.,
$V$ is a coherent reflexive $\nbigo_{X_{\hyperk}}(\ast H_{\hyperk})$-module
with an integrable connection
$\nabla:V\to V\otimes\Omega_{X_{\hyperk}}^1$.
The following theorem is proved in \cite{kedlaya2}.
If $\hyperk=\cnum$,
it is also proved in \cite{Mochizuki-wild}.

\begin{thm}
There exists a projective birational morphism
$\rho_{\hyperk}:X'_{\hyperk}\to X_{\hyperk}$
such that
(i) $X'_{\hyperk}$ is smooth,
(ii) $H'_{\hyperk}=\rho_{\hyperk}^{-1}(H_{\hyperk})$ is normal crossing,
(iii) $X'_{\hyperk}\setminus H'_{\hyperk}\simeq
X_{\hyperk}\setminus H_{\hyperk}$,
(iv) $\rho_{\hyperk}^{\ast}(V,\nabla)$ is good
on $(X'_{\hyperk},H'_{\hyperk})$.
\hfill\qed
\end{thm}

\begin{thm}
\label{thm;25.12.8.100}
There exists the unique irregularity of
a meromorphic Lagrangian cover $\vecI(V,\nabla)$
such that the following holds.
\begin{itemize}
 \item Let $\rho:X'_{\hyperk}\to X_{\hyperk}$
       be a projective morphism of
       smooth $\hyperk$-varieties
       such that
       (i) $H'_{\hyperk}=\rho^{-1}(H_{\hyperk})$
       is normal crossing,
       (ii) $X'_{\hyperk}\setminus H'_{\hyperk}\simeq
       X_{\hyperk}\setminus H_{\hyperk}$.
       Then,
       $\rho^{\ast}(V,\nabla)$
       is good on $(X'_{\hyperk},H_{\hyperk'})$
       if and only if
       $\rho^{\ast}\vecI(V,\nabla)$
       is etale locally good on $(X'_{\hyperk},H'_{\hyperk})$.
       In case,
       we have
       $\vecnbigi(\rho^{\ast}(V,\nabla))
       =
       \vecnbigi(\rho^{\ast}\vecI(V,\nabla))$.
\end{itemize}
\end{thm}
\pf
Let us consider the case
where $\hyperk\subset\cnum$.
We obtain $X_{\cnum}$, $H_{\cnum}$
and $(V_{\cnum},\nabla)$
from $X_{\hyperk}$, $H_{\hyperk}$
and $(V_{\hyperk},\nabla)$.
There exists a Lagrangian meromorphic cover
$\Sigma(V_{\cnum})\subset T^{\ast}(X_{\cnum}\setminus H_{\cnum})$
associated with $(V_{\cnum},\nabla)$
as in \S\ref{subsection;25.12.8.30}.
Let
$\rho_{\hyperk}:(X'_{\hyperk},H'_{\hyperk})\to (X_{\hyperk},H_{\hyperk})$
be such that $\rho_{\hyperk}^{\ast}(V_{\hyperk},\nabla)$
is good on $(X'_{\hyperk},H'_{\hyperk})$.
Let $\rho_{\cnum}:(X'_{\cnum},H'_{\cnum})\to (X_{\cnum},H_{\cnum})$
be the induced morphism.
By Corollary \ref{cor;25.12.8.31},
$\rho_{\cnum}^{\ast}(\Sigma(V_{\cnum}))$ is good on
$(X'_{\cnum},H_{\cnum}')$,
and we have
$\vecnbigi(\rho_{\cnum}^{\ast}(\Sigma(V_{\cnum})))
=\vecnbigi(\rho_{\cnum}^{\ast}(V_{\cnum},\nabla))$.
There exists
$\Sigma_{\hyperk}\subset T^{\ast}(X_{\hyperk}\setminus H_{\hyperk})$
such that
$\vecnbigi\bigl(
 (\rho_{\hyperk}^{\ast}\Sigma_{\hyperk})
 \times \cnum
 \bigr)
=\vecnbigi(\rho_{\cnum}^{\ast}\Sigma (V_{\cnum}))$
by Proposition \ref{prop;25.10.14.10}.
By Proposition \ref{prop;25.12.8.40},
Corollary \ref{cor;25.12.8.31}
and Proposition \ref{prop;25.10.13.10},
the claim holds in this case.

In the general case,
there exists an algebraically closed field
$\hyperk_0$ of characteristic $0$
such that
(i) $\hyperk_0$ is an algebraically closure of
a field finitely generated over $\rnum$,
(ii) $\hyperk_0\subset\hyperk$,
(iii) there exists $(X_{\hyperk_0},H_{\hyperk_0},V_{\hyperk_0},\nabla)$
which induces $(X_{\hyperk},H_{\hyperk},V_{\hyperk},\nabla)$.
By the previous consideration,
there exists
the irregularity of meromorphic Lagrangian cover
$\vecI(V_{\hyperk_0},\nabla)$
associated with $(V_{\hyperk_0},\nabla)$.
We set
$\vecI(V_{\hyperk},\nabla)=
\vecI(V_{\hyperk_0},\nabla)\times \hyperk$.
Then,
By Proposition \ref{prop;25.12.8.40},
Corollary \ref{cor;25.12.8.31}
and Proposition \ref{prop;25.10.13.10},
the claim holds in this case.
\hfill\qed

\subsubsection{Comparison with the pole orders of the connections}

Let $X$ be a normal $\hyperk$-variety
with a hypersurface $H=\bigcup_{i\in\Lambda}H_i$.
There exists a closed subset $Z\subset H$
with $\codim_X(Z)\geq 2$
such that
(i) $X^{\circ}=X\setminus Z$ is smooth,
(ii) $H^{\circ}=H\setminus Z$ is a normal crossing hypersurface of $X^{\circ}$.
We also set $H_i^{\circ}=H_i\setminus Z$.

Let $(V,\nabla)$ be a meromorphic integrable connection
on $(X,H)$.
Note that there exists a closed subset
$Z\subset Z'\subset H$ 
with $\codim_XZ'\geq 2$
such that
$(V,\nabla)_{|X\setminus Z'}$
is good on $(X\setminus Z',H\setminus Z')$.

\begin{prop}
\label{prop;25.10.15.51}
For a tuple of non-negative integers
$\vecm=(m_i\,|\,i\in\Lambda)$,
the following conditions are equivalent.
\begin{itemize}
\item There exists a reflexive lattice $V_0\subset V$
such that
$\nabla(V_{0|X^{\circ}})\subset
      V_{0|X^{\circ}}\otimes
      \Omega^1_{X^{\circ}}(\log H^{\circ})(\sum m_iH^{\circ}_i)$.
 \item
      Let $Z\subset Z'$ be any closed subset of $H$
      with $\codim_XZ'\geq 2$
      such that
      $(V,\nabla)_{|X\setminus Z'}$
      is good on $(X\setminus Z',H\setminus Z')$.
      Let $V_0\subset V$ be any coherent reflexive $\nbigo_X$-module
      such that (i) $V_0(\ast H)=V$,
      (ii) $V_{0|X\setminus Z'}$
      is a good lattice of
      $(V,\nabla)_{|X\setminus Z'}$.
      Then, we have
      $\nabla(V_{0|X^{\circ}})
      \subset
      V_{0|X^{\circ}}
      \otimes\Omega^1_{X^{\circ}}(\log H^{\circ})(\sum m_iH_i^{\circ})$.
 \item For any $\Sigma\in\vecI(V,\nabla)$,
       the closure of $\Sigma$ in
       $T^{\ast}X^{\circ}(\log H^{\circ})(\sum m_iH_i^{\circ})$
       is proper over $X^{\circ}$.
\end{itemize}
The conditions are independent of the choice of $Z$.
\end{prop}
\pf
The second condition implies the first.
Let us prove that the first condition implies the second.
It is enough to consider the case $\dim X=1$.
Moreover, it is enough to consider formal and unramified case.
Let $V$ be a finite dimensional
$\hyperk(\!(z)\!)$-vector space
with a connection $\nabla$.
Suppose that there is
a finite subset $\nbigi\subset z^{-1}\hyperk[z^{-1}]$
and a decomposition
\[
 (V,\nabla)
 =\bigoplus_{\gminia\in\nbigi}
 (V_{\gminia},\nabla_{\gminia})
\]
such that $\nabla_{\gminia}-d\gminia$ are regular.
Let $E=\bigoplus E_{\gminia}$ be an unramifiedly good lattice.
Let $N_1=\max\{\deg_{z^{-1}}\gminia\,|\,\gminia\in\nbigi\}$.
Take $\gminia\in\nbigi$ such that
$\deg_{z^{-1}}\gminia=N_1$.
Let $G$ be any lattice of $V$.
There exists $N_2\in\seisuu_{\geq 0}$ such that 
$z\nabla_{\del_z}(G)\subset z^{-N_2}G$.
We set $G_{\gminia}=G\cap E_{\gminia}$.
We have
$z\nabla_{\del_z}G_{\gminia}
\subset
z^{-N_2}G_{\gminia}$.
Then, it is easy to check that
$N_2\geq N_1$.
Thus, we obtain that the second claim implies the first.

Suppose that the second condition is satisfied.
Because
$\vecnbigi(\Sigma)=\vecnbigi(V,\nabla)$ on
$(X\setminus Z',H\setminus Z')$,
we obtain that the closure
$\Sigmabar$ of
$\Sigma$ in
$T^{\ast}(X\setminus Z')
(\log H^{\circ})(\sum m_iH_i^{\circ})$
is proper over $X\setminus Z'$.
By Lemma \ref{lem;25.10.15.30},
the third condition is satisfied.
Conversely, if the third condition is satisfied,
$\nabla(V_0)\subset
V_0\otimes\Omega^1_{X^{\circ}}(\log H^{\circ})(\sum m_iH_i^{\circ})$
holds on $X\setminus Z'$.
By the Hartogs property,
the second condition is satisfied.
\hfill\qed

\subsection{Universal integrable meromorphic relative connections}

\subsubsection{Universal meromorphic integrable relative connection of a given bundle}

Let $\nbigs$ be a $\hyperk$-variety.
Let $\nbigx\to\nbigs$ be a smooth projective morphism of
$\hyperk$-varieties
with relatively very ample line bundle $\nbigo_{\nbigx}(1)$.
Let $\nbigh\subset\nbigx$ be a hypersurface of $\nbigx$
which is flat over $\nbigs$.
For any $\hyperk$-variety $\nbigs_1$ over $\nbigs$,
we set $\nbigx_{\nbigs_1}=\nbigx\times_{\nbigs}\nbigs_1$
and $\nbigh_{\nbigs_1}=\nbigh\times_{\nbigs}\nbigs_1$.
For any coherent $\nbigo_{\nbigx}$-module $\nbigf$ flat over $\nbigs$,
let $\nbigf_{\nbigs_1}$ denotes the pull back of $\nbigf$
to $\nbigx_{\nbigs_1}$.

\begin{prop}
\label{prop;24.8.18.11}
Let $\nbigv$ be a locally free
$\nbigo_{\nbigx}$-module of finite rank.
For any positive integer $N$,
there exists a complex variety over $\nbigs_1$ over $\nbigs$
such that the following holds.
\begin{itemize}
\item
The locally free $\nbigo_{\nbigx_{\nbigs_1}}$-module
$\nbigv_{\nbigs_1}$ is equipped with
     an integrable meromorphic relative connection
\[
     \nabla:\nbigv_{\nbigs_1}
     \to
     \nbigv_{\nbigs_1}
     \otimes
     \Omega^1_{\nbigx/\nbigs}(N\nbigh)_{\nbigs_1}.
\]
    
 \item Let $s\in\nbigs$ be any closed point.
       Let $\nabla:\nbigv_{s}
       \to\nbigv_{s}\otimes\Omega^1_{\nbigx/\nbigs}(N\nbigh)_s$
       be a meromorphic integrable connection.
       Then, there exists a closed point
       $s_1\in \nbigs_1$ over $s$
       such that
       $(\nbigv_{s},\nabla)
       \simeq
       (\nbigv_{\nbigs_1},\nabla)_{s_1}$,
       where
       $(\nbigv_{\nbigs_1},\nabla)_{s_1}$
       denotes the restriction of
       $(\nbigv_{\nbigs_1},\nabla)$
       to $(\nbigx_{\nbigs_1})_{s_1}$.       
\end{itemize}
\end{prop}
\pf
By considering a covering by affine open subsets of $\nbigs$,
we may assume that
(i) $\nbigs$ is affine,
(ii) there exist mutually distinct smooth hypersurfaces $\nbigh'_i$ 
$(i=1,\ldots,m)$ obtained as sections of $\nbigo_{\nbigx}(\ell)$
for some $\ell\geq 1$,
(iii) $\bigcap_{i=1}^m\nbigh_i'=\emptyset$.
We set $\nbigh'=\bigcup \nbigh_i'$.
Let $N'$ be a sufficiently large positive integer.
Let $p_{\nbigs}:\nbigx\to\nbigs$ denote the projection.
We may assume that
\[
R^jp_{\nbigs\ast}\Bigl(
 \End(\nbigv)
 \otimes
 \Omega^1_{\nbigx/\nbigs}(N\nbigh+N'\nbigh_i')
 \Bigr)
 =0
 \quad(i=1,\ldots,m),
 \quad
 R^jp_{\nbigs\ast}\Bigl(
 \End(\nbigv)
 \otimes
 \Omega^1_{\nbigx/\nbigs}(N\nbigh+N'\nbigh')
 \Bigr)
 =0
\]
for any $j>0$,
and 
$p_{\nbigs\ast}\Bigl(
\End(\nbigv)
 \otimes
 \Omega^1_{\nbigx/\nbigs}(N\nbigh+N'\nbigh')
 \Bigr)$
and
$p_{\nbigs\ast}\Bigl(
\End(\nbigv)
\otimes
\Omega^1_{\nbigx/\nbigs}(N\nbigh+N'\nbigh'_i)
\Bigr)$
$(i=1,\ldots,m)$
are locally free.
We may also assume that
$R^jp_{\nbigs\ast}\Bigl(
 \End(\nbigv)
 \otimes
 \Omega^2_{\nbigx/\nbigs}(2N\nbigh+2N'\nbigh')
 \Bigr)=0$ for any $j>0$,
and 
\[
 p_{\nbigs\ast}\Bigl(
 \End(\nbigv)
 \otimes
 \Omega^2_{\nbigx/\nbigs}(2N\nbigh+2N'\nbigh')
 \Bigr)
\]
is locally free.

\begin{lem}
\label{lem;25.10.7.20}
If $N'$ is sufficiently large,
there exist meromorphic relative connections
\[
\nabla^{(i)}_0:\nbigv\to
 \nbigv\otimes
 \Omega_{\nbigx/\nbigs}^1(N'\nbigh_i'),
 \quad
 \nabla^{(0)}_0:
 \nbigv\to
 \nbigv\otimes
 \Omega_{\nbigx/\nbigs}^1(N'\nbigh'),
\]
which are not necessarily integrable.
\end{lem}
\pf
Because $\nbigx\setminus\nbigh_i'$ is affine,
there exists a sufficiently large $M$
such that $\nbigv(\ast \nbigh_i')$ is a direct summand of 
$\nbigo_{\nbigx}(\ast \nbigh_i')^{\oplus M}$.
Let
$\pi_i:
\nbigo_{\nbigx}(\ast \nbigh_i')^{\oplus M}
\to \nbigv(\ast\nbigh_i')$
denote the projection
with respect to a decomposition into a direct sum.
Let
$j_i:\nbigv(\ast \nbigh_i')
\to
\nbigo_{\nbigx}(\ast\nbigh_i')^{\oplus M}$
denote the inclusion.
Let $\nabla^{\prime(i)}_0$ be the trivial connection
of the trivial bundle 
$\nbigo_{\nbigx}(\ast\nbigh_i')^{\oplus M}$.
We obtain a connection
$\nabla^{(i)}_0=\pi_i\circ\nabla^{\prime(i)}_0\circ j_i$
of $\nbigv(\ast\nbigh_i')$.
Similarly,
we can construct $\nabla^{(0)}_0$.
\hfill\qed

\vspace{.1in}

Let $\nabla^{(i)}_0$ and $\nabla^{(0)}_0$
be meromorphic relative connections
as in Lemma \ref{lem;25.10.7.20}.
We regard 
$p_{\nbigs\ast}\Bigl(
\End(\nbigv)\otimes
\Omega_{\nbigx/\nbigs}^1(N\nbigh+N'\nbigh')
\Bigr)$
as an affine variety $\nbigs^{(0)}$ over $\nbigs$.
There exists the universal section $\Psi^{(0)}$ of
\[
 \End(\nbigv_{\nbigs^{(0)}})\otimes
 \Omega^1_{\nbigx/\nbigs}(N\nbigh+N'\nbigh')_{\nbigs^{(0)}}
\]
on $\nbigx_{\nbigs^{(0)}}$.
By using $\nabla^{(0)}_0$ and $\Psi^{(0)}$,
we obtain the meromorphic relative connection
\[
       \nablatilde^{(0)}:\nbigv_{\nbigs^{(0)}}
       \to
       \nbigv_{\nbigs^{(0)}}
       \otimes
       \Omega^1_{\nbigx/\nbigs}(N\nbigh+N'\nbigh')_{\nbigs^{(0)}}.
\]
The following holds.
\begin{itemize}
 \item Let $s\in\nbigs$ be any closed point.
       Let $\nabla:\nbigv_s\to
       \nbigv_s\otimes\Omega^1_{\nbigx/\nbigs}(N\nbigh+N'\nbigh')_s$
       be a meromorphic connection.
       Then, there exists
       $s'\in
       \nbigs^{(0)}$ over $s$
       such that
       $
       (\nbigv_{\nbigs^{(0)\prime}},\nablatilde^{(0)})_{s'}
       \simeq
       (\nbigv_s,\nabla)$.       
\end{itemize}
Similarly,
we regard
$p_{\nbigs\ast}\Bigl(
\End(\nbigv)\otimes
\Omega_{\nbigx/\nbigs}^1(N\nbigh+N'\nbigh'_i)
\Bigr)$
as an affine variety $\nbigs^{(i)}$ over $\nbigs$.
There exists the universal section $\Psi^{(i)}$ of
\[
 \End(\nbigv_{\nbigs^{(i)}})\otimes
 \Omega^1_{\nbigx/\nbigs}(N\nbigh+N'\nbigh'_i)_{\nbigs^{(i)}}
\]
on $\nbigx_{\nbigs^{(i)}}$.
By using $\nabla^{(i)}_0$ and $\Psi^{(i)}$,
we obtain the meromorphic relative connection
\[
       \nablatilde^{(i)}:\nbigv_{\nbigs^{(i)}}
       \to
       \nbigv_{\nbigs^{(i)}}
       \otimes
       \Omega^1_{\nbigx/\nbigs}(N\nbigh+N'\nbigh'_i)_{\nbigs^{(i)}}.
\]
The following holds.
\begin{itemize}
 \item Let $s\in\nbigs$ be any closed point.
       Let $\nabla:\nbigv_s\to
       \nbigv_s\otimes\Omega^1_{\nbigx/\nbigs}(N\nbigh+N'\nbigh_i')_s$
       be a meromorphic connection.
       Then,
       there exists
       $s'\in\nbigs^{(i)}$
       over $s$
       such that
       $(\nbigv_{\nbigs^{(i)\prime}},\nablatilde^{(i)})_{s'}
       \simeq
       (\nbigv_s,\nabla)$.       
      
\end{itemize}

We obtain the section
$\nablatilde^{(i)}-\nabla_0^{(0)}
=(\nabla^{(i)}_0-\nabla^{(0)}_0)
+\Psi^{(i)}$ of
$\End(\nbigv_{\nbigs^{(i)}})
 \otimes
 \Omega^1_{\nbigx/\nbigs}(N\nbigh+N\nbigh')_{\nbigs^{(i)}}$.
We obtain the induced closed embeddings
$\nbigs^{(i)}\to\nbigs^{(0)}$
such that the pull back of
$(\nbigv_{\nbigs^{(0)}},\nablatilde^{(0)})$
are isomorphic to
$(\nbigv_{\nbigs^{(i)}},\nablatilde^{(i)})$.
Let $p_{\nbigs^{(0)}}:\nbigx_{\nbigs^{(0)}}\to\nbigs^{(0)}$
denote the projection.
As the curvature,
we obtain
the section $(\nablatilde^{(0)})^2$
of the locally free sheaf
\[
 p_{\nbigs^{(0)}\ast}\Bigl(
 \End(\nbigv_{\nbigs^{(0)}})
 \otimes
 \Omega^2_{\nbigx/\nbigs}(2N\nbigh+2N'\nbigh')_{\nbigs^{(0)}}
 \Bigr).
\]
Let $\nbigs_1\subset\nbigs^{(0)}$ be the subscheme 
as the intersection of
$\nbigs^{(i)}$ $(i=1,\ldots,m)$
and the zero of $(\nablatilde^{(0)})^2$.
Let $(\nbigv_{\nbigs_1},\nabla)$ be
the pull back of
$(\nbigv_{\nbigs^{(0)}},\nablatilde^{(0)})$.
Then, it has the desired property.
\hfill\qed

\subsubsection{Appendix: Universal extension}

Let $\nbigv_i$ be a locally free $\nbigo_{\nbigx}$-modules of finite rank
with an integrable meromorphic relative connection
\[
 \nabla_i:\nbigv_i\to
 \nbigv_i\otimes
 \Omega_{\nbigx/\nbigs}^1(N\nbigh).
\]

\begin{prop}
\label{prop;24.8.19.1}
There exist a $\hyperk$-variety $\nbigs_1$ over $\nbigs$,
a locally free
$\nbigo_{\nbigx_{\nbigs_1}}$-module $\nbigv_{3}$
equipped with a relative integrable meromorphic connection
\[
 \nabla:\nbigv_3\to
 \nbigv_3
 \otimes
 \Omega^1_{\nbigx/\nbigs}(N\nbigh)_{\nbigs_1}
\]
and the exact sequence of integrable meromorphic connections
\begin{equation}
\label{eq;24.8.18.31}
       0\lrarr
       (\nbigv_{1},\nabla)_{\nbigs_1}
       \lrarr
       (\nbigv_3,\nabla)
       \lrarr
       (\nbigv_{2},\nabla)_{\nbigs_1}
       \lrarr 0,
\end{equation} 
such that the following holds. 
\begin{itemize}
 \item Let $s\in\nbigs$.
      Let $V$ be a locally free $\nbigo_{\nbigx_s}$-module
       with an integrable meromorphic connection
       $\nabla:V\to V\otimes\Omega^1_{\nbigx/\nbigs}(N\nbigh)_s$
       with an exact sequence
\begin{equation}
\label{eq;24.8.18.30}
       0\lrarr (V_1,\nabla)\lrarr (V,\nabla)
       \lrarr (V_2,\nabla)\lrarr 0
\end{equation}
       such that
       $(V_i,\nabla)\simeq
       (\nbigv_i,\nabla)_s$.
       Then,
       there exists
       $s'\in\nbigs_1$ over $s$
       the pull back of {\rm(\ref{eq;24.8.18.31})}
       to $s'$
       is isomorphic to {\rm(\ref{eq;24.8.18.30})}.
\end{itemize}
\end{prop}
\pf
By setting
$\nbigc^j=\nhom(\nbigv_2,\nbigv_1)
 \otimes\Omega^{j}_{\nbigx/\nbigs}(jN\nbigh)$,
we obtain the complex
$\nbigc^{\bullet}$ on $\nbigx$.
By taking an affine open covering of $\nbigs$,
we may assume that $\nbigs$ is smooth and affine.
We may also assume that there exist
mutually distinct hypersurfaces $\nbigh_i$ $(i=1,\ldots,m)$ of $\nbigx$
such that
(i) $\nbigh_i$ are smooth over $\nbigs$
and obtained as the zero of
sections of $\nbigo_{\nbigx}(\ell)$ for some $\ell>0$,
(ii) $\bigcap \nbigh_i=\emptyset$.
We set $\nbigh(I)=\bigcup_{i\in I} \nbigh_i$
for any $I\subset \{1,\ldots,m\}$.
For any non-empty $I$,
we may assume that
$R^kp_{\nbigs\ast}\bigl(
 \nbigc^j\otimes \nbigo((j+1)\nbigh(I))
 \bigr)=0$ $(k>0)$,
and that
$p_{\nbigs\ast}\bigl(
\nbigc^j\otimes \nbigo((j+1)\nbigh(I))
\bigr)$
are locally free.

We set
$U=\bigoplus_{i=1}^m \cnum e_i$.
For an ordered subset
$I=\{i_1,\ldots,i_p\}\subset \{1,\ldots,m\}$,
we set
$e_I=e_{i_1}\wedge\cdots\wedge e_{i_p}$.
We set $U_I=\cnum e_I\subset \bigwedge^{|I|}U$.
The multiplication of $e_i$
induces 
$\varphi_i:U_I\to U_{I\sqcup\{i\}}$.
Together with the natural inclusion
$p_{\nbigs\ast}\bigl(
\nbigc^j\otimes \nbigo((j+1)\nbigh(I))
\bigr)
\to
p_{\nbigs\ast}\bigl(
\nbigc^j\otimes \nbigo((j+1)\nbigh(I\cup\{i\}))
\bigr)$,
we obtain
\[
 p_{\nbigs\ast}\bigl(
\nbigc^j\otimes \nbigo((j+1)\nbigh(I))
\bigr)
\otimes U_{I}
\lrarr
p_{\nbigs\ast}\bigl(
\nbigc^j\otimes \nbigo((j+1)\nbigh(I\cup\{i\}))
\bigr)
\otimes U_{I\cup\{i\}}.
\]
We set
\[
 \nbigctilde^{j,p}
=\bigoplus_{|I|=p}
 p_{\nbigs\ast}\bigl(
\nbigc^j\otimes \nbigo((j+1)\nbigh(I))
\bigr)
\otimes U_{I}.
\]
We obtain
$\nbigctilde^{j,p}
 \to
 \nbigctilde^{j,p+1}$.
We obtain the morphisms
$p_{\nbigs\ast}\bigl(
\nbigc^j\otimes \nbigo((j+1)\nbigh(I))
\bigr)
\to 
p_{\nbigs\ast}\bigl(
\nbigc^{j+1}\otimes \nbigo((j+2)\nbigh(I))
\bigr)$
induced by the integrable meromorphic relative connections,
and hence
$\nbigctilde^{j,p}\to \nbigctilde^{j+1,p}$.
Thus, we obtain a double complex
$\nbigctilde^{\bullet,\bullet}$.
Let $\nbigchat^{\bullet}$
denote the total complex.
We may regard
$\nbigchat^i$ are vector bundles
on $\nbigs_0$,
and
we have the morphisms
$\nbigchat^i\to\nbigchat^{i+1}$.

We regard $\nbigchat^1$ as a variety $\nbigs^{(0)}$ over $\nbigs$.
Let $\pi_1:\nbigs^{(0)}\to\nbigs$ denote the projection.
We obtain the section $\Psi$ of
$\pi_1^{\ast}(\nbigchat^2)$ on $\nbigs_0$
induced by $\nbigchat^1\to\nbigchat^2$.
Let $\nbigs_1\subset\nbigs_0$
denote the subscheme obtained as $\Psi^{-1}(0)$.

For $i\neq j$,
we have the morphisms
\[
 \psi_{i,j}:
 \nbigv_{2,\nbigs_1}
 \lrarr
 \nbigv_{1,\nbigs_1}
 \otimes
 \nbigo(\nbigh_i\cup \nbigh_j)_{\nbigs_1}.
\]
For $i$,
we have the morphisms
\[
 \phi_{i}:
 \nbigv_{2,\nbigs_1}
 \lrarr
 \nbigv_{1,\nbigs_1}
 \otimes
 \bigl(
 \Omega^1_{\nbigx/\nbigs}(N\nbigh)
 \otimes
 \nbigo(2\nbigh_i)
 \bigr)_{\nbigs_1}.
\]
They satisfy the cocycle condition:
\[
\psi_{i,j}+\psi_{j,k}=\psi_{i,k},
\quad
\nabla_1\circ \psi_{i,j}-\psi_{i,j}\circ\nabla_2
+\phi_i-\phi_j=0,
\quad
\nabla_1\circ\phi_i+\phi_i\circ\nabla_2=0.
\]
We obtain the isomorphisms
\begin{equation}
\label{eq;25.10.8.1}
 \left(
 \begin{array}{cc}
  \id_{\nbigv_{1}} & \psi_{i,j} \\
  0 & \id_{\nbigv_2}
 \end{array}
 \right)
 :
\Bigl(
 (\nbigv_{1,\nbigs_1}\oplus
 \nbigv_{2,\nbigs_1})_{|(\nbigx\setminus\nbigh_j)_{\nbigs_1}}
 \Bigr)_{|(\nbigx\setminus(\nbigh_i\cup\nbigh_j))_{\nbigs_1}}
 \simeq
 \Bigl(
 (\nbigv_{1,\nbigs_1}\oplus\nbigv_{2,\nbigs_1})
 _{|(\nbigx\setminus\nbigh_i)_{\nbigs_1}}
 \Bigr)_{|(\nbigx\setminus(\nbigh_i\cup\nbigh_j))_{\nbigs_1}}.
\end{equation}
They satisfy the cocycle conditions.
By patching
$(\nbigv_{1,\nbigs_1}\oplus
 \nbigv_{2,\nbigs_1})_{|(\nbigx\setminus\nbigh_j)_{\nbigs_1}}$
via the isomorphisms (\ref{eq;25.10.8.1}),
we obtain a vector bundle $\nbigv_3$
on $\nbigs_1$.
The integrable meromorphic relative connections $\nabla$ of
$(\nbigv_{i})_{|\nbigs_3\times (X\setminus H_i)}$
and $\phi_i$ induces
an integrable meromorphic relative connections of
$(\nbigv_{1,\nbigs_1}\oplus
 \nbigv_{2,\nbigs_1})_{|(\nbigx\setminus\nbigh_j)_{\nbigs_1}}$.
By patching them,
we obtain 
a integrable meromorphic relative connection of $\nbigv_3$.
There exists the exact sequence (\ref{eq;24.8.18.31}).
The desired conditions are satisfied by the construction.
\hfill\qed

\subsection{Estimate of maximal slopes in terms of slopes}
\label{subsection;25.10.15.12}

\subsubsection{Slopes}

Let $\hyperk$ be an algebraically closed field of characteristic $0$.
Let $X$ be an $n$-dimensional smooth connected projective
$\hyperk$-variety
with an ample line bundle $\nbigo_X(1)$.
We refer \cite{Fulton}
for the intersection theory and 
the characteristic classes of coherent $\nbigo_X$-modules
in the Chow group of $X$.
We set $\omega=c_1(\nbigo_X(1))$.
For any torsion-free coherent $\nbigo_X$-module $V$,
we set
$\deg_{\omega}(V)=\int_Xc_1(V)\omega^{n-1}$
and
$\mu_{\omega}(V)
=\deg_{\omega}(V)/\rank(V)$.
There exists the Harder-Narasimhan filtration $\nbigf^{HN}$ of $V$
with respect to the slope $\mu_{\omega}$,
i.e.,
there exists a unique decreasing filtration $\nbigf^{HN}$ of $V$
indexed by $\rnum$ such that
(i) $\nbigf^{HN}_{a}(V)=0$ and $\nbigf^{HN}_{a}(V)=V$
for any sufficiently large $a>0$,
(ii) $\Gr^{\nbigf^{HN}}_{a}(V)$ are $\mu_{\omega}$-semistable sheaves
with $\mu\bigl(\Gr^{{\nbigf}^{HN}}_{a}(V)\bigr)=a$.
We set
\[
 \mu_{\omega,\max}(V)
 =\max\bigl\{
 a\,\big|\,\Gr^{\nbigf^{HN}}_{a}(V)\neq 0
 \bigr\},
 \quad
 \mu_{\omega,\min}(V)
 =\min\bigl\{
 a\,\big|\,\Gr^{\nbigf^{HN}}_{a}(V)\neq 0
 \bigr\}.
\]

\begin{lem}
\label{lem;26.1.3.1}
Let $V_i$ $(i=1,2)$ be torsion-free coherent $\nbigo_X$-modules
with a non-zero morphism $\varphi:V_1\to V_2$.
Suppose that $V_1$ is $\mu_{\omega}$-semistable.
Then, we obtain $\mu_{\omega}(V_1)\leq \mu_{\omega,\max}(V_2)$.
\end{lem}
\pf
Let $\nbigf^{HN}(V_2)$ be a Harder-Narasimhan filtration of $V_2$.
There exists $b(0)$ such that
$\varphi(V_1)\subset \nbigf^{HN}_{b(0)}(V_2)$
and $\varphi(V_1)\not\subset \nbigf^{HN}_{<b(0)}(V_2)$.
We obtain a non-zero morphism
$\varphi_1:V_1\to \Gr^{\nbigf}_{b(0)}(V_2)$.
Because
$V_1$ and $\Gr^{\nbigf}_{b(0)}(V_2)$
are $\mu_{\omega}$-semistable,
we obtain
\[
 \mu_{\omega}(V_1)
 \leq
 \mu_{\omega}(\Image \varphi_1)
 \leq
 \mu_{\omega}(\Gr^{\nbigf}_{b(0)}(V_2))
 \leq
 \mu_{\omega,\max}(V_2).
\]
Thus, we are done.
\hfill\qed

\begin{lem}
\label{lem;26.1.3.11}
Let $V$ be a coherent $\nbigo_X$-module
with an increasing filtration $F_j(V)$ $(j=0,\ldots,m)$
by coherent $\nbigo_X$-submodules
such that
(i) $F_0(V)=0$ and $F_m(V)=V$,
(ii) $\Gr^F_j(V)$ are torsion-free.
Then, we have 
\[
 \min_{1\leq j\leq m}\mu_{\omega}(\Gr^{\nbigf}_j(V))
 \leq
 \mu_{\omega}(V)
 \leq
 \max_{1\leq j\leq m}\mu_{\omega}(\Gr^{\nbigf}_j(V)).
\]
\end{lem}
\pf
If $m=1$, the claim is trivial.
Suppose $m\geq 2$.
Either one of the following holds.
\begin{itemize}
 \item $\mu_{\omega}(\nbigf_{m-1}(V))
       =\mu_{\omega}(V)=\mu_{\omega}(\Gr^{\nbigf}_m(V))$.
 \item $\mu_{\omega}(\nbigf_{m-1}(V))
       \leq
       \mu_{\omega}(V)\leq
       \mu_{\omega}(\Gr^{\nbigf}_m(V))$.
 \item $\mu_{\omega}(\nbigf_{m-1}(V))
       \geq
       \mu_{\omega}(V)\geq
       \mu_{\omega}(\Gr^{\nbigf}_m(V))$.
\end{itemize}
Then, we obtain the claim of the lemma by an induction on $m$.
\hfill\qed

\subsubsection{Estimate of maximal slopes in terms of slopes}

For a sufficiently large integer $M_1$,
there exists a monomorphism
$\Omega_X^1\to \nbigo_X(M_1)^{n}$.
We set $M_2=\deg_{\omega}(\nbigo_X(M_1))$.
Let $H$ be a normal crossing hypersurface.
Let $N$ be any positive integer.
\begin{lem}
Let $V_i$ $(i=1,2)$
be torsion-free coherent $\mu_{\omega}$-semistable $\nbigo_X$-modules
with a non-zero morphism
$V_1\to V_2\otimes\Omega^1_X(NH)$.
Then,
we obtain
$\mu_{\omega}(V_1)\leq \mu_{\omega}(V_2)+M_2+N\deg_{\omega}(H)$ 
\end{lem}
\pf
There exists a non-zero morphism
$V_1\to V_2\otimes\nbigo_X(M_1)$.
Because $V_1$ and $V_2\otimes\nbigo_X(M_1)$
are $\mu_{\omega}$-semistable,
we obtain
$\mu_{\omega}(V_1)\leq
\mu_{\omega}(V_2\otimes\nbigo_X(M_1))
=\mu_{\omega}(V_2)+M_2$.
\hfill\qed

\begin{prop}
\label{prop;24.8.18.2}
Let $V$ be a locally free $\nbigo_X$-module
equipped with a meromorphic connection
$\nabla:V\to V\otimes\Omega^1_X(NH)$
such that
$(V,\nabla)$ is irreducible,
i.e.,
there is no saturated $\nbigo_X$-submodule $V'$ of $V$
such that $\nabla(V')\subset V'\otimes\Omega_X^1(NH)$
with $0<\rank(V')<\rank(V)$.
Then, we obtain
\[
       \mu_{\omega,\max}(V)-
       \mu_{\omega,\min}(V)
       \leq
       \rank(V)\cdot \bigl(
       M_2+N\deg_{\omega}(H)
 \bigr).
\]       
As a result, we obtain
\[
 \mu_{\omega,\max}(V)
 \leq
 \mu_{\omega}(V)
+\rank(V)\cdot \bigl(M_2+N\deg_{\omega}(H)\bigr).
\]
\end{prop}
\pf
Let $(V,\nabla)$ be as in the statement of the proposition.
If $V$ is a $\mu_{\omega}$-semisimple sheaf,
$\mu_{\omega,\max}(V)-\mu_{\omega,\min}(V)=0$.
Suppose that $V$ is not $\mu_{\omega}$-semistable.
We set
$V_{a}=\Gr^{\nbigf^{HN}}_{a}(V)$.
Let $a_{0}>a_1>\cdots>a_m$
be the set
$\{a\,|\,V_{a}\neq 0\}$.
We have $a_0=\mu_{\omega,\max}(V)$
and $a_{m}=\mu_{\omega,\min}(V)$.
We set $i(0)=0$.

\begin{lem}
\label{lem;24.8.18.1}
There exists
a strictly increasing sequence
$i(0)<i(1)<\cdots <i(q)=m$
and
$k(p)\leq i(p-1)$ $(p=1,\ldots,q)$
such that the following holds:
\begin{equation}
\label{eq;25.10.7.1}
 a_{k(p)}\leq a_{i(p)}+M_2+N\deg_{\omega}(H).
 \end{equation}
\end{lem}
\pf
Because $(V,\nabla)$ is irreducible,
the induced morphism of $\nbigo_{X}$-modules
$\nbigf_{a_0}\to (V/\nbigf_{a_0})\otimes \Omega^1_X(NH)$
is non-trivial.
There exist $0<i(1)\leq m$
and a non-zero morphism
$V_{a_0}\to
V_{a_{i(1)}}\otimes
\Omega^1_X(NH)$.
We set
$k(1)=0$.
We obtain
\[
 a_{k(1)}\leq a_{i(1)}+M_2+N\deg_{\omega}H.
\]
Suppose we have already constructed
$0=i(0)<\cdots<i(s)<m$
and $k(p)\leq i(p-1)$ $(p=1,\ldots,s)$
such that (\ref{eq;25.10.7.1}) holds.
Because $(V,\nabla)$ is irreducible,
the induced morphism of $\nbigo_X$-modules
$\nbigf_{a_{i(s)}}\to
(V/\nbigf_{a_{i(s)}})\otimes \Omega^1_X(NH)$
is non-trivial.
There exist 
$k(s+1)\leq i(s)<i(s+1)$
and a non-zero morphism
$V_{a_{k(s+1)}}\to
V_{a_{i(s+1)}}\otimes
\Omega^1_X(NH)$.
It implies that
$a_{k(s+1)}\leq a_{i(s+1)}+M_2+N\deg_{\omega}(H)$.
Hence, we obtain the claim of Lemma \ref{lem;24.8.18.1}
by an induction.
\hfill\qed

\vspace{.1in}
For $p=1,\ldots,q$, we obtain
\[
 a_{i(p-1)}
 \leq
 a_{k(p)}
 \leq
 a_{i(p)}+M_2+N\deg_{\omega}(H).
\]
We obtain
$a_0\leq a_m+q(M_2+N\deg_{\omega}(H))\leq
a_m+\rank(V)(M_2+N\deg_{\omega}(H))$.
Thus, we obtain Proposition \ref{prop;24.8.18.2}
in the case that $V$ is not $\mu_{\omega}$-semistable.
\hfill\qed

\begin{cor}
\label{cor;26.1.3.20}
Let $V$ be a locally free $\nbigo_X$-module
with a meromorphic connection
$\nabla:V\to V\otimes\Omega^1_X(NH)$
and an increasing filtration $F_j(V)$ $(j=0,\ldots,m)$
by $\nbigo_X$-submodules
 such that
(i) $F_0(V)=0$ and $F_m(V)=V$,
(ii) $\nabla(F_j(V))\subset F_j(V)\otimes\Omega^1_X(NH)$,
 (iii) $\Gr^{F}_j(V):=F_j(V)/F_{j-1}(V)$
$(j=1,\ldots,m)$  are locally free $\nbigo_X$-modules.
Then, we obtain 
\begin{equation}
\label{eq;26.1.3.2}
 \mu_{\omega,\max}(V)
 \leq
 \max_{1\leq j\leq m}\mu_{\omega}(\Gr^{F}_j(V))
+\rank(V)(M_2+N\deg_{\omega}(H)).
\end{equation}
If moreover there exist constants $A_i$ $(i=1,2)$ such that 
$A_1\leq \mu_{\omega}(\Gr^{\nbigf}_j(V))\leq A_2$
$(j=1,\ldots,m)$,
we obtain
\begin{equation}
\label{eq;26.1.3.10}
 \mu_{\omega,\max}(V)
 \leq
 \mu_{\omega}(V)
 +\rank(V)(M_2+N\deg_{\omega}(H))
 +A_2-A_1.
\end{equation}
\end{cor}
\pf
By Proposition \ref{prop;24.8.18.2},
we obtain
\begin{equation}
\label{eq;26.1.3.3}
 \mu_{\omega,\max}(\Gr^{F}_j(V))
 \leq 
 \mu_{\omega}(\Gr^{F}_j(V))
+\rank\Gr^{F}_j(V)(M_2+N\deg_{\omega}(H)).
\end{equation}
Let $\nbigf^{HN}\Gr^F_j(V)$ denote the Harder-Narasimhan filtration
of $\Gr^F_j(V)$.
Let $\nbigf^{HN}(V)$ denote the Harder-Narasimhan filtration of $V$.
We set $a(0)=\mu_{\omega,\max}(V)$.
There exists $j(0)$ such that
$\nbigf^{HN}_{a(0)}(V)\subset F_{j(0)}(V)$
and
$\nbigf^{HN}_{a(0)}(V)\not\subset F_{j(0)-1}(V)$.
We obtain the non-zero morphism
$\nbigf^{HN}_{a(0)}(V)\to \Gr^{F}_{j(0)}(V)$.
By Lemma \ref{lem;26.1.3.1},
we obtain
$a(0)\leq \mu_{\omega,\max}(\Gr^F_{j(0)}(V))$.
Then, we obtain (\ref{eq;26.1.3.2}) from (\ref{eq;26.1.3.3}).
We obtain (\ref{eq;26.1.3.10})
from (\ref{eq;26.1.3.2})
and Lemma \ref{lem;26.1.3.11}.
\hfill\qed

\subsection{Hilbert polynomials}

\subsubsection{Hilbert polynomials of torsion-free sheaves}
\label{subsection;25.12.10.1}

Let $\hyperk$ be an algebraically closed field of characteristic $0$.
Let $X$ be an $n$-dimensional projective smooth $\hyperk$-variety
with an ample line bundle $\nbigo_X(1)$.
For any torsion-free $\nbigo_X$-module $\nbigf$,
let $P_{\nbigf}\in\rnum[t]$ denote the Hilbert polynomial of $\nbigf$,
i.e.,
$P_{\nbigf}(m)$
equals
the Euler number of $\nbigf\otimes\nbigo_X(m)$
for any $m\in\seisuu$,
i.e.,
\[
 P_{\nbigf}(m)
=\int_X \td(X)\ch(\nbigf)\ch(\nbigo_X(m)).
\]
There exist integers $a_i(\nbigf)$ $(i=0,\ldots,n)$
such that
\[
 P_{\nbigf}(m)
 =\sum_{i=0}^n
 a_i(\nbigf)
 \left(
 \begin{array}{c}
  m+n-1 \\ n-i
 \end{array}
 \right).
\]
We have $a_0(\nbigf)=\rank(\nbigf)\int_Xc_1(\nbigo_X(1))^{n}$.

\subsubsection{Maruyama's boundedness Theorem}

Let $\nbigx_S\to S$ be a smooth projective morphism of
smooth complex algebraic varieties
with a relatively ample line bundle $\nbigo(1)$.
Let $r$ be a positive integer,
let $a_1,a_2\in\seisuu$.
Let $M>0$.
Let $\nbigc(\nbigx_s,r,a_1,a_2,M)$
denote the category of
torsion-free $\nbigo_{\nbigx_s}$-sheaves $E$
satisfying the following conditions:
\[
 \mu_{\max}(E)<\mu(E)+M,
 \quad
 \rank(E)=r,
 \quad
 |a_1(E)|\leq a_1,
 \quad
 a_2(E)\geq a_2.
\]
\begin{prop}[Maruyama
\mbox{\cite{Maruyama}}]
\label{prop;25.12.14.110}
The family $(\nbigc(\nbigx_s,r,a_1,a_2,M)\,|\,s\in S(\cnum))$
is bounded in the following sense.
\begin{itemize}
 \item There exists a morphism of smooth complex varieties
       $\nbigs\to S$.
 \item There exists a coherent
       $\nbigo_{\nbigx_{\nbigs}}$-module $\nbige$
       flat over $\nbigs$.
 \item For any $b\in S(\cnum)$
       and $E\in\nbigc(\nbigx_s,r,a_1,a_2,M)$,
       there exists
       $\btilde\in \nbigs(\cnum)$ over $b$
       such that
       $\nbige_{\btilde}\simeq E$.
       \hfill\qed
\end{itemize}
 \end{prop}

\begin{rem}
This was generalized by Langer
to the positive characteristic case.
See {\rm\cite{Langer-2022}}.
\hfill\qed
\end{rem}

\subsection{Estimates for some lattices}

Let $\hyperk$ be an algebraically closed field of characteristic $0$.
Let $X$ be an $n$-dimensional projective smooth $\hyperk$-variety
with an ample line bundle $\nbigo_X(1)$.
We set $\omega=c_1(\nbigo_X(1))$.
Let $H=\bigcup_{j\in\Lambda}H_j$
be a simple normal crossing hypersurface of $X$.
By the intersection theory \cite{Fulton},
we obtain the following numbers
associated with $(X,\nbigo_X(1))$:
\begin{equation}
\label{eq;25.10.15.3}
\int_X\omega^n,\quad
\int_Xc_1(X)\omega^{n-1},\quad
\int_{X}c_1(X)^2\omega^{n-2},\quad
\int_{X}c_2(X)\omega^{n-2},\quad
\end{equation}
We also obtain the following additional numbers
associated with $(X,H,\nbigo_X(1))$:
\begin{equation}
\label{eq;25.10.15.4}
\int_Xc_1(X)[H_i]\omega^{n-2},\quad
\int_X[H_i]\omega^{n-1},\quad
\int_X[H_i]\cdot [H_j]\omega^{n-2}.
\end{equation}

\begin{prop}
\label{prop;25.10.15.11}
There exist
a constant $A>0$ depending only on 
$\int_X[H_i]\omega^{n-1}$ $(i\in\Lambda)$,
and a constant $C(r)>0$
depending on a positive number $r$
and the numbers {\rm(\ref{eq;25.10.15.3}, \ref{eq;25.10.15.4})}
such that the following holds.
\begin{itemize}
 \item Let $(V,\nabla)$ be any good meromorphic flat bundle on $(X,H)$
       with $\rank(V)\leq r$.
       Let $\nbigf$ be a filtration of $(V,\nabla)$.
       Then, there exists a good lattice $V_0$ of $V$
       such that
       (i) $\Gr^{\nbigf}_j(V_0)$ are locally free,
       (ii) $0\leq \mu_{\omega}(\Gr^{\nbigf}_jV_0)\leq A$,
       (iii) 
       $|a_1(\Gr^{\nbigf}_j(V_0))|
       +|a_2(\Gr^{\nbigf}_j(V_0))|\leq C(r)$,
       (iv) $\Gr^{\nbigf}_j(V_0)$ generates
       $\Gr^{\nbigf}_j(V)$ as a $\nbigd_X$-module,
       i.e.,
       $\nbigd_X\cdot\Gr^{\nbigf}_j(V_0)=\Gr^{\nbigf}_j(V)$.
\end{itemize}
\end{prop}
\pf
Let us consider the case $\hyperk=\cnum$.
Let $V^{\DM}\subset V$ denote the good Deligne-Malgrange lattice
of $(V,\nabla)$.
By the Riemann-Roch formula,
we have
\begin{multline}
P_{V^{\DM}}(m)
=m^n\rank(V)
 \int_X\omega^n
 +m^{n-1}\int_X
 \bigl(
 c_1(V^{\DM})+\td_1(X)\rank(V^{\DM})
 \bigr)\omega^{n-1}
\\
+m^{n-2}\int_{X}\Bigl(
 \ch_2(V^{\DM})+c_1(V^{\DM})\td_1(X)
 +\rank(V^{\DM})\td_2(X)
 \Bigr)\omega^{n-2}
 +O(m^{n-3}).
\end{multline}
We may assume that $\nbigo_X(1)$ is very ample.
There exists a smooth surface $Z$
obtained as the complete intersection of
$(n-2)$-sections of $\nbigo_X(1)$
such that
$H_Z=Z\cap H$ is normal crossing in $Z$.
The restriction
$(V^{\DM},\nabla)_{|Z}$
is the good Deligne-Malgrange lattice
of $(V,\nabla)_{|Z}$.
By Proposition \ref{prop;25.10.9.40},
there exists $C_1>0$
depending only on 
$\int_X[H_i]\omega^{n-1}$
such that 
\[
 \left|
\int_Xc_1(V^{\DM})\omega^{n-1}
\right|
\leq \rank(V)C_1.
\]
There also exists $C_2>0$
depending only on
$\rank(V)$,
$\int_X[H_i]\omega^{n-1}$,
$\int_X[H_i]c_1(X)\omega^{n-2}$
and 
$\int_X[H_i]\cdot [H_j]\omega^{n-2}$
such that
\[
\left|
\int_Xc_1(V^{\DM})c_1(X)\omega^{n-2}
\right|
+
\left|
\int_Xc_1(V^{\DM})^2\omega^{n-2}
\right|
+
\left|
\int_Xc_2(V^{\DM})\omega^{n-2}
\right|
\leq C_2.
\]
Then, we obtain the desired estimate in the case $\hyperk=\cnum$.

Let us study the general case.
It is enough to consider the case
$\hyperk\subset\cnum$.
We consider the splitting $\tau:\hyperk/\seisuu\to\hyperk$
induced by the bijection
$\{a\in\hyperk\,|\,0\leq \Re(a)<1\}\simeq \hyperk/\seisuu$.
By \cite[Theorem 5.3.4]{kedlaya2},
there exists the good Deligne-Malgrange lattice
$V_0\subset V$
associated with the choice $\tau$.
It induces the good Deligne-Malgrange lattice of
the complexification $(V_{\cnum},\nabla)$
on $X\times\cnum$.
We also note that
$\Gr^{\nbigf}_j(V_0)$
are good Deligne-Malgrange lattices
of $\Gr^{\nbigf}_j(V,\nabla)$.
Hence,
$V_0$ satisfies the desired conditions.
\hfill\qed

\section{Boundedness as meromorphic objects}

\subsection{Families of formal meromorphic flat bundles
and formal $\nbigd$-modules}

\subsubsection{Families of good formal meromorphic flat bundles}

Let $R$ be a regular Henselian local ring over $\hyperk$.
Let $\gminim$ be the maximal ideal.
We have $\hyperk\simeq R/\gminim$.

We set $R_1=R[\![x_1,\ldots,x_n]\!]$
and $\Rtilde_1=R_1(\ast x_1\ldots x_{\ell})$.
Let $\nbige$ be a finitely generated reflexive $\Rtilde_1$-module
with an integrable connection
$\nabla:\nbige\to\nbige\otimes\Omega^1_{R_1/R}$.
Such $(\nbige,\nabla)$
is called a formal meromorphic integrable connection
over $\Rtilde_1$ relative to $R$.
A lattice $E$ of $\nbige$ means
a finitely generated $R_1$-submodule
such that $E\otimes\Rtilde_1=\nbige$.

\begin{df}
A lattice $E$ of $\nbige$ is called logarithmic with respect to $\nabla$
if 
\[
 \nabla E
 \subset
 \bigoplus_{i=1}^{\ell}
 E\otimes
 \bigl(dx_i/x_i\bigr)
 \oplus
  \bigoplus_{i=\ell+1}^{n}
 E\otimes dx_i.
\]
We say that $(\nbige,\nabla)$ is regular
if there exists a logarithmic lattice. 
\hfill\qed
\end{df}

For any $f\in \Rtilde_1$,
we have $df\in \Omega^1_{\Rtilde_1/R}$.

\begin{df}
\label{df;26.1.24.1}
A lattice $E$ of $\nbige$ is called unramifiedly good
if there exist a good set $\nbigi\subset \Rtilde_1/R_1$
and a decomposition 
$(E,\nabla)=\bigoplus_{\gminia\in\nbigi}(E_{\gminia},\nabla_{\gminia})$
such that 
$(E_{\gminia},\nabla_{\gminia}-d\gminiatilde)$
are logarithmic.
Here, $\gminiatilde\in\Rtilde_1$ denote lifts of $\gminia$.
We say that $(\nbige,\nabla)$ is unramifiedly good
 if there exists an unramifiedly good lattice.
\hfill\qed
\end{df}

For any positive integer $e$,
we set $\Rtilde_1^{(e)}=\Rtilde_1[x_1^{1/e},\ldots,x_{\ell}^{1/e}]$.
We say
$\nbige^{(e)}=\nbige\otimes\Rtilde_1^{(e)}$.
It is naturally equipped with an integrable connection
$\nabla^{(e)}$ relative to $R$.
Let $\Gal_e$ denote the Galois group of
the extension $\Rtilde_1^{(e)}/\Rtilde_1$.
We have the natural $\Gal_e$-action on
$\nbige^{(e)}$.

\begin{df}
We say that
$(\nbige,\nabla)$ is good 
if $(\nbige^{(e)},\nabla^{(e)})$ is unramifiedly good
for some $e\in\seisuu_{>0}$.
\hfill\qed
\end{df}

\begin{df}
When $\nbige$ is good,
we say that a lattice $E$ of $\nbige$ is good if 
there exists a $\Gal_e$-invariant unramifiedly good lattice
$E'$ of $\nbige^{(e)}$
such that $E$ is the descent of $E'$. 
\hfill\qed
\end{df}

\subsubsection{The $V\nbigd_{R_1/R}$-submodules and
$\nbigd_{R_1/R}$-submodules associated with lattices}

We set $\nbigd_{R_1/R}=R_1\langle\del_1,\ldots,\del_{n} \rangle$
and
$V\nbigd_{R_1/R}=R_1\langle x_1\del_{1},
\ldots,x_{\ell}\del_{\ell},\del_{\ell+1},
\ldots,\del_n\rangle$.
Any meromorphic integrable connection $(\nbige,\nabla)$
over $R_1$ relative to $R$
is naturally regarded as a $V\nbigd_{R_1/R}$-module
or a $\nbigd_{R_1/R}$-module.

Let $E$ be a lattice of $(\nbige,\nabla)$.
We obtain the $V_{\nbigd_{R_1/R}}$-submodule
$V\nbigd_{R_1/R}\cdot E\subset \nbige$.
We set
\[
 \nbigd_{R_1/R}(E)
 :=\nbigd_{R_1/R}\otimes_{V\nbigd_{R_1/R}}
 \bigl(
 V\nbigd_{R_1/R}\cdot E
 \bigr).
\]
For $E\subset E'$,
we obtain $\nbigd_{R_1/R}$-homomorphism
$\nbigd_{R_1/R}(E)\to\nbigd_{R_1/R}(E')$.

\begin{prop}
\label{prop;26.1.24.5}
Suppose that $E$ is good.
\begin{itemize}
 \item  $V\nbigd_{R_1/R}\cdot E$
	and $\nbige/(V\nbigd_{R_1/R}\cdot E)$ are flat over $R$.
	There exists a natural isomorphism
\begin{equation}
\label{eq;26.1.26.1}
	(V\nbigd_{R_1/R}\cdot E)\otimes_R\hyperk
	\simeq
	V\nbigd_{\hyperk[\![x_1,\ldots,x_n]\!]/\hyperk}\cdot
	 (E\otimes_R\hyperk).
\end{equation}
 \item	 $\nbigd_{R_1/R}(E)$ is flat over $R$.
      There exists the following natural isomorphism
\[
 \nbigd_{R_1/R}(E)\otimes_R\hyperk
 \simeq
 \nbigd_{\hyperk[\![x_1,\ldots,x_{n}]\!]/\hyperk}(E\otimes\hyperk).
\]
\end{itemize}
\end{prop}
\pf
The second and third claims follow from the first.
Let us study the first.
Let us consider the case where $E$ is unramifiedly good.
There exist a good set $\nbigi$ and a decomposition
$E=\bigoplus E_{\gminia}$
as in Definition \ref{df;26.1.24.1}.
For each $\gminia$,
let $I(\gminia)$
denote the set of $1\leq i\leq \ell$
such that $\ord(\gminia)_i\neq 0$.
It is easy to see that
$V\nbigd_{R_1/R}\cdot E_{\gminia}
=E_{\gminia}(\ast x^{I(\gminia)})$.
We may assume that 
$I(\gminia)=\{1,\ldots,k\}$ for some $k$.
Let $\nbigi(k)$ denote the set of $\gminia\in\nbigi$
such that 
$I(\gminia)=\{1,\ldots,k\}$.
We set
$E_k=\bigoplus_{\gminia\in\nbigi(k)}
 E_{\gminia}$.
We obtain the decomposition
\begin{equation}
\label{eq;26.1.24.2}
 E=\bigoplus_{k=1}^{\ell}E_k.
\end{equation}
Then, we have
\begin{equation}
 V\nbigd_{R_1/R}E
=\bigoplus_{k=1}^{\ell}
 E_k(\ast x_1\cdots x_k).
\end{equation}
From this description,
we obtain that
$V\nbigd_{R_1/R}E$ and
$\nbige/V\nbigd_{R_1/R}E$ are flat over $R$.
We also obtain the isomorphism (\ref{eq;26.1.26.1}).

Suppose that $E$ is good.
There exists a $\Gal_e$-invariant unramifiedly good lattice
$E'\subset\nbige^{(e)}$
such that $E$ is the Galois descent of $E'$.
We have the decomposition
\begin{equation}
\label{eq;26.1.24.3}
 E'=\bigoplus_{k=1}^{\ell} E'_k
\end{equation}
as in (\ref{eq;26.1.24.2}).
The decomposition (\ref{eq;26.1.24.3})
is $\Gal_e$-equivariant,
Let $E_k$ denote the Galois descent of $E_k'$.
We obtain the decomposition
\begin{equation}
 E=\bigoplus_{k=1}^{\ell} E_k.
\end{equation}
It is easy to see that
$V\nbigd_{R_1/R}\cdot E
=\bigoplus_{k=1}^{\ell}E_k(\ast x_1\cdots x_k)$.
As a result,
we obtain Proposition \ref{prop;26.1.24.5}.
\hfill\qed

 \subsubsection{Complement}

Let $\nbige$ be a free $\Rtilde_1$-module.
Let $E$ be an $R_1$-lattice of $\nbige$.
Let $\nbigi\subset \Rtilde_1/R_1$ be a good set of irregular values.
Let $\hyperk_1$ be an algebraic closure of the fractional field of $R$.
We set $\gbigr_1=\hyperk_1[\![x_1,\ldots,x_n]\!]$
and $\gbigrtilde_1=\gbigr_1(\ast x_1\cdots x_{\ell})$.
By the extensions $R_1\to \gbigr_1$
and $\Rtilde_1\to \gbigrtilde_1$,
we obtain
$\gbige=\gbigr_1\otimes E$
and $\gbigetilde=\gbigrtilde_1\otimes\nbige$.
It is standard to obtain the following lemma.
\begin{lem}
\label{lem;26.1.21.40}
Suppose that there exists the decomposition
$\gbige=\bigoplus_{\gminia\in\nbigi}
 \gbige_{\gminia}$
such that
\[
(\nabla-d\gminiatilde\id)\gbige_{\gminia}
 \subset
 \bigoplus_{i=1}^{\ell}
 \gbige_{\gminia}\otimes (dx_i/x_i)
 \oplus
  \bigoplus_{i=\ell+1}^{n}
 \gbige_{\gminia}\otimes dx_i.
\] 
Then, there exists a decomposition
$E=\bigoplus_{\gminia\in\nbigi}E_{\gminia}$
such that  
\[
(\nabla-d\gminiatilde\id)E_{\gminia}
 \subset
 \bigoplus_{i=1}^{\ell}
 E_{\gminia}\otimes (dx_i/x_i)
 \oplus
  \bigoplus_{i=\ell+1}^{n}
 E_{\gminia}\otimes dx_i.
\] 
\hfill\qed
\end{lem}

\subsection{Good families of meromorphic flat connections}

\subsubsection{Families of meromorphic flat connections}
\label{subsection;25.1.21.30}

Let $\hyperk$ be an algebraically closed field of characteristic $0$.
Let $U=\Spec(S)$ be a smooth affine $\hyperk$-scheme.
Let $X$ be a normal $\hyperk$-variety with a morphism
$X\to U$ such that
(i) $X$ is projective and flat over $U$,
(ii) each geometric fiber of $X\to U$ is irreducible and normal.
Let $H\subset X$ be a hypersurface flat over $S$.
We assume that $X^{\circ}=X\setminus H$ is smooth over $U$.
For any geometric point $u:\Spec(\hyperk_1)\to U$,
we set $(X_u,H_u)=(X,H)\times_U \Spec(\hyperk_1)$.

Let $j:X^{\circ}\to X$ denote the inclusion.
We have the relative exterior derivative
$d:j_{\ast}\nbigo_{X^{\circ}}
\to j_{\ast}(\Omega^1_{X^{\circ}/U})$.

Let $\nbige$ be a reflexive coherent $\nbigo_X(\ast H)$-module
with an integrable relative connection
$\nabla:\nbige\to\nbige\otimes j_{\ast}\Omega^1_{X^{\circ}/U}$
such that $\nbige$ is flat over $S$.
Such $(\nbige,\nabla)$ is called
a meromorphic flat connection on $(X,H)$ relative to $U$.
If $\nbige$ is a locally free $\nbigo_X(\ast H)$-module,
$(\nbige,\nabla)$ is called a meromorphic flat bundle
on $(X,H)$ relative to $U$.

For each geometric point $u:\Spec(\hyperk_1)\to U$,
we obtain the coherent $\nbigo_{X_u}(\ast H_u)$-module $\nbige_u$
equipped with an integrable connection $\nabla_u$.

\begin{lem}
$\nbige_{|X^{\circ}}$ is a locally free $\nbigo_{X^{\circ}}$-module.
\end{lem}
\pf
Because $\nbige_{u|X_u\setminus H_u}$ are flat
$\nbigo_{X_u\setminus H_u}$-modules for any closed point,
we obtain that $\nbige_{|X^{\circ}}$ is flat,
i.e., locally free.
\hfill\qed

\subsubsection{Families of good meromorphic flat bundles}

Suppose that $X$ is smooth
and that $X\to U$ is smooth.
We also assume that
(i) $H$ is normal crossing relative to $S$,
(ii) the monodromy is trivial.

Let $x\in H$ be a closed point.
Let $n=\dim X$.
There exists a parameter system
$(x_1,\ldots,x_n)$ such that
the ideal of $H$ at $x$ is generated by
$\prod_{i=1}^{\ell}x_i$ for some $\ell$.

Let $(\nbigv,\nabla)$ be a meromorphic flat bundle on $(X,H)$
relative to $U$.
We set 
$\nbigv_{|\xhat}=\nbigv\otimes\nbigohat_{X,x}$.
We obtain the induced integrable relative connection
$\nabla:
 \nbigv_{|\xhat}
 \to
 \nbigv_{|\xhat}
 \otimes\Omega^1_{X/U}$.
\begin{df}
$(\nbigv,\nabla)$ is called regular
(resp. good, unramifiedly good) at
a closed point $x\in H$
if 
$(\nbigv_{|\xhat},\nabla)$ is 
regular
(resp. good, unramifiedly good). 
We say that $(\nbigv,\nabla)$ is
good (resp. regular, unramifiedly good) on $(X,H)$
if $(\nbigv,\nabla)_{|\xhat}$ is
good (resp. regular, unramifiedly good) at any closed point $x\in H$. 
\hfill\qed
\end{df}

\subsubsection{Etale coordinate systems}

An etale coordinate neighbourhood
of $x\in X$ relative to $U$
is a tuple
$(\nbigu,x',\varphi,\psi)$
of a smooth affine $\hyperk$-variety $\nbigu$,
a closed point $x'\in \nbigu$,
an etale morphism $\varphi:(\nbigu,x')\to (X,x)$,
and an etale morphism
$\psi:\nbigu\to U\times \Spec(\hyperk[x_1\ldots,x_n])$
satisfying the following conditions.
\begin{itemize}
 \item $H_{\nbigu}:=\varphi^{-1}(H)$ equals
       $\psi^{-1}\bigl(
       \bigcup_{i=1}^{\ell}
       \{z_i=0\}
       \bigr)$
       for some $1\leq \ell\leq n$.
 \item $\psi(x')\in U\times(0,\ldots,0)$.
 \item $\psi^{-1}(U\times\{z_i=0\})$
       are connected.
\end{itemize}
For $e\in\seisuu_{>0}$,
there exists the ramified covering
\[
U\times\Spec(\hyperk[x_1^{1/e},\ldots,x_{\ell}^{1/e},x_{\ell+1},\ldots,x_n])
\to
U\times\Spec(\hyperk[x_1,\ldots,x_n]).
\]
Let $\nbigu^{(e)}$ denote the fiber product of
$\nbigu$ and
$U\times\Spec(\hyperk[x_1^{1/e},\ldots,x_{\ell}^{1/e},x_{\ell+1},\ldots,x_n])$
over $U\times\Spec(\hyperk[x_1,\ldots,x_n])$.
Let $\varphi_{\nbigu,e}:\nbigu^{(e)}\to X$ denote the induced morphism.
We set
$H_{\nbigu}^{(e)}=\varphi_{\nbigu,e}^{-1}(H_{\nbigu})$.

\begin{lem}
Let $(\nbige,\nabla)$ be a good meromorphic flat bundle
 on $(X,H)$ relative to $U$.
Let $\{(\nbigu_i,\varphi_i,\psi_i)\}_{i=1,\ldots,m}$
be an etale covering
of $X$ by etale coordinate systems.
Then,
there exists $e\in\seisuu_{>0}$
such that
$\varphi_{\nbigu_i,e}^{\ast}(\nbige,\nabla)_{|\nbigu_i}$
is unramifiedly good
on  $(\nbigu_i^{(e)},H_{\nbigu_i}^{(e)})$.
\hfill\qed
\end{lem}

\subsubsection{Good lattices}

Let $(\nbige,\nabla)$ be a good meromorphic flat bundle
on $(X,H)$ relative to $U$.
Let $E\subset\nbige$ be a locally free $\nbigo_X$-submodule
such that $\nbige=E(\ast H)$.

\begin{df}
Suppose that $(\nbige,\nabla)$ is unramifiedly good.
Then, $E$ is called unramifiedly good if 
the following holds for any $x\in H$.
\begin{itemize}
 \item We have the decomposition
       $E_{|\xhat}=
       \bigoplus_{\gminia\in\nbigi(\nbige,\nabla,x)}
       \Ehat_{x,\gminia}$
       such that 
       $\Ehat_{x,\gminia}=(\nbigehat_{x,\gminia}\cap E_{|\xhat})$,
       and
       $(\nabla-d\gminiatilde)\Ehat_{x,\gminia}
       \subset
       \Ehat_{x,\gminia}
       \otimes\Omega^1_{X/U}(\log H)$.
       \hfill\qed
\end{itemize}
\end{df}

\begin{df}
Let $\{(\nbigu_i,\varphi_i,\psi_i)\}_{i=1,\ldots,m}$
and $e\in\seisuu_{>0}$
be as above.
We say $E$ is a good lattice of $(\nbige,\nabla)$
if there exist unramifiedly good lattices
$E_{e,i}$ of $\varphi_e^{\ast}(\nbige,\nabla)$
such that 
(i) $E_{e,i}$ is invariant under the Galois group of
$\nbigu_i^{(e)}/\nbigu_i$,
(ii) the Galois descent of $E_{e,i}$ equals $\varphi_i^{\ast}E$. 
\hfill\qed
\end{df}

\subsubsection{Decompositions into the regular part and the irregular part}

Let $(V,\nabla)$ be a good meromorphic flat bundle on $(X,H)$
relative to $U$.
Let $E$ be a good lattice.
Let $x\in H_k$.
We have the natural map
\[
 \rho_k:
 \nbigi(V,\nabla,x)
 \to
 \nbigohat^{(H,e)}_X(\ast H)_x
 \big/
 \nbigohat^{(H,e)}_X(\ast H_{\neq k})_x
\]
Let $\nbigi(V,\nabla,x)_{k,0}=\rho_k^{-1}(0)$.
We set
\[
(E_{|\xhat})^{\reg,k}
 =\bigoplus_{\gminia\in \nbigi(V,\nabla,x)_{k,0}}
  (E_{|\xhat}\cap V^{(H,e)}_{\xhat,\gminia}),
  \quad\quad
(E_{|\xhat})^{\irreg,k}
=\bigoplus_{\gminia\not\in \nbigi(V,\nabla,x)_{k,0}}
 (E_{|\xhat}\cap V^{(H,e)}_{\xhat,\gminia}).
\]
We obtain the decomposition
\begin{equation}
\label{eq;26.2.21.32}
 E_{|\xhat}=
 (E_{|\xhat})^{\reg,k}
 \oplus
 (E_{|\xhat})^{\irreg,k}.
\end{equation}
It is easy to check the following
\begin{equation}
\label{eq;26.2.21.34}
 \nabla (E_{|\xhat})^{\reg,k}
  \subset
  (E_{|\xhat})^{\reg,k}
  \otimes
  \Bigl(
  \Omega^1_{X/U}(\log H_k)
  +
  \Omega^1_{X/U}(\ast H_{\neq k})
  \Bigr),
  \quad\quad
  \nabla (E_{|\xhat})^{\irreg,k}
  \subset
  (E_{|\xhat})^{\irreg,k}
 \otimes
  \Omega^1_{X/U}(\ast H).
\end{equation}
The following proposition is similar to
Proposition \ref{prop;26.2.21.20}.

\begin{prop}
\label{prop;26.2.21.30}
There exists the decomposition
\begin{equation}
\label{eq;26.2.21.33}
 E_{|\Hhat_k}
=E_{\Hhat_k}^{\reg}\oplus
E_{\Hhat_k}^{\irreg}
\end{equation}
which induces {\rm(\ref{eq;26.2.21.32})}
for any $x\in H_k$.
 Moreover, we have
 $\nabla(E_{\Hhat}^{\reg})
 \subset
 E_{\Hhat}^{\reg}
 \otimes
 \Bigl(
 \Omega^1_{X/U}(\log H_k)
+\Omega^1_{X/U}(\ast H_{\neq k})
 \Bigr)$
and 
$\nabla(E_{\Hhat}^{\irreg})
\subset
E_{\Hhat}^{\irreg}
\otimes
\Omega^1_{X/U}(\ast H)
\Bigr)$. 
\hfill\qed
\end{prop}

By (\ref{eq;26.2.21.34}),
we obtain the endomorphism
$\Res_{H_k}(\nabla)$
of $(E^{\reg}_{\Hhat_k|H_k})_{|\xhat}$.

\begin{cor}
\label{cor;26.2.21.40}
Let $E$ be good lattice.
We obtain the residue endomorphism
$\Res_{H_k}(\nabla)$
on $(E_{\Hhat_k}^{\reg})_{|H_k}$
which induces $\Res_{H_k}(\nabla)$
on $(E^{\reg}_{\Hhat_k|H_k})_{|\xhat}$
for any $x\in H_k$.
\hfill\qed
\end{cor}

\subsubsection{Condition for lattices to be unramifiedly good}

Let $(V,\nabla)$ be a meromorphic flat bundle on $(X,H)$
relative to $U$ with a lattice $E$.
Let $\Sigma\subset T^{\ast}((X\setminus H)/U)$ be
a unramifiedly good meromorphic Lagrangian cover.
Assume that $U$ is irreducible.
Let $u:\Spec(\hyperk_1)\to U$ be
a geometric point obtained as an algebraic closure
of the generic point.

\begin{lem}
\label{lem;26.1.21.41}
Assume the following conditions.
\begin{itemize}
 \item $(\nbige_u,\nabla_u)$ is unramifiedly good,
       and that $E_u$ is unramifiedly good.
 \item $\Sigma_{u}$ is a representative of
       $\vecI(\nbige_u,\nabla_u)$.
\end{itemize}
Then, the following holds.
\begin{itemize}
 \item $(\nbige,\nabla)$ is unramifiedly good, and 
       $E$ is unramifiedly good.
\end{itemize}
\end{lem}
\pf
Let $x\in H$ be any closed point.
Let us prove that $E$ is unramifiedly good at $x$.
It is enough to consider the case where
there exists a section $U\to X$
whose image $U_1$ contains $x$.
We may assume that $X$ is affine, $X=\Spec(R)$.
Let $R_{U_1}$ denote the completion of $R$ along $U_1$.
We have
$\nbigi\subset R_{U_1}(\ast H)/R_{U_1}$
and a decomposition
$\Sigma=\bigsqcup_{\gminia\in\nbigi}
 \Sigma_{\gminia}$
such that
$\kappa_{d\gminiatilde}^{\ast}\Sigma_{\gminia}$
are logarithmic.
By using Lemma \ref{lem;26.1.21.40},
we can prove that there exists a decomposition
$E_{|\xhat}=\bigoplus_{\gminia\in \nbigi}\Ehat_{x,\gminia}$
such that 
$(\nabla-d\gminiatilde)E_{x,\gminia}
\subset
E_{x,\gminia}\otimes \Omega^1_{X/U}(\log H)$.
\hfill\qed

\subsubsection{Refinements of parameter spaces}

Let $U,X,H$ be as in \S\ref{subsection;25.1.21.30}.
Let $(\nbige,\nabla)$ be a meromorphic integrable connection
on $(X,H)$ relative to $U$.

\begin{prop}
\label{prop;26.1.21.42}
There exist
a smooth $\hyperk$-variety $\Utilde$
with a surjection $\Utilde\to U$,
and a projective birational morphism
$\varphi:\Xhat_{\Utilde}\to X_{\Utilde}$,
such that 
the following holds.
\begin{itemize}
 \item $\Xhat_{\Utilde}$ is smooth over $\Utilde$,
       and $\varphi^{-1}(H_{\Utilde})$ is normal crossing
       relative to $\Utilde$
       with the trivial monodromy.
 \item $\varphi^{\ast}(\nbige,\nabla)$ is good.
 \item There exists a good lattice $E$ of
       $\varphi^{\ast}(\nbige,\nabla)$.
\end{itemize} 
\end{prop}
\pf
It is enough to consider the case where
$\hyperk$ is an algebraic closure of
a field finitely generated over $\rnum$.
Let $T\subset U$ be any locally closed irreducible smooth subset.
Let $u:\hyperk_1\to T$
be an algebraic closure of the generic point of $T$.
We obtain a meromorphic flat connection
$(\nbige_u,\nabla_u)$ on $(X_u,H_u)$.
There exists a projective birational morphism
$\rho_u:\Xhat_u\to X_u$
such that
(i) $\Hhat_u=\rho_u^{-1}(H_u)$ is simple normal crossing,
(ii) $\Xhat_u\setminus \Hhat_u\simeq X_u\setminus H_u$,
(iii) $\rho_u^{\ast}(\nbige_u,\nabla_u)$ is good on
$(\Xhat_u,\Hhat_u)$.
Let $\Sigma(\nbige_u,\nabla_u)$
be a representative of
the Lagrangian irregularity  $\vecI(\nbige_u,\nabla_u)$.
We fix an embedding $\hyperk_1\subset\cnum$.
There exists a $\cnum$-Deligne-Malgrange lattice
$\nbigp^{\DM}_0(\nbige_u)$
of
$(\nbige_u,\nabla_u)$.
There exist an etale covering of $\Xhat_u$ by
etale coordinate systems
$(\nbigu_{u,i},\varphi_{u,i},\psi_{u,i})$.
There exist unramifiedly good lattices
$\nbigp^{\DM}_{0}\bigl(
\varphi_{u,i,e}^{\ast}(\nbige_u)\bigr)$
of $\varphi_{u,i,e}^{\ast}(\nbige_u,\nabla)$
on $\nbigu_{u,i}^{(e)}$.

There exist a Zariski open subset $T'$ of $T$,
a finite etale covering $\Ttilde'\to T'$,
a model $(\Xhat_{T'},\rho_{T'})$ of $(\Xhat_u,\rho_u)$,
models $(\nbigu_{T',i},\varphi_{T',i},\psi_{T',i})$
of $(\nbigu_{u,i},\varphi_{u,i},\psi_{u,i})$,
a model $\nbigp^{\DM}_{0}(\nbige_{T'})$ of
$\nbigp^{\DM}_{0}(\nbige_u)$,
a model $\Sigma_{T'}$ of $\Sigma(\nbige_u,\nabla_u)$,
and models
$\nbigp^{\DM}_{0}\bigl(
\varphi_{T',i,e}^{\ast}(\nbige_{T'})\bigr)$.
By Lemma \ref{lem;26.1.21.41},
we obtain that
$\nbigp^{\DM}_{0}\bigl(
\varphi_{T',i,e}^{\ast}(\nbige_{T'})\bigr)$
are unramifiedly good.
Hence, we obtain that
$\nbigp^{\DM}_{0}(\nbige_{T'})$ is good.

By using the previous argument
and a Noetherian induction,
we obtain the claim of Proposition \ref{prop;26.1.21.42}.
\hfill\qed

\subsection{Boundedness of good meromorphic flat bundles
with bounded irregularity}

\subsubsection{Meromorphic good systems of ramified irregular values}

Let $X$ be a smooth projective $\hyperk$-variety.
Let $H$ be a normal crossing hypersurface of $X$.
\begin{df}
Let $\vecnbigi=(\nbigi_x\,|\,x\in H)$
be a tuple of
good sets of ramified irregular values
$\nbigi_x\subset
\nbigohat^{(H)}_{X,x}(\ast H)/\nbigohat^{(H)}_{X,x}$.
It is called meromorphic 
if there exists a meromorphic Lagrangian cover
$\Sigma\subset T^{\ast}(X\setminus H)$
such that $\vecnbigi=\vecnbigi(\Sigma)$.
\hfill\qed
\end{df}

Let $\vecnbigi$ be a meromorphic
tuple of good sets of irregular values.
Let $\Sigma$ be a meromorphic Lagrangian cover
such that $\vecnbigi(\Sigma)=\vecnbigi$.
Let $\vecm\in\seisuu_{\geq 0}^{\Lambda}$,
where $\Lambda$ denotes the index set of
the irreducible decomposition
$H=\bigcup_{j\in \Lambda}H_j$.
We say
$\ord(\vecnbigi)\geq \vecm$
if the closure of
$\Sigma$ in $T^{\ast}X(\log H)\bigl(\sum m_iH_i\bigr)$
is proper over $X$.

\begin{lem}
Suppose that $\ord(\vecnbigi)\geq \vecm$.
Let $(V,\nabla)$ be a good meromorphic flat bundle on $(X,H)$
such that $\vecnbigi(V)\subset\vecnbigi$.
Let $V_0$ be any good lattice of $(V,\nabla)$.
Then, we obtain
$\nabla(V_0)\subset V_0\otimes\Omega^1_X(\log H)\bigl(\sum m_iH_i\bigr)$.
\end{lem}
\pf
In the unramified case,
we can check it by direct computation.
The general case is reduced to the unramified case.
\hfill\qed

\subsubsection{Boundedness of good meromorphic flat bundles
as meromorphic objects}
\label{subsection;25.12.14.121}

Let $X_S\to S$ be a smooth projective morphism of $\hyperk$-varieties
such that each closed fiber is connected.
Let $H_S\subset X_S$ be a hypersurface
which is normal crossing relative to $S$
(see \S\ref{subsection;25.12.7.40}).
For any $\hyperk$-variety $S'$ over $S$,
we set $X_{S'}=X_S\times_SS'$
and $H_{S'}=H_S\times_SS'$.
For any closed point $s\in S$,
let $\nbigc^{\good}(X_s,H_s,r,\vecm)$
be a family of good meromorphic flat bundles
$(V,\nabla)$ of rank $r$ on $(X_s,H_s)$
such that $\ord(\vecnbigi)\geq \vecm$.

\begin{thm}
\label{thm;25.10.15.40}
There exist a $\hyperk$-variety $\nbigs$ over $S$,
a locally free $\nbigo_{X_{\nbigs}}$-module $E_{\nbigs}$,
an integrable connection
\[
 \nabla_{\nbigs}:
 E_{\nbigs}\to
 E_{\nbigs}\otimes\Omega^1_{X_{\nbigs}/\nbigs}(\ast H_{\nbigs})
\]
relative to $\nbigs$ such that the following holds.
\begin{itemize}
 \item For any $s\in S$,
       $(V,\nabla)\in
       \nbigc^{\good}(X_s,H_s,r,\vecm)$
       and a $\cnum$-Deligne-Malgrange lattice $V_0$,
       there exist $\stilde\in\nbigs$ over $s$
       such that
       $(V_0,\nabla)\simeq
       (E_{\nbigs},\nabla_{\nbigs})_{\stilde}$.
\item $(E_{\nbigs}(\ast H_{\nbigs}),\nabla_{\nbigs})$
       is a good meromorphic flat bundle
       on $(X_{\nbigs},H_{\nbigs})$
       relative to $\nbigs$.
\end{itemize}
In this sense,
the family $(\nbigc^{\good}(X_s,H_s,r,\vecm)\,|\,s\in S(\hyperk))$
is bounded as meromorphic objects.
\end{thm}
\pf
Let $N\in\seisuu_{>0}$
such that
$\Omega_{X_s}^1(\log H_s)(\sum m_iH_{s,i})
\subset \Omega_{X_s}^1(NH_s)$.
Let $M_2>0$ be as in \S\ref{subsection;25.10.15.12}
for some $s\in S(\hyperk)$.
For $s\in S(\hyperk)$,
let $\nbigg_s(r)$
be the family of locally free sheaves $E$ of rank $r$ on $X_s$
satisfying the following conditions.
\begin{itemize}
 \item $|a_1(E)|+|a_2(E)|\leq C(r)$,
       where
       $a_i(E)$ denote the constants
       as in \S\ref{subsection;25.12.10.1},
       and $C(r)$ denotes a constant
       as in Proposition \ref{prop;25.10.15.11}.
 \item $\mu_{\max}(E)\leq \mu(E)+r(M_2+N\deg_{\omega}(H_s))+A$,
       where $A$ denotes the constant as in Proposition \ref{prop;25.10.15.11}.
       (See \S\ref{subsection;25.10.15.12} for
       $\mu_{\max}(E)$, $\mu(E)$ and $M_2$.)
\end{itemize}
By the boundedness theorem of Maruyama \cite{Maruyama}
(see Proposition \ref{prop;25.12.14.110}),
there exist
a morphism of $\hyperk$-varieties 
$\nbigs_1\to S$
and
a locally free sheaf
$\nbige$
of rank $r$
on $X_{\nbigs_1}$
such that
for any $s\in S(\hyperk)$
and any $E\in \nbigg_s(r)$,
there exists $s_1\in\nbigs_1(\hyperk)$ over $s$
such that $\nbige_s\simeq E$.
By Proposition \ref{prop;24.8.18.11},
there exists a morphism of
$\hyperk$-varieties $\nbigs_2\to\nbigs_1$
and a universal integrable meromorphic connection
\[
\nabla^u:\nbige_{\nbigs_2}\to
\nbige_{\nbigs_2}\otimes
\Omega^1_{X_{\nbigs_2/\nbigs_2}}\bigl(NH_{\nbigs_2}\bigr)
\]
relative to $\nbigs_2$
such that the following holds.
\begin{itemize}
 \item For any $s_1\in \nbigs_1(\hyperk)$
       and an integrable meromorphic connection
       $\nabla:\nbige_{s_1}
       \to \nbige_{s_1}
       \otimes\Omega^1_{X_{s_1}}(NH_{s_1})$,
       there exists
       $s_2\in\nbigs_2$ over $s_1$
       such that
       $(\nbige_{\nbigs_2},\nabla^u)_{s_2}\simeq
       (\nbige_{s_1},\nabla)$.
\end{itemize}
By refining $\nbigs_2$,
we may assume that
$(\nbige_{\nbigs_2}(\ast H_{\nbigs_2}),\nabla^u)$
is good relative to $\nbigs_2$.

For any $s\in S(\hyperk)$ and
$(V,\nabla)\in \nbigc^{\good}(X_s,H_s,r,\vecm)$,
there exists
a $\cnum$-Deligne-Malgrange lattice $V_0$ of $V$
as in Proposition \ref{prop;25.10.15.11}.
By Corollary \ref{cor;26.1.3.20},
$V_0$ is contained in $\nbigg_s(r)$.
By the construction,
there exists
$s_2\in \nbigs_2$
such that
$(V_0,\nabla)\simeq
(\nbige,\nabla^u)_{s_2}$.
\hfill\qed

\begin{rem}
Later, we shall study
the boundedness of the family 
$(\nbigc^{\good}(X_s,H_s,r,\vecm)\,|\,s\in S(\hyperk))$
as $\nbigd$-modules
(Theorem {\rm\ref{thm;26.2.7.31}}).
 \hfill\qed
\end{rem}

\subsection{Boundedness of meromorphic flat bundles with bounded irregularity}
\label{subsection;26.2.7.32}

Let $X$ be a normal projective $\hyperk$-variety
with a hypersurface $H$
such that $X\setminus H$ is smooth.
There exists a closed subset $Z\subset H$
with $\codim_X(Z)\geq 2$
such that
(i) $X^{\circ}=X\setminus Z$ is smooth,
(ii) $H^{\circ}=H\setminus Z$ is a normal crossing hypersurface of $X^{\circ}$.
We also set $H_i^{\circ}=H_i\setminus Z$.

For $\vecm\in\seisuu_{>0}^{\Lambda}$,
let $\nbigc(X,H,r,\vecm)$ be the family of
meromorphic flat bundles $(V,\nabla)$ of rank $r$
on $(X,H)$
such that the following holds.
(See Proposition \ref{prop;25.10.15.51} for the equivalent conditions.)
\begin{itemize}
 \item There exists a reflexive lattice $V_0\subset V$
       such that
       $\nabla(V_{0|X^{\circ}})\subset
       V_0\otimes\Omega_{X^{\circ}}^1(\log H^{\circ})(\sum m_iH^{\circ}_i)$.
\end{itemize}

\begin{thm}
\label{thm;25.12.14.30}
There exist a smooth $\hyperk$-variety $\nbigs$,
a smooth $\hyperk$-variety $\Xhat_{\nbigs}$ over $\nbigs$
 with a morphism $\rho_{\nbigs}:
 \Xhat_{\nbigs}\to X\times\nbigs$  over $\nbigs$,
a hypersurface $\Hhat_{\nbigs}\subset \Xhat_{\nbigs}$
normal crossing relative to $\nbigs$,
a good meromorphic flat bundle $(\nbigv,\nabla)$ on
$(\Xhat_{\nbigs},\Hhat_{\nbigs})$ relative to $\nbigs$,
and a lattice $\Ehat_{\nbigs}\subset \nbigv$ 
such that 
the following holds.
\begin{itemize}
 \item $\rho^{-1}_{\nbigs}(H\times\nbigs)=\Hhat_{\nbigs}$,
       and $\Xhat_{\nbigs}\setminus\Hhat_{\nbigs}\simeq
       (X\setminus H)\times\nbigs$.
 \item For any $(V,\nabla)\in \nbigc(X,H,r,\vecm)$,
       there exists $s\in\nbigs$
       such that
       $\rho_s^{\ast}(V,\nabla)
       \simeq
       (\nbigv,\nabla)_s$
       and that
       $\Ehat_s$
       is a $\cnum$-Deligne-Malgrange lattice of
       $\rho_s^{\ast}(V,\nabla)$.
\end{itemize}
\end{thm}
\pf
Let $N>m_i+1$ for any $i$.
By Proposition \ref{prop;25.10.18.10},
there exist a smooth $\hyperk$-variety $S$,
a Zariski closed subset
$\Sigma_S\subset S\times T^{\ast}(X\setminus H)$
and a projective morphism
$\rho_S:\Xhat_S\to S\times X$,
such that the following holds.
\begin{itemize}
 \item $\Sigma_S\to S\times (X\setminus H)$
       is a meromorphic Lagrangian cover over $S\times(X,H)$,
       and $\Sigma_S$ is flat over $S$.
       The number of any fiber is less than $r$.
 \item For any any meromorphic Lagrangian cover
       $\Sigma\subset T^{\ast}(X\setminus H)$
       over $(X^{\circ},H^{\circ})$
       such that
       the closure of $\Sigma$
       in $T^{\ast}X^{\circ}(NH^{\circ})$ is proper over $X^{\circ}$,
       there exists $s\in S$
       such that
       $\Sigma=\Sigma_S\times_Ss=:\Sigma_s$.
 \item $\Xhat_S$ is smooth projective over $S$,
       and
       $\Hhat_S=\rho_S^{-1}(S\times H)$
       is normal crossing relative to $S$.
 \item For each $s\in S$,
       $\rho_s^{\ast}(\Sigma_s)$
       is good on
       $(\Xhat_s,\Hhat_s)=
       (\Xhat_S,\Hhat_S)\times_Ss$,
       where $\rho_s:\Xhat_s\to X$
       denotes the morphism induced by $\rho_S$.
\end{itemize}
For any smooth $\hyperk$-variety $T$ over $S$,
we set
$(\Xhat_T,\Hhat_T)=(\Xhat_S,\Hhat_S)\times_ST$.

Let $S'$ be a connected component of $S$.
Let
$\Hhat_{S'}
=\bigcup_{i\in\Lambda(S')}
\Hhat_{S',i}$
denote the irreducible decomposition.
For each $s\in S'$,
we obtain the tuple $\vecnbigi(\Sigma_s)$
on $(\Xhat_s,\Hhat_s)$.
There exists $\vecm(S')\in\seisuu_{\geq 0}^{\Lambda(S')}$
such that
$\ord(\vecnbigi(\Sigma_s))\geq \vecm(S')$
for any $s\in S'$.

By Theorem \ref{thm;25.10.15.40},
there exist
a morphism of $\hyperk$-varieties
$S_1\to S$,
a locally free
$\nbigo_{\Xhat_{S_1}}(\ast \Hhat_{S_1})$-module
$\nbigvtilde^u_{S_1}$
with an integrable meromorphic connection
$\nabla:
 \nbigvtilde^u_{S_1}
 \to
 \nbigvtilde^u_{S_1}\otimes
 \Omega^1_{\Xhat_{S_1}/S_1}$
relative to $S_1$,
and a lattice $\Ehat_{S_1}$ of $\nbigvtilde^u_{S_1}$
such that the following holds.
\begin{itemize}
 \item $(\nbigvtilde^u_{S_1},\nabla)$
       is good relative to $S_1$.
 \item Let $s'\in S'\subset S$ be any closed point.
       Let $(V,\nabla)$
       be a good meromorphic flat bundle
       on $(\Xhat_s,\Hhat_s)$
       such that
       $\vecnbigi(V)\geq \vecm(S')$.
       Then,
       there exists $s_1\in S_1$ over $s$
       such that
       $(\nbigvtilde^u_{S_1},\nabla)_{|s_1}\simeq
       \rho_{s_1}^{\ast}(V,\nabla)$,
       and that
       $\Ehat_{s_1}$ is a $\cnum$-DM lattice of
       $\rho_{s_1}^{\ast}(V,\nabla)$.       
\end{itemize}
By refining $S_1$,
we may also assume the following.
\begin{itemize}
 \item The meromorphic Lagrangian irregularity
       of $(\nbigvtilde^u_{S_1},\nabla)$
       equals
       $\vecI(\Sigma_{S_1})$.
\end{itemize}

Let $(V,\nabla)\in\nbigc(X,H,r,\vecm)$.
Let $\vecI(V,\nabla)$ be the associated meromorphic Lagrangian irregularity.
Take any $\Sigma(V,\nabla)\in\vecI(V,\nabla)$.
By Proposition \ref{prop;25.10.15.51},
the closure of $\Sigma(V,\nabla)$
in $T^{\ast}(X^{\circ})(\log H^{\circ})(\sum m_iH^{\circ}_i)$
is proper over $X^{\circ}$.
There exists $s\in S$
such that $\Sigma(V,\nabla)=\Sigma_{s}$,
and $\rho_s^{\ast}\Sigma(V,\nabla)$ is good on
$(\Xhat_s,\Hhat_s)$.
By Theorem \ref{thm;25.12.8.100},
$\rho_s^{\ast}(V,\nabla)$ is good on $(\Xhat_s,\Hhat_s)$,
and $\vecnbigi(\rho_s^{\ast}(V,\nabla))=\vecnbigi(\Sigma(V))$.
Hence,
there exists $s_1\in S_1$ over $s$
such that
$(\nbigvtilde^u_{S_1},\nabla)_{|s_1}\simeq\rho_{s_1}^{\ast}(V,\nabla)$,
and that
$\Ehat_{s_1}$ is a $\cnum$-Deligne-Malgrange lattice of
$\rho_{s_1}^{\ast}(V,\nabla)$.
\hfill\qed

\begin{cor}
\label{cor;26.2.23.2}
There exists a smooth $\hyperk$-variety $\nbigs_1$
and a coherent torsion-free $\nbigo_{\nbigs_1\times X}$-module
$E_{\nbigs_1}$
with a meromorphic integrable connection 
\[
 \nabla_{\nbigs_1}:E_{\nbigs_1}\to
 E_{\nbigs_1}\otimes\Omega^1_{\nbigs_1\times X/\nbigs_1}
 (\ast (H\times\nbigs_1))
\]
relative to $\nbigs_1$ such that the following holds.
\begin{itemize}
 \item $E_{\nbigs_1}$ is a reflexive $\nbigo_{\nbigs_1\times X}$-module
       flat over $\nbigs_1$.
 \item For any $(V,\nabla)\in \nbigc(X,H,r,\vecm)$,
       there exist $s\in\nbigs_1$
       and a $\cnum$-Deligne-Malgrange lattice $V_0$
       of $(V,\nabla)$
       such that
       $(E_{\nbigs_1},\nabla_{\nbigs_1})_{s_1}
       \simeq (V_0,\nabla)$.
\end{itemize} 
In this sense, the family $\nbigc(X,H,r,\vecm)$ is bounded.
\end{cor}
\pf
By refining $\nbigs$,
we may assume that
$\rho_{\nbigs\ast}\Ehat_{\nbigs}$
and
$(\rho_{\nbigs\ast}\Ehat_{\nbigs})^{\lor}$
are flat over $\nbigs$.
We set
$E_{\nbigs}=
\bigl(
(\rho_{\nbigs\ast}\Ehat_{\nbigs})^{\lor}
\bigr)^{\lor}$.
By refining $\nbigs$,
we may assume that
$(E_{\nbigs})
\simeq
\bigl((\rho_{\nbigs\ast}\Ehat_{\nbigs})^{\lor}\bigr)^{\lor}$.
Moreover, we may assume
$(E_{\nbigs})_s
\simeq
\bigl((\rho_{\nbigs\ast}\Ehat_{\nbigs})^{\lor}_s\bigr)^{\lor}$
for any $s\in\nbigs(\hyperk)$.
In particular,
$(E_{\nbigs})_s$ are reflexive.

We have the meromorphic integrable relative connection
$\nabla_{\nbigs}$ of $E_{\nbigs}$
induced by $(\Ehat_{\nbigs},\nabla)$.
If $(\Ehat_{\nbigs})_s$ is a good $\cnum$-Deligne-Malgrange lattice
of $((\Ehat_{\nbigs})_s(\ast \Hhat_s),\nabla)$
for some $s\in \nbigs(\hyperk)$,
there exists $Z\subset H$ such that
$\codim_X(Z)\geq 2$
such that
$(E_{\nbigs})_{s|X\setminus Z}$ is
a $\cnum$-Deligne-Malgrange lattice
of $((E_{\nbigs})_s(\ast H),\nabla)_{|X\setminus Z}$.
Because $(E_{\nbigs})_s$ is reflexive,
$(E_{\nbigs})_s$ is a $\cnum$-DM lattice
of $((E_{\nbigs})_s(\ast H),\nabla)$. 
\hfill\qed

\section{Families of $\nbigd$-modules}

\subsection{Families of $\nbigd$-modules}

\subsubsection{$\nbigd_{X_S/S}$-modules}

Let $\hyperk$ be an algebraically closed field of characteristic $0$.
Let $S$ be a $\hyperk$-scheme.
Let $X_S$ be a quasi-projective smooth scheme over $S$
of relative dimension $n$,
i.e., $n=\rank\Omega^1_{X_S/S}$.
Let $\nbigd_{X_S/S}$ denote
the sheaf of differential operators
on $X_S$ relative to $S$.
The sheaf of algebras $\nbigd_{X_S/S}$
is a coherent sheaf of rings
in the sense of \cite[Definition 1.4.8]{Hotta-Takeuchi-Tanisaki}.

We set $\Omega_{X_S/S}=\bigwedge^n\Omega_{X_S/S}^1$.
For any left $\nbigd_{X_S/S}$-module $M$,
we obtain a right $\nbigd_{X_S/S}$-module
$\Omega_{X_S/S}\otimes_{\nbigo_{X_S}} M$.
Conversely,
for any right $\nbigd_{X_S/S}$-module $N$,
we obtain a left $\nbigd_{X_S/S}$-module
$N\otimes_{\nbigo_{X_S}}(\Omega_{X_S/S})^{-1}$.
These procedures induce an equivalence
between left $\nbigd_{X_S/S}$-modules
and right $\nbigd_{X_S/S}$-modules.
In this paper,
we shall consider left $\nbigd_{X_S/S}$-modules.
The following lemma is clear.
\begin{lem}
The category
of $\nbigd_{X_S/S}$-modules 
$\Mod(\nbigd_{X_S/S})$ is an abelian category.
The full subcategory of 
 coherent
 $\nbigd_{X_S/S}$-modules
$\Mod_c(\nbigd_{X_S/S})$ is
an abelian subcategory,
and  closed under extensions.
Similarly, the full subcategory of
quasi-coherent $\nbigd_{X_S/S}$-modules
$\Mod_{qc}(\nbigd_{X_S/S})$ is
an abelian subcategory,
and closed under extensions.
 \hfill\qed
\end{lem}

Let $D(\nbigd_{X_S/S})$ denote
the derived category of
complexes of $\nbigd_{X_S/S}$-modules.
Let $D_c(\nbigd_{X_S/S})$
(resp. $D_{qc}(\nbigd_{X_S/S})$)
denote
the derived category of
cohomologically coherent
(resp. quasi-coherent)
complexes of $\nbigd_{X_S/S}$-modules.
Let $D^b(\nbigd_{X_S/S})\subset D(\nbigd_{X_S/S})$,
$D^b_c(\nbigd_{X_S/S})\subset D_c(\nbigd_{X_S/S})$
and $D^b_{qc}(\nbigd_{X_S/S})\subset D_{qc}(\nbigd_{X_S/S})$
denote the full subcategories of
cohomologically bounded objects.

\begin{df}
A $\nbigd_{X_S/S}$-module is called
$S$-flat if it is flat over $\nbigo_S$. 
An object of $D(\nbigd_{X_S/S})$
is called $S$-flat if 
it is cohomologically $S$-flat.
\hfill\qed
\end{df}

\begin{rem}
The analytic theory for 
family of $\nbigd$-modules 
has been developed in 
{\rm\cite{Fiorot-Monteiro-Fernandes-Sabbah1,
Fiorot-Monteiro-Fernandes-Sabbah2}}.
\hfill\qed
\end{rem}

\subsubsection{Functoriality in the parameter direction}

Let $g:S'\to S$ be any morphism of $\hyperk$-schemes.
We set $X_{S'}=X_S\times_SS'$.
The induced morphism 
$X_{S'}\to X_S$ is also denoted by $g$.
We have the natural morphism of sheaves of algebras
$g^{-1}(\nbigd_{X_S/S})\to \nbigd_{X_{S'}/S'}$.

\paragraph{Inverse image}
We assume that $S$ is a smooth $\hyperk$-variety.
For any $\nbigd_{X_S/S}$-module $M$,
\[
g^{\ast}(M)
:=\nbigo_{X_{S'}/S'}\otimes_{g^{-1}\nbigo_{X_S}}g^{-1}(M)\simeq
\nbigd_{X_{S'}/S'}\otimes_{g^{-1}(\nbigd_{X_S/S})}
g^{-1}(M)
\]
is naturally a $\nbigd_{X_{S'}/S'}$-module.
We obtain
$Lg^{\ast}:D^b_c(\nbigd_{X_S/S})\to D_c^b(\nbigd_{X_{S'}/S'})$
by setting
\[
 Lg^{\ast}(M^{\bullet})
 =\nbigd_{X_{S'}/S'}\otimes^L_{g^{-1}(\nbigd_{X_{S}/S})}
g^{-1}(M^{\bullet}).
\]
It also induces
$Lg^{\ast}:D^b_{qc}(\nbigd_{X_S/S})\to D^b_{qc}(\nbigd_{X_{S'}/S'})$.

\begin{lem}
If $M$ is $S$-flat,
we have $Lg^{\ast}(M)\simeq g^{\ast}(M)$.
\end{lem}
\pf
Let $\hyperk_{X_{S'}}$
denote the sheaf of locally constant
$\hyperk$-valued functions on $X_{S'}$.
It is naturally a sheaf of algebras.
Let $D(\hyperk_{X_{S'}})$
denote the derived category of
$\hyperk_{X_{S'}}$-complexes.
There exists the natural forgetful functor
$a:D(\nbigd_{X_{S'}/S'})\to D(\hyperk_{X_{S'}})$.
It is enough to prove the vanishing of
the $j$-th cohomology sheaves
of $a(Lg^{\ast}(M))$ for any $j<0$.

Let $\pi:X_{S'}\to S'$ denote the projection.
There exist the natural isomorphisms
$\nbigo_{X_{S'}}
\simeq\pi^{-1}(\nbigo_{S'})\otimes_{\pi^{-1}g^{-1}\nbigo_S}
g^{-1}(\nbigo_{X_S})$,
and 
\[
 \nbigd_{X_{S'}/S'}
\simeq\pi^{-1}(\nbigo_{S'})\otimes_{\pi^{-1}g^{-1}\nbigo_S}
g^{-1}(\nbigd_{X_S/S}).
\]
Let $\nbigp\to \nbigo_{S'}$
be a $g^{-1}(\nbigo_{S})$-flat resolution.
We obtain the quasi-isomorphism of 
the complex of right $g^{-1}\nbigd_{X_S/S}$-modules
\[
 \pi^{-1}(\nbigp)\otimes_{\pi^{-1}g^{-1}\nbigo_{S}}
 g^{-1}(\nbigd_{X_S/S})
 \to
 \nbigd_{X_{S'}/S'}.
\]
It is a $g^{-1}(\nbigd_{X_S/S})$-flat resolution of
$\nbigd_{X_{S'}/S'}$.
Then, $a(Lg^{\ast}(M))$
is represented by
\[
\Bigl(
 \pi^{-1}(\nbigp)
 \otimes_{\pi^{-1}g^{-1}\nbigo_{S}}
 g^{-1}(\nbigd_{X_{S}/S})
 \Bigr)
 \otimes_{g^{-1}\nbigd_{X_{S}/S}}g^{-1}(M)
 \simeq
 \pi^{-1}\nbigp\otimes_{\pi^{-1}g^{-1}(\nbigo_{S})}g^{-1}(M).
\]
Because $M$ is $S$-flat,
we obtain
the vanishing of
the $j$-th cohomology sheaves
of $a(Lg^{\ast}(M))$ for any $j<0$.
\hfill\qed

\begin{rem}
If $M$ is $S$-flat,
$g^{\ast}(M)$ is often denoted by $M_{S'}$.
\hfill\qed
\end{rem}

\paragraph{Direct image}

We assume that both $S$ and $S'$ are smooth $\hyperk$-varieties.
For any $\nbigd_{X_{S'}/S'}$-module $N$,
we obtain
the $\nbigd_{X_S/S}$-module $g_{X\ast}(N)$.
For any $\nbigd_{X_{S'}/S'}$-complex $N^{\bullet}$,
by taking a $\nbigd_{X_{S'}/S'}$-injective resolution
$N^{\bullet}\to \nbigj$,
we obtain a complex of $\nbigd_{X_S/S}$-modules
$g_{X\ast}(\nbigj)$.
In this way,
we obtain
$Rg_{\ast}:
D^b(\nbigd_{X_{S'}/S'})\to D^b(\nbigd_{X_S/S})$.
We assume that both $S$ and $S'$ are
smooth $\hyperk$-varieties.
It also induces
$Rg_{\ast}:
D^b_{qc}(\nbigd_{X_{S'}/S'})\to D^b_{qc}(\nbigd_{X_S/S})$.
If $g$ is proper,
we also obtain
$D^b_{c}(\nbigd_{X_{S'}/S'})\to D^b_{c}(\nbigd_{X_S/S})$.

\subsubsection{Good filtrations}

Assume that $S$ is a Noetherian $\hyperk$-scheme.
We have the natural filtration
$F(\nbigd_{X_S/S})$
by the order of differential operators
with which $\nbigd_{X_S/S}$
is a sheaf of filtered rings
in the sense of
\cite[A.1(c)]{kashiwara_text}.
Let $M$ be a coherent $\nbigd_{X_S/S}$-module.
The notion of good filtration for
coherent $\nbigd_{X_S/S}$-module is defined
in the standard way.
Namely, it is an increasing filtration $F(M)$
by coherent $\nbigo_{X_S}$-submodules
such that the following holds.
\begin{itemize}
 \item $\bigcup F_j(M)=M$.
       There exist $j_0(x)$ such that
       $F_{j_0}(M)=0$.
 \item $F_j\nbigd_{X_S/S}\cdot F_k(M)\subset F_{j+k}(M)$.
 \item $\bigoplus F_j(M)$ is coherent over
       $\bigoplus F_j(\nbigd_{X_S/S})$.
\end{itemize}
As in the case of $\nbigd_X$-modules,
such a good filtration always exists.
Moreover, if $F$ and $F'$ are two such filtrations,
there exists $i_0>0$ such that
$F'_{j-i_0}(M)\subset F_j(M)\subset F'_{j+i_0}(M)$
for any $j\in\seisuu$
(See \cite[Theorem 2.1.3]{Hotta-Takeuchi-Tanisaki}.)
By taking the graduation with respect to the filtration,
we obtain
the coherent $\Gr^F(\nbigd_{X_S/S})$-module $\Gr^F(M)$.
Let $T^{\ast}(X_S/S)$ denote the relative cotangent bundle
of $X_S$ over $S$.
We obtain 
a coherent $\nbigo_{T^{\ast}(X_S/S)}$-module $\Gr^F(M)^{\sim}$
from $\Gr^F(M)$.
The support is independent of
the choice of a good filtration
as in the case of $\nbigd_X$-modules
(see \cite[Theorem 2.2.1]{Hotta-Takeuchi-Tanisaki}).
The support is denoted by $\Ch(M)$, called
the characteristic variety of $M$.

\begin{df}
A good filtration $F(M)$ is called $S$-flat if
$\Gr^F(M)$ is flat over $\nbigo_S$.
\hfill\qed
\end{df}

\begin{lem}
Suppose that $M$ has an $S$-flat good filtration $F$.
\begin{itemize}
 \item The $\nbigd_{X_S/S}$-module
       $M$ and
       the $\nbigo_{X_S}$-modules $F_j(M)$ are flat over $S$.
 \item Let $\iota:S'\to S$ be an immersion
       of a locally closed subset.
       The induced morphism
       $X_{S'}\to X_S$ is also denoted by $\iota$.
       Then,
       $\iota^{\ast}F(M)=(\iota^{\ast}F_j(M))$
       is an $S'$-flat good filtration of
       the $\nbigd_{X_{S'}/S'}$-module
       $\iota^{\ast}(M)$.
       We also have
       $\iota^{\ast}\Gr^F(M)
       =\Gr^{\iota^{\ast}F}(\iota^{\ast}M)$.       
\hfill\qed
\end{itemize}
\end{lem}

\begin{lem}
\label{lem;25.12.10.40}
Let $M$ be any coherent $\nbigd_{X_S/S}$-module.
For any good filtration $F$ of $M$,
there exists a non-empty open subset $S'$ of $S$
such that
the induced good filtration $F$ of $M_{|X_{S'}}$
is $S'$-flat.
In particular, $M_{|X_{S'}}$ is $S'$-flat.
\end{lem}
\pf
By \cite[Theorem 6.9.1]{EGA-IV-2},
there exists a non-empty open subset $S'$ of $S$
such that $\Gr^F(M)_{|X_{S'}}$ is flat over $S'$.

\hfill\qed

\begin{cor}
\label{cor;26.1.26.20}
For any morphism $\varphi:M_1\to M_2$ of
coherent $\nbigd_{X_S/S}$-modules,
There exists a non-empty Zariski open subset $S'$
such that
$\Ker\varphi_{|X_{S'}}$, $\Image\varphi_{|X_{S'}}$
and $\Cok\varphi_{|X_{S'}}$ 
have $S'$-flat good filtrations.
Here,
$\varphi_{|X_{S'}}$ denotes
the induced morphism
$M_{1|X_{S'}}\to M_{2|X_{S'}}$.
\hfill\qed
\end{cor}

\begin{cor}
For any $S$-flat coherent $\nbigd_{X_S/S}$-module $M$,
there exists a decomposition $S=\bigsqcup_{j\in \Gamma}S_j$
such that
each $M_{S_j}$ has an $S_j$-flat good filtration.
For any morphism of $S$-flat coherent $\nbigd_{X_S/S}$-modules
$\varphi:M_1\to M_2$,
there exists a decomposition
$S=\bigsqcup S_j$ such that
$M_{1,S_j}$,
$M_{2,S_j}$
$\Ker\varphi_{S_j}$,
$\Image\varphi_{S_j}$
and $\Cok\varphi_{S_j}$
have $S$-flat good filtrations.
Here, $\varphi_{S_j}:M_{1,S_j}\to M_{2,S_j}$
denotes the induced morphism.
\hfill\qed
\end{cor}

Recall that a geometric point of $S$
is a morphism $\eta:\Spec(\hyperk_1)\to S$,
where $\hyperk_1$ is an algebraically closed field.
The induced $\hyperk_1$-variety is denoted by $X_{\eta}$.
The induced $\nbigd_{X_{\eta}}$-module
is denoted by $M_{\eta}$.

\begin{cor}
Suppose that $S$ is irreducible.
Let $\eta$ be a geometric point
obtained as an algebraic closure of 
the generic point of $S$.
Let $M$ be a coherent $\nbigd_{X_S/S}$-module
such that $M_{\eta}=0$.
Then, there exists a dominant quasi-finite etale morphism
of $\hyperk$-schemes $f:S'\to S$  
such that 
(i) $\eta$ factors through $f$,
(ii) $f^{\ast}M=0$.
\end{cor}
\pf
It is enough to consider the case
where $S$, $X_S$ and $T^{\ast}(X_S/S)$ are affine.
We may assume that $M$ has an $S$-flat good filtration $F$.
It induces a good filtration $F_{\eta}$ of $M_{\eta}$.
Because $\Gr^{F_{\eta}}(M_{\eta})=0$,
we obtain $\eta^{\ast}\Gr^F(M)=0$.
Because $T^{\ast}(X_S/S)$ is affine,
there exist sections $a_1,\ldots,a_m$
generating $\Gr^F(M)$
over $\nbigo_{T^{\ast}(X_S/S)}$.
There exists
a dominant quasi-finite etale morphism
$f:S'\to S$ such that 
$f^{\ast}(m_i)$ are $0$
in $f^{\ast}\Gr^F(M)$.
It implies 
$f^{\ast}\Gr^F(M)=0$
and $f^{\ast}(M)=0$.
\hfill\qed

\subsubsection{Holonomic $\nbigd_{X_S/S}$-modules}

Let $M$ be an $S$-flat coherent $\nbigd_{X_S/S}$-module.

\begin{df}
\label{df;25.12.12.1}
$M$ is called holonomic if there exists a
decomposition $\bigsqcup_{j\in\Gamma}S_j$
by locally closed irreducible subsets
such that 
(i) each $M_{S_j}$ has an $S_j$-flat good filtration,
(ii) $\dim \Ch(M_{S_j})=\dim X_{S_j}$.
\hfill\qed
\end{df}

\begin{lem}
A coherent $\nbigd_{X_S/S}$-module
$M$ is holonomic if and only if
the following holds.
\begin{description}
 \item[(A)] 
       For any geometric point $\eta$ of $S$,
       $M_{\eta}$ is a holonomic $\nbigd_{X_{\eta}}$-module.
\end{description}
\end{lem}
\pf
There exists a decomposition $S=\bigsqcup S_j$
such that each $M_{S_j}$ has an $S_j$-flat good filtration $F$.
Let $\eta:\Spec(\hyperk_1)\to S$ be any geometric point of $S$.
There exists $S_j$ such that
$\eta$ factors through $S_j$.
The $S_j$-flat good filtration $F$ of $M_{S_j}$ induces
a good filtration $F_{\eta}$ of $M_{\eta}$
for which we have
$\Gr^{F_{\eta}}(M_{\eta})
=\Gr^F(M_{S_j})\otimes_{\nbigo_{S_j}}\hyperk_1$.
Then, the claim is clear.
\hfill\qed

\subsubsection{Model of morphisms}

Suppose that $S$ is irreducible.
Let $M_i$ $(i=1,2)$ be coherent $\nbigd_{X_S/S}$-modules.
\begin{lem}
Let $\eta:\Spec\hyperk_1\to S$ be a geometric point
obtained as an algebraic closure of
the generic point of $S$.
For any morphism $\varphi:M_{1,\eta}\to M_{2,\eta}$,
there exist
a dominant quasi-finite etale morphism of $\hyperk$-varieties $f:S'\to S$
and a morphism of $\nbigd_{X_{S'}/S'}$-modules
$\varphi_{S'}:f^{\ast}(M_1)\to f^{\ast}(M_2)$
such that
(i) $\eta$ factors through $f$,
(ii) $\varphi$ is induced by $\varphi_{S'}$,
(iii) $f(S')$ contains a Zariski open subset of $S$.
 \hfill\qed
\end{lem}

\subsection{Basic functoriality}

\subsubsection{Pull back in the fiber direction}

Let $S$ be a smooth $\hyperk$-variety.
Let $f:X_S\to Y_S$ be a morphism $\hyperk$-varieties over $S$.
We set
$\nbigd_{(X_S\to Y_S)/S}:=
f^{\ast}\nbigd_{Y_S/S}$,
which is naturally
a left $\nbigd_{X_S/S}$-module
and a right $f^{-1}\nbigd_{Y_S/S}$-module.
For any left $\nbigd_{Y_S/S}$-module $M$,
we set
\[
Lf^{\ast}(M)=
\nbigd_{(X_S\to Y_S)/S}
\otimes^L_{f^{-1}\nbigd_{Y_S/S}}
f^{-1}(M).
\]
It induces
$Lf^{\ast}:D^b(\nbigd_{Y_S/S})\to
D^b(\nbigd_{X_S/S})$.
It also induces
$Lf^{\ast}:D^b_{qc}(\nbigd_{Y_S/S})\to
D^b_{qc}(\nbigd_{X_S/S})$
as in the case of $\nbigd$-modules.
(See \cite[Proposition 1.5.8]{Hotta-Takeuchi-Tanisaki}.
See also \cite[Remark 1.5.10]{Hotta-Takeuchi-Tanisaki}.)

\subsubsection{Direct images in the fiber direction}

Let $f:X_S\to Y_S$ be a projective morphism.
We set
\[
 \nbigd_{(Y_S\larr X_S)/S}=
 \Omega_{X_S/S}^{-1}\otimes_{f^{-1}\nbigo_{Y_S}}
 f^{\ast}(\nbigd_{Y_S/S}\otimes\Omega_{Y_S/S}).
\]
For any $\nbigd_{X_S/S}$-complex $M^{\bullet}$,
we obtain
\[
 f_{\dagger}(M^{\bullet}):=
 Rf_{\ast}\Bigl(
 \nbigd_{(Y_S\larr X_S)/S}
 \otimes^L_{\nbigd_{X_S/S}}M^{\bullet}
 \Bigr)
\in D^b(\nbigd_{Y_S/S}).
\]

The following is an analogue of
\cite[Theorem 2.5.1]{Hotta-Takeuchi-Tanisaki}.
\begin{prop}
We obtain the functors
 $f_{\dagger}:D^b_{\sharp}(\nbigd_{X_S/S})\to
 D^b_{\sharp}(\nbigd_{Y_S/S})$
$(\sharp=c,qc)$.
\hfill\qed
\end{prop}

For any locally closed smooth subset $T\subset S$,
let $f_T:X_T\to Y_T$ denote the induced morphism.
Let $\iota_X:X_T\to X_S$
and $\iota_Y:Y_T\to Y_S$ denote the inclusions.
\begin{lem}
For $M\in D^b_c(\nbigd_{Y_S/S})$,
there exists the natural isomorphism
\[
 f_{S\dagger}\iota_{X\ast}(M)\simeq
 \iota_{Y\ast}f_{T\dagger}(M).
\]
\end{lem}
\pf
There exists the following natural isomorphism
\[
 \nbigd_{(Y_S\larr X_S)/S}
 \otimes^L_{\nbigd_{X_S/S}}
 \iota_{X\ast}M
\lrarr
 \iota_{X\ast}
 \bigl(
 \iota_X^{-1}
 \nbigd_{(Y_S\larr X_S)/S}
 \otimes^L_{\iota_X^{-1}\nbigd_{X_S/S}}
 M
 \bigr)
 \lrarr
 \iota_{X\ast}
 \bigl(
 \nbigd_{(Y_T\larr X_T)/T}
 \otimes^L_{\nbigd_{X_T/T}}
 M
 \bigr).
\]
It is easy to check that
the composition is an isomorphism
in $D^b(f_S^{-1}(\nbigd_{Y_S/S}))$.
We obtain the following isomorphism
in $D^b(\nbigd_{Y_S/S})$:
\[
 f_{S\dagger}(\iota_{X\ast}M)
 \simeq
 f_{S\ast}\circ\iota_{X\ast}\bigl(
  \nbigd_{(Y_T\larr X_T)/T}
 \otimes^L_{\nbigd_{X_T/T}}
 M
 \bigr)
 \simeq
 \iota_{Y\ast}f_{T\dagger}(M).
\]
\hfill\qed

\begin{prop}
\label{prop;25.12.10.41}
Let $M$ be an $S$-flat coherent $\nbigd_{X_S/S}$-module.
Then, 
there is a decomposition  $S=\bigsqcup_{i\in\Lambda} S_i$
by locally closed smooth subsets
such that 
$f^k_{S_i\dagger}(M_{S_i})$ $(k\in\seisuu)$
are $S_i$-flat coherent $\nbigd_{X_{S_i}/S_i}$-modules.
For any locally closed smooth subset $T\subset S_i$,
there exist the natural isomorphisms
 $f^k_{S_i\dagger}(M_{S_i})_{T}
 \simeq
 f^k_{T\dagger}(M_{T})$.
\end{prop}
\pf
Let $S'\subset S$ be any locally closed smooth subset.
By Lemma \ref{lem;25.12.10.40},
there exists a non-empty open subset $S''\subset S'$
such that
$f^k_{S''\dagger}(M_{S''})$ are $S''$-flat.

Let $T$ be a locally closed smooth subset of $S$.
Let $\iota_X:X_T\to X_S$ and $\iota_Y:Y_T\to Y_S$
denote the induced inclusions.
Let
$\nbigj$ be a cohomologically bounded
$\nbigd_{X_S/S}$-injective complex.
Let
$\iota_X^{-1}\nbigj^{\bullet}\to \nbigj_1$
be an $\iota_X^{-1}\nbigd_{X_S/S}$-injective complex.
There exist the following natural morphisms:
\[
 \iota_Y^{-1}f_{S\ast}(\nbigj)
 \to
 f_{T\ast}\iota_X^{-1}(\nbigj)
 \to
 f_{T\ast}\nbigj_1.
\]
Hence, there exists a natural morphism
\begin{multline}
 \iota_Y^{-1}
 \Bigl(
 Rf_{S\ast}\bigl(
 \nbigd_{(Y_S\larr X_S)/S}
 \otimes^L_{\nbigd_{X_S/S}}
 M
 \bigr)
 \Bigr)
 \lrarr
 Rf_{T\ast}\iota_X^{-1}
 \bigl(
 \nbigd_{(Y_S\larr X_S)/S}
 \otimes^L_{\nbigd_{X_S/S}}
  M
  \bigr)
\\
 \simeq
   Rf_{T\ast}
 \bigl(
 \iota_X^{-1}\nbigd_{(Y_S\larr X_S)/S}
 \otimes^L_{\iota_X^{-1}\nbigd_{X_S/S}}
 \iota_X^{-1} M
  \bigr).
\end{multline}
Let $g:\iota_Y\circ f_T=f_S\circ\iota_X$.
There exists the natural isomorphism
\[
 f_T^{-1}(\nbigd_{Y_T/T})
 \otimes_{g^{-1}\nbigd_{Y_S/S}}
 \iota_X^{-1}(\nbigd_{(Y_S\larr X_S)/S})
 \simeq
 f_T^{-1}(\nbigd_{Y_T/T})
 \otimes^L_{g^{-1}\nbigd_{Y_S/S}}
 \iota_X^{-1}(\nbigd_{(Y_S\larr X_S)/S})
 \simeq
 \nbigd_{(Y_T\larr X_T)/T}.
\]
Hence, there exists the following natural morphism
in $D^b_c(\nbigd_{Y_T/T})$:
\begin{equation}
\label{eq;26.2.3.2}
 L\iota_Y^{\ast}\bigl(
 f_{S\dagger}(M)
 \bigr)
 \lrarr
 f_{T\dagger}
 L\iota_X^{\ast}(M)
\end{equation}
Let us check that (\ref{eq;26.2.3.2}) is a quasi-isomorphism.
It is enough to check the case
where 
$X=Y\times \proj^m$ for a positive integer $m$,
$Y$ is affine,
and $f$ is the projection.
It is enough to prove that
$\iota_{Y\ast}
 L\iota_Y^{\ast}\bigl(
f_{S\dagger}(M)
 \bigr)
 \to
\iota_{Y\ast}f_{T\dagger}
 L\iota_X^{\ast}(M)$
is a quasi-isomorphism.
Let $\nbigp_Y\to \iota_{Y\ast}\nbigo_Y$
be a finite resolution by locally free $\nbigo_X$-modules
of finite ranks.
We have
\[
 \iota_{Y\ast}
 L\iota_Y^{\ast}\bigl(
f_{S\dagger}(M)
 \bigr)
 \simeq
 f_{S\dagger}(M\otimes f_S^{\ast}\nbigp_Y)
 \simeq
 f_{S\dagger}\bigl(
 \iota_{X\ast}L\iota_X^{\ast}(M)
 \bigr)
 \simeq
 \iota_{Y\ast}
 f_{T\dagger}(L\iota_X^{\ast}M).
\]

By using a Noetherian induction,
we obtain the claims of Proposition \ref{prop;25.12.10.41}.
\hfill\qed

\begin{cor}
Suppose that $M$ is $S$-flat,
and that $Y=S$.
Then, there is a decomposition $S=\bigsqcup S_i$
by locally closed smooth subsets
such that  
$f^j_{S_i\dagger}(M_{S_i})$
are locally free $\nbigo_{S_i}$-modules.
Moreover,
for any closed point $s\in S_i$,
there exist the natural isomorphisms
$f^j_{s\dagger}(M_{S_i})_{|s}
 \simeq
 H^{\dim X+j}(X_s,M_{s}\otimes \Omega^{\bullet})$.
\hfill\qed
\end{cor}

\subsubsection{Duality}

Assume that $S$ is smooth quasi-projective over $\hyperk$.
Let $n$ denote the relative dimension of $X_S$ over $S$.
For any $\nbigd_{X_S/S}$-complex $M^{\bullet}$,
we set
\[
 \DD(M^{\bullet})=\nrhom_{\nbigd_{X_S/S}}
 (M^{\bullet},
 \nbigd_{X_S/S}\otimes\Omega_{X_S/S}^{-1})[n]
 \in
 D(\nbigd_{X_S/S}).
\]
It induces
a contravariant functor
$\DD:D^b_c(\nbigd_{X_S/S})\to D^b_c(\nbigd_{X_S/S})$.

\begin{prop}
\label{prop;25.12.13.10}
Let $M$ be a holonomic $S$-flat 
$\nbigd_{X_S/S}$-module.
Then, 
there exists a decomposition
$S=\bigsqcup_{j\in \Gamma}S_j$
by locally closed smooth subsets
such that
$\DD(M_{S_j})$ are 
holonomic $\nbigd_{X_{S_j}/S_j}$-modules
flat over $S_j$. 
For any locally closed smooth subset $T\subset S_j$,
there exists the natural isomorphism
$\DD(M_{S_j})_{T}\simeq\DD(M_T)$.
\end{prop}
\pf
Let $S'\subset S$ be any smooth irreducible locally closed subset.
Let $\nbigp\to M_{S'}$
be a finite $\nbigd_{X_{S'}/S'}$-projective resolution.
Then,
$\DD(M_{S'})$ is represented by
\[
\nhom_{\nbigd_{X_{S'}/S'}}
\bigl(
\nbigp,\nbigd_{X_{S'/S'}}\otimes\Omega_{X_{S'}/S'}^{-1}
\bigr)[n].
\]
Let $\eta:\Spec(\hyperk_1)\to S'$ be a geometric point
obtained as an algebraic closure of the generic point of $S'$.
The induced morphism
$X_{\eta}\to X_{S'}$ is also denoted by $\eta$.
We obtain the induced locally projective $\nbigd_{X_{\eta}}$-resolution
$\eta^{\ast}\nbigp\to M_{\eta}$.
There exists the following natural isomorphism:
\begin{equation}
\label{eq;26.2.3.1}
\eta^{\ast}
 \nhom_{\nbigd_{X_{S'}}}\bigl(
 \nbigp,
 \nbigd_{X_{S'}/S'}
 \otimes
 \Omega_{X_{S'}/S'}^{-1}
 \bigr)[n]
\simeq
 \nhom_{\nbigd_{X_{\eta}}}\bigl(
 \eta^{\ast}\nbigp,
  \nbigd_{X_{\eta}}
 \otimes
 \Omega_{X_{\eta}}^{-1}
 \bigr)[n].
\end{equation}
We obtain $(\DD(M_{S'}))_{\eta}=\DD(M_{\eta})$,
which is a holonomic $\nbigd_{X_{\eta}}$-module.
Hence, there exists
a Zariski open subset $S''\subset S'$
such that
$\DD(M_{S''})$ is a holonomic
$S''$-flat $\nbigd_{X\times S''/S''}$-module.

Let $T$ be a locally closed subset of $S$.
Let $\iota:X_T\to X_{S}$ denote the inclusion.
We obtain the induced locally projective
$\nbigd_{X_T/T}$-resolution
$\iota^{\ast}\nbigp\to M_T$.
There exists the isomorphism
\[
 \iota^{\ast}
 \nhom_{\nbigd_{X_S/S}}
 \bigl(
 \nbigp,\nbigd_{X_S/S}
 \otimes\Omega_{X_S/S}^{-1}
 \bigr)[n]
 \simeq
  \nhom_{\nbigd_{X_T/T}}
 \bigl(
 \iota^{\ast}\nbigp,\nbigd_{X_T/T}
 \otimes\Omega_{X_T/T}^{-1}
 \bigr)[n].
\]
We obtain
$L\iota^{\ast}\DD(M)
\simeq
\DD(M_T)$.

Then, we obtain the first claim of the proposition
by a Noetherian induction.
\hfill\qed

\subsubsection{Kashiwara's equivalence}

Let $\iota:X_S\to Y_S$ be a closed embedding
of smooth quasi-projective $\hyperk$-varieties over $S$.
Let $\Mod_{\sharp}^{X_S}(\nbigd_{Y_S/S})$ $(\sharp=c,qc)$
denote the full subcategory of
$\Mod_{\sharp}(\nbigd_{Y_S/S})$
consisting of objects whose supports are contained in $\iota(X_S)$.
Let $D^{b,X_S}_{\sharp}(\nbigd_{Y_S/S})$ $(\sharp=c,qc)$
denote the full subcategory of
$D^b_{\sharp}(\nbigd_{Y_S/S})$
consisting of objects whose cohomological supports
are contained in $\iota(X_S)$.
The following proposition
is similar to
\cite[Theorem 1.6.1, Corollary 1.6.2]{Hotta-Takeuchi-Tanisaki}.
\begin{prop}
\label{prop;26.1.30.1}
The direct image functors induce the equivalences
\[
 i_{\dagger}:
 \Mod_{\sharp}(\nbigd_{X_S/S})
 \simeq
 \Mod_{\sharp}^{X_S}(\nbigd_{Y_S/S}),
 \quad
 i_{\dagger}:
 D^b_{\sharp}(\nbigd_{X_S/S})
 \simeq
 D_{\sharp}^{b,X_S}(\nbigd_{Y_S/S})
 \quad\quad(\sharp=c,qc).
\]
It preserves the $S$-flat condition and the holonomicity condition.
The quasi-inverse are induced by
\[
Li^{\ast}[\dim X-\dim Y]:
D_{\sharp}^{b,X_S}(\nbigd_{Y_S/S})\to
D^b_{\sharp}(\nbigd_{X_S/S}).
\]
\hfill\qed
\end{prop}

\subsubsection{Pull back of quasi-coherent $\nbigd_{Y/S}$-modules by
closed immersions}

Let $\iota:X_S\to Y_S$ be a closed embedding
of smooth quasi-projective $\hyperk$-varieties over $S$.
We regard $X_S$ as a closed subvariety of $Y_S$.
Let $H_j$ $(j=1,\ldots,m)$ be closed hypersurfaces of $Y_S$
such that $X_S=\bigcap_{j=1}^m H_j$.

We set $U=\bigoplus_{j=1}^m\hyperk e_j$.
For an ordered subset $J=(j_1,\ldots,j_k)\subset \{1,\ldots,m\}$,
we set
$e_J=e_{j_1}\wedge\cdots\wedge e_{j_k}\in \bigwedge^kU$.
We also set
$H(J)=\bigcup_{j\in J}H_j$.
We formally set $e_{\emptyset}=1$
and $H(\emptyset)=\emptyset$.

Let $M$ be a quasi-coherent $\nbigd_{X_S/S}$-module.
We set 
\[
 \nbigc^k(M,\vecH)
 =\bigoplus_{|J|=k}
 M(\ast H(J))\otimes e_J.
\]
For any $J\subset \{1,\ldots,m\}$
and $i\in \{1,\ldots,m\}\setminus J$,
there exists the morphism
$M(\ast H(J))\otimes e_J
\to M(\ast H(J\cup\{i\}))e_i\wedge e_{J}$
induced by
$M(\ast H(J))\to M(\ast H(J\cup\{i\}))$
and the multiplication of $e_i$.
The morphisms induce 
$d:\nbigc^k(M,\vecH)\to \nbigc^{k+1}(M,\vecH)$
such that $d\circ d=0$.
Because $(\nbigc^{\bullet}(M,\vecH),d)_{|Y_S\setminus X_S}$
is acyclic,
$(\nbigc^{\bullet}(M,\vecH),d)
\in D^{b,X_S}_{qc}(\nbigd_{Y_S/S})$.

\begin{lem}
There exists the natural morphism 
$C^{\bullet}(M,\vecH)\simeq 
 \iota_{\dagger}
 L\iota^{\ast}(M)[\dim X-\dim Y]$
in $D^{b,X_S}_{qc}(\nbigd_{Y_S/S})$.
\end{lem}
\pf
There exists the natural morphism of complexes
$C^{\bullet}(M,\vecH)\to M$.
It induces
\[
 L\iota^{\ast}\nbigc^{\bullet}(M,\vecH)[\dim X-\dim Y]
 \simeq
 L\iota^{\ast}(M)[\dim X-\dim Y].
\]
By Proposition \ref{prop;26.1.30.1},
there exists a natural isomorphism
\[
 \nbigc^{\bullet}(M,\vecH)[\dim X-\dim Y]
 \simeq
 \iota_{\dagger}\bigl(
 L\iota^{\ast}\nbigc^{\bullet}(M,\vecH)[\dim X-\dim Y]
 \bigr)
 \simeq
 \iota_{\dagger}\bigl(
 L\iota^{\ast}(M)[\dim X-\dim Y]
 \bigr).
\]
\hfill\qed

\subsection{The restriction of extensions}

\subsubsection{Adjunction}

Let $S$ be a smooth $\hyperk$-variety.
Let $X_S$ be a $\hyperk$-variety which is smooth projective over $S$.
Let $M$ be any coherent $\nbigd_{X_S/S}$-module.
Let $T\subset S$ be a closed complex submanifold.
Let $\iota:T\to S$ denote the inclusion.
Let $N_T$ be any $\nbigd_{X_T/T}$-module.

\begin{lem}
\label{lem;25.12.13.20}
There exists the following natural isomorphism
in $D^b(p^{-1}\nbigo_S)$:
\begin{equation}
\label{eq;25.12.13.1}
 \iota_{\ast}\nrhom_{\nbigd_{X_T/T}}
 (L\iota^{\ast}M,N_T)
\simeq 
 \nrhom_{\nbigd_{X_S/S}}(M,\iota_{\ast}N_T)
\end{equation}
\end{lem}
\pf
Let $\nbigj_0$ be a $\nbigd_{X_T/T}$-injective resolution of
$N_T$.
Let $\nbigj_1$ be a $\nbigd_{X_S/S}$-injective resolution of
$\iota_{\ast}N_T$.
There exists a quasi-isomorphism of
$\nbigd_{X_S/S}$-complexes
$\iota_{\ast}\nbigj_0\to \nbigj_1$
extending the identity of $\iota_{\ast}N_T$.
Let $P\to M$ be a locally $\nbigd_{X_S/S}$-projective resolution
of $M$.
Then, $\nrhom_{\nbigd_{X_S/S}}(M,\iota_{\ast}N_T)$
is represented by
$\nhom_{\nbigd_{X_S/S}}(P,\nbigj_1)$,
and 
$\iota_{\ast}\nrhom_{\nbigd_{X_T/T}}
 (L\iota^{\ast}M,N_T)$
is represented by
$\iota_{\ast}\nhom_{\nbigd_{X_T/T}}
 \bigl(
 \iota^{\ast}P,
 \nbigj_0
 \bigr)$.
There exists the following natural isomorphism and quasi-isomorphism
\begin{equation}
\label{eq;25.12.13.22}
\iota_{\ast}\nhom_{\nbigd_{X_T/T}}
 \bigl(
 \iota^{\ast}P,
 \nbigj_0
 \bigr)
 \simeq
 \nhom_{\nbigd_{X_S/S}}(P,\iota_{\ast}\nbigj_0)
 \lrarr
  \nhom_{\nbigd_{X_S/S}}(P,\nbigj_1).
\end{equation}
We obtain (\ref{eq;25.12.13.1}).
\hfill\qed

\vspace{.1in}
Let $p:X_S\to S$ denote the projection.

\begin{lem}
\label{lem;25.12.14.4}
There exists the following morphisms 
in $D(p^{-1}\nbigo_S)$
such that $b\circ a$ equals the morphism
{\rm(\ref{eq;25.12.13.1})}:
\begin{equation}
\label{eq;25.12.13.23}
\begin{CD}
 \iota_{\ast}
 \nrhom_{\nbigd_{X_T/T}}\bigl(
 L\iota^{\ast}M,N_T
 \bigr)
  @>{a}>>
 \nrhom_{\nbigd_{X_S/S}}\bigl(
 \iota_{\ast}L\iota^{\ast}M,\iota_{\ast}N_T
 \bigr)
 @>{b}>>
 \nrhom_{\nbigd_{X_S/S}}\bigl(
 M,\iota_{\ast}N_T
 \bigr)
\end{CD}
\end{equation}
\end{lem}
\pf
We use the notation in the proof of
Lemma \ref{lem;25.12.13.20}.
We have the following commutative diagram:
\[
 \begin{CD}
  \nhom_{\nbigd_{X_S/S}}
  \bigl(\iota_{\ast}\iota^{\ast}P,\iota_{\ast}\nbigj_0\bigr)
  @>{c_1}>>
  \nhom_{\nbigd_{X_S/S}}
  \bigl(
  \iota_{\ast}\iota^{\ast}P,\nbigj_1
  \bigr)
  \\
  @V{c_2}V{\simeq}V @V{c_3}VV \\
  \nhom_{\nbigd_{X_S/S}}
  \bigl(P,\iota_{\ast}\nbigj_0\bigr)
  @>{c_4}>{\rm qis}>
  \nhom_{\nbigd_{X_S/S}}
  \bigl(
  P,\nbigj_1\bigr).
 \end{CD}
\]
Here, $c_2$ is an isomorphism,
and $c_4$ is a quasi-isomorphism.
We have the following natural isomorphism
\begin{equation}
\label{eq;25.12.13.21}
 \iota_{\ast}
 \nhom_{\nbigd_{X_T/T}}
 \bigl(\iota^{\ast}P,\nbigj_0
 \bigr)
 \simeq
 \nhom_{\nbigd_{X_S/S}}
 \bigl(\iota_{\ast}\iota^{\ast}P,\iota_{\ast}\nbigj_0\bigr).
\end{equation}
The composition of (\ref{eq;25.12.13.21}), $c_2$ and $c_4$
is the quasi-isomorphism (\ref{eq;25.12.13.22}).
Note that
$\nrhom_{\nbigd_{X_S/S}}\bigl(
 \iota_{\ast}L\iota^{\ast}M,\iota_{\ast}N_T
 \bigr)$
is represented by
$\nhom_{\nbigd_{X_S/S}}
\bigl(
\iota_{\ast}\iota^{\ast}P,\nbigj_1
\bigr)$.
Hence,
we obtain the morphisms in (\ref{eq;25.12.13.23}).
\hfill\qed

\subsubsection{Tensor product with sheaves on the parameter space}

Let $M$ be any coherent $\nbigd_{X_S/S}$-module.
Let $N$ be any $\nbigd_{X_S/S}$-module.
\begin{lem}
\label{lem;25.12.14.1}
For any $\nbigo_S$-module $F$,
there exists the natural isomorphism
\begin{equation}
\label{eq;25.12.13.2}
 \nrhom_{\nbigd_{X_S/S}}(M,N)\otimes^L_{p^{-1}\nbigo_S}p^{-1}F
 \simeq
 \nrhom_{\nbigd_{X_S/S}}(M,N\otimes^L_{p^{-1}\nbigo_S}p^{-1}F ).
\end{equation}
\end{lem}
\pf
Let $P\to M$ be a finite locally projective $\nbigd_{X_S/S}$-resolution.
Let $N\to \nbigj_2$ be a $\nbigd_{X_S/S}$-injective resolution.
Let $L\to F$ be any $\nbigo_S$-flat resolution.
Then,
$\nrhom_{\nbigd_{X_S/S}}(M,N)\otimes^L_{p^{-1}\nbigo_S}p^{-1}F$
is represented by
\begin{equation}
\label{eq;25.12.13.30}
 \nhom_{\nbigd_{X_S/S}}(P,\nbigj_2)\otimes_{p^{-1}\nbigo_S}
 p^{-1}L
 \simeq
 \nhom_{\nbigd_{X_S/S}}(P,\nbigj_2\otimes_{p^{-1}\nbigo_S}p^{-1}L).
\end{equation}
Let $N\otimes_{p^{-1}\nbigo_S}p^{-1}L\to\nbigj_3$ be 
a $\nbigd_{X_S/S}$-injective resolution.
Then,
$\nrhom_{\nbigd_{X_S/S}}(M,N\otimes^L_{p^{-1}\nbigo_S}F)$
is represented by
$\nhom_{\nbigd_{X_S/S}}(P,\nbigj_3)$.
There exists a quasi-isomorphism of
$\nbigd_{X_S/S}$-complexes
$\nbigj_2\otimes_{p^{-1}\nbigo_S}
p^{-1}L
\lrarr
\nbigj_3$.
It induces a quasi-isomorphism
\begin{equation}
\label{eq;25.12.13.31}
 \nhom_{\nbigd_{X_S/S}}(P,
 \nbigj_2\otimes_{p^{-1}\nbigo_S}p^{-1}L)
  \lrarr
  \nhom_{\nbigd_{X_S/S}}(P,\nbigj_3).
\end{equation}
We obtain (\ref{eq;25.12.13.2})
from (\ref{eq;25.12.13.30}) and (\ref{eq;25.12.13.31}).
\hfill\qed

\subsubsection{Change of parameter spaces}

By Lemma \ref{lem;25.12.13.20}
and Lemma \ref{lem;25.12.14.1},
we obtain
\begin{multline}
 \iota_{\ast}
 \nrhom_{\nbigd_{X_T/T}}
 \bigl(
 L\iota^{\ast}M,L\iota^{\ast}N
 \bigr)
 \simeq
 \nrhom_{\nbigd_{X_S/S}}(M,\iota_{\ast}L\iota^{\ast}N)
 \simeq
 \nrhom_{\nbigd_{X_S/S}}
 \bigl(M,
 N\otimes^L_{p^{-1}\nbigo_{S}}p^{-1}\iota_{\ast}\nbigo_{T}
 \bigr)
 \\
 \simeq
  \nrhom_{\nbigd_{X_S/S}}
 \bigl(M,
 N
 \bigr)
 \otimes^L_{p^{-1}\nbigo_{S}}p^{-1}\iota_{\ast}\nbigo_{T}.
\end{multline}
As a result,
we obtain
\[
 \iota_{\ast}Rp_{\ast}
  \nrhom_{\nbigd_{X_T/T}}
 \bigl(
 L\iota^{\ast}M,L\iota^{\ast}N
 \bigr)
 \simeq
 Rp_{\ast}
 \nrhom_{\nbigd_{X_S/S}}
 \bigl(M,N \bigr)
 \otimes^L_{\nbigo_{S}}\iota_{\ast}\nbigo_{T}.
\]

In particular, we obtain the following.
\begin{cor}
\label{cor;25.12.13.13}
We assume that
$M$ and $N$ are flat over $S$.
Then, there exist the natural isomorphisms
\[
 \iota_{\ast}R^kp_{\ast}
  \nrhom_{\nbigd_{X_T/T}}
 \bigl(
 \iota^{\ast}M,\iota^{\ast}N
 \bigr)
 \simeq
 R^kp_{\ast}
 \nrhom_{\nbigd_{X_S/S}}
 \bigl(M,N
 \bigr)
 \otimes^L_{\nbigo_{S}}\iota_{\ast}\nbigo_{T}.
\]
\hfill\qed
\end{cor}

\subsubsection{Commutative diagrams}

We impose the following.
\begin{condition}
\label{condition;25.12.14.2}
We assume that
(i) $S$ is affine,
(ii) $M$ and $N$ are flat over $S$.
\hfill\qed
\end{condition}

We set $N_T=\iota^{\ast}N$.
Let $\iota_{\ast}N_T\to\nbigj_1$
be a $\nbigd_{X_S/S}$-injective resolution.
Let $N\to\nbigj_2$
be a $\nbigd_{X_S/S}$-injective resolution.
There exists a morphism of
$\nbigd_{X_S/S}$-complexes
$\nbigj_2\to\nbigj_1$
extending $N\to \iota_{\ast}N_T$.

Let $L(S)\to \nbigo_S$ and $L(T)\to\iota_{\ast}\nbigo_T$
be a finite resolution by locally free $\nbigo_S$-modules
of finite rank
with the following commutative diagram:
\[
 \begin{CD}
  L(S) @>>> L(T) \\
  @VVV @VVV \\
  \nbigo_S@>>>\iota_{\ast}\nbigo_T.
 \end{CD}
\]
Note that
$N\otimes_{p^{-1}\nbigo_S} p^{-1}L(T)\to 
\nbigj_2\otimes_{p^{-1}\nbigo_S} p^{-1}L(T)$
is a $\nbigd_{X_S/S}$-injective resolution.
There exists the following commutative diagram of
quasi-isomorphisms:
\[
 \begin{CD}
  N\otimes_{p^{-1}\nbigo_S} p^{-1}L(T)
  @>>>
  \nbigj_2\otimes_{p^{-1}\nbigo_S} p^{-1}L(T)
  \\
  @VVV @VVV\\
  \iota_{\ast}N_T @>>>
  \nbigj_1.
 \end{CD}
\]
\begin{lem}
The following diagram is commutative
in $D(\nbigd_{X_S/S})$:
\begin{equation}
\label{eq;25.12.13.40}
\begin{CD}
 \nbigj_2\otimes_{p^{-1}\nbigo_S}p^{-1}L(S)
 @>>>
 \nbigj_2\otimes_{p^{-1}\nbigo_S}p^{-1}L(T)
 \\
 @VVV @VVV \\
 \nbigj_2
 @>>>\nbigj_1.
\end{CD}
\end{equation}
\end{lem}
\pf
Because the following diagram is commutative,
\[
 \begin{CD}
 N\otimes_{p^{-1}\nbigo_S}p^{-1}L(S)
 @>>>
 N\otimes_{p^{-1}\nbigo_S}p^{-1}L(T)
 \\
 @VVV @VVV \\
 N
  @>>>
  \iota_{\ast}N_T,
 \end{CD}
\]
we obtain the commutativity of (\ref{eq;25.12.13.40}).
\hfill\qed

\vspace{.1in}
Let $\alpha:M\to\nbigj_2[k]$ be a morphism of
$\nbigd_{X_S/S}$-complexes.
We obtain
$\alpha_1:M\to \nbigj_1[k]$
as the composition of
\[
 M\stackrel{\alpha}{\lrarr} \nbigj_2[k]
  \lrarr
 \nbigj_1[k].
\]
We also obtain
$\alpha_2:M\otimes_{p^{-1}\nbigo_S}p^{-1}L(S)\to
\nbigj_2[k]\otimes_{p^{-1}\nbigo_S}p^{-1}\nbigo_T$
as the composition of
the following morphisms:
\[
 M\otimes p^{-1}L(S)
 \stackrel{\alpha\otimes\id}{\lrarr}
 \nbigj_2[k]\otimes p^{-1}L(S)
 \lrarr
 \nbigj_2[k]\otimes p^{-1}L(T).
\]

\begin{lem}
We have $\alpha_1=\alpha_2$
in $\Ext^k_{\nbigd_{X_S/S}}(M,\iota_{\ast}N_T)$.
\end{lem}
\pf
We have the following commutative diagram
in $D^b(\nbigd_{X_S/S})$:
\[
 \begin{CD}
  M\otimes_{p^{-1}\nbigo_S} p^{-1}L(S)
  @>{\alpha\otimes\id}>>
  \nbigj_2[k]\otimes_{p^{-1}\nbigo_S} p^{-1}L(S)
  @>>>
  \nbigj_2[k]\otimes_{p^{-1}\nbigo_S} p^{-1}L(T)
  \\
  @VV{qis}V @VV{qis}V @VV{qis}V\\
  M @>{\alpha}>>
  \nbigj_2[k]
  @>>>
  \nbigj_1[k].
 \end{CD}
\]
Hence, we obtain the claim of the lemma.
\hfill\qed

\vspace{.1in}

There exists the following commutative diagram:
\[
 \begin{CD}
  \nhom_{\nbigd_{X_S/S}}
  (M,\nbigj_2[k])
  @>>>
  \nhom_{\nbigd_{X_S/S}}
  (M\otimes p^{-1}L(S),\nbigj_2[k]\otimes p^{-1}L(S))
  \\
  @VVV @VVV
  \\
  \nhom_{\nbigd_{X_S/S}}
  (M\otimes p^{-1}L(T),\nbigj_2[k]\otimes p^{-1}L(T))
  @>>>
  \nhom_{\nbigd_{X_S/S}}
  (M\otimes p^{-1}L(S),\nbigj_2[k]\otimes p^{-1}L(T)).
 \end{CD}
\]
From $\alpha$,
we obtain
\[
\alpha_3
=\alpha\otimes\id
:M\otimes_{p^{-1}\nbigo_S} p^{-1}L(T)
\to
\nbigj_2[k]\otimes_{p^{-1}\nbigo_S} p^{-1}L(T).
\]
The following lemma is clear.
\begin{lem}
\label{lem;26.1.31.1}
$\alpha_2$ equals
the image of $\alpha_3$
via 
$\Ext^k_{\nbigd_{X_S/S}}(\iota_{\ast}M_T,\iota_{\ast}N_T)
\to
 \Ext^k_{\nbigd_{X_S/S}}(M,\iota_{\ast}N_T)$.
\hfill\qed
\end{lem}

\begin{cor}
\label{cor;25.12.14.3}
The morphism
$\Ext^k_{\nbigd_{X_S/S}}(M,N)
\to \Ext^k_{\nbigd_{X_S/S}}(M,\iota_{\ast}N_T)$
factors through as follows:
\[
\Ext^k_{\nbigd_{X_S/S}}(M,N)
\to
\Ext^k_{\nbigd_{X_S/S}}(\iota_{\ast}M_T,\iota_{\ast}N_T)
\to \Ext^k_{\nbigd_{X_S/S}}(M,\iota_{\ast}N_T).
\]
\hfill\qed
\end{cor}

\subsubsection{The restriction of extensions}

We assume Condition \ref{condition;25.12.14.2}.
Let $T\subset S$ be any smooth closed subset of $S$.
The inclusion maps $T\to S$ and $X_T\to X_S$
are denoted by $\iota$.
For any $\nbigd_{X_SS/S}$-module $G$ flat over $S$,
we set $G_T=\iota^{\ast}G$.
By Corollary \ref{cor;25.12.13.13},
we obtain
\[
\Ext^1_{\nbigd_{X_S/S}}(M,N)
\to
\Ext^1_{\nbigd_{X_S/S}}(M,\iota_{\ast}N_T)
\simeq
\Ext^1_{\nbigd_{X_T/T}}(M_{T},N_{T}).
\]
For any $\alpha\in\Ext^1_{\nbigd_{X_S/S}}(M,N)$,
the image is denoted by
$\alpha_T\in\Ext^1_{\nbigd_{X_T/T}}(M_{T},N_{T})$.

\begin{prop}
\label{prop;25.12.14.5}
If $\alpha$ corresponds to
an extension
$0\to N\to \Mtilde\to M\to 0$,
then 
$\alpha_T$ corresponds to 
the induced extension 
$0\to N_{T}\to \Mtilde_T\to M_{T}\to 0$.
Note that $\Mtilde$ is flat over $S$.
\end{prop}
\pf
We have the following factorization
as in Corollary \ref{cor;25.12.14.3}:
\[
\Ext^1_{\nbigd_{X_S/S}}(M,N)
\to
\Ext^1_{\nbigd_{X_S/S}}(\iota_{\ast}M_T,\iota_{\ast}N_T)
\to
\Ext^1_{\nbigd_{X_T/T}}(M_{T},N_{T}).
\]
We have $\alpha'_T\in
\Ext^1_{\nbigd_{X_S/S}}(\iota_{\ast}M_T,\iota_{\ast}N_T)$
as the image of $\alpha$.
By Lemma \ref{lem;26.1.31.1},
it is enough to prove that
$\alpha'_T$ corresponds to
the extension induced by
$0\to \iota_{\ast}M_T\to \iota_{\ast}(\Mtilde_T)\to
\iota_{\ast}(N_T)\to 0$.

Let $\alphatilde:M\to\nbigj_2[1]$
be a morphism representing $\alpha$.
We also have the natural morphism of $\nbigd_{X_S/S}$-complexes 
$M\to \Cone(\Mtilde\to M)$,
which also represents $\alpha$.
There exists the following commutative diagram
in $D^b_c(\nbigd_{X_S/S})$:
\[
 \begin{CD}
  M@>{\simeq}>> M\\
  @V{\alphatilde}VV @VVV \\
  \nbigj_2[1]
  @>{\simeq}>> \Cone(\Mtilde\to M).
 \end{CD}
\]
Let
$\alphatilde'_T:
M\otimes_{\pi^{-1}\nbigo_S}\pi^{-1}L(T)
\to
\nbigj_2[1]\otimes_{\pi^{-1}\nbigo_S} \pi^{-1}L(T)$
be the morphism induced by $\alphatilde$,
which represents $\alpha'_T$.
We obtain the following commutative diagram
in $D^b_c(\nbigd_{X_S/S})$:
\[
 \begin{CD}
  M\otimes \pi^{-1}L(T)
  @>{\simeq}>>
  M\otimes\pi^{-1}L(T)
  @>{\simeq}>>
  \iota_{\ast}(M_s)
  \\
  @V{\alphatilde'_T}VV
 @VVV @VVV \\
  \nbigj_2[1] \otimes \pi^{-1}L(T) 
  @>{\simeq}>>
  \Cone(\Mtilde\to M)\otimes \pi^{-1}L(T)
  @>{\simeq}>>
  \iota_{\ast}\Cone(\Mtilde_{s}\to M_s).
 \end{CD}
\]
Note that
$\iota_{\ast}(M_s)\to
\iota_{\ast}\Cone(\Mtilde_{s}\to M_s)$
corresponds to
the extension
$0\to \iota_{\ast}M_T\to \iota_{\ast}(\Mtilde_T)\to
\iota_{\ast}(N_T)\to 0$.
Then, we obtain the claim of Proposition
\ref{prop;25.12.14.5}.
\hfill\qed

\subsection{$\nbigd$-modules associated with good lattices
in the absolute case}

\subsubsection{The $\nbigd$-modules associated with lattices}

Let $X$ be a smooth $\hyperk$-variety
with a hypersurface $H$.
Let $(\nbige,\nabla)$ be a meromorphic flat connection
on $(X,H)$.

Let $E$ be a lattice of a meromorphic integrable connection
$(\nbige,\nabla)$ on $(X,H)$.
We obtain the $V\nbigd_{X}$-submodule
$V\nbigd_X\cdot E\subset \nbige$.
We set
\[
 \nbigd_X(E):=
 \nbigd_{X}\otimes_{V\nbigd_X}
 (V\nbigd_X\cdot E).
\]
There exists a natural isomorphism
$\nbigd_X(E)(\ast H)=\nbige$.

For a morphism of meromorphic integrable connections
$(\nbige',\nabla)\to(\nbige,\nabla)$.
Let $E'\to E$ be a morphism of lattices
which induces the morphism $\nbige'\to\nbige$.
Then, 
we obtain $\nbigd_X(E')\to \nbigd_X(E)$.

\subsubsection{The $\nbigd$-modules associated with good lattices}

Suppose that $H$ is simple normal crossing,
and that 
$(\nbige,\nabla)$ is a good meromorphic flat bundle
on $(X,H)$.
We set
$\nbige^{\lor}=\nhom_{\nbigo_X(\ast H)}(\nbige,\nbigo_X(\ast H))$
which is equipped with the induced integrable connection $\nabla$.
We obtain the $\nbigd$-module
$\nbige(!H)=\DD(\nbige^{\lor})$.

Let $H=\bigcup_{i\in\Lambda} H_j$ be the irreducible decomposition.
For $I\subset \Lambda$,
we set $H(I)=\bigcup_{j\in I}H_j$.
We set
\[
\nbige(!H(I)\ast H(I^c))
=\nbige(!H)\otimes\nbigo_X(\ast H(I^c)),
\]
where $I^c=\Lambda\setminus I$.
For $I\supset J$,
there exists the natural morphism
$\nbige(!H(I)\ast H(I^c))
\to
\nbige(!H(J)\ast H(J^c))$.

\begin{prop}
\label{prop;26.2.5.10}
Let $E$ be a good lattice of $(\nbige,\nabla)$.
Assume that for each $i\in\Lambda$
one of the following holds. 
\begin{description}
 \item[(a)] $\Sp(E,\nabla,H_i)^{\reg}\cap \seisuu_{> 0}=\emptyset$.
 \item[(b)] $\Sp(E,\nabla,H_i)^{\reg}\cap \seisuu_{\leq 0}=\emptyset$.
\end{description}
Let $\Lambda(a)$ (resp. $\Lambda(b)$)
denote the set of $i\in\Lambda$
satisfying (a) (resp. (b)).
We set $H(a)=\bigcup_{i\in\Lambda(a)}H_i$
and $H(b)=\bigcup_{i\in\Lambda(b)}H_i$.
Then, there exists the natural isomorphism
$\nbigd_X(E)\simeq
 \nbige(!H(b)\ast H(a))$.
\hfill\qed
\end{prop}

Let $\nbigp_{\ast}^{\DM}(\nbige)$
be a good $\cnum$-Deligne-Malgrange filtered bundle
as in \S\ref{subsection;26.2.4.1}.
\begin{prop}
For $\veca\in \rnum_{<0}^{I}\times\rnum_{\geq 0}^{I^c}$,
there exists
a natural isomorphism
\[
 \nbigd_X\bigl(\nbigp^{\DM}_{\veca}(\nbige)\bigr)
 \simeq
 \nbige\bigl(!H(I)\ast H(I^c)\bigr).
\]
Let $I\supset J$.
Let $\veca(I)\in\rnum_{<0}^I\times\rnum_{\geq 0}^{I^c}$
and $\veca(J)\in\rnum_{<0}^{J}\times\rnum_{\geq 0}^{J^c}$
such that $\veca(I)_i=\veca(J)_i$ for $J\sqcup I^c$.
Then, there exists the following commutative diagram:
\[
  \begin{CD}
  \nbigd_X\bigl(\nbigp^{\DM}_{\veca(I)}(\nbige)\bigr)
  @>{\simeq}>>
  \nbige\bigl(!H(I)\ast H(I^c)\bigr)\\
  @VVV @VVV \\
  \nbigd_X\bigl(\nbigp^{\DM}_{\veca(J)}(\nbige)\bigr)
  @>{\simeq}>>
  \nbige(!H(J)\ast H(J^c)).
 \end{CD}
\]
\hfill\qed
\end{prop}

\begin{prop}
For any good lattice $E_1$ of $\nbige$
such that 
$\nbigp^{\DM}_{(0,\ldots,0)}(\nbige)\subset E_1$,
we have 
$\nbigd_X(E_1)\simeq \nbige$.
For any good lattice $E_2$ of $\nbige$
such that 
$\nbigp_{<(0,\ldots,0)}(\nbige)\supset E_2$,
we have 
$\nbigd_X(E_2)\simeq \nbige(!H)$.
\hfill\qed
\end{prop}

\subsection{Families of $\nbigd$-modules associated with good lattices}

\subsubsection{Families of $\nbigd$-modules associated with lattices}

Let $S$ be a smooth $\hyperk$-variety.
Let $X$ be a $\hyperk$-variety 
with a smooth projective morphism $X\to S$
such that each geometric fiber is irreducible.
Let $H\subset X$ be a hypersurface of $X$
which is normal crossing relative to $S$
with trivial monodromy.
Let $H=\bigcup_{i\in\Lambda}H_i$
denote the irreducible decomposition.

Let $(\nbige,\nabla)$ be
a meromorphic integrable connection on $(X,H)$ relative to $S$.
For any lattice $E$ of $\nbige$,
we obtain the $V\nbigd_{X/S}$-submodule
$V\nbigd_{X/S}\cdot E\subset\nbige$.
We obtain $\nbigd_{X/S}$-module
\[
\nbigd_{X/S}(E)=
\nbigd_{X/S}\otimes_{V\nbigd_{X/S}}
\bigl(
V\nbigd_{X/S}\cdot E
\bigr).
\]

\subsubsection{Families of $\nbigd$-modules associated with good lattices}

Suppose that $E$ is good.
Let $u$ be any closed point of $S$.
Let $x\in H$ be a closed point over $u$.
We set
$\Rhat_u=\nbigohat_{U,u}$
and
$\Rhat_x=\nbigohat_{X,x}$.

\begin{lem}
\label{lem;26.2.5.1}
There exist the natural isomorphisms
$(V\nbigd_{X/S}\cdot E)_{|\xhat}
 \simeq
 V\nbigd_{\Rhat_x/\Rhat_u}\cdot E_{|\xhat}$
and 
$\nbigd_X(E)_{|\xhat}
\simeq
\nbigd_{\Rhat_x/\Rhat_u}(E_{|\xhat})$.
\end{lem}
\pf
The image of the natural morphism
\begin{equation}
\label{eq;26.1.26.10}
 \Rhat_x\otimes_{R_x}
(V\nbigd_{X/S}\cdot E)_x
\lrarr \nbige_{|\xhat}
\end{equation}
equals 
$V\nbigd_{\Rhat_x/\Rhat_u}\cdot E_{|\xhat}$.
Because $\Rhat_x$ is flat over $R_x=\nbigo_{X,x}$,
the morphism (\ref{eq;26.1.26.10}) is injective.
Hence, the natural morphism
$\Rhat_x\otimes_{R_x}
(V\nbigd_{X/S}\cdot E)_x
\to
 V\nbigd_{\Rhat_x/\Rhat_u}\cdot E_{|\xhat}$
is an isomorphism.
We obtain the second isomorphism from the first. 
\hfill\qed

\begin{lem}
$V\nbigd_{X/S}\cdot E$
is flat over $S$.
For any closed point $u\in S$,
the natural morphism
$(V\nbigd_{X/S}\cdot E)_u
\to
V\nbigd_{X_u}\cdot E_u$
is an isomorphism.
\end{lem}
\pf
Let $x\in H$ be a closed point over $u$.
Let us prove that
$(V\nbigd_{X/S}\cdot E)_x$ is flat over
$R_u=\nbigo_{U,u}$.
Let $I$ be any ideal of $R_u$.
It is enough to prove that
$I\otimes (V\nbigd_{X/S}\cdot E)_x
\to (V\nbigd_{X/S}\cdot E)_x$
is injective
for the $S$-flatness.
We set $\Ihat=I\cdot \Rhat_u\simeq I\otimes\Rhat_u$.
Because 
$V\nbigd_{\Rhat_x/\Rhat_u}\cdot
E_{|\xhat}
\subset \nbige_{|\xhat}$
is flat over $\Rhat_u$
by Proposition \ref{prop;26.1.24.5},
$\Ihat\otimes_{\Rhat_u}
\bigl(
V\nbigd_{\Rhat_x/\Rhat_u}\cdot
E_{|\xhat}
\bigr)
\lrarr
\bigl(
V\nbigd_{\Rhat_x/\Rhat_u}\cdot
E_{|\xhat}
\bigr)$
is injective.

By Lemma \ref{lem;26.2.5.1},
there exists the isomorphism
$\Rhat_x\otimes_{R_x}
(V\nbigd_{X/S}\cdot E)_x
\simeq
V\nbigd_{\Rhat_x/\Rhat_u}\cdot
E_{|\xhat}$.
We also have
\[
\Rhat_x\otimes_{R_x}
\bigl(
I\otimes_{R_u}
(V\nbigd_{X/S}\cdot E)_x\bigr)
\simeq
\Rhat_x\otimes_{R_x}
\bigl(
IR_x\otimes_{R_x}
(V\nbigd_{X/S}\cdot E)_x\bigr)
\simeq
I\Rhat_x\otimes_{\Rhat_x}
(V\nbigd_{\Rhat_x/\Rhat_u}\cdot E_{|\xhat})
\simeq
\Ihat\otimes_{\Rhat_u}
 V\nbigd_{\Rhat_x/\Rhat_u}\cdot E_{|\xhat}.
\]
There exists the following commutative diagram:
\[
 \begin{CD}
  \Rhat_x\otimes_{R_x}
  \bigl(
  I\otimes_{R_u}
(V\nbigd_{X/S}\cdot E)_x\bigr)
  @>>>
  \Rhat_x\otimes_{R_x}
(V\nbigd_{X/S}\cdot E)_x \\
  @V{\simeq}VV @V{\simeq}VV \\
\Ihat\otimes_{\Rhat_u}
\bigl(
V\nbigd_{\Rhat_x/\Rhat_u}\cdot
E_{|\xhat}
\bigr)
 @>>>
V\nbigd_{\Rhat_x/\Rhat_u}\cdot
E_{|\xhat}.
 \end{CD}
\]
Because $\Rhat_x$ is faithfully flat over $R_x$,
we obtain that
$I\otimes_{R_u}
(V\nbigd_{X/S}\cdot E)_x
\lrarr
(V\nbigd_{X/S}\cdot E)_x$
is injective.
Thus, we obtain that
$V\nbigd_{X/S}\cdot E$ is $S$-flat.

Let $\Rhat_{1,x}=\nbigohat_{X_u,x}$.
We have the following commutative diagram:
\[
  \begin{CD}
   \Rhat_{1,x}
   \otimes
 (V\nbigd_{X/S}\cdot E)_u
 @>>>
   \Rhat_{1,x}
   \otimes
   V\nbigd_{X_u}E_u
   \\
   @V{\simeq}VV @V{\simeq}VV \\
   \hyperk\otimes_{\Rhat_u}
   \bigl(
   V\nbigd_{\Rhat_x/\Rhat_u}\cdot E_{|\xhat}
   \bigr)
   @>{\simeq}>>
   V\nbigd_{\Rhat_{1,x}/\hyperk}\cdot E_{u|\xhat}.
  \end{CD}
\]
We obtain that 
$(V\nbigd_{X/S}\cdot E)_u
\simeq
V\nbigd_{X_u}E_u$.
\hfill\qed

\begin{cor}
$\nbigd_{X/S}(E)$ is $S$-flat,
and there exists the natural morphism
$\nbigd_{X/S}(E)_u
\simeq
\nbigd_{X_u}(E_u)$.
\hfill\qed
\end{cor}

Let $\varphi:S'\to S$ be a morphism of smooth $\hyperk$-varieties.
We obtain
$(X_{S'},H_{S'})=(X,H)\times_SS'$
with the induced morphism
$\varphi:(X_{S'},H_{S'})\to (X,H)$.
We obtain
$(\nbige_{S'},\nabla)=
\varphi^{\ast}(\nbige,\nabla)$
on $(X_{S'},H_{S'})$.
For any lattice $E$ of $\nbige$,
we obtain the induced lattice $E_{S'}$.

\begin{prop}
Suppose that $E$ is a good lattice of $(\nbige,\nabla)$.
Then, there exists a natural isomorphism
\[
 \varphi^{\ast}\nbigd_{X/S}(E)
 \simeq
 \nbigd_{X_{S'}/S'}(E_{S'}).
\] 
\hfill\qed
\end{prop}

Let $\eta:\Spec(\hyperk_1)\to S$
be a geometric point.
\begin{prop}
Suppose that $E$ is a good lattice of $(\nbige,\nabla)$.
Then, there exists a natural isomorphism
\[
 \nbigd_{X/S}(E)_{\eta}
 \simeq
 \nbigd_{X_{\eta}}(E_{\eta}).
\] 
\hfill\qed
\end{prop}

\subsubsection{Exceptional subsets}

Suppose that $S$ is irreducible.
Let $\eta(0):\Spec(\hyperk_1)\to S$
be the geometric point obtained as
an algebraic closure of the generic point of $S$.
Let $\Lambda=I\sqcup I^c$ be a decomposition.
Let $E$ be a good lattice such that
$\nbigd_{X_{\eta(0)}}(E_{\eta(0)})
 \simeq
\nbige_{\eta(0)}(!H_{\eta(0)}(I)\ast H_{\eta(0)}(I^c))$.
By Proposition \ref{prop;26.2.5.10},
we have
$\Sp(E_{\eta(0)},\nabla,i)^{\reg}
\cap \seisuu_{>0}=\emptyset$ for $i\in I^c$,
and
$\Sp(E_{\eta(0)},\nabla,i)^{\reg}
\cap \seisuu_{\leq 0}=\emptyset$ for $i\in I$.
By Corollary \ref{cor;26.2.21.40},
for $m\in\seisuu$ and $i\in \Lambda$,
we obtain the Zariski closed subset
\[
 W(E,i,m)=\{u\,|\,m\in\Sp(E_u,\nabla)\}\subset S.
\]
We set
$A=\bigl\{
 (i,m)\in I^c\times\seisuu_{>0}
 \bigr\}
 \cup
 \bigl\{
 (i,m)\in I\times\seisuu_{\leq 0}
 \bigr\}$.
We obtain the following countable union of Zariski closed subsets:
\[
 W(E)=\bigcup_{(i,m)\in A}W(E,i,m)\subsetneq S.
\]
By Proposition \ref{prop;26.2.5.10},
we obtain the following.
\begin{prop}
\label{prop;26.3.2.10}
Let $\eta'$ be any geometric point of $S$
such that the image of $\eta$
is not contained in $W(E)$.
Then,
$\nbigd_{X_{\eta'}}(E_{\eta'})
\simeq
\nbige_{\eta'}(!H_{\eta'}(I)\ast H_{\eta'}(I^c))$.
\hfill\qed
\end{prop}

\begin{lem}
\label{lem;26.2.6.20}
Let $E'$ be another good lattice such that
$\nbigd_{X_{\eta(0)}}(E'_{\eta(0)})
 \simeq
\nbige_{\eta(0)}\bigl(!H_{\eta(0)}(I)\ast H_{\eta(0)}(I^c)\bigr)$. 
Then, there exists a Zariski closed subset
$W_1\subset W(E)\cap W(E')$
such that 
$\nbigd_{X/S}(E)_{S\setminus W_1}
\simeq
\nbigd_{X/S}(E')_{S\setminus W_1}$.
\end{lem}
\pf
There exists a positive integer $N$
such that the following holds.
\[
\Sp(E_{\eta(0)},\nabla,i)^{\reg}
\subset
\Bigl\{
\alpha+n\,\Big|\,
\alpha\in \Sp(E'_{\eta(0)},\nabla,i)^{\reg},\,\,
n\in\seisuu,\,\,|n|\leq N
\Bigr\},
\]
\[
\Sp(E'_{\eta(0)},\nabla,i)^{\reg}
\subset
\Bigl\{
\alpha+n\,\Big|\,
\alpha\in \Sp(E_{\eta(0)},\nabla,i)^{\reg},\,\,
n\in\seisuu,\,\,|n|\leq N
\Bigr\}.
\]
Then, the claim is clear.
\hfill\qed

\subsubsection{The associated $\nbigd_{X/S}$-modules and morphisms}
\label{subsection;26.2.5.20}

Suppose that $S$ is irreducible.
Let $H=\bigcup_{j\in\Lambda} H_j$ be the irreducible decomposition.
For any $\vecm\in\seisuu^{\Lambda}$,
we set $\vecm H=\sum m_iH_i$.

Let $(\nbige,\nabla)$ be a good meromorphic flat bundle
on $(X,H)$ relative to $S$.
Let $E$ be a good lattice of $(\nbige,\nabla)$.
We assume the following.
\begin{condition}
\mbox{{}}\label{condition;26.2.6.2}
For any geometric point $\eta(0)$
obtained as an algebraic closure of
the generic point of $S$,
we have
$E_{\eta(0)}\subset
\nbigp^{\DM}_{<0}(\nbige_{\eta(0)})$
for a good $\cnum$-Deligne-Malgrange filtered bundle
$\nbigp^{\DM}_{\ast}(\nbige_{\etabar(0)})$
\hfill\qed
\end{condition}

Take $\vecm\in\seisuu_{>0}^{\Lambda}$
satisfying the following.
\begin{condition}
\mbox{{}}\label{condition;26.2.6.3}
$\nbigp^{\DM}_{(1,\ldots,1)}(\nbige_{\eta(0)})
\subset
E_{\eta(0)}(\vecm H_{\eta(0)})$
for a good $\cnum$-Deligne-Malgrange filtered bundle
$\nbigp^{\DM}_{\ast}(\nbige_{\etabar(0)})$.
\hfill\qed
\end{condition}
We set 
$W(E,\vecm)=W(E)\cup W(E(\vecm H))$.

Let $\nbigb(E,\vecm)\subset S(\hyperk)$
denote the set of closed points $u$
such that the following holds.
\begin{itemize}
 \item $E_u\subset \nbigp^{\DM}_{<0}(\nbige_u)$
       and
       $\nbigp^{\DM}_{(1,\ldots,1)}(\nbige_u)
       \subset E_u(\vecm H_u)$
       for a good $\cnum$-Deligne-Malgrange filtered bundle
       $\nbigp^{\DM}_{\ast}(\nbige_u)$
       of $(\nbige_u,\nabla)$.
\end{itemize}

The following lemma is clear.
\begin{lem}
$W(E,\vecm)\cap \nbigb(E,\vecm)=\emptyset$.
\hfill\qed
\end{lem}

For any $I\subset \Lambda$,
we set $\vecm(I)\in\seisuu_{\geq 0}^{\Lambda}$
such that
$\vecm(I)_i=m_i$ $(i\in I)$ and $\vecm(I)_i=0$ $(i\not\in I)$.
We obtain the $\nbigd_{X/S}$-module
$\nbigd_{X/S}(E(\vecm(I)\cdot H))$.

\begin{lem}
Let $\eta$ be a geometric point such that
the image of $\eta$ is not contained in $W(E,\vecm)$.
Then, there exists a natural isomorphism
$\nbigd_{X/S}\bigl(E(\vecm(I)\cdot H)\bigr)_{\eta}
\simeq
\nbige_{\eta}\bigl(!H_{\eta}(I^c)\ast H_{\eta}(I)\bigr)$.
\end{lem}
\pf
By the construction,
the following holds
for any $u\not\in W(E,\vecm)$.
\begin{itemize}
 \item $\Sp\bigl(E_u(\vecm(I)H_u),\nabla,i\bigr)\cap
       \seisuu_{>0}=
       \Sp\bigl(E_u(\vecm H_u),\nabla,i\bigr)\cap
       \seisuu_{>0}=
       \emptyset$
       for any $i\in I$.
 \item $\Sp\bigl(E_u(\vecm(I)H_u),\nabla,i\bigr)\cap
       \seisuu_{\leq 0}=
       \Sp\bigl(E_u,\nabla,i\bigr)\cap
       \seisuu_{\leq 0}
       =\emptyset$
       for any $i\in I^c$.
\end{itemize}
Then, the claim of the lemma follows.
\hfill\qed

\vspace{.1in}

For $I\subset J$,
there exists the natural morphism
\begin{equation}
\label{eq;26.1.25.1} 
 \nbigd_{X/S}\bigl(E(\vecm(I)\cdot H)\bigr)
 \lrarr
\nbigd_{X/S}\bigl(E(\vecm(J)\cdot H)\bigr).
\end{equation}

By Corollary \ref{cor;26.1.26.20},
we obtain the following lemma.
\begin{lem}
There exists a non-empty open subset
$S'\subset S$
such that
the cokernel of the morphism
\[
\nbigd_{X/S}\bigl(E(\vecm(I)\cdot H)\bigr)_{|X_{S'}}
 \lrarr
\nbigd_{X/S}\bigl(E(\vecm(J)\cdot H)\bigr)_{|X_{S'}}
\]
on $X_{S'}$ is flat over $S'$ 
for any pair $I\subset J$.
\hfill\qed
\end{lem}

Let $E'$ be another good lattice of $(\nbige,\nabla)$
satisfying Condition \ref{condition;26.2.6.2}.
Take $\vecm'\in\seisuu^{\Lambda}$
satisfying Condition \ref{condition;26.2.6.3} for $E'$.

\begin{prop}
There exists a Zariski closed subset
$W_1\subset
W(E,\vecm)\cup W(E',\vecm')$
such that
\[
 \nbigd_{X/S}\bigl(E(\vecm(I)\cdot H)\bigr)_{S\setminus W_1}
 \simeq
 \nbigd_{X/S}\bigl(E'(\vecm'(I)\cdot H)\bigr)_{S\setminus W_1}.
\]
In particular,
 $W_1\cap\bigl(
 \nbigb(E,\vecm)\cap \nbigb(E',\vecm')\bigr)=\emptyset$.
\end{prop}
\pf
It follows from Lemma \ref{lem;26.2.6.20}.
\hfill\qed

\vspace{.1in}
We also have the following lemma for exceptional subsets.

\begin{lem}
If $E'\subset E$ and $E(\vecm H)\subset E'(\vecm'H)$ are satisfied,
then
we have $W(E',\vecm')\subset W(E,\vecm)$
and $\nbigb(E,\vecm)\subset \nbigb(E',\vecm')$.
\hfill\qed
\end{lem}
 
Let us give a complement.

\begin{lem}
\label{lem;26.3.3.10}
Let $E_0$ be a lattice of $\nbige$
which is not necessarily good.
Assume that $E\subset E_0(-H)$
and $E_0(H)\subset E(\vecm H)$.
Let $\nbiga(E_0)$
denote the set of closed points $u$ of $S$
such that $E_{0,u}$ is 
a good $\cnum$-Deligne-Malgrange lattice
of $(\nbige_u,\nabla_u)$.
Then,
we have $\nbiga(E_0)\subset\nbigb(E,\vecm)$.
In particular,
$W(E,\vecm)\cap \nbiga(E_0)=\emptyset$.
\hfill\qed
\end{lem}

\subsubsection{Push-forward}

We continue to assume that $S$ is irreducible.
Let $Y$ be a $\hyperk$-variety
with a smooth projective morphism $Y\to S$.
Let $H_Y$ be a hypersurface of $Y$ flat over $S$,
which is not necessarily normal crossing.
Let $H_Y=\bigcup_{j\in \Lambda_Y}H_{Y,j}$
be the irreducible decomposition.

Let $f:X\to Y$ be a morphism of projective $\hyperk$-varieties over $S$
such that $H=f^{-1}(H_Y)$
and that $f:X\setminus H\to Y\setminus H_Y$
is a closed immersion.
We have the induced map
$f_{\ast}:\Lambda\to \Lambda_Y$.

Let $(\nbige,\nabla)$ be a good meromorphic flat bundle
on $(X,H)$ relative to $S$.
Let $E$ and $\vecm$
be as in \S\ref{subsection;26.2.5.20}.

\begin{lem}
There exists an open subset $S'\subset S$
such that
$f_{\dagger}^j\Bigl(
 \nbigd_{X/S}\Bigl(
 E\bigl(
 \vecm(f_{\ast}^{-1}(I_Y))H\bigr)
 \Bigr)
\Bigr)_{|Y_{S'}}=0$ $(j\neq 0)$
for any $I_Y\subset\Lambda_Y$. 
For any geometric point $\eta$ of $S'$
whose image is not contained
in $W(E,\vecm)$,
there exists the natural isomorphism
\begin{equation}
\label{eq;26.2.6.1}
 f_{\dagger}^0\Bigl(
 \nbigd_{X/S}\Bigl(
 E\bigl(
 \vecm(f_{\ast}^{-1}(I_Y))H\bigr)
 \Bigr)
 \Bigr)_{\eta}
 \simeq
 f_{\dagger}^0\bigl(
 \nbige_{\eta}(!H_{\eta})
 \bigr)\bigl(\ast H_{Y,\eta}(I_Y)\bigr).
\end{equation}
\end{lem}
\pf
Let $\etabar$ denote an algebraic closure of
a generic point of $S$.
We have
$f_{\dagger}^j\Bigl(
\nbigd_{X_{\etabar}}\bigl(E_{\etabar}(\vecm(I_1)H_{\etabar})\bigr)
\Bigr)_{|Y_{\etabar}}=0$ for $j\neq 0$,
where $I_1=f_{\ast}^{-1}(I_Y)$.
Then, we obtain the first claim.
The second claim follows from Proposition \ref{prop;26.3.2.10}.
\hfill\qed

\begin{cor}
There exists a decomposition
$S=\bigsqcup S_j$
by locally closed smooth subsets
 such that
\[
 f^k_{S_j\dagger}\Bigl(
 \nbigd_{X_{S_j}/S_j}\Bigl(
 E_{S_j}\bigl(
 \vecm(f_{\ast}^{-1}(I_Y))H_{S_j}\bigr)
 \Bigr)
\Bigr)=0\quad(k\neq 0)
\] 
for any $I_Y\subset\Lambda_Y$. 
For any geometric point $\eta$ of $S_j$
whose image is not contained
in $W(E,\vecm)\cap S_j$,
there exists the natural isomorphism
{\rm(\ref{eq;26.2.6.1})}. 
\hfill\qed
\end{cor}

\subsubsection{Refinements}
\label{subsection;26.2.7.1}

Let $K\subset X$ be a hypersurface
such that $K$ and $H\cup K$ are flat over $S$
and that $\codim (H\cap K)\geq 2$.
We do not assume that $K\cup H$ is normal crossing
relative to $S$.

Let $X'$ be a $\hyperk$-variety
with a smooth projective morphism $X'\to S$.
Let $H'$ be a hypersurface of $X$
normal crossing relative to $S$.
Let $H'=\bigcup_{j\in\Lambda'}H'_j$
be the irreducible decomposition.
Let $f:X'\to X$ be a morphism over $S$
such that $f^{-1}(H\cup K)=H'$
and that $X'\setminus H'\simeq X\setminus (H\cup K)$.
Note that $f^{-1}(H)$ is also normal crossing relative to $S$.

Let $(\nbige,\nabla)$ be
a good meromorphic flat bundle on $(X,H)$
relative to $S$.
Let $E$ be a good lattice of $(\nbige,\nabla)$
satisfying Condition \ref{condition;26.2.6.2}.
Take $\vecm$ satisfying Condition \ref{condition;26.2.6.3}.

We obtain a good meromorphic flat bundle
$(\nbige',\nabla')=f^{\ast}(\nbige,\nabla)$
on $(X',H')$ relative to $S$.
We also obtain the good lattice $E'=f^{\ast}(E)(-H')$
satisfying Condition \ref{condition;26.2.6.2}
for $(\nbige',\nabla')$.
Take $\vecm'\in\seisuu_{>0}^{\Lambda'}$
such that
$f^{\ast}(E(\vecm))(H')\subset E'(\vecm'H')$.
It satisfies Condition \ref{condition;26.2.6.3} for $E'$.

For any decomposition $I\sqcup I^c$,
we obtain
the $\nbigd_{X/S}$-module
$\nbigd_{X/S}\bigl(E(\vecm(I)H)\bigr)$.
There exists $I_1\subset\Lambda'$
such that
$H'(I_1)=f^{-1}(H(I)\cup K)$.
There exists $I_2\subset\Lambda'$
such that
$H'(I_2)=f^{-1}(H(I))$.
We obtain
\[
\nbigd_{X'/S}\bigl(E'(\vecm'(I_2)H')\bigr)
\lrarr
\nbigd_{X'/S}\bigl(E'(\vecm'(I_1)H')\bigr).
\]
There exist the following natural morphisms:
{\footnotesize
\begin{equation}
\label{eq;26.2.6.30}
\begin{CD}
 f^0_{\eta(0)\dagger}
  \nbigd_{X'_{\eta(0)}}(E'_{\eta(0)}(\vecm'(I_2) H'_{\eta(0)}))
 @>>>
  \nbigd_{X_{\eta(0)}}(E_{\eta(0)}(\vecm H_{\eta(0)}))
 @>>>
  f^0_{\eta(0)\dagger}
 \nbigd_{X'_{\eta(0)}}(E'_{\eta(0)}(\vecm'(I_1) H'_{\eta(0)}))
 \\
 @V{\simeq}VV @V{\simeq}VV @V{\simeq}VV \\
 f^0_{\eta(0)\dagger}\Bigl(
 \nbige'_{\eta(0)}
 \bigl(!H'_{\eta(0)}\bigr)\Bigr)\bigl(\ast H_{\eta(0)}(I)\bigr)
 @>>>
 \nbige_{\eta(0)}\bigl(!H_{\eta(0)}(I^c)\ast H_{\eta(0)}(I)\bigr)
 @>>>
 f^0_{\eta(0)\dagger}
 \nbige'_{\eta(0)}
 \Bigl(
 \bigl(!H'_{\eta(0)}\bigr)
 \Bigr)\bigl(\ast (H_{\eta(0)}(I)\cup K_{\eta(0)})\bigr).
\end{CD}
\end{equation}
}

\begin{prop}
\label{prop;26.2.7.10}
There exist a Zariski open subset $S'\subset S$
and morphisms
\[
 f_{\dagger}^0
 \Bigl(\nbigd_{X'/S}\bigl(
  E'(\vecm'(I_2)H')
 \bigr)
 \Bigr)_{S'}
 \stackrel{a}{\lrarr}
 \nbigd_{X/S}\bigl(
 E(\vecm(I)H)
 \bigr)_{S'}
 \stackrel{b}{\lrarr}
 f_{\dagger}^0
 \Bigl(\nbigd_{X'/S}\bigl(
  E'(\vecm'(I_1)H')
 \bigr)
 \Bigr)_{S'}
\]
which induces {\rm(\ref{eq;26.2.6.30})}. 
Moreover the following holds.
\begin{itemize}
 \item The morphism $a$ is an epimorphism,
       the morphism $b$ is a monomorphism,
       and the cokernel of $b$ is flat over $S'$.
 \item Let $\eta$ be a geometric point of $S'$
       whose image is not contained in
       $W(E,\vecm)\cup W(E',\vecm')$.
       Then,
       there exists the following commutative diagram:
{\small
\[
\begin{CD}
f_{\dagger}^0
 \Bigl(\nbigd_{X'/S}\bigl(
  E'(\vecm'(I_2)H')
 \bigr)
 \Bigr)_{\eta}
 @>>>
 \nbigd_{X/S}\bigl(
 E(\vecm(I)H)
 \bigr)_{\eta}
 @>>>
 f_{\dagger}^0
 \Bigl(\nbigd_{X'/S}\bigl(
  E'(\vecm'(I_1)H')
 \bigr)
 \Bigr)_{\eta}
 \\
 @V{\simeq}VV @V{\simeq}VV @V{\simeq}VV \\
 f_{\dagger}^0\Bigl(
 \nbige'_{\eta}(!H'_{\eta})
 \Bigr)(\ast H_{\eta}(I))
 @>>>
 \nbige_{\eta}\bigl(!H_{\eta}(I^c)\ast H_{\eta}(I)\bigr)
 @>>>
 f_{\dagger}^0\Bigl(
 \nbige'_{\eta}(!H'_{\eta})
 \Bigr)
 \bigl(\ast (H_{\eta}(I)\cup K_{\eta})\bigr).
\end{CD}
\]     }  
\end{itemize} 
\end{prop}
\pf
There exists a Zariski open subset $S'_1\subset S$
such that the following holds.
\begin{itemize}
 \item
      $f^j_{\dagger}\Bigl(
      \nbigd_{X'/S}\bigl(
      E'(\vecm'(I_1)H')
      \bigr)_{S'_1}
      \Bigr)=0$
      and
      $f^j_{\dagger}\Bigl(
      \nbigd_{X'/S}\bigl(
      E'(\vecm'(I_2)H')
      \bigr)_{S'_1}
      \Bigr)=0$ for $j\neq 0$.
 \item
      The cokernel of the natural morphism
\begin{equation}
\label{eq;26.2.6.32}
  f_{\dagger}^0
 \Bigl(\nbigd_{X'/S}\bigl(
  E'(\vecm'(I_2)H')
 \bigr)_{S'_1}
 \Bigr)
 \lrarr
 f_{\dagger}^0
 \Bigl(\nbigd_{X'/S}\bigl(
  E'(\vecm'(I_1)H')
 \bigr)_{S'_1}
 \Bigr)
\end{equation}
is flat over $S'_1$.      
\end{itemize}
Let $\nbigm$ denote the image of (\ref{eq;26.2.6.32}).
Note that $\nbigm_{\eta(0)}$
equals the image of
the middle term in (\ref{eq;26.2.6.30}).
Moreover, for any geometric point $\eta$
whose image is not contained in $W(E',\vecm')$,
there exist the following commutative diagram:
{\small
\[
\begin{CD}
f_{\dagger}^0
 \Bigl(\nbigd_{X'/S}\bigl(
  E'(\vecm'(I_2)H')
 \bigr)
 \Bigr)_{\eta}
 @>>>
 \nbigm_{\eta}
 @>>>
 f_{\dagger}^0
 \Bigl(\nbigd_{X'/S}\bigl(
  E'(\vecm'(I_1)H')
 \bigr)
 \Bigr)_{\eta}
 \\
 @V{\simeq}VV @V{\simeq}VV @V{\simeq}VV \\
 f_{\dagger}^0\Bigl(\nbige'_{\eta}\bigl(!H'_{\eta}(I_2^c)\bigr)\Bigr)
 (\ast H_{\eta}(I))
 @>>>
 \nbige_{\eta}(!H_{\eta}(I^c)\ast H_{\eta}(I))
 @>>>
 f_{\dagger}^0\Bigl(
 \nbige'_{\eta}(!H'_{\eta})\Bigr)
 \bigl(\ast (H_{\eta}(I)\cup K_{\eta})\bigr).
\end{CD}
\]     }  
By shrinking $S_1'$,
we may assume that
$\nbigm_{S_1'}\simeq
\nbigd_{X/S}\bigl(E(\vecm(I)H)\bigr)_{S_1'}$.
\hfill\qed

\begin{cor}
There exist a Zariski closed subset
$W\subset W(E,\vecm)\cup W(E',\vecm')$,
a decomposition
$S\setminus W=\bigsqcup S_j$,
and morphisms
\[
 f_{S_j\dagger}^0
 \Bigl(\nbigd_{X'/S}\bigl(
  E'(\vecm'(I_2)H')
 \bigr)_{S_j}
 \Bigr)
\stackrel{a_j}{\lrarr}
 \nbigd_{X/S}\bigl(
 E(\vecm(I)H)
 \bigr)_{S_j}
\stackrel{b_j}{\lrarr}
 f_{S_j\dagger}^0
 \Bigl(\nbigd_{X'/S}\bigl(
  E'(\vecm'(I_1)H')
 \bigr)_{S_j}
 \Bigr)
\]
such that the following holds.
\begin{itemize}
 \item The morphisms $a_j$ are epimorphisms,
       the morphisms $b_j$ are monomorphisms,
       and the cokernel of $b_j$ are flat over $S_j$.
 \item Let $\eta$ be a geometric point of $S_j$
       whose image is not contained in
       $W(E,\vecm)\cup W(E',\vecm')$.
       Then,
       there exist the following commutative diagram:
{\small
\[
\begin{CD}
f_{S_j\dagger}^0
 \Bigl(\nbigd_{X'/S}\bigl(
  E'(\vecm'(I_2)H')
 \bigr)_{S_j}
 \Bigr)_{\eta}
 @>>>
 \nbigd_{X/S}\bigl(
 E(\vecm(I)H)
 \bigr)_{\eta}
 @>>>
 f_{S_j\dagger}^0
 \Bigl(\nbigd_{X'/S}\bigl(
  E'(\vecm'(I_1)H')
 \bigr)
 \Bigr)_{\eta}
 \\
 @V{\simeq}VV @V{\simeq}VV @V{\simeq}VV \\
 f_{S_j\dagger}^0\Bigl(
 \nbige'_{\eta}(!H'_{\eta})
 \Bigr)
 \bigl(\ast H_{\eta}(I)\bigr)
 @>>>
 \nbige_{\eta}\bigl(!H_{\eta}(I^c)\ast H_{\eta}(I)\bigr)
 @>>>
 f_{S_j\dagger}^0\Bigl(
 \nbige'_{\eta}(!H_{\eta})
 \Bigr)
 \bigl(
 \ast (H_{\eta}(I)\cup K)
 \bigr).
\end{CD}
\]     }  
\end{itemize} 
\hfill\qed
\end{cor}

Let us give a complement.
\begin{lem}
$\nbigb(E,\vecm)
\subset \nbigb(E',\vecm')$.
As a result,
we have
$\nbigb(E,\vecm)\cap(W(E,\vecm)\cup W(E,\vecm'))=\emptyset$.
\end{lem}
\pf
We have
$f^{\ast}(E)(-H)_u\subset
f^{\ast}(\nbigp_{<0}^{\DM}\nbige_u)(-H)
\subset
\nbigp_{<0}^{\DM}(f^{\ast}\nbige_u)$.
We also have
\[
\nbigp^{\DM}_{(1,\ldots,1)}(\nbige'_u)(H)
\subset
f^{\ast}(\nbigp^{\DM}_{(1,\ldots,1)}\nbige_u)(H)
\subset
f^{\ast}(E(\vecm H))(H)
\subset E'(\vecm'H').
\]
Then, we obtain the claim of the lemma.
\hfill\qed

\subsection{Almost minimal extensions associated with cells}

\subsubsection{Cells and almost minimal extensions}

Let $X$ be a smooth $\hyperk$-variety.
Let $Z$ be a  closed subset of $X$.
Let $H$ be a hypersurface of $X$.

Let $S$ be a smooth $\hyperk$-variety
which is not necessarily irreducible.
Let $Z_S$ be a smooth $\hyperk$-variety projective over $S$
with a normal crossing hypersurface
$H_{Z_S}$ relative to $S$.
Let $\rho:Z_S\to S\times X$ be a projective morphism over $S$
such that
(i) $H_{Z_S}=\rho^{-1}(H)$,
(ii) $Z_S\setminus H_{Z_S}\simeq S\times(Z\setminus H)$.
Let $(\nbige,\nabla)$ be a good meromorphic flat bundle
on $(Z_S,H_{Z_S})$ relative to $S$
with a lattice $E_0$.

Let $H_{Z_S}=\bigcup_{j\in\Lambda}H_{Z_S,j}$
denote the irreducible decomposition.
Let $E$ be a good lattice of $(\nbige,\nabla)$
satisfying Condition \ref{condition;26.2.6.2}.
Take $\vecm\in\seisuu_{>0}^{\Lambda}$
satisfying Condition \ref{condition;26.2.6.3}.

Let $H=\bigcup_{j\in\Gamma}H_j$ be the irreducible decomposition.
For any $I\subset \Gamma$,
we obtain $\rho^{\ast}(I)\subset \Lambda$
such that
$\rho^{-1}(S\times H(I))=\bigcup_{j\in \rho^{\ast}(I)}H_{Z_S,j}$.

\begin{lem}
\label{lem;26.2.7.21}
There exist a Zariski closed subset
$W\subset W(E,\vecm)$,
a morphism of smooth $\hyperk$-varieties
$\Stilde\to S$
such that the following holds:
\begin{itemize}
 \item The image of $\Stilde\to S$ is $S\setminus W$.
 \item For any $I\subset\Gamma$,
       $\rho_{\Stilde\dagger}\Bigl(
       \nbigd_{Z_{\Stilde}/\Stilde}
       \bigl(E(\vecm(\rho^{\ast}I)H_{Z_S})
       \bigr)_{\Stilde}
       \Bigr)$
       are $\Stilde$-flat,
       and 
       $\rho^k_{\Stilde\dagger}\Bigl(
       \nbigd_{Z_S/S}\bigl(E(\vecm(\rho^{\ast}I)H_{Z_S})\bigr)_{\Stilde}
       \Bigr)=0$ $(k\neq 0)$.
       
 \item Let $\eta$ be any geometric point of $\Stilde$
       whose image is not contained in
       $\Stilde\times_SW(E,\vecm)$.
       Then, there exist the following natural isomorphisms
\[
       \rho_{\Stilde\dagger}\Bigl(
       \nbigd_{Z_S/S}\bigl(E(\vecm(\rho^{\ast}I)H_{Z_S})\bigr)_{\Stilde}
       \Bigr)_{\eta}
       \simeq
       \rho_{\eta\dagger}\Bigl(
       \nbige_{\eta}\bigl(!H_{Z_S,\eta}
       \bigr)
       \Bigr)(\ast H(I)_{\eta}).
\]
       Moreover,
       the following diagrams are commutative
       for $I\subset J$:
\[
 \begin{CD}
  \rho_{\Stilde\dagger}\Bigl(
  \nbigd_{Z_S/S}(E(\vecm(\rho^{\ast}I)H_{Z_S}))_{\Stilde}
  \Bigr)_{\eta}
  @>{\simeq}>>
 \rho_{\eta\dagger}\Bigl(
  \nbige_{\eta}\bigl(!H_{Z_S,\eta}\bigr)
  \Bigr)(\ast H(I)_{\eta})
 \\
 @VVV @VVV \\
  \rho_{\Stilde\dagger}\Bigl(
  \nbigd_{Z_S/S}(E(\vecm(\rho^{\ast}J)H_{Z_S}))_{\Stilde}
  \Bigr)_{\eta}
  @>{\simeq}>>
  \rho_{\eta\dagger}\Bigl(
  \nbige_{\eta}\bigl(!H_{Z_S,\eta}\bigr)
  \Bigr)(\ast H(J)_{\eta}).
 \end{CD} 
\]
 \item The cokernel of
       $\rho_{\Stilde\dagger}
       \Bigl(\nbigd_{Z_S/S}\bigl(E(\vecm(\rho^{\ast}I)H_{Z_S})\bigr)
       _{\Stilde}\Bigr)
       \to
       \rho_{\Stilde\dagger}\Bigl(
       \nbigd_{Z_S/S}\bigl(E(\vecm(\rho^{\ast}J)H_{Z_S})\bigr)_{\Stilde}
       \Bigr)$
       are flat over $\Stilde$.
\hfill\qed
\end{itemize}
\end{lem}

By replacing $S$ with $\Stilde$,
we may assume that the claims of Lemma \ref{lem;26.2.7.21}
even in the case $\Stilde=S$.
As a consequence,
there exist the following natural isomorphisms:
\begin{multline}
\Cok\Bigl(
       \rho_{\dagger}
       \Bigl(\nbigd_{Z_S/S}\bigl(E(\vecm(\rho^{\ast}I)H_{Z_S})\bigr)\Bigr)
       \to
       \rho_{\dagger}\Bigl(
       \nbigd_{Z_S/S}\bigl(E(\vecm(\rho^{\ast}J)H_{Z_S})\bigr)
       \Bigr)
       \Bigr)
       \simeq \\
 \Cok\Bigl(
       \rho_{\dagger}
       \Bigl(\nbigd_{Z_S/S}\bigl(E(\vecm(\rho^{\ast}I)H_{Z_S})\bigr)\Bigr)
       \to
       \rho_{\dagger}\Bigl(
       \nbigd_{Z_S/S}\bigl(E(\vecm(\rho^{\ast}J)H_{Z_S})\bigr)
       \Bigr)
 \Bigr),
\end{multline}
\begin{multline}
 \Ker\Bigl(
       \rho_{\dagger}
       \Bigl(\nbigd_{Z_S/S}\bigl(E(\vecm(\rho^{\ast}I)H_{Z_S})\bigr)\Bigr)
       \to
       \rho_{\Stilde\dagger}\Bigl(
       \nbigd_{Z_S/S}\bigl(E(\vecm(\rho^{\ast}J)H_{Z_S})\bigr)
       \Bigr)
       \Bigr)
       \simeq \\
 \Ker\Bigl(
       \rho_{\dagger}
       \Bigl(\nbigd_{Z_S/S}\bigl(E(\vecm(\rho^{\ast}I)H_{Z_S})\bigr)\Bigr)
       \to
       \rho_{\dagger}\Bigl(
       \nbigd_{Z_S/S}\bigl(E(\vecm(\rho^{\ast}J)H_{Z_S})\bigr)
       \Bigr)
       \Bigr),
\end{multline}
\begin{multline}
 \Image\Bigl(
       \rho_{\dagger}
       \Bigl(\nbigd_{Z_S/S}\bigl(E(\vecm(\rho^{\ast}I)H_{Z_S})\bigr)\Bigr)
       \to
       \rho_{\dagger}\Bigl(
       \nbigd_{Z_S/S}\bigl(E(\vecm(\rho^{\ast}J)H_{Z_S})\bigr)
       \Bigr)
       \Bigr)
       \simeq \\
 \Image\Bigl(
       \rho_{\dagger}
       \Bigl(\nbigd_{Z_S/S}\bigl(E(\vecm(\rho^{\ast}I)H_{Z_S})\bigr)\Bigr)
       \to
       \rho_{\dagger}\Bigl(
       \nbigd_{Z_S/S}\bigl(E(\vecm(\rho^{\ast}J)H_{Z_S})\bigr)
       \Bigr)
       \Bigr).
\end{multline}

\begin{df}
For a tuple $\nbigu=(S,Z_S,\rho,\nbige,\nabla,E,\vecm)$
such that
the claim of Lemma {\rm\ref{lem;26.2.7.21}}
 holds with $\Stilde=S$,
let
$\nbigl(\nbigu)$
denote the image of
\[
       \rho_{\dagger}
       \Bigl(\nbigd_{Z_S/S}\bigl(E\bigr)\Bigr)
       \to
       \rho_{\dagger}\Bigl(
       \nbigd_{Z_S/S}\bigl(E(\vecm H_{Z_S})\bigr)
       \Bigr).
\]
For any $K=\bigcup_{i\in I}H_i\subset H$,
let 
$\nbigl^{\ast K}(\nbigu)$
denote the image of
\[
       \rho_{\dagger}
 \Bigl(\nbigd_{Z_S/S}\bigl(E(\vecm(\rho^{\ast}I)
 H_{Z_S})\bigr)\Bigr)
       \to
       \rho_{\dagger}\Bigl(
       \nbigd_{Z_S/S}\bigl(E(\vecm H_{Z_S})\bigr)
       \Bigr).
\]
Note that $\nbigl(\nbigu)$
and $\nbigl^{\ast K}(\nbigu)$ are $S$-flat.
\hfill\qed
\end{df}

Note that $E(\vecm H)^{\lor}$ be a good lattice of
$(\nbige^{\lor},\nabla)$
satisfying Condition \ref{condition;26.2.6.2},
and that
$\vecm$ satisfies Condition \ref{condition;26.2.6.3}
for $(\nbige^{\lor},\nabla)$ and $E(\vecm H)^{\lor}$.
Hence, from a tuple $\nbigu=(S,Z_S,\rho,\nbige,\nabla,E,\vecm)$
we obtain another tuple 
$\nbigu^{\lor}=(S,Z_S,\rho,\nbige^{\lor},\nabla,
E(\vecm H)^{\lor},\vecm)$.

\begin{df}
A tuple $\nbigu=(S,Z_S,\rho,\nbige,\nabla,E,\vecm)$
is called a cell if the following holds.
\begin{itemize}
 \item 
The claim of Lemma {\rm\ref{lem;26.2.7.21}}
holds $\nbigu=(S,Z_S,\rho,\nbige,\nabla,E,\vecm)$
and
$\nbigu^{\lor}=(S,Z_S,\rho,\nbige^{\lor},\nabla,
E(\vecm H)^{\lor},
\vecm)$
with 
$\Stilde=S$.
 \item $\DD \nbigl(\nbigu)$
       and $\DD\nbigl(\nbigu^{\lor})$
       are $S$-flat.
\end{itemize}
We set
$\nbigb(\nbigu)=\nbigb(E,\vecm)$.
Note that
$\nbigb(\nbigu)=\nbigb(\nbigu^{\lor})$.
\hfill\qed
\end{df}

\begin{lem}
\label{lem;26.3.3.1}
Let $\nbigu=(S,Z_S,\rho,\nbige,\nabla,E,\vecm)$
be a cell.
Let $\eta$ be any geometric point of $S$
whose image is not contained in
$W(E,\vecm)$.
Then, there exists the natural isomorphism
$(\DD\nbigl(\nbigu))_{\eta}
 \simeq
 \nbigl(\nbigu^{\lor})_{\eta}$.
\hfill\qed
\end{lem}

\begin{lem}
\label{lem;26.3.3.2}
There exist a Zariski closed subset
$W\subset W(E,\vecm)\cup W(E(\vecm H)^{\lor},\vecm)$
and a natural isomorphism
$\DD(\nbigl(\nbigu))_{S\setminus W}
 \simeq
 \nbigl(\nbigu^{\lor})_{S\setminus W}$. 
\hfill\qed
\end{lem}

\subsubsection{Change of hypersurfaces}

Let $H'$ be a hypersurface of $X$
such that $H\subset H'$.
By refining $S$,
we assume that there exists
a projective morphism $f:Z_S'\to Z_S$ over $S$
such that
$H'_{Z_S'}=(\rho\circ f)^{-1}(H')$
is simply normal crossing relative to $S$.
We set $\rho':=\rho\circ f$.

We set
$(\nbige',\nabla)=f^{\ast}(\nbige,\nabla)$ on $(Z'_S,H'_{Z_S'})$.
We obtain a good lattice 
$E'=f^{\ast}(E)(-H'_{Z_S'})$
of $(\nbige',\nabla)$
satisfying Condition \ref{condition;26.2.6.2}
for $(\nbige',\nabla)$.
Take $\vecm'$
such that
$f^{\ast}(E(\vecm H_{Z_S}))(H'_{Z_S'})
\subset E'(\vecm' H'_{Z_S'})$,
which satisfies
Condition \ref{condition;26.2.6.3}
for $(\nbige',\nabla)$ and $E'$.
By refining $S$ if necessary,
we obtain a cell
$\nbigu'=
(S,Z_S',\rho',\nbige',\nabla',E',\vecm')$.

\begin{lem}
\label{lem;26.2.7.20}
There exist a Zariski closed subset
$W\subset W(E,\vecm)\cup W(E',\vecm')$
and a natural isomorphism
$\nbigl(\nbigu)_{S\setminus W}
 \simeq
 \nbigl(\nbigu')_{S\setminus W}$.
We also have
$\nbigb(E,\vecm)\subset \nbigb(E',\vecm')$
and hence
 $\nbigb(E,\vecm)\cap
 \bigl(W(E,\vecm)\cup W(E',\vecm')\bigr)=\emptyset$.
In particular,
$W\cap \nbigb(E,\vecm)=\emptyset$.
\end{lem}
\pf
It follows from Proposition \ref{prop;26.2.7.10}.
\hfill\qed

\subsection{Approximation of universal extensions}

\subsubsection{Localizability over subsets}

Let $S$ be a smooth $\hyperk$-variety
with a subset $A\subset S(\hyperk)$.
Let $M_S$ be an $S$-flat holonomic $\nbigd_{S\times X/S}$-module.

\begin{df}
\label{df;26.3.3.1}
We say that $M_S$ is localizable over $A$
if for any hypersurface $H$
there exist a quasi-finite morphism of smooth $\hyperk$-varieties
$\Stilde\to S$,
a subset $\Atilde\subset \Stilde(\hyperk)$,
a holonomic $\nbigd_{\Stilde\times X/\Stilde}$-module
$M_{\Stilde,1}$, 
$\nbigd_{\Stilde\times X/\Stilde}$-homomorphisms
\[
 M_{\Stilde}\to M_{\Stilde,1}\to M_{\Stilde}(\ast (\Stilde\times H)),
\] 
and a countable union of Zariski closed subsets
$W(M_{\Stilde,1})\subset\Stilde$
such that the following holds.
\begin{itemize}
 \item $\Atilde$ is mapped onto $A$.
 \item $W(M_{\Stilde,1})\cap \Atilde=\emptyset$.
 \item Let $\eta$ be any geometric point of $\Stilde$
       whose image is not contained in
       $W(M_{\Stilde,1})$.
       Then, the induced morphism
       $(M_{\Stilde,1})_{\eta}
       \to
       M_{\Stilde}(\ast (\Stilde\times H))_{\eta}$
       is an isomorphism.       
\hfill\qed
\end{itemize}
\end{df}

\begin{rem}
Such $\Stilde$ and $M_{\Stilde,1}$ are not unique.
\hfill\qed
\end{rem}

The following lemma is clear.
\begin{lem}
Suppose that $M_S$ is localizable over $A$.
Let $T\to S$ be any locally closed irreducible subset.
We set $A_T=A\cap T(\hyperk)$.
Then, $M_T$ is localizable over $A_T$.
\hfill\qed 
\end{lem}

\begin{lem}
Let $0\to M'_S\to M_S\to M''_S\to 0$
be an exact sequence of
$S$-flat holonomic $\nbigd_{S\times X/S}$-modules.
If two of $M'_S$, $M_S$ and $M''_S$ are localizable over $A$,
the remaining one is also localizable over $A$.
\end{lem}
\pf
Let $T\subset S$ be any locally closed irreducible subset.
We set $A_T=A\cap T(\hyperk)$.
It is enough to consider the case $A_T\neq\emptyset$.
Let $\eta$ be a geometric point
obtained as an algebraic closure of the generic point of $T$.
There exist a quasi-finite morphism
$\Ttilde\to T$
and 
an exact sequence
of $\Ttilde$-flat holonomic $\nbigd_{\Ttilde\times X/\Ttilde}$-modules
$0\to N'_{\Ttilde}\to N_{\Ttilde}\to N''_{\Ttilde}\to 0$
which induces 
$0\to M'_{\eta}(\ast H_{\eta})\to
M_{\eta}(\ast H_{\eta})
\to M''_{\eta}(\ast H_{\eta})\to 0$.
We may assume that 
there exist morphisms of exact sequences
$(0\to M'_{\Ttilde}\to M_{\Ttilde}\to M''_{\Ttilde}\to 0)
\to
(0\to N'_{\Ttilde}\to N_{\Ttilde}\to N''_{\Ttilde}\to 0)
\to
(0\to M'_{\Ttilde}\to M_{\Ttilde}\to M''_{\Ttilde}\to 0)
(\ast \Ttilde\times H)$.
Let $A_{\Ttilde}\subset \Ttilde(\hyperk)$
be the inverse image of $A_T$.

Let us show that $M_S$ is localizable over $A$
by assuming that $M'_S$ and $M''_S$ are localizable over $A$.
The other claims can be proved similarly.
We may assume that there exist
$M'_{\Ttilde}\to M'_{\Ttilde,j}\to M'_{\Ttilde}(\ast (\Ttilde\times H))$
and
$M''_{\Ttilde}\to M''_{\Ttilde,j}\to M''_{\Ttilde}(\ast (\Ttilde\times H))$
as in the definition.
We may also assume that
there exist isomorphisms
$M'_{\Ttilde,j}\simeq N'_{\Ttilde}$
and
$M''_{\Ttilde,j}\simeq N''_{\Ttilde}$.
We set $W(M_{\Ttilde,j})_{\Ttilde}
=\bigl(W(M'_{\Ttilde,j})\cup W(M''_{\Ttilde,j})\bigr)\times_S\Ttilde$.
Let $\eta$ be a geometric point of $\Ttilde$
whose image is not contained in $W_{\Ttilde}$.
Because
$(M'_{\Ttilde,j})_{\eta}\simeq
(M'_{\Ttilde})_{\eta}(\ast H_{\eta})$
and
$(M''_{\Ttilde,j})_{\eta}\simeq
(M''_{\Ttilde})_{\eta}(\ast H_{\eta})$,
we obtain that
$(M_{\Ttilde,j})_{\eta}\simeq
(M_{\Ttilde})_{\eta}(\ast H_{\eta})$.
Hence, we obtain that $M$ is also localizable over $A$
by using a Noetherian induction.
\hfill\qed

\begin{lem}
Let $\kappa:M\to M'$ be a morphism of $S$-flat holonomic
$\nbigd_{S\times X/S}$-modules.
Assume that $M$ and $M'$ are localizable over $A$,
and that $\Cok\kappa$ is $S$-flat.
Then,  $\Image\kappa$, $\Ker\kappa$ and $\Cok(\kappa)$
are also localizable over $A$.
\end{lem}
\pf
There exist
$(\Stilde_i,\Atilde_i,(M_i)_{\Stilde_i,1},W(M_{\Stilde_i,1}))$
as in Definition \ref{df;26.3.3.1}.
By refining $\Stilde_i$,
we may assume that $\Stilde_1=\Stilde_2=:\Stilde$
and $\Atilde_1=\Atilde_2=:\Atilde$.
We may also assume that
there exists the following commutative diagram:
\[
 \begin{CD}
 (M_1)_{\Stilde}
 @>>> (M_1)_{\Stilde,1}
 @>>>
  (M_1)_{\Stilde}(\ast (\Stilde\times H))
  \\
  @VVV @V{\kappa_{\Stilde,1}}VV @VVV \\
  (M_2)_{\Stilde}
 @>>> (M_2)_{\Stilde,1}
 @>>>
  (M_2)_{\Stilde}(\ast (\Stilde\times H)).
 \end{CD}
\]
Moreover, we may assume that the cokernel of 
$\kappa_{\Stilde,1}$
is $\Stilde$-flat.
Then, the claim is clear.
\hfill\qed

\begin{lem}
For $\nbigu=(S,Z_S,\rho,\nbige,\nabla,E,\vecm)$
be any cell,
$\nbigl(\nbigu)$ and $\DD\nbigl(\nbigu)$ are localizable over
$\nbigb(\nbigu)$.
\end{lem}
\pf
It follows from
Lemma \ref{lem;26.3.3.2} and
Lemma \ref{lem;26.2.7.20}.
\hfill\qed

\subsubsection{Strong localizability}

Let $M$ be an $S$-flat holonomic $\nbigd_{S\times X/S}$-module
with $A\subset S(\hyperk)$.

\begin{df}
We say that $M$ is strongly localizable
if the following holds 
 for any smooth $\hyperk$-variety $Y$
 any cell
 $\nbigu_Y=(S_1,\rho_1,\nbige,\nabla,E,\vecm)$
 on $Y$.
\begin{itemize}
 \item The $\nbigd_{(S_1\times S)\times(Y\times X)
       /(S_1\times S)}$-modules
       $\nbigl\bigl(\nbigu_Y\bigr)
       \boxtimes M$
       and $\DD
       \nbigl\bigl(\nbigu_Y\bigr)
	\boxtimes M$
	are localizable over $\nbigb(\nbigu_Y)\times A$.
\hfill\qed
\end{itemize}
\end{df}

\begin{lem}
Let $0\to M'\to M\to M''\to 0$
be an exact sequence of $S$-flat holonomic $\nbigd_{S\times X/S}$-modules.
If two of them are strongly localizable,
then the rest is also strongly localizable. 
\hfill\qed
\end{lem}

\begin{lem}
Let $\varphi:M\to M'$ be a morphism of
$S$-flat holonomic $\nbigd_{S\times X/S}$-modules.
Suppose that $\Cok\varphi$ is $S$-flat,
and that $M$ and $M'$ are strongly localizable.
Then, $\Cok(\varphi)$, $\Image(\varphi)$ and $\Ker(\varphi)$
are strongly localizable.
\hfill\qed 
\end{lem}

\subsubsection{Approximation of the space of extensions}

Let $S$ be a smooth $\hyperk$-variety.
Let $A\subset S(\hyperk)$ be a subset.
Let $M$ be an $S$-flat holonomic $\nbigd_{S\times X/S}$-module
which is strongly localizable over $A$.

Let $\nbigu=(T,Z,H,\rho,\nbige,E)$ be a cell on $X$.
We obtain
the $(T\times S)$-flat holonomic
$\nbigd_{(T\times S)\times(X\times X)/T\times S}$-module
$\DD \nbigl(\nbigu)\boxtimes M$.
It is localizable over $\nbigb(\nbigu)\times A$.

Let $\Delta:X\to X\times X$ denote the diagonal embedding.
It induces
$(T\times S)\times X\to (T\times S)\times X\times X$.
Let $H_j$ $(j=1,\ldots,m)$
be a tuple of hypersurfaces
such that $\Delta(X)=\bigcap_{i=1}^{\ell} H_i$.

There exist $U\to T\times S$
and $W\subset U$
such that the following holds.
\begin{itemize}
 \item $W$ is a countable union of
       closed Zariski subsets of $U$
       such that
       $W\cap\bigl(
       U\times_{T\times S}(\nbigb(\nbigu)\times A)
       \bigr)=\emptyset$.
 \item For any $J\subset \{1,\ldots,m\}$,
       there exist
       coherent $\nbigd_{T\times S\times X\times X/T\times S}$-modules
       $(\DD\nbigl(\nbigu)\boxtimes M)_{U,J}$
       with morphisms
       $(\DD\nbigl(\nbigu)\boxtimes M)_{U,J}
       \to 
       (\DD\nbigl(\nbigu)\boxtimes M)_{U,J'}$
       for $J\subset J'$
       such that
\[
\begin{CD}
 \bigl((\DD \nbigl(\nbigu)\boxtimes M)_{U,J}\bigr)_{\eta}
 @>>>
 \bigl((\DD \nbigl(\nbigu)\boxtimes M)_{U,J'}\bigr)_{\eta}
 \\
 @V{\simeq}VV @V{\simeq}VV \\
 (\DD \nbigl(\nbigu)\boxtimes M)_{\eta}(\ast H(J)_{\eta})
 @>>>
 (\DD \nbigl(\nbigu)\boxtimes M)_{\eta}(\ast H(J')_{\eta})
\end{CD}
\]       
       for any geometric point $\eta$ of $U$
       whose image is not contained in $W$.
\end{itemize}

Let $U$ denote the $\hyperk$-vector space
with the base $e_1,\ldots,e_{\ell}$.
For an ordered subset $J=(j_1,\ldots,j_{m})$,
we set $e_J=e_{j_1}\wedge\cdots\wedge e_{j_m}$.
We set
\[
 \nbigc^{\app,k}\bigl(\nbigl(\nbigu),M\bigr)
 =\bigoplus_{|J|=k}
 \bigl(\DD \nbigl(\nbigu)\boxtimes M\bigr)_{U,J}
 \otimes e_J.
\]
We obtain the complex
$\nbigc^{\app,\bullet}(\nbigl(\nbigu),M)$.
We also set
\[
 \nbigc^k\bigl(\nbigl(\nbigu),M\bigr)
 =\bigoplus_{|J|=k}
 (\DD \nbigl(\nbigu)\boxtimes M)(\ast H(J))\otimes e_J.
\]
We obtain the complex 
$\nbigc^{\bullet}(\nbigl(\nbigu),M)$.
There exists the natural morphism of
complexes.
\[
 \nbigc^{\app,\bullet}(\nbigl(\nbigu),M)
 \to
 \nbigc^{\bullet}(\nbigl(\nbigu),M).
\]
For any geometric point $\eta$
whose image is not contained in $W$,
we have
\[
 \nbigc^{\app,\bullet}(\nbigl(\nbigu),M)_{\eta}
\simeq
 \nbigc^{\bullet}(\nbigl(\nbigu),M)_{\eta}.
\]

Let $\pi:U\times X\to U$ denote the projection.
We obtain the following morphisms
\[
 R^k\pi_{\ast}
 \nbigc^{\app,\bullet}(\nbigl(\nbigu),M)
 \lrarr
 R^k\pi_{\ast}
 \nbigc^{\bullet}(\nbigl(\nbigu),M).
\]
By refining $U$,
we may assume that
$R^k\pi_{\ast}\nbigc^{\app\bullet}(\nbigl(\nbigu),M)$
are locally free $\nbigo_U$-modules of finite rank.

For any geometric point $\eta$
whose image is not contained in $W$,
we have
\[
\Bigl(
R^k\pi_{\ast}\nbigc^{\app\bullet}(\nbigl(\nbigu),M)
\Bigr)_{\eta}
\simeq
 \Bigl(
 R^k\pi_{\ast}\nbigc^{\app\bullet}(\nbigl(\nbigu),M)
 _{\eta}
 \Bigr)
 \simeq
  \Bigl(
 R^k\pi_{\ast}\nbigc^{\bullet}(\nbigl(\nbigu)_{\eta},M_{\eta})
 \Bigr)
 \simeq
 \Ext^k\Bigl(
 \nbigl(\nbigu)_{\eta},M_{\eta}
 \Bigr).
\]

Let $E^1$ be the vector bundle over $U$
corresponding to
$R^1\pi_{\ast}\nbigc^{\app,\bullet}(\nbigl(\nbigu),M)$.
Let $p:E^1\to U$ denote the projection.
Let $\pi_{1}:E^1\times X\to E^1$ denote the projection.
We obtain 
\[
 p^{\ast}
 R^1\pi_{\ast}\nbigc^{\app,\bullet}(\nbigl(\nbigu),M)
 \to 
 p^{\ast}
 R^1\pi_{\ast}\nbigc^{\bullet}(\nbigl(\nbigu),M)
 \simeq
 R^1\pi_{1\ast}\nbigc^{\bullet}(\nbigl(\nbigu)_{E^1},M_{E^1})
\]
on $E^1$.
The universal section of
$p^{\ast}R^1\pi_{\ast}\nbigc^{\app,\bullet}(\nbigl(\nbigu),M)$
induces the section $s_1$
of
$R^1\pi_{1\ast}\nbigc^{\bullet}(\nbigl(\nbigu)_{E^1},M_{E^1})$
on $E^1$.
Because $E^1$ is an affine scheme,
it induces an element
\[
 \alpha(s_1)
 \in
 \Ext^1_{\nbigd_{E^1\times X/E^1}}
 \Bigl(
 \nbigl(\nbigu)_{E^1},
 M_{E^1}
 \Bigr).
\]
We obtain the extension of
$\nbigd_{E^1\times X/E^1}$-modules
corresponding to $\alpha(s_1)$:
\[
 0\lrarr M_{E^1}
 \lrarr
 \Mtilde_{E^1}
 \lrarr
 \nbigl(\nbigu)_{E^1}
 \lrarr 0.
\]

Let $a\in A\subset S(\hyperk)$
and $b\in \nbigb(\nbigu)\subset T(\hyperk)$.
Let $u\in\Utilde(\hyperk)$ be over $(a,b)$.
There exists an isomorphism
\[
 E^1_{|u}
 \simeq
 \Ext^1_{\nbigd_X}(\nbigl(\nbigu)_b,M_a).
\]
We obtain the following proposition from
Proposition \ref{prop;25.12.14.5}.
\begin{prop}
\label{prop;26.1.31.11}
For $x\in E^1_{|u}$,
the extension 
\[
 0\lrarr M_{E^1|x\times X}
 \lrarr
 \Mtilde_{E^1|x\times X}
 \lrarr
 \nbigl(\nbigu)_{E^1|x\times X}
 \lrarr 0
\]
corresponds to 
$x\in \Ext^1_{\nbigd_X}(\nbigl(\nbigu)_b,M_a)$.
\hfill\qed
\end{prop}

\section{Boundedness for some families of $\nbigd$-modules}

\subsection{Boundedness condition}

Let $\hyperk$ be an algebraically closed field of characteristic $0$.
Let $S$ be a $\hyperk$-variety.
Let $S(\hyperk)$ denote the set of closed points.
Let $\pi:X_S\to S$ be a smooth projective morphism
of relative dimension $n$.
For each $s\in S(\hyperk)$,
we set $X_s=\{s\}\times_SX_S$.
Let $\Mod(\nbigd_{X_s})$ denote the category of
$\nbigd_{X_s}$-modules.

\begin{df}
\label{df;26.2.21.50}
Let $\nbigc_s$ $(s\in S(\hyperk))$
be full subcategories of $\Mod(\nbigd_{X_s})$.
We say that the family $(\nbigc_s\,|\,s\in S(\hyperk))$ is bounded
if there exist the following objects.
\begin{itemize}
 \item A $\hyperk$-variety $\nbigs$ over $S$.
       We set
       $X_{\nbigs}=\nbigs\times_SX_S$.
 \item A coherent $\nbigd_{X_{\nbigs}/\nbigs}$-module
       $\nbigm$
       which is flat over $\nbigs$.
 \item For any $s\in S(\hyperk)$ and 
       any object $M$ in $\nbigc_s$,
       there exists a closed point
       $\stilde\in \nbigs(\hyperk)$ over $s$
       such that
       $\nbigm_s\simeq M$.
\hfill\qed
\end{itemize}
\end{df}

\subsubsection{Boundedness of meromorphic flat connections}

We have the following boundedness.
\begin{thm}
\label{thm;26.2.7.31}
Let $S$, $X_S$, $H_S$ be as in
{\rm\S\ref{subsection;25.12.14.121}}.
Then, the family $(\nbigc^{\good}(X_s,H_s,r,\vecm)\,|\,s\in S(\hyperk))$
is bounded in the sense of Definition {\rm\ref{df;26.2.21.50}}.
\end{thm}
\pf
Let $\nbigs$ and $(E_{\nbigs},\nabla_{\nbigs})$
be as in Theorem \ref{thm;25.10.15.40}.
We set
$\nbige_{\nbigs}=E_{\nbigs}(\ast H_{\nbigs})$.
We obtain the good meromorphic flat bundle
$(\nbige_{\nbigs},\nabla_{\nbigs})$ relative to $\nbigs$.
By refining $\nbigs$,
we may assume that there exists a good lattice
of $(\nbige_{\nbigs},\nabla_{\nbigs})$.
Let $\nbigs_1$ be an irreducible component of $\nbigs$.
We obtain
$(X_{\nbigs_1},H_{\nbigs_1})
=(X_{\nbigs},H_{\nbigs})\times_{\nbigs}\nbigs_1$.
As the restriction,
we obtain
$(\nbige_{\nbigs_1},\nabla_{\nbigs_1})$
and the induced lattice $E_{\nbigs_1}$.
Let $H_{\nbigs_1}=\bigcup_{j\in \Lambda_1}H_{\nbigs_1,j}$
be the irreducible decomposition.
Let $E^{(1)}_{\nbigs_1}$
be a good lattice of
$(\nbige_{\nbigs_1},\nabla_{\nbigs_1})$
satisfying Condition \ref{condition;26.2.6.2}.
We take $\vecm_1\in\seisuu_{>0}^{\Lambda_1}$
satisfying Condition \ref{condition;26.2.6.3}
for $E^{(1)}_{\nbigs_1}$.
We also assume the following condition.
\begin{itemize}
 \item $E^{(1)}_{\nbigs_1}
       \subset E_{\nbigs_1}(-H_{\nbigs_1})$
       and
       $E_{\nbigs_1}(H_{\nbigs_1})\subset
       E^{(1)}_{\nbigs_1}(\vecm_1H_{\nbigs_1})$.
       It implies
       $\nbiga(E_{\nbigs_1})
       \subset\nbigb(E^{(1)}_{\nbigs_1},\vecm_1)$.
       (See Lemma {\rm\ref{lem;26.3.3.10}} for $\nbiga(E_{\nbigs_1})$.)
\end{itemize}
We obtain the coherent $\nbigd_{X_{\nbigs_1}/\nbigs_1}$-module
$\nbigd_{X_{\nbigs_1}/\nbigs_1}\bigl(
 E^{(1)}_{\nbigs_1}(\vecm_1H_{\nbigs_1})
 \bigr)$
which is flat over $\nbigs_1$.
Applying this construction to each irreducible component of $\nbigs$,
we obtain the claim of the theorem.
\hfill\qed

\begin{thm}
\label{thm;26.2.23.1}
Let $X$, $H$, $r$, $\vecm$ be as in
{\rm\S\ref{subsection;26.2.7.32}}.
Then, the family $\nbigc(X,H,r,\vecm)$
is bounded 
in the sense of Definition {\rm\ref{df;26.2.21.50}}.
\end{thm}
\pf
Let $\nbigs$, $\Xhat_{\nbigs}$,
$\rho_{\nbigs}:\Xhat_{\nbigs}\to \nbigs\times X$,
$H_{\nbigs}$, 
a good meromorphic flat bundle
$(\nbigv,\nabla)$
on $(\Xhat_{\nbigs},\Hhat_{\nbigs})$
relative to $\nbigs$,
and a lattice $\Ehat_{\nbigs}$
as in Theorem \ref{thm;25.12.14.30}.
By refining $\nbigs$,
we may assume that
$(\nbigv,\nabla)$ has a good lattice.
As in the proof of Theorem \ref{thm;26.2.7.31},
we construct
a coherent $\nbigd_{\Xhat_{\nbigs}/\nbigs}$-module $\nbigm$
such that the following holds.
\begin{itemize}
 \item 
       $\nbigm_{s}\simeq \nbigv_s$
       for any $u\in \nbiga(E_{\nbigs})$.
\end{itemize}
By considering
$\rho_{\nbigs\dagger}\nbigm$,
we obtain the claim of the theorem.
\hfill\qed

\subsection{Boundedness of holonomic $\nbigd$-modules with dominated
characteristic cycles}

Let $X$ be a smooth projective $\hyperk$-varieties.
Let $\Lambda=\sum m_i\Lambda_i$,
where $m_i\in\seisuu_{>0}$,
and $\Lambda_i$ are complex Lagrangian cone in $T^{\ast}X$.
For any holonomic $\nbigd_X$-module $M$,
we have the characteristic cycle $\CC(M)$.
Let $\nbigc(X,\Lambda)$
denote the family of holonomic $\nbigd_X$-modules $M$
such that $\CC(M)\leq \sum m_i\Lambda_i$,
i.e.,
$\CC(M)=\sum n_i\Lambda_i$ with $n_i\leq m_i$.
\begin{thm}
\label{thm;25.12.14.20}
The family $\nbigc(X,\Lambda)$ is bounded
in the sense of Definition {\rm\ref{df;26.2.21.50}}.
\end{thm}

\subsubsection{Boundedness of minimal extensions}

Let $Z\subset X$ be an irreducible closed subset,
which is not necessarily normal.
Let $H$ be a hypersurface of $X$
such that $Z\setminus H$ is smooth.
We set $H_Z=H\cap Z$.
Let $\iota_Z:Z\setminus H_Z\to X\setminus H$ denote the immersion.

A meromorphic connection on $(Z,H_Z)$ means
an algebraic $\nbigo_{Z\setminus H_Z}$-module $V$
with an integrable connection $\nabla$.
It induces a $\nbigd_{X\setminus H}$-module
$\iota_{Z\dagger}(V)$.
The minimal extension is denoted by
$(\iota_{Z\dagger}V)_{\min}$.

Let $\Mero_{\min}(Z,H_Z,\Lambda)$
denote the category of $\nbigd_X$-modules
obtained as 
$(\iota_{Z\dagger}V)_{\min}$
for a meromorphic flat connection $V$ on $(Z,H_Z)$
such that
\[
 \CC\bigl(
 (\iota_{Z\dagger}V)_{\min}
 \bigr)
 \leq\Lambda.
\]

\begin{prop}
\label{prop;26.1.31.10}
There exists a cell
$\nbigu=(S,Z_S,\rho,\nbige,\nabla,E,\vecm)$
such that the following holds.
\begin{itemize}
 \item For any
       $(\iota_{Z\dagger}V)_{\min}
       \in \Mero_{\min}(Z,H_Z,\Lambda)$,
       there exists
       $a \in \nbigb(\nbigu)$ such that
       $(\iota_{Z\dagger}V)_{\min}
       \simeq
       \nbigl(\nbigu)_a$.
\end{itemize}
\end{prop}
\pf
Let $\varphi:\Ztilde\to Z$ denote the normalization.
Note that $\varphi$ is a finite map.
Let $\Ztilde^{\sm}$ denote the smooth part.
We set $H_{\Ztilde}=\varphi^{-1}(H_Z)$.

Let $H=\bigcup_{i\in\Lambda}H_i$ denote
the irreducible decomposition.
Let $H^{\circ}$ denote the smooth part of $H$.
We set $H_i^{\circ}=H_i\cap H^{\circ}$
which is Zariski open in $H_i$.
There exists a Zariski open subset
$H^{\sankaku}_{Z,i}\subset H_{Z,i}^{\circ}$
such that
$\varphi^{-1}(H^{\sankaku}_{Z,i})\to H_{Z,i}^{\sankaku}$
is etale.

Let $(V,\nabla)$ be a meromorphic flat bundle
on $(\Ztilde,H_{\Ztilde})$.
Let $r$ denote the multiplicity of
$T^{\ast}_{Z^{\circ}}X$ in $\Lambda$.
We obtain
$\rank(V)\leq r$.

Let $V_{\min}$ denote the minimal extension
on $\Ztilde^{\sm}$.
We set
$V_!=\DD\bigl(\DD(V)(\ast H_{\Ztilde})\bigr)$.
There exists the natural morphism
$V_{!}\to V$,
and $V_{\min}$ equals the image.
Let $C(V)$ and $K(V)$ denote the cokernel and the kernel.
Their supports are contained in $H_{\Ztilde}$.
On $H^{\sankaku}_{Z,i}$,
we have
$\varphi^j_{\dagger}(C(V))=\varphi^j_{\dagger}(K(V))=0$
unless $j=0$.
Hence,
around $H^{\sankaku}_{Z,i}$,
we have
$\varphi_{\dagger}(V_{\min})=
\bigl(
 \varphi_{\dagger}(V)
 \bigr)_{\min}$.

\begin{lem}
\label{lem;25.12.14.60}
Let $\CC_{i}(V_{\min})$
denote the multiplicity of $T^{\ast}H_{\Ztilde_i}$
in $\CC(V_{\min})$.
Let $k_i$ denote the multiplicity of
$T^{\ast}_{H_i}X$ in $\Lambda$.
We obtain
$\CC_i(V_{\min})\leq k_i$.
\end{lem}
\pf
It is enough to study the case $\hyperk=\cnum$.
We can work in the complex analytic setting.
By \cite[Proposition 9.4.2, Theorem 11.3.3]{Kashiwara-Schapira},
around $H^{\sankaku}_{Z,i}$,
we obtain
\[
\varphi_{\ast}(\CC(V_{\min})) 
=\CC\bigl(
(\varphi_{\ast}V)_{\min}
\bigr).
\]
Here $\varphi_{\ast}$ for Lagrangian cone is
defined as in
\cite[Definition 9.3.3]{Kashiwara-Schapira}.

Let $P$ be a general point of $H_{\Ztilde,i}$.
There exists a neighbourhood $U$ of $\varphi(P)$
such that for the irreducible decomposition
$Z\cap U=\bigcup Z_{P,j}$,
we have $Z_{P,j}\cap Z_{P,k}=H_{Z,i}$ $(j\neq k)$.
Also, $\varphi^{-1}(Z_j)\to Z_j$ are homeomorphisms.
Note that
$\varphi_{\ast}(T^{\ast}_{H_{\Ztilde,i}}\Ztilde)
 =T^{\ast}_{H_{Z,i}}X$
because $\varphi_{\ast}\cnum_{H_{\Ztilde,i}}=\cnum_{H_{Z,i}}$.
Hence, we obtain the claim of Lemma \ref{lem;25.12.14.60}.
\hfill\qed

\vspace{.1in}

The following lemma is easy to see.
\begin{lem}
There exist $N_i$ depending only on $k_i$ and $r$
such that the following holds.
\begin{itemize}
 \item  Let $\Ztilde^{\circ}$ denote the smooth part of $\Ztilde$.
	Let $\iota_{Z\dagger}(V)_{\min}\in
	\Mero(Z,H_Z,\Lambda)$.
        For any $\Sigma\in\vecI(V)$,
	the closure of $\Sigma$
	in $T^{\ast}\Ztilde^{\circ} (\sum N_iH_i)$
	is proper over $\Ztilde^{\circ}$.
	\hfill\qed
\end{itemize}
\end{lem}

Then, we obtain
the claim of Proposition \ref{prop;26.1.31.10}.
\hfill\qed

\subsubsection{Proof of Theorem \ref{thm;25.12.14.20}}

\begin{df}
An ordered partition of $\Lambda$ is
a decomposition
$\Lambda=\sum \Lambda^{(i)}$,
where $\Lambda^{(i)}=\sum m_j^{(i)}\Lambda_j$
for some non-negative integers $m_j^{(i)}$.
Let $\nbigp(\Lambda)$ denote the set of ordered partition.
\hfill\qed
\end{df}

We obtain the following proposition
from Proposition \ref{prop;26.1.31.10}
and Proposition \ref{prop;26.1.31.11}.

\begin{prop}
Let $\vecLambda=(\Lambda^{(i)})\in\nbigp(\Lambda)$.
There exists a complex variety $S(\vecLambda)$
 and a coherent $\nbigd_{X\times S/S}$-module
 $\nbigm(\vecLambda)$
with a filtration $\nbigf$ by $\nbigd$-submodules
such that the following holds.
\begin{itemize}
 \item $(\nbigm(\vecLambda),\nbigf)$ is flat over $S$,
       and $\nbigm(\vecLambda)$ is holonomic.
 \item Let $M$ be a holonomic $\nbigd_X$-module
       with a Jordan-H\"older decomposition $\nbigf$
       such that
       $\Gr^{\nbigf}_j(M)\leq \Lambda^{(j)}$.
       Then, there exists
       $s\in S$
       such that
       $(M,\nbigf)\simeq
       (\nbigm(\vecLambda),\nbigf)_s$.
       \hfill\qed
\end{itemize}
\end{prop}

Let $M$ be a holonomic $\nbigd$-module on $X$
such that $\CC(M)\leq \Lambda$.
There exists a Jordan-H\"older filtration $\nbigf(M)$.
We have
$\CC(M)=\sum \CC(\Gr^{\nbigf}_i(M))$.
There exists partition $\vecLambda\in\nbigp(\Lambda)$
such that $\CC(\Gr^{\nbigf}_i(M))\leq \Lambda^{(i)}$.
There exists $s\in \nbigs(\vecLambda)$
such that
$M\simeq \nbigm(\vecLambda)_{s}$.
Therefore,
$\nbigs=\bigsqcup_{\vecLambda\in\nbigp(\Lambda)} \nbigs(\vecLambda)$
with the induced $\nbigd_{X\times \nbigs/\nbigs}$-module
satisfies the desired property in
Theorem \ref{thm;25.12.14.20}.
\hfill\qed

\subsubsection{Appendix: Constraint on the characteristic cycles}

Let $M$ be a coherent $\nbigd_{S\times X/S}$-module
which is flat over $S$ and holonomic.
Let $\Lambda_i\subset T^{\ast}X$ $(i=1,\ldots,\ell)$
be complex Lagrangian cones.
Let $m_i$ $(i=1,\ldots,\ell)$ be positive integers.
We set $\Lambda=\sum m_i\Lambda_i$.

\begin{lem}
\label{lem;26.2.24.1}
There exists a constructible subset $S'\subset S$
such that 
$\CC(M_s)\leq \Lambda$ if and only if $s\in S'$. 
\end{lem}
\pf
Let $T\subset S$ be a locally closed smooth subvariety.
There exists a Zariski open subset $T'\subset T$
such that
$M_{T'}$ has a good filtration flat over $T'$.
We obtain the characteristic variety $\Ch(M_{T'})$
and the characteristic cycle $\CC(M_{T'})$
by using the good filtration.

Let $\hyperk_1$ be an algebraic closure of
the generic point of $T$.
Let $\eta:\Spec\hyperk_1\to T$
be the induced morphism.
We obtain $M_{\eta}$ on $X_{\hyperk_1}$.
We obtain a good filtration of $M_{\eta}$
from the good filtration of $M_{T'}$.
We have
$\Ch(M_{\eta})=\Ch(M_{T'})_{\eta}$.
If $\Ch(M_{\eta})\not\subset |\Lambda|\times_{\hyperk}\hyperk_1$,
there exists a Zariski open subset $T'\subset T$
such that
$\Ch(M_s)\not\subset \bigcup \Lambda_{i,\hyperk_1}$
 for any $s\in T'$.
If $\Ch(M_{\eta})=\bigcup \Lambda_{i,\hyperk_1}$,
we have
$\CC(M_{T'})=\sum n_i(T'\times\Lambda_i)$
for some $n_i\in\seisuu_{\geq 0}$.
We have $n_i\leq m_i$ if and only if
$\CC(M_{\eta})\leq \sum m_i \Lambda_{i,\hyperk_1}$.
Then, we obtain the claim of the lemma
by using a Noetherian induction.
\hfill\qed

\end{document}